# MULTIRATE METHODS FOR ORDINARY DIFFERENTIAL EQUATIONS

MICHAEL GÜNTHER* AND ADRIAN SANDU†

**Abstract.** This survey provides an overview of state-of-the art multirate schemes, which exploit the different time scales in the dynamics of a differential equation model by adapting the computational costs to different activity levels of the system. We start the discussion with the straightforward approach based on interpolating and extrapolating the slow–fast coupling variables; the multirate Euler scheme, used as a base example, falls into this class. Next we discuss higher order multirate schemes that generalize classical singlerate linear multistep, Runge-Kutta, and extrapolation methods.

**1. Introduction.** Throughout this paper we consider dynamical systems described by the initial-value problem (IVP) of ordinary differential equations (ODEs)

$$\dot{\mathbf{y}} = \mathbf{f}(\mathbf{y}), \quad t \in [t_0, t_\mathrm{f}], \quad \mathbf{y}(t_0) = \mathbf{y}_0, \quad \mathbf{y}(t) \in \mathbb{R}^d, \tag{1.1}$$

with $\mathbf{f} : \mathbb{R}^d \to \mathbb{R}^d$ Lipschitz continuous in $\mathbf{y}$ to obtain a unique solution. For simplicity, and without loss of generality, we consider here autonomous dynamics.

Many dynamical systems (1.1) display a multiple time scale dynamics, with some components of the system evolving at faster pace, and others evolving at a slower pace. Different time scales may be associated with different activity levels of various components (e.g., fast signals in active transistors and slow voltage changes in latent components of an integrated circuit), or with different processes that drive the dynamics (e.g., fast chemical reactions and slow tracer transport driving pollutant concentrations in the atmosphere).

*Multirate time discretization schemes leverage the different time scales in the dynamics of a differential equation model by adapting the step sizes, and therefore the computational costs, to different activity levels of the system. Specifically, faster parts are solved with smaller time steps, and slower parts with larger time steps. The goal of multirating is to provide solutions of the target accuracy while considerably improving the overall computational efficiency when compared to traditional singlerate methods.*

The remainder of this paper is organized as follows. Section 2 discusses multiscale dynamics, and the benefits and challenges of multirate numerical integration. The historical – and intuitive – approach to constructing multirate schemes based on interpolating and extrapolating the coupling variables between subsystems is discussed in Section 3. As an intrinsic example for such schemes multirate Euler schemes are discussed in full detail in Section 4, as well as linear multistep schemes in Section 6. Multirate Runge-Kutta schemes in Section 5, especially the class of GARK schemes, provide another approach on defining multirate schemes besides the idea of interpolating and extrapolating. Another example are linearly-implicit GARK schemes, an extension of linearly-implicit Runge-Kutta schemes, which we discuss in Section 7.

*Bergische Universität Wuppertal, School of Mathematics and Sciences, IMACM. D-42097 Wuppertal, Germany (guenther@uni-wuppertal.de).

†Computational Science Laboratory, Department of Computer Science, Virginia Tech, Blacksburg, VA 24061 (sandu@cs.vt.edu).

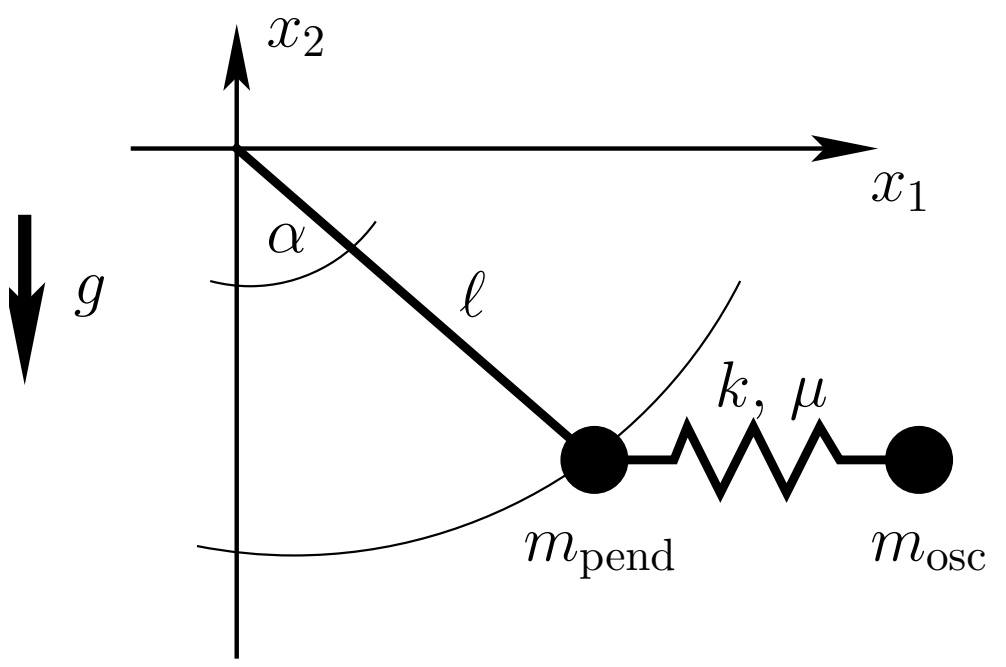

Fig. 1: Mathematical pendulum coupled to an oscillator (taken after [3]).

Multirate generalizations of extrapolation schemes are investigated in Section 8 and multirate infinitesimal schemes, again considered in the GARK framework, in Section 9. Finally, conclusions are drawn in Section 10.

## 2. Multiscale dynamics and its efficient numerical solution.

**2.1. Multiscale dynamics with partitioned components evolving at different timescales.** An ODE system (1.1) where different components display a multiscale behavior, can be partitioned into two subsystems $\mathbf{y}^{\{\text{s}\}} \in \mathbb{R}^{d^{\{\text{s}\}}}$ (slow part) and $\mathbf{y}^{\{\text{F}\}} \in \mathbb{R}^{d^{\{\text{F}\}}}$ (fast part), with $d^{\{\text{s}\}} + d^{\{\text{F}\}} = d$ and usually $d^{\{\text{F}\}} \ll d^{\{\text{s}\}}$, according to their dynamics,

$$\begin{bmatrix} \dot{\mathbf{y}}^{\{\text{s}\}} \\ \dot{\mathbf{y}}^{\{\text{F}\}} \end{bmatrix} = \begin{bmatrix} \mathbf{f}^{\{\text{s}\}}(\mathbf{y}^{\{\text{s}\}}, \mathbf{y}^{\{\text{F}\}}) \\ \mathbf{f}^{\{\text{F}\}}(\mathbf{y}^{\{\text{s}\}}, \mathbf{y}^{\{\text{F}\}}) \end{bmatrix}, \qquad \begin{bmatrix} \mathbf{y}^{\{\text{s}\}}(t_0) \\ \mathbf{y}^{\{\text{F}\}}(t_0) \end{bmatrix} = \begin{bmatrix} \mathbf{y}_0^{\{\text{s}\}} \\ \mathbf{y}_0^{\{\text{F}\}} \end{bmatrix}, \tag{2.1}$$

with $\mathbf{f}^{\{\text{s}\}} : \mathbb{R}^{d^{\{\text{s}\}}} \times \mathbb{R}^{d^{\{\text{F}\}}} \to \mathbb{R}^{d^{\{\text{s}\}}}$, $\mathbf{f}^{\{\text{F}\}} : \mathbb{R}^{d^{\{\text{s}\}}} \times \mathbb{R}^{d^{\{\text{F}\}}} \to \mathbb{R}^{d^{\{\text{F}\}}}$ and $\mathbf{f}^{\{\text{s}\}}$ and $\mathbf{f}^{\{\text{F}\}}$ are assumed to be Lipschitz continuous in both $\mathbf{y}^{\{\text{s}\}}$ and $\mathbf{y}^{\{\text{F}\}}$ with corresponding Lipschitz constants $L^{\{\text{s},\text{s}\}}, L^{\{\text{s},\text{F}\}}, L^{\{\text{F},\text{s}\}}$, and $L^{\{\text{F},\text{F}\}}$, respectively. Such a splitting naturally arises when different systems are coupled in a multiphsics setting, for example. Note that in this case the splitting is fixed and does not change with respect to time. Here the multirate approach is to discretize the evolution of slow components $\mathbf{y}^{\{\text{s}\}}$ using large step sizes, and employ small step sizes to resolve fast components $\mathbf{y}^{\{\text{F}\}}$ such that they are approximated accurately enough.

Example 1 (Coupled pendulum). *As an example of a system with component-wise multiscale behavior we consider a mathematical pendulum of constant length $\ell$ that is coupled to a damped oscillator with a horizontal degree of freedom, which is discussed in [3] and illustrated in Figure 1. The system consists of two rigid bodies: the first mass $m_{\text{pend}}$ is connected to a second mass $m_{\text{osc}}$ by a soft spring.*

*The minimal set of coordinates $\mathbf{q} = [\alpha, x_1]^\top$ uniquely describes the position of both bodies. The equations of motion given by the Euler-Lagrange equation:*

$$\begin{pmatrix} m_{\text{pend}}\,\ell & 0 \\ 0 & m_{\text{osc}} \end{pmatrix} \ddot{\mathbf{q}} = \begin{pmatrix} -m_{\text{pend}}\, g\, \sin(\alpha) + \cos(\alpha)\, F \\ -F \end{pmatrix} =: \mathbf{f}(\mathbf{q}, \dot{\mathbf{q}}),$$

*with the spring-damper force $F$ modeled as*

$$F = k\,(x_1 - \ell\, \sin(\alpha)) + \mu\,(\dot{x}_1 - \ell\, \dot{\alpha}\, \cos(\alpha)),$$

*where $k$ denotes the spring stiffness and $\mu$ is the coefficient of friction. This two-dimensional second order system of differential equations can be easily transformed into a four-dimensional ODE system (1.1) by introducing the derivatives $\dot{\mathbf{q}}$ as additional unknowns. The time evolution of the two components is illustrated in Figure 2. We see that $\alpha(t)$ oscillates quickly, while $x_1(t)$ oscillates slowly.*

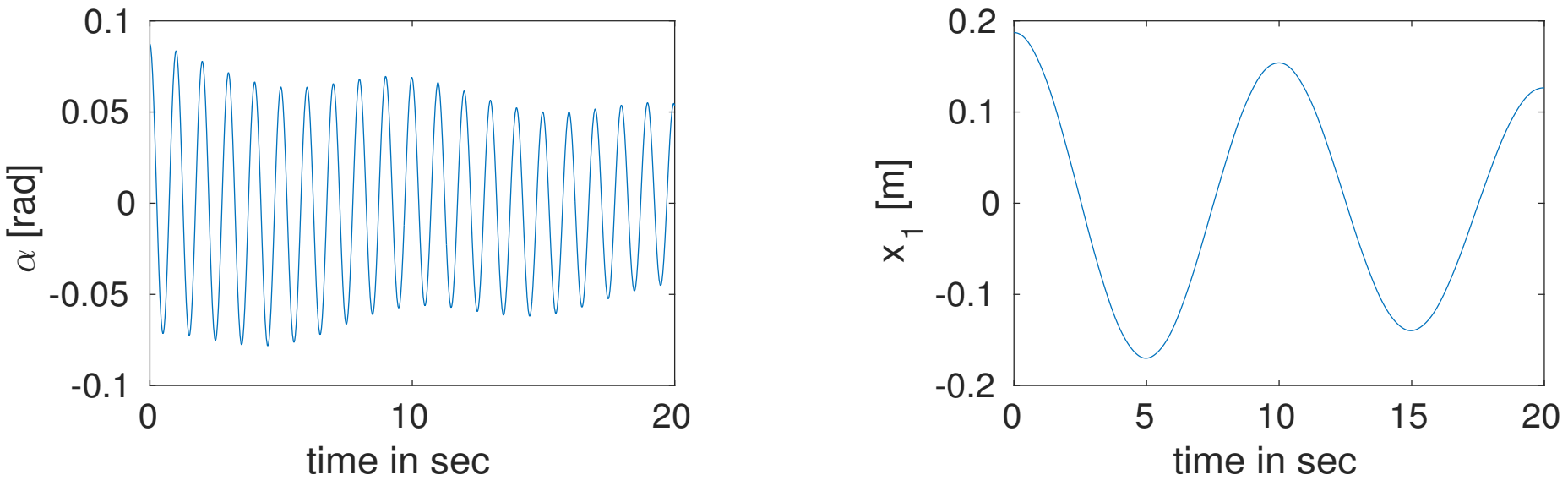

(a) Fast angle $\alpha(t)$. (b) Slow coordinate $x_1(t)$.

Fig. 2: Time evolution of components of the two-body problem.

**2.2. Multiscale dynamics driven by multiple physical processes with different timescales.** The right-hand-side $\mathbf{f}(\mathbf{y})$ of (1.1) represents the cumulative effect of all physical processes that drive the evolution of the system. A second cause of multiscale dynamics is due to the system being driven by multiple physical processes with different dynamics. In this case the right hand side $\mathbf{f}(\mathbf{y})$ can be split into two components of different activity levels $\mathbf{f}^{\{\mathrm{S}\}}(\mathbf{y})$ and $\mathbf{f}^{\{\mathrm{F}\}}(\mathbf{y})$, and we end up with the additively split initial-value problem

$$\dot{\mathbf{y}} = \mathbf{f}^{\{\mathrm{S}\}}(\mathbf{y}) + \mathbf{f}^{\{\mathrm{F}\}}(\mathbf{y}), \quad t \in [t_0, t_{\mathrm{f}}], \quad \mathbf{y}(t_0) = \mathbf{y}_0. \tag{2.2}$$

Here $\mathbf{f}^{\{\mathrm{F}\}}$ defines a fast changing and $\mathbf{f}^{\{\mathrm{S}\}}$ a slowly changing process. If $\mathbf{f}^{\{\mathrm{F}\}}$ is inexpensive to evaluate, but $\mathbf{f}^{\{\mathrm{S}\}}$ is expensive to evaluate, which is the case in many applications such as lattice quantum-chromoynamics [63], then the multirate potential to improve computational efficiency lies in evaluating the slow and expensive part $\mathbf{f}^{\{\mathrm{S}\}}$ less often than the fast and cheap part $\mathbf{f}^{\{\mathrm{F}\}}$.

Example 2 (Coupled pendulum, continued). *Recall the rigid body from Example 1. The right-hand side is also a source of multirate behavior. Specifically, the right hand side $\mathbf{f}(\mathbf{q}, \dot{\mathbf{q}})$ can be split into two components $\mathbf{f}(\mathbf{q}, \dot{\mathbf{q}}) = \mathbf{f}^{\{\mathrm{S}\}}(\mathbf{q}, \dot{\mathbf{q}}) + \mathbf{f}^{\{\mathrm{F}\}}(\mathbf{q})$ with*

$$\mathbf{f}^{\{\mathrm{S}\}}(\mathbf{q}, \dot{\mathbf{q}}) = \begin{pmatrix} \frac{\cos(\alpha) F}{m_{\mathrm{pend}} \ell} \\ -\frac{F}{m_{\mathrm{osc}}} \end{pmatrix} \quad \textit{and} \quad \mathbf{f}^{\{\mathrm{F}\}}(\mathbf{q}) = \begin{pmatrix} -g \sin(q_l) \\ 0 \end{pmatrix}$$

*of different activity. As seen in Figure 3, $\mathbf{f}^{\{\mathrm{F}\}}$ defines a fast changing force, whereas $\mathbf{f}^{\{\mathrm{S}\}}$ is slowly changing. If the evaluation of the spring-damper force $F$ becomes expensive, the multirate potential is obvious.*

**2.3. Equivalence of partitioned and split systems.** Note that the partitioned and the additively split systems are formulations that can be transformed into

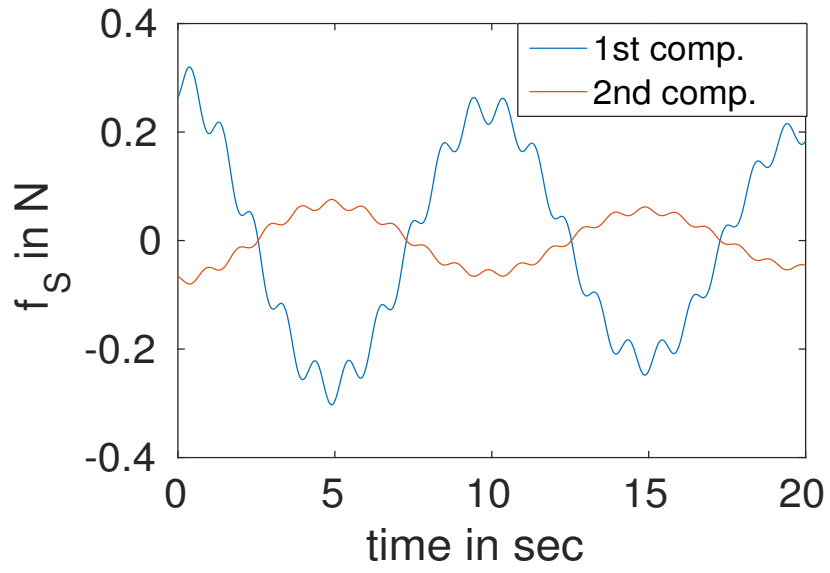

(a) Components of slow force $\mathbf{f}^{\{S\}}(\mathbf{y}(t))$.

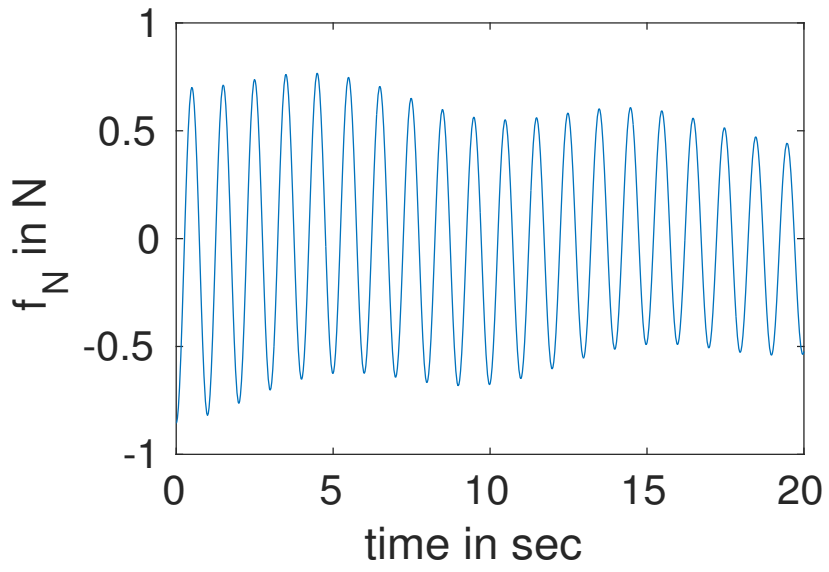

(b) Fast force (1st component of $\mathbf{f}^{\{F\}}(\mathbf{y}(t))$.

Fig. 3: Time evolution of slow and fast forces acting on the two-body system.

one another [67]. For example (2.1) can be written in the form (2.2) as follows:

$$\text{(2.3a)} \qquad \mathbf{y}(t) \equiv \begin{bmatrix} \mathbf{y}^{\{S\}}(t) \\ \mathbf{y}^{\{F\}}(t) \end{bmatrix},$$

$$\text{(2.3b)} \qquad \dot{\mathbf{y}} = \begin{bmatrix} \dot{\mathbf{y}}^{\{S\}} \\ \dot{\mathbf{y}}^{\{F\}} \end{bmatrix} = \begin{bmatrix} \mathbf{f}^{\{S\}}(\mathbf{y}) \\ 0 \end{bmatrix} + \begin{bmatrix} 0 \\ \mathbf{f}^{\{F\}}(\mathbf{y}) \end{bmatrix}, \quad t \in [t_0, t_f].$$

Similarly, (2.2) can be written in the form (2.1) as follows:

$$\text{(2.4a)} \qquad \mathbf{y}(t) \equiv \mathbf{y}^{\{S\}}(t) + \mathbf{y}^{\{F\}}(t) - \mathbf{y}_0, \quad t \in [t_0, t_f],$$

$$\text{(2.4b)} \qquad \begin{bmatrix} \dot{\mathbf{y}}^{\{S\}} \\ \dot{\mathbf{y}}^{\{F\}} \end{bmatrix} = \begin{bmatrix} \mathbf{f}^{\{S\}}(\mathbf{y}^{\{S\}} + \mathbf{y}^{\{F\}} - \mathbf{y}_0) \\ \mathbf{f}^{\{F\}}(\mathbf{y}^{\{S\}} + \mathbf{y}^{\{F\}} - \mathbf{y}_0) \end{bmatrix}, \quad \begin{bmatrix} \mathbf{y}^{\{S\}}(t_0) \\ \mathbf{y}^{\{F\}}(t_0) \end{bmatrix} = \begin{bmatrix} \mathbf{y}_0 \\ \mathbf{y}_0 \end{bmatrix}.$$

We note that when a split system is rewritten in partitioned form (2.4) the number of variables doubles, and the full solution is the sum of individual solutions [67].

REMARK 1 (Arbitrary number of partitions). *The equivalency between additively and component partitioned forms has been shown for two partitions. This equivalency holds for an arbitrary number* N *of partitions.*

We will use the partitioned system (2.1) as our default formulation in this article. Since the two formulations are equivalent it suffices to perform the accuracy and stability analyses for one form.

**2.4. Numerical solutions of multiscale systems.** We now discuss the challenges faced by traditional time discretizations when applied to multiscale dynamics.

**2.4.1. For a given level of accuracy step sizes are limited by the fast-changing components.** Consider a system (2.1) with a solution involving slow-varying and fast-varying components:

$$\text{(2.5a)} \qquad \begin{bmatrix} \dot{\mathbf{y}}^{\{S\}} \\ \dot{\mathbf{y}}^{\{F\}} \end{bmatrix} = \begin{bmatrix} \mathbf{A}^{\{S,S\}} & \mathbf{v} \\ \mathbf{u}^T & \lambda^{\{F,F\}} \end{bmatrix} \begin{bmatrix} \mathbf{y}^{\{S\}} - \mathbf{v}\cos(\omega t) \\ \mathbf{y}^{\{F\}} - \cos(\mathtt{m}\,\omega\, t) \end{bmatrix} - \begin{bmatrix} \omega \sin(\omega t)\,\mathbf{v} \\ \mathtt{m}\,\omega \sin(\mathtt{m}\,\omega\, t) \end{bmatrix},$$

$$\text{(2.5b)} \qquad t_0 = 0, \quad \begin{bmatrix} \mathbf{y}^{\{S\}}(0) \\ \mathbf{y}^{\{F\}}(0) \end{bmatrix} = \begin{bmatrix} \mathbf{v} \\ 1 \end{bmatrix},$$

(2.5c) $$\mathbf{u}, \mathbf{v}, \mathbf{y}^{\{\mathrm{S}\}}(t) \in \mathbb{R}^d, \quad \mathbf{y}^{\{\mathrm{F}\}}(t) \in \mathbb{R}, \quad \mathbf{A}^{\{\mathrm{S},\mathrm{S}\}}\,\mathbf{v} = \lambda^{\{\mathrm{S},\mathrm{S}\}}\,\mathbf{v}.$$

The system (2.5) has a solution involving slow and fast components:

$$\mathbf{y} = \begin{bmatrix} \mathbf{y}^{\{\mathrm{S}\}} \\ \mathbf{y}^{\{\mathrm{F}\}} \end{bmatrix} = \begin{bmatrix} \cos(\omega\, t)\,\mathbf{v} \\ \cos(\mathtt{m}\,\omega\, t) \end{bmatrix}, \quad \mathtt{m} > 1.$$

Application of forward Euler method with step size $H$ over $[t_0, t_1 = t_0 + H]$ leads to the following local truncation error:

$$\begin{bmatrix} \boldsymbol{\delta}_1^{\{\mathrm{S}\}} \\ \boldsymbol{\delta}_1^{\{\mathrm{F}\}} \end{bmatrix} = \begin{bmatrix} \mathbf{y}_1^{\{\mathrm{S}\}} - \mathbf{y}^{\{\mathrm{S}\}}(t_1) \\ \mathbf{y}_1^{\{\mathrm{F}\}} - \mathbf{y}^{\{\mathrm{F}\}}(t_1) \end{bmatrix} = \frac{H^2}{2} \begin{bmatrix} \ddot{\mathbf{y}}^{\{\mathrm{S}\}}(t_a) \\ \ddot{\mathbf{y}}^{\{\mathrm{F}\}}(t_a) \end{bmatrix} = -\frac{H^2}{2} \begin{bmatrix} \omega^2\, \cos(\omega\, t_a)\mathbf{v} \\ \mathtt{m}^2\, \omega^2\, \cos(\mathtt{m}\,\omega\, t_a) \end{bmatrix}.$$

To ensure a local truncation error smaller than a desired tolerance `tol` in each component we need to restrict step sizes as follows:

$$H \le H^{\{\mathrm{S}\}}_{\max} = \frac{\sqrt{2\,\mathtt{tol}}}{\omega} \quad \Rightarrow \quad |\boldsymbol{\delta}_1^{\{\mathrm{S}\}}| \le \mathtt{tol}; \quad H \le H^{\{\mathrm{F}\}}_{\max} = \frac{\sqrt{2\,\mathtt{tol}}}{\mathtt{m}\,\omega} \quad \Rightarrow \quad |\boldsymbol{\delta}_1^{\{\mathrm{F}\}}| \le \mathtt{tol}.$$

In order to achieve the desired level of accuracy using a traditional single-step implementation of forward Euler, the step size is restricted by the fastest component since $H^{\{\mathrm{F}\}}_{\max} = H^{\{\mathrm{S}\}}_{\max}/\mathtt{m}$. We make the following remarks.

1. Integration with a step size $H \le H^{\{\mathrm{F}\}}_{\max}$ leads to a slow local truncation error $|\boldsymbol{\delta}_1^{\{\mathrm{S}\}}| \le \mathtt{tol}/\mathtt{m}^2$ much smaller than required.
2. Let $c^{\{\mathrm{F}\}}$ and $c^{\{\mathrm{S}\}}$ be the costs of evaluating the fast and slow components time derivatives $\mathbf{f}^{\{\mathrm{F}\}}$ and $\mathbf{f}^{\{\mathrm{S}\}}$, respectively. In (2.5) the dimension of the slow subsystem is much larger than the dimension of the fast subsystem, therefore $c^{\{\mathrm{S}\}} \gg c^{\{\mathrm{F}\}}$, which is the case of interest in multirating. Integration over $[t_0, t_\mathrm{f}]$ requires $N^{\{\mathrm{F}\}}$ time step to keep local truncation errors below `tol`, at a cost:

   (2.6)
   $$N^{\{\mathrm{F}\}} = \left\lceil \frac{t_\mathrm{f} - t_0}{H^{\{\mathrm{F}\}}} \right\rceil = \left\lceil \frac{\mathtt{m}\,\omega\,(t_\mathrm{f} - t_0)}{\sqrt{2\,\mathtt{tol}}} \right\rceil, \quad \mathrm{Cost}^{\mathrm{singlerate}} = N^{\{\mathrm{F}\}} \left( c^{\{\mathrm{F}\}} + c^{\{\mathrm{S}\}} \right).$$

   The slow component is integrated with a higher accuracy than needed, at the expense of a computational cost $c^{\{\mathrm{S}\}} \cdot \mathcal{O}(\mathtt{m})$ that is higher than needed.

**2.4.2. For numerical stability step sizes of explicit schemes are limited by the fastest dynamics.** Consider a decoupled version of system (2.5) with

$$\mathbf{u} = \mathbf{v} = \mathbf{0}_{d\times 1}, \quad \mathbf{A}^{\{\mathrm{S},\mathrm{S}\}} = \lambda^{\{\mathrm{S},\mathrm{S}\}}\,\mathbf{I}_{d\times d}, \quad \lambda^{\{\mathrm{F},\mathrm{F}\}}, \lambda^{\{\mathrm{S},\mathrm{S}\}} < 0, \quad |\lambda^{\{\mathrm{F},\mathrm{F}\}}| = \mathtt{m}\,|\lambda^{\{\mathrm{S},\mathrm{S}\}}|.$$

The system components $\mathbf{y}^{\{\mathrm{S}\}}$ and $\mathbf{y}^{\{\mathrm{F}\}}$ evolve with fast and slow dynamics driven by the corresponding eigenvalues $\lambda^{\{\mathrm{S},\mathrm{S}\}}$ and $\lambda^{\{\mathrm{F},\mathrm{F}\}}$, respectively. Application of forward Euler method with step size $H$ has the following numerical stability restriction:

$$H \le H^{\{\mathrm{S}\}}_{\max} = \frac{2}{|\lambda^{\{\mathrm{S},\mathrm{S}\}}|}; \qquad H \le H^{\{\mathrm{F}\}}_{\max} = \frac{2}{|\lambda^{\{\mathrm{F},\mathrm{F}\}}|} = \frac{H^{\{\mathrm{S}\}}_{\max}}{\mathtt{m}}.$$

Numerical integration using a traditional singlerate forward Euler implementation requires the step size to be restricted by the fastest dynamics since $H^{\{\mathrm{F}\}}_{\max} = H^{\{\mathrm{S}\}}_{\max}/\mathtt{m}$. Again, the slow component is integrated with lower step size than needed for slow numerical stability, at the expense of a computational cost that is higher than needed.

**2.4.3. Solution for computational efficiency: treat the fast and the slow components differently.** The discussion above indicates that, for computational efficiency, one needs a multimethod approach that discretizes the slow and fast components differently. The multimethod selection strategy depends on the goals of the simulation, as shown in Figure 4.

- If only an accurate slow component is of interest, then one damps out the fast scales using an implicit-explicit (IMEX) approach: integrate the fast component implicitly, and the slow component explicitly with appropriate slow step size. This alleviates the stability limitation on the step size due to the fast component, discussed in Subsection 2.4.2. Specifically, the implicit fast discretization can accommodate large step sizes (appropriate for the slow component) without numerical instability.
- If the goal is to obtain both slow and fast accurate solutions then we employ a multirate strategy. The slow component is solved with a large step $H$, and the fast component with multiple small steps $h = H/\mathtt{m}$. This strategy alleviates the accuracy limitation on the step size due to the fast component, discussed in Subsection 2.4.1. In the example discussed in Subsection 2.4.1 discretize the slow component with $N^{\{\mathrm{S}\}}$ large steps of size $H$:

$$H \equiv H^{\{\mathrm{S}\}} \le \frac{\sqrt{2\,\mathtt{tol}}}{\omega} \quad \Rightarrow \quad |\boldsymbol{\delta}^{\{\mathrm{S}\}}| \le \mathtt{tol}; \quad N^{\{\mathrm{S}\}} = \left\lceil \frac{t_\mathrm{f} - t_0}{H^{\{\mathrm{S}\}}} \right\rceil = \left\lceil \frac{\omega\,(t_\mathrm{f} - t_0)}{\sqrt{2\,\mathtt{tol}}} \right\rceil,$$

  and the fast component with $N^{\{\mathrm{F}\}}$ small time steps of size $h$

$$h = \frac{H}{\mathtt{m}} \le \frac{\sqrt{2\,\mathtt{tol}}}{\mathtt{m}\,\omega} \quad \Rightarrow \quad |\boldsymbol{\delta}^{\{\mathrm{F}\}}| \le \mathtt{tol}, \quad N^{\{\mathrm{F}\}} = \mathtt{m}\,N^{\{\mathrm{S}\}},$$

  where $N^{\{\mathrm{F}\}}$ is given in (2.6). The cost of the multirate forward Euler approach is

$$\text{Cost}^{\text{multirate}} = N^{\{\mathrm{S}\}}\,c^{\{\mathrm{S}\}} + N^{\{\mathrm{F}\}}\,c^{\{\mathrm{F}\}} = N^{\{\mathrm{F}\}}\left(\frac{1}{\mathtt{m}}\,c^{\{\mathrm{S}\}} + c^{\{\mathrm{F}\}}\right) \tag{2.7}$$

  is to be compared with the cost of the singlerate approach (2.6). The fast evaluation cost remains the same, while the slow evaluation cost decreases by a factor of $\mathtt{m}$. For the case $c^{\{\mathrm{S}\}} \gg c^{\{\mathrm{F}\}}$ the multirate approach realizes important cost savings without degrading the solution accuracy. The singlerate approach (2.6) spends $\mathcal{O}(\mathtt{m}\,c^{\{\mathrm{S}\}})$ cost to evaluate the expensive slow function, and computes the slow component more accurately than required, at a cost considerably higher than needed.

Application of a multimethod approach that discretizes the slow and fast components differently faces a number of challenges, discussed next.

**2.4.4. Clean separation of time scales.** A prerequisite for efficient use of multirate schemes is that different components are characterized by different time scales, and no scale dynamics significantly impacts other scales.

In general, it can be said that with strong coupling between the subsystems, the fast time constant generally affects the slow subsystem, which will therefore no longer be "purely slow", therefore reducing the multirate potential.

With weak coupling, on the other hand, the fast time constant (if at all) will only enter the slow system as a slight perturbation, the magnitude of which will be negligible compared to the slow waveforms. The partitioning of the right-hand side in Figure 3 provides an example of this, where the coupling strength is given by the stiffness and the friction of the spring.

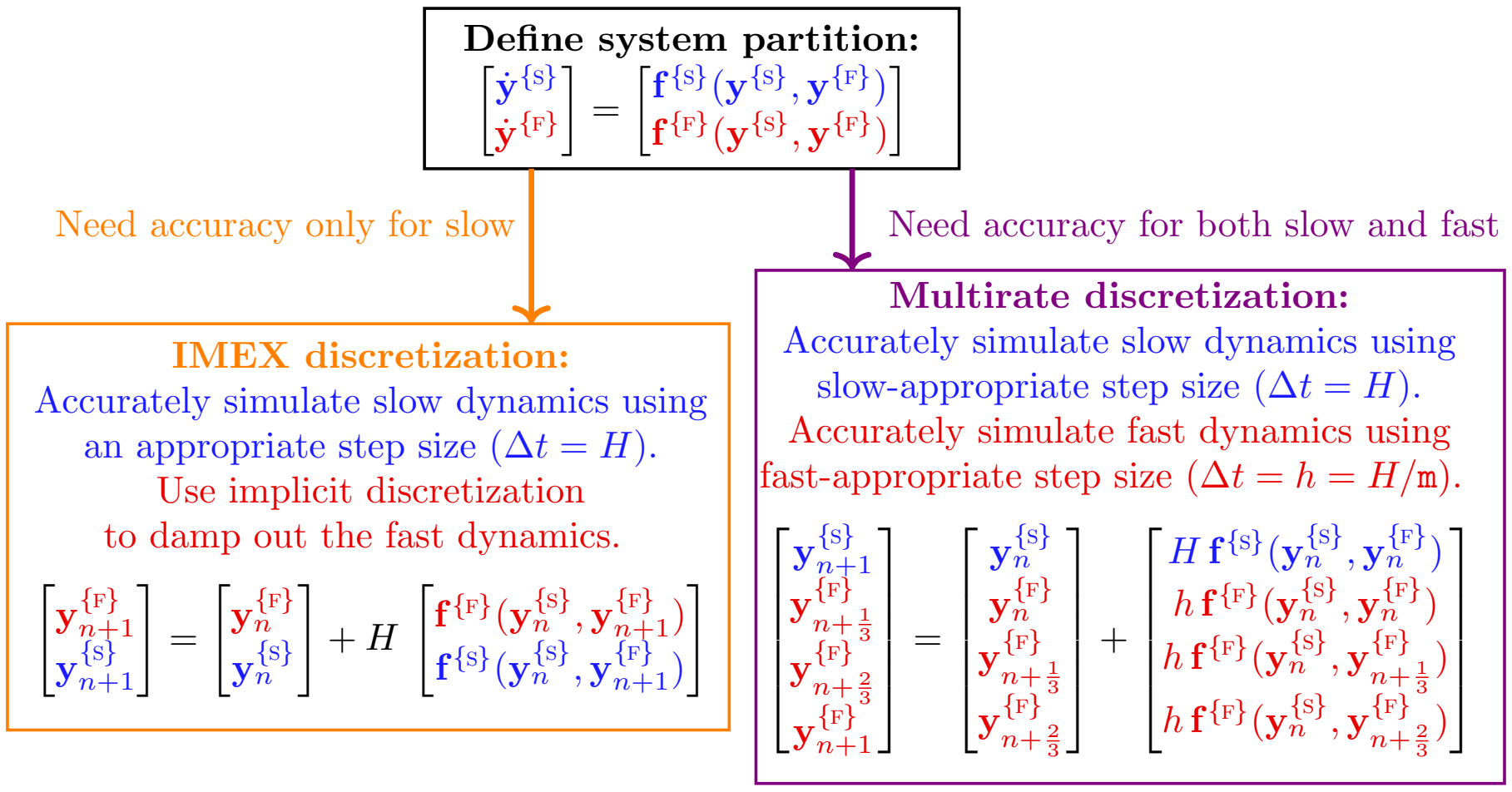

Fig. 4: Strategies to select an efficient multimethod that discretizes differently the fast and the slow components

Consider as an example for this behavior the linear system

$$\begin{bmatrix} \dot{\mathbf{y}}^{\{S\}}(t) \\ \dot{\mathbf{y}}^{\{F\}}(t) \end{bmatrix} = \underbrace{\begin{bmatrix} \lambda^{\{S,S\}} & \lambda^{\{F,S\}} \\ \lambda^{\{S,F\}} & \lambda^{\{F,F\}} \end{bmatrix}}_{A} \begin{bmatrix} \mathbf{y}^{\{S\}}(t) \\ \mathbf{y}^{\{F\}}(t) \end{bmatrix}, \quad \begin{bmatrix} \mathbf{y}^{\{S\}}(0) \\ \mathbf{y}^{\{F\}}(0) \end{bmatrix} = \begin{bmatrix} \mathbf{y}^{\{S\}}_0 \\ \mathbf{y}^{\{F\}}_0 \end{bmatrix}.$$

Application of forward Euler method with step size $H$ over $[t_0, t_1 = t_0 + H]$ leads to the following local truncation error:

$$\begin{bmatrix} \boldsymbol{\delta}^{\{S\}}_1 \\ \boldsymbol{\delta}^{\{F\}}_1 \end{bmatrix} = -\frac{H^2}{2}\, A^2 \exp\{A\tau\}\, \mathbf{y}_0 = -\frac{H^2}{2} \begin{bmatrix} {\lambda^{\{S,S\}}}^2\,(1 + \delta(\rho\mu + \mu + 1)) \\ {\lambda^{\{F,F\}}}^2 \left(1 + \rho\left(\frac{\delta}{\mu} + (1 + \frac{1}{\mu})\right)\right) \end{bmatrix}$$

for some $\tau \in [0, H]$, where we start from the initial value $\mathbf{y}_0 = \exp(-A\tau)\,[1,1]^{\mathrm{T}}$ and use the notation

$$\mu := \frac{\lambda^{\{F,F\}}}{\lambda^{\{S,S\}}}, \quad |\mu| \gg 1, \quad \delta := \frac{\lambda^{\{F,S\}}}{\lambda^{\{S,S\}}}, \quad \rho := \frac{\lambda^{\{S,F\}}}{\lambda^{\{F,F\}}},$$

for the ratio of the time scales $\mu$, the coupling from the fast into the slow part $\delta$ and the coupling from the slow into the fast part $\rho$.

To ensure a local truncation error smaller than a desired tolerance tol in each component we need to restrict step sizes as follows:

$$\begin{aligned} H &\leq H^{\{S\}}_{\max} = \frac{\sqrt{2\,\mathtt{tol}}}{|\lambda^{\{S,S\}}|}\, \frac{1}{\sqrt{|1 + \delta(1 + \mu + \rho\mu)|}} \quad \Rightarrow \quad |\boldsymbol{\delta}^{\{S\}}_1| \leq \mathtt{tol}; \\ H &\leq H^{\{F\}}_{\max} = \frac{\sqrt{2\,\mathtt{tol}}}{|\mu||\lambda^{\{S,S\}}|}\, \frac{1}{\sqrt{|1 + \rho(1 + \frac{1+\delta}{\mu})|}} \quad \Rightarrow \quad |\boldsymbol{\delta}^{\{F\}}_1| \leq \mathtt{tol}. \end{aligned}$$

In order to achieve the desired level of accuracy using forward Euler, the step size is restricted by the fastest component. If we have only a weak coupling, i.e., $|\delta| \ll 1/|\mu|$ and $|\rho| \ll |1|$, the ratio of the step sizes is given by the ratio of the time constants (see Subsection 2.4.1): $H_{\max}^{\{F\}} \approx H_{\max}^{\{S\}}/|\mu|$. However, even in the case of only a one-sided coupling the multirate potential may be lost: for example, with $\rho = 0$ and $\delta = \mu - 1$ we get $H_{\max}^{\{F\}} = H_{\max}^{\{S\}}$. The slow component is no longer slow, as the coupling allows the fast mode to pollute the slow component. Any step size ratio of $H_{\max}^{\{F\}}$ is possible, depending on the coupling:

$$\frac{H_{\max}^{\{F\}}}{H_{\max}^{\{S\}}} \to \begin{cases} 0, & \delta \text{ fixed and } \rho \to -\mu/(\mu+1+\delta), \\ \infty, & \rho \text{ fixed and } \delta \to -1/(1+\mu+\rho\mu). \end{cases}$$

**2.4.5. Challenge: partitioning the system into components with different time scales.** In order to apply a computationally efficient multimethod approach on has to first partition the system into fast and slow parts, as shown in Figure 4.

The partitioning raises two main practical questions for multirate schemes. First, how can a given problem be split into (coupled) subproblems with separated time scales? Secondly, can this partition be done automatically? If the problem is given by a coupled problem consisting of subsystems with different time scales, as described above, then partitioning is naturally given by dividing the problem into the different subsystems. Otherwise, the answer to these questions is somewhat more subtle. If the partition is not given by the (coupled) problem itself, an algorithm is required to determine the partitions automatically. Various approaches have been discussed in the literature here (stiffness, nonlinearity, comparison between extrapolated and numerical solution, component-wise step-size prediction, to name just a few), but dynamic partitioning involves large overheads that may exceed the computational gains obtained by multirating.

**2.4.6. Challenge: building numerical schemes that correctly capture the coupling between subsystems.** While the discussion above focuses on individual fast and slow systems, in problems of interest the fast and slow subsystems evolve together in a coupled manner. Evolving different components with different time steps raises the question of how to treat the coupling: at intermediate points one has computed either the slow or the fast component, while the other is missing. How can we coordinate the partial solutions with different time steps such that the overall solution is stable and accurate? This review will discuss in detail the construction of numerical schemes that correctly capture the couplings, such that the overall solution meets the desired stability and accuracy requirements.

**2.4.7. When should I use a multirate approach, and when not?.** Multirate schemes involve overheads associated with partitioning the system and with the more complex book-keeping of information. Multirate schemes are effective when computational gains exceed these overheads.

As a rule of thumb, multirate schemes are particularly efficient when applied to a problem with a fixed, problem-defined partition into a slow and a fast part, where evaluation of the slow part requires a high computational effort, while evaluation of the fast part is inexpensive. The computational gains come from evaluating the expensive slow part less often.

Should the system be characterized by more than two separated time scales, then mutirate schemes defined for two time scales can be applied in a hierarchical manner to take advantage of the multiple time scales.

**3. Interpolation and extrapolation based multirate schemes.** A natural way to define multirate schemes for the partitioned system (2.1) is based on extrap-olating and interpolating the coupling variables. Specifically, the coupling variables in (2.1) are replaced over a macro-step from $t_n$ to $t_{n+1} = t_n + H$ by approximations $\tilde{\mathbf{y}}^{\{\mathrm{S}\}}(t)$ and $\tilde{\mathbf{y}}^{\{\mathrm{F}\}}(t)$ obtained via interpolation or extrapolation from previously computed solution quantities. This leads to the modified system:

$$\text{(3.1a)}\quad \dot{\mathbf{y}}^{\{\mathrm{S}\}} = \mathbf{f}^{\{\mathrm{S}\}}\left(\mathbf{y}^{\{\mathrm{S}\}}, \tilde{\mathbf{y}}^{\{\mathrm{F}\}}(t)\right) =: \tilde{\mathbf{f}}^{\{\mathrm{S}\}}(t, \mathbf{y}^{\{\mathrm{S}\}}),$$

$$\text{(3.1b)}\quad \dot{\mathbf{y}}^{\{\mathrm{F}\}} = \mathbf{f}^{\{\mathrm{F}\}}\left(\tilde{\mathbf{y}}^{\{\mathrm{S}\}}(t), \mathbf{y}^{\{\mathrm{F}\}}\right) =: \tilde{\mathbf{f}}^{\{\mathrm{F}\}}(t, \mathbf{y}^{\{\mathrm{F}\}}), \quad t \in [t_n, t_{n+1}] \subseteq [t_0, t_\mathrm{f}].$$

One step of a multirate scheme solves the modified system (3.1) over $[t_n, t_{n+1}] \subseteq [t_0, t_\mathrm{f}]$ using a large single step of macro step size $H$ for the subsystem $\mathbf{y}^{\{\mathrm{S}\}}$, and $\mathtt{m} \in \mathbb{N}$ steps of (micro step) size $h = H/\mathtt{m}$ for $\mathbf{y}^{\{\mathrm{F}\}}$. During the solve, the respective coupling variables need to be evaluated at the intermediate steps and stages. Depending on the sequence of computation of the unknowns $\mathbf{y}^{\{\mathrm{S}\}}$ and $\mathbf{y}^{\{\mathrm{F}\}}$, one distinguishes the following versions of extra-/and interpolation techniques as suggested in [22]:

1. *Slowest-first approach.* One starts with a slow first step to solve (3.1a) and advance slow variables $\mathbf{y}^{\{\mathrm{S}\}}$ from $t_n$ to $t_{n+1}$. The modified slow equation (3.1a) uses extrapolated values $\tilde{\mathbf{y}}^{\{\mathrm{F}\}}$ based on information available at $t_n$ to evaluate the fast coupling variable in the current macrostep. Next, we perform $\mathtt{m}$ micro-steps to solve (3.1b) and advance the fast variables $\mathbf{y}^{\{\mathrm{F}\}}$ from $t_n$ to $t_{n+1}$. The modified fast equation (3.1b) uses interpolated values $\tilde{\mathbf{y}}^{\{\mathrm{S}\}}$ based on the information available from the recently computed slow solution in the macro step.
2. *Fastest-first approach.* One starts with $\mathtt{m}$ micro steps to solve (3.1b), using extrapolated slow values $\tilde{\mathbf{y}}^{\{\mathrm{S}\}}$ based on information available at $t_n$. Next, one macro step is performed to solve the slow dynamics (3.1a) from $t_n$ to $t_{n+1}$, using interpolated fast values $\tilde{\mathbf{y}}^{\{\mathrm{F}\}}$ based on available information from the recently computed fast solution.
3. *Fully-decoupled approach.* The modified fast (3.1b) and slow equations (3.1a) variables are integrated in parallel using in both cases extrapolated coupling variable based on information available at $t_n$, the beginning of the current macro step [6].

We note that the above list of possible couplings is not exhaustive.

Remark 2. *The restriction that the extrapolation can only be based on the information at $t_n$ can be relaxed to employ the data of the preceding macro step $[t_n - H, t_n]$. In fact, one can encode such an information (e.g., as a spline model) and update and transport it from macro step to macro step [62].*

For ODE systems, all variants of extrapolation/interpolation-based multirate schemes have convergence order $p$ (in the final asymptotic phase) provided that the following hold [5]:

1. the basic integration schemes (i.e., the schemes used to solve the slow and the fast subsystems with given coupling data, respectively) have order $p$, and
2. the extrapolation/interpolation schemes are of approximation order $p - 1$.

Theorem 3.1 (Consistency of slowest-first multirate schemes [5]). *Consider the coupled ODE (2.1) where $\mathbf{f}^{\{\mathrm{S}\}}$ and $\mathbf{f}^{\{\mathrm{F}\}}$ are Lipschitz continuous. Construct the modified system (3.1) using extrapolation and interpolation procedures of order $p - 1$ for both the slow $\tilde{\mathbf{y}}^{\{\mathrm{S}\}}$ and the fast $\tilde{\mathbf{y}}^{\{\mathrm{F}\}}$ coupling variables. Furthermore, we apply two basic integration schemes of order $p$ to solve (3.1): one to solve (3.1a) for $\mathbf{y}^{\{\mathrm{S}\}}$ using one macro-step of size $H$, and a second one to solve (3.1b) for $\mathbf{y}^{\{\mathrm{F}\}}$ using $\mathtt{m} \in \mathbb{N}$* *steps of size $h = H/\mathtt{m}$. If the basic integration schemes are applied in a slowest-first approach, the resulting slowest-first multirate scheme has order $p$.*

*Proof.* Consider the IVP system (2.1) with exact initial condition at time $t_n$, $\mathbf{y}_n^{\{S\}} = \mathbf{y}^{\{S\}}(t_n)$ and $\mathbf{y}_n^{\{F\}} = \mathbf{y}^{\{F\}}(t_n)$. The unique solution of (2.1) over $[t_n, t_{n+1}]$ is $(\mathbf{y}^{\{S\}}(t),\, \mathbf{y}^{\{F\}}(t))$. Consider also the modified system (3.1) with exact initial condition $\mathbf{y}_n^{\{S\}} = \mathbf{y}^{\{S\}}(t_n)$ and $\mathbf{y}_n^{\{F\}} = \mathbf{y}^{\{F\}}(t_n)$, and denote its unique solution by $(\hat{\mathbf{y}}^{\{S\}}(t),\, \hat{\mathbf{y}}^{\{F\}}(t))$.

In the slowest-first approach, the coupling variable $\tilde{\mathbf{y}}^{\{F\}}$ is given by extrapolation of order $p-1$, with a constant $L_F > 0$, using information on $\mathbf{y}^{\{F\}}$ in $[0, t_n]$:

$$\|\mathbf{y}^{\{F\}}(t) - \tilde{\mathbf{y}}^{\{F\}}(t)\| \le L_F \cdot H^p + \mathcal{O}(H^{p+1}) \quad \text{for any } t \in [t_n,\, t_{n+1}]. \tag{3.2}$$

The interpolated values $\tilde{\mathbf{y}}^{\{S\}}$ of order $p-1$ are based on the slow numerical approximations already computed in $[t_n, t_{n+1}]$ such that it holds (with a constant $L_S > 0$):

$$\|\mathbf{y}^{\{S\}}(t) - \tilde{\mathbf{y}}^{\{S\}}(t)\| \le L_S \cdot H^p + \mathcal{O}(H^{p+1}) \quad \text{for any } t \in [t_n,\, t_{n+1}]. \tag{3.3}$$

Now, we apply the two basic integration schemes of order $p$ in multirate fashion to the decoupled model (3.1) and denote the numerical solution $(\mathbf{y}_{n+1}^{\{S\}},\, \mathbf{y}_{n+1}^{\{F\}})$.

Then, the distance between the multirate approximation and the exact solution can be estimated as

(3.4)

$$\begin{pmatrix} \|\mathbf{y}_{n+1}^{\{S\}} - \mathbf{y}^{\{S\}}(t_{n+1})\| \\ \|\mathbf{y}_{n+1}^{\{F\}} - \mathbf{y}^{\{F\}}(t_{n+1})\| \end{pmatrix} \le \begin{pmatrix} \|\mathbf{y}_{n+1}^{\{S\}} - \hat{\mathbf{y}}^{\{S\}}(t_{n+1})\| \\ \|\mathbf{y}_{n+1}^{\{F\}} - \hat{\mathbf{y}}^{\{F\}}(t_{n+1})\| \end{pmatrix} + \begin{pmatrix} \|\hat{\mathbf{y}}^{\{S\}}(t_{n+1}) - \mathbf{y}^{\{S\}}(t_{n+1})\| \\ \|\hat{\mathbf{y}}^{\{F\}}(t_{n+1}) - \mathbf{y}^{\{F\}}(t_{n+1})\| \end{pmatrix}.$$

a) The numerical schemes of order $p$ give for the first term on the right-hand side:

$$\begin{pmatrix} \|\mathbf{y}_{n+1}^{\{S\}} - \hat{\mathbf{y}}^{\{S\}}(t_{n+1})\| \\ \|\mathbf{y}_{n+1}^{\{F\}} - \hat{\mathbf{y}}^{\{F\}}(t_{n+1})\| \end{pmatrix} \le \begin{pmatrix} c_S H^{p+1} + \mathcal{O}(H^{p+2}) \\ c_F H^{p+1} + \mathcal{O}(H^{p+2}) \end{pmatrix}$$

employing constants $c_S$, $c_F > 0$ (for leading errors). The fast error estimate has the order of a local truncation error due to the fixed number $\mathtt{m}$ of micro-steps.
b) We now estimate the second term on the right-hand side of (3.4). For the slow part and fast part we get, using Lipschitz continuity of $\mathbf{f}^{\{S\}}$ and $\mathbf{f}^{\{F\}}$, the extrapolation and interpolation estimates (3.2) and (3.3), and finally applying Gronwall's lemma,

$$\begin{aligned} \|\hat{\mathbf{y}}^{\{S\}}(t_{n+1}) - \mathbf{y}^{\{S\}}(t_{n+1})\| &\le L^{\{S,F\}} L_F\, e^{L^{\{S,S\}} H}\, H^{p+1} + \mathcal{O}(H^{p+2}), \\ \|\hat{\mathbf{y}}^{\{F\}}(t_{n+1}) - \mathbf{y}^{\{F\}}(t_{n+1})\| &\le L^{\{F,S\}} L_S e^{L^{\{F,F\}} H}\, H^{p+1} + \mathcal{O}(H^{p+2}). \end{aligned}$$

For more details see [5].
c) Finally, we need to form the total error as the sum of both error terms in the right-hand side (3.4):

$$\begin{aligned} \|\mathbf{y}_{n+1}^{\{S\}} - \mathbf{y}^{\{S\}}(t_{n+1})\| &\le \left(c_S + L^{\{S,F\}} L_F e^{L^{\{S,S\}} H}\right) H^{p+1} + \mathcal{O}(H^{p+2}), \\ \|\mathbf{y}_{n+1}^{\{F\}} - \mathbf{y}^{\{F\}}(t_{n+1})\| &\le \left(c_F + L^{\{F,S\}} L_S e^{L^{\{F,F\}} H}\right) H^{p+1} + \mathcal{O}(H^{p+2}). \end{aligned}$$

This concludes the proof. ∎

Remark 3. *If one uses extrapolation and interpolation schemes of order $p$ instead of $p-1$, one has to replace the term $L_S H^p$ by $L_S H^{p+1}$, and $L_F H^p$ by $L_F H^{p+1}$ which yields the estimate*

$$\begin{aligned}
\|\mathbf{y}_{n+1}^{\{S\}} - \mathbf{y}^{\{S\}}(t_{n+1})\| &\leq c_S H^{p+1} + \mathcal{O}(H^{p+2}), \\
\|\mathbf{y}_{n+1}^{\{F\}} - \mathbf{y}^{\{F\}}(t_{n+1})\| &\leq c_F H^{p+1} + \mathcal{O}(H^{p+2}),
\end{aligned}$$

*that is, the multirate error is dominated by the error of the base numerical integration schemes.*

The proof can be slightly adapted to verify the convergence result for the fastest-first and fully-decoupled approach as well:

Corollary 3.2 (Consistency of fastest-first [5] and fully-decoupled [6] multirate schemes). *The convergence result of Theorem 3.1 remains valid if the slowest-first approach is replaced by the fastest-first or fully-decoupled one.*

Remark 4. *For the basic integration schemes employed in Theorem 3.1 and corollary 3.2 we can choose either (a) one-step integration schemes, or (b) two zero-stable multistep schemes with appropriate initializations.*

Remark 5 (Extrapolation schemes for one-step schemes). *To keep the characteristics of one-step schemes, the extrapolation and interpolation must be based on information in the current macro step size $[t_n, t_{n+1}]$. Extrapolation of order 0 (or 1) can be easily obtained from the initial solution data at $t_n$ (and the derivative information at $t_n$ provided by the ODE). This allows directly the construction of multirate one-step methods of order 1 (or 2).*

*Higher-order multirate schemes require higher-order extra-/inter-po-la-tors, which can be constructed only if information of previous time steps is used. Generally, this may turn a one-step scheme into a multi-step scheme, and raise questions concerning stability. However, if the extrapolation is computed sequentially in a spline-oriented fashion, the modified functions $\tilde{\mathbf{f}}^{\{S\}}$ and $\tilde{\mathbf{f}}^{\{F\}}$ are the same for all time intervals inside $[t_0, t_f]$, then the multirate scheme based on extrapolation / interpolation can still be considered as a one-step scheme applied to the modified ODE equations [62].* □

This convergence result shows that any one-step integration scheme of order $p$ can be used as a base scheme for a multirate scheme of order $p$ provided that the interpolation/extrapolation procedures are at least of order $p-1$. Note that one has to use a nonautonomous implementation of the base scheme, as the coupling variables enter the modified equations (3.1) as time-dependent inputs. An example of such schemes is given by multirate ROW schemes derived in [26].

We will apply the convergence results derived in this section to Euler schemes and multirate linear multistep schemes, which are discussed in the following chapters 4 and 6, respectively.

**4. A canonical example: multirate Euler schemes.** As a canonical example we discuss how the basic Euler integration schemes can be generalized to multirate methods. We discuss the convergence and stability properties of different variants, which drastically differ with respect to efficiency and computational costs.

**4.1. Multirate forward Euler method.** Let us first start to develop multirate forward Euler schemes (MRFE) to exploit multirate behavior given by the different dynamics in the components (component partitioning).

**4.1.1. Method formulation.** Consider the partitioned initial-value problem (2.1). Since the rates of change are different we seek to discretize the slow components with a large step size $H$ and the fast components with a small step size $h = H/\mathtt{m}$.

The discrete approximation times $t_n, t_{n+1}, \ldots$ are given by the large step size $H$. The intermediate approximation points corresponding to the small step size $h$ are considered as fractions of the full steps. Thus we will use the notation:

$$t_{n+1} = t_n + H \quad \text{and} \quad t_{n+\ell/\mathtt{m}} = t_n + \ell\, h, \quad \ell = 0, \ldots, \mathtt{m}. \tag{4.1}$$

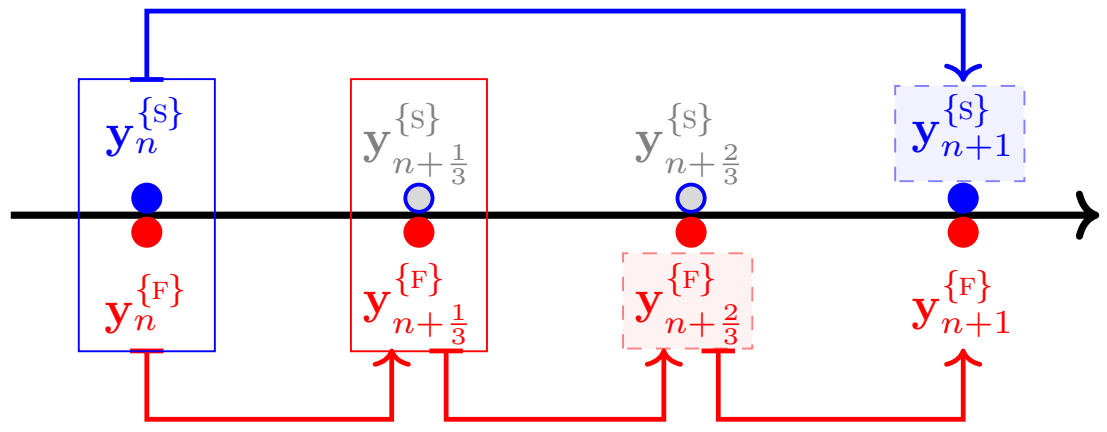

Fig. 5: Cartoon of MRFE method for $\mathtt{m} = 3$. The values computed by Euler steps are represented by filled discs. Empty circles represent the slow values not computed by Euler steps, and that need to be approximated.

The slow MRFE step (4.2a) uses the known values $(\mathbf{y}_n^{\{\mathrm{S}\}}, \mathbf{y}_n^{\{\mathrm{F}\}})$ to compute $\mathbf{y}_{n+1}^{\{\mathrm{S}\}}$:

$$\text{Slow:} \quad \begin{cases} t_{n+1} = t_n + H, \\ \mathbf{y}_{n+1}^{\{\mathrm{S}\}} = \mathbf{y}_n^{\{\mathrm{S}\}} + H\, \mathbf{f}^{\{\mathrm{S}\}}(\mathbf{y}_n^{\{\mathrm{S}\}}, \mathbf{y}_n^{\{\mathrm{F}\}}). \end{cases} \tag{4.2a}$$

Step (4.2a) is illustrated by the blue arrow in Figure 5.

The MRFE scheme takes $\mathtt{m}$ successive small steps (4.2b) to advance the fast component:

$$\text{Fast:} \quad \begin{cases} h := H/\mathtt{m}; \qquad t_{n+(\ell+1)/\mathtt{m}} = t_n + (\ell+1))\, h, \quad \ell = 0, \ldots, \mathtt{m}-1; \\ \mathbf{y}_{n+(\ell+1)/\mathtt{m}}^{\{\mathrm{F}\}} = \mathbf{y}_{n+\ell/\mathtt{m}}^{\{\mathrm{F}\}} + h\, \mathbf{f}^{\{\mathrm{F}\}}(\mathrm{y}_{n+\ell/\mathtt{m}}^{\{\mathrm{S}\}}, \mathbf{y}_{n+\ell/\mathtt{m}}^{\{\mathrm{F}\}}). \end{cases} \tag{4.2b}$$

Steps (4.2b) are illustrated by the red arrows in Figure 5.

Remark 6 (Computational costs). *One step of MRFE (4.2) requires one slow function evaluation and* $\mathtt{m}$ *fast function evaluations, while the regular FE with step size* $h$ *requires* $\mathtt{m}$ *evaluations of both functions. This leads to considerable computational savings whenever the cost to evaluate* $\mathbf{f}^{\{\mathrm{S}\}}$ *greatly exceeds the cost to evaluate* $\mathbf{f}^{\{\mathrm{F}\}}$*.*

As shown in Figure 5, the first fast step uses the known values $(\mathbf{y}_n^{\{\mathrm{S}\}}, \mathbf{y}_n^{\{\mathrm{F}\}})$ to compute $\mathbf{y}_{n+1/3}^{\{\mathrm{F}\}}$. The second fast step needs both $(\mathrm{y}_{n+1/3}^{\{\mathrm{S}\}}, \mathbf{y}_{n+1/3}^{\{\mathrm{F}\}})$ to compute $\mathbf{y}_{n+2/3}^{\{\mathrm{F}\}}$. However, $\mathrm{y}_{n+1/3}^{\{\mathrm{S}\}}$ is *not* computed explicitly by the slow scheme (4.2a). Thus the intermediate slow variables $\mathrm{y}_{n+\ell/\mathtt{m}}^{\{\mathrm{S}\}}$ for $\ell = 1, \ldots, \mathtt{m}-1$ need to be approximated; different approximation strategies lead to different flavors of the multirate scheme.

In the *slowest–first strategy* one starts with solving (4.2a) to obtain the slow solution $\mathbf{y}_{n+1}^{\{\mathrm{S}\}}$. The intermediate slow values are then obtained by interpolation:

$$\mathrm{y}_{n+\ell/\mathtt{m}}^{\{\mathrm{S}\}} = \begin{cases} \mathbf{y}_n^{\{\mathrm{S}\}}, & \text{constant interpolation}, \\ \frac{\mathtt{m}-\ell}{\mathtt{m}} \mathbf{y}_n^{\{\mathrm{S}\}} + \frac{\ell}{\mathtt{m}} \mathbf{y}_{n+1}^{\{\mathrm{S}\}}, & \text{linear interpolation}. \end{cases} \tag{4.2c}$$

In the *fastest–first strategy* one starts with solving (4.2b) to advance the fast components to $\mathbf{y}^{\{\mathrm{F}\}}_{n+1}$. The intermediate slow values are then obtained as:

$$\text{(4.2d)}\qquad \mathrm{y}^{\{\mathrm{S}\}}_{n+\ell/\mathtt{m}} = \begin{cases} \mathbf{y}^{\{\mathrm{S}\}}_n, & \text{constant extrapolation,} \\ \mathbf{y}^{\{\mathrm{S}\}}_n + \frac{\ell}{\mathtt{m}}\, H\, \mathbf{f}^{\{\mathrm{S}\}}(\mathbf{y}^{\{\mathrm{S}\}}_n, \mathbf{y}^{\{\mathrm{F}\}}_n), & \text{linear Taylor series extrapolation.} \end{cases}$$

Note that the fastest–first and slowest–first multirate forward Euler methods coincide. The results of Section 3 show that the multirate forward Euler scheme has convergence order one, provided that we use at least constant extrapolation.

#### 4.1.2. Linear stability analysis of multirate forward Euler.

*A two-dimensional test problem.* Dahlquist's test problem needs to be generalized to account for two time scale problems in order to investigate the stability of multirate schemes. Following the analysis done by Kværnø [36], we consider the following generalized linear test problem:

$$\text{(4.3)}\qquad \begin{bmatrix} \dot{\mathbf{y}}^{\{\mathrm{S}\}}(t) \\ \dot{\mathbf{y}}^{\{\mathrm{F}\}}(t) \end{bmatrix} = \underbrace{\begin{bmatrix} \lambda^{\{\mathrm{S},\mathrm{S}\}} & \lambda^{\{\mathrm{F},\mathrm{S}\}} \\ \lambda^{\{\mathrm{S},\mathrm{F}\}} & \lambda^{\{\mathrm{F},\mathrm{F}\}} \end{bmatrix}}_{A} \begin{bmatrix} \mathbf{y}^{\{\mathrm{S}\}}(t) \\ \mathbf{y}^{\{\mathrm{F}\}}(t) \end{bmatrix} =: \begin{bmatrix} \mathbf{f}^{\{\mathrm{S}\}}(\mathbf{y}^{\{\mathrm{S}\}}(t), \mathbf{y}^{\{\mathrm{F}\}}(t)) \\ \mathbf{f}^{\{\mathrm{F}\}}(\mathbf{y}^{\{\mathrm{S}\}}(t), \mathbf{y}^{\{\mathrm{F}\}}(t)) \end{bmatrix}.$$

The test system (4.3) is assumed to have the following properties:

1. The system has real coefficients,

$$\text{(4.4a)}\qquad \lambda^{\{\mathrm{S},\mathrm{S}\}},\ \lambda^{\{\mathrm{F},\mathrm{F}\}},\ \lambda^{\{\mathrm{S},\mathrm{F}\}},\ \lambda^{\{\mathrm{F},\mathrm{S}\}} \in \mathbb{R}.$$

2. The system has two underlying time scales associated with the slow ($\lambda^{\{\mathrm{S},\mathrm{S}\}}$) and the fast ($\lambda^{\{\mathrm{F},\mathrm{F}\}}$) dynamics. Each subsystem is stable when run in decoupled mode, therefore the two diagonal terms are negative

$$\text{(4.4b)}\qquad \lambda^{\{\mathrm{S},\mathrm{S}\}} < 0, \qquad \lambda^{\{\mathrm{F},\mathrm{F}\}} < 0.$$

3. The dynamics is characterized by the following coefficients:

$$\text{(4.4c)}\qquad \text{scale ratio:} \qquad \mu = \frac{|\lambda^{\{\mathrm{F},\mathrm{F}\}}|}{|\lambda^{\{\mathrm{S},\mathrm{S}\}}|},$$

$$\text{(4.4d)}\qquad \text{coupling coefficient:} \qquad \mathtt{k} = \frac{\lambda^{\{\mathrm{F},\mathrm{S}\}}\, \lambda^{\{\mathrm{S},\mathrm{F}\}}}{\lambda^{\{\mathrm{F},\mathrm{F}\}}\, \lambda^{\{\mathrm{S},\mathrm{S}\}}}.$$

   To sample the fast part accurately enough one may require that $\mathtt{m} \geq \mu$. The system is weakly coupled for $|\mathtt{k}| \ll 1$.

4. The coupled system (4.3) is assumed to be stable. The coupling between these two components is represented by $\lambda^{\{\mathrm{F},\mathrm{S}\}}$ and $\lambda^{\{\mathrm{S},\mathrm{F}\}}$. The system (4.3) is stable if the real part of the eigenvalues of $A$ is negative. The characteristic polynomial of this matrix is:

$$p(z) = (z - \lambda^{\{\mathrm{S},\mathrm{S}\}})(z - \lambda^{\{\mathrm{F},\mathrm{F}\}}) - \lambda^{\{\mathrm{S},\mathrm{F}\}}\, \lambda^{\{\mathrm{F},\mathrm{S}\}} = z^2 - \mathrm{tr}(A) z + \det(A).$$

   According to the continuous Routh-Hurwitz criterion all the roots are in the left half plane iff all the coefficients of the characteristic polynomial are positive:

$$\text{(4.4ea)}\qquad \mathrm{tr}(A) < 0 \qquad \Leftrightarrow \qquad \lambda^{\{\mathrm{S},\mathrm{S}\}} + \lambda^{\{\mathrm{F},\mathrm{F}\}} < 0,$$

(4.4eb) $\det(A) > 0 \quad \Leftrightarrow \quad \lambda^{\{\mathrm{S,S}\}}\,\lambda^{\{\mathrm{F,F}\}} > \lambda^{\{\mathrm{S,F}\}}\,\lambda^{\{\mathrm{F,S}\}} \quad \Leftrightarrow \quad \mathtt{k} < 1.$

The trace relation (4.4ea) is fulfilled since the separate dynamics of the two subsystems are stable (4.4b). The determinant relation (4.4eb) can be expressed in terms of the coupling coefficient (4.4d); note that the linear test system is stable for $\mathtt{k} \to -\infty$.

We introduce the following variables:

$$z^{\{\mathrm{S,S}\}} = H\,\lambda^{\{\mathrm{S,S}\}}, \quad z^{\{\mathrm{F,F}\}} = H\,\lambda^{\{\mathrm{F,F}\}}, \quad z^{\{\mathrm{F,S}\}} = H\,\lambda^{\{\mathrm{F,S}\}}, \quad z^{\{\mathrm{S,F}\}} = H\,\lambda^{\{\mathrm{S,F}\}}.$$

From the properties of the test problem (4.4b) and (4.4e) we have that:

$$z^{\{\mathrm{S,S}\}} < 0, \qquad z^{\{\mathrm{F,F}\}} < 0, \qquad \text{and} \qquad z^{\{\mathrm{S,S}\}}\,z^{\{\mathrm{F,F}\}} > z^{\{\mathrm{S,F}\}}\,z^{\{\mathrm{F,S}\}}. \tag{4.5}$$

*Multirate forward Euler with constant interpolation.* The multirate forward Euler applied to the test problem (4.3) reads:

$$\mathbf{y}^{\{\mathrm{S}\}}_{n+1} = (1 + z^{\{\mathrm{S,S}\}})\,\mathbf{y}^{\{\mathrm{S}\}}_{n} + z^{\{\mathrm{F,S}\}}\,\mathbf{y}^{\{\mathrm{F}\}}_{n} \tag{4.6a}$$

$$\mathbf{y}^{\{\mathrm{F}\}}_{n+(\ell+1)/\mathtt{m}} = (1 + z^{\{\mathrm{F,F}\}}/\mathtt{m})\,\mathbf{y}^{\{\mathrm{F}\}}_{n+\ell/\mathtt{m}} + (z^{\{\mathrm{S,F}\}}/\mathtt{m})\,\mathbf{y}^{\{\mathrm{S}\}}_{n}, \quad \ell = 0,\dots,\mathtt{m}-1. \tag{4.6b}$$

From equation (4.6b) we have

$$\mathbf{y}^{\{\mathrm{F}\}}_{n+1} = (1 + z^{\{\mathrm{F,F}\}}/\mathtt{m})^{\mathtt{m}}\,\mathbf{y}^{\{\mathrm{F}\}}_{n} + \frac{(1 + z^{\{\mathrm{F,F}\}}/\mathtt{m})^{\mathtt{m}} - 1}{z^{\{\mathrm{F,F}\}}/\mathtt{m}}\,(z^{\{\mathrm{S,F}\}}/\mathtt{m})\,\mathbf{y}^{\{\mathrm{S}\}}_{n}.$$

Therefore

$$\begin{bmatrix} \mathbf{y}^{\{\mathrm{S}\}}_{n+1} \\ \mathbf{y}^{\{\mathrm{F}\}}_{n+1} \end{bmatrix} = \underbrace{\begin{bmatrix} \mathrm{R}^{\mathrm{S}}_{\mathrm{FE}} & z^{\{\mathrm{F,S}\}} \\ (\mathrm{R}^{\mathrm{F}}_{\mathrm{FE}} - 1)(z^{\{\mathrm{S,F}\}}/z^{\{\mathrm{F,F}\}}) & \mathrm{R}^{\mathrm{F}}_{\mathrm{FE}} \end{bmatrix}}_{\mathrm{R}^{\mathrm{CON}}_{\mathrm{MRFE}}} \cdot \begin{bmatrix} \mathbf{y}^{\{\mathrm{S}\}}_{n} \\ \mathbf{y}^{\{\mathrm{F}\}}_{n} \end{bmatrix}, \tag{4.7}$$

where we have used the stability functions $\mathrm{R}^{\mathrm{S}}_{\mathrm{FE}} = 1 + z^{\{\mathrm{S,S}\}}$ and $\mathrm{R}^{\mathrm{F}}_{\mathrm{FE}} = (1 + z^{\{\mathrm{F,F}\}}/\mathtt{m})^{\mathtt{m}}$ of the underlying forward Euler schemes.

Definition 4.1 (Linear stability). *The multirate forward Euler method is linearly stable if all the eigenvalues of the transfer matrix $\mathrm{R}^{\mathrm{CON}}_{\mathrm{MRFE}}$ have absolute values smaller than or equal one, and the eigenvalues of absolute value one are simple.*

The stability can be easily analyzed with the help of the following result.

Lemma 4.2 (Routh-Hurwitz discrete stability criterion [32, 45]). *The matrix $\mathrm{R} \in \mathbb{R}^{2\times 2}$ has a spectral radius bounded by one if and only if the following three conditions hold:*

$$a) \qquad 1 + \mathrm{tr}(\mathrm{R}) + \det(\mathrm{R}) > 0, \tag{4.8a}$$

$$b) \qquad \det(\mathrm{R}) < 1, \qquad \text{and} \tag{4.8b}$$

$$c) \qquad 1 - \mathrm{tr}(\mathrm{R}) + \det(\mathrm{R}) > 0. \tag{4.8c}$$

Using the stability functions $\mathrm{R}^{\mathrm{S}}_{\mathrm{FE}} = 1 + z^{\{\mathrm{S,S}\}}$ and $\mathrm{R}^{\mathrm{F}}_{\mathrm{FE}} = (1 + z^{\{\mathrm{F,F}\}}/\mathtt{m})^{\mathtt{m}}$ of the underlying forward Euler schemes and the coupling coefficient $\mathtt{k}$, we get

$$\mathrm{tr}(\mathrm{R}^{\mathrm{CON}}_{\mathrm{MRFE}}) = \mathrm{R}^{\mathrm{S}}_{\mathrm{FE}} + \mathrm{R}^{\mathrm{F}}_{\mathrm{FE}},$$
$$\det(\mathrm{R}^{\mathrm{CON}}_{\mathrm{MRFE}}) = \mathrm{R}^{\mathrm{S}}_{\mathrm{FE}}\mathrm{R}^{\mathrm{F}}_{\mathrm{FE}} - \mathtt{k}(1 - \mathrm{R}^{\mathrm{S}}_{\mathrm{FE}})(1 - \mathrm{R}^{\mathrm{F}}_{\mathrm{FE}}).$$

Remark 7 (Stability of multirate forward Euler scheme applied to a one-way coupled system). *Note that in case of one-way coupling ($z^{\{\mathrm{F,S}\}} = 0$ or $z^{\{\mathrm{S,F}\}} = 0 \Rightarrow \mathtt{k} = 0$) the stability of the multirate scheme is given by the stability of the forward Euler scheme applied to each component with the corresponding steps $H$ and $h$, respectively.*

Theorem 4.3 (Linear stability of multirate forward Euler schemes with constant interpolation). *Assume that the step sizes $H$ and $h$ are sufficiently small such that both base schemes, applied to the decoupled slow and fast systems, are linearly stable. Then the multirate forward Euler scheme applied to the two-dimensional system (4.3) with constant interpolation is stable for*

$$\mathtt{k} \in \left( \frac{\mathrm{R}^{\mathrm{S}}_{\mathrm{FE}}\,\mathrm{R}^{\mathrm{F}}_{\mathrm{FE}} - 1}{(1-\mathrm{R}^{\mathrm{S}}_{\mathrm{FE}})(1-\mathrm{R}^{\mathrm{F}}_{\mathrm{FE}})}, \frac{(1+\mathrm{R}^{\mathrm{S}}_{\mathrm{FE}})(1+\mathrm{R}^{\mathrm{F}}_{\mathrm{FE}})}{(1-\mathrm{R}^{\mathrm{S}}_{\mathrm{FE}})(1-\mathrm{R}^{\mathrm{F}}_{\mathrm{FE}})} \right) \cap (-\infty, 1),$$

*and unstable otherwise.*

*Proof.* One checks the three criteria (4.8) discussed in Lemma 4.2.

Remark 8. *As*

$$\frac{\mathrm{R}^{\mathrm{S}}_{\mathrm{FE}}\,\mathrm{R}^{\mathrm{F}}_{\mathrm{FE}} - 1}{(1-\mathrm{R}^{\mathrm{S}}_{\mathrm{FE}})(1-\mathrm{R}^{\mathrm{F}}_{\mathrm{FE}})} < 0 < \frac{(1+\mathrm{R}^{\mathrm{S}}_{\mathrm{FE}})(1+\mathrm{R}^{\mathrm{F}}_{\mathrm{FE}})}{(1-\mathrm{R}^{\mathrm{S}}_{\mathrm{FE}})(1-\mathrm{R}^{\mathrm{F}}_{\mathrm{FE}})}$$

*holds, stability holds in a neighborhood of $\mathtt{k} = 0$, i.e., there exists some $\epsilon > 0$ depending on $\mathrm{R}^{\mathrm{S}}_{\mathrm{FE}}$ and $\mathrm{R}^{\mathrm{F}}_{\mathrm{FE}}$ such that the multirate forward Euler schemes are stable for all $\mathtt{k} \in [-\epsilon, \epsilon]$. However, $\epsilon$ becomes arbitrarily small for $\mathrm{R}^{\mathrm{S}}_{\mathrm{FE}}, \mathrm{R}^{\mathrm{F}}_{\mathrm{FE}} \to -1$, for example.*

**4.2. Multirate backward Euler method.** As unconditional stability requires implicit schemes, we now derive multirate backward Euler methods (MRBE) to aim at unconditional stability for system (4.3) for strong coupling values $\mathtt{k} \ll -1$.

**4.2.1. Method formulation.** Consider the partitioned IVP (2.1). We apply the backward Euler method to the slow and fast components. If the classical backward Euler method is used then a system of $d = d^{\{\mathrm{S}\}} + d^{\{\mathrm{F}\}}$ nonlinear equations is solved at each time step. Assuming that a Newton-like method is employed, and that a direct linear algebra solver is used, then the cost per time step is $\mathcal{O}(d^3)$. The cost of advancing the entire system from $t_n$ to $t_{n+1}$ using $\mathtt{m}$ backward Euler steps with a small time step $h$ is therefore:

$$\text{cost of traditional backward Euler} \sim \mathcal{O}(\mathtt{m}\,(d^{\{\mathrm{S}\}} + d^{\{\mathrm{F}\}})^3). \tag{4.9}$$

In a multirate approach the slow and the fast subsystems are solved with different time steps as shown in Fig. 6. The goal is to perform the integration at a cost smaller than (4.9), without sacrificing accuracy.

There are several different possibilities to couple the two subsystems. These approaches are discussed next.

*The fully-coupled approach.* In the spirit of single-rate backward Euler method we take backward steps for both the slow and the fast components:

$$\mathbf{y}^{\{\mathrm{S}\}}_{n+1} = \mathbf{y}^{\{\mathrm{S}\}}_{n} + H\,\mathbf{f}^{\{\mathrm{S}\}}(\mathbf{y}^{\{\mathrm{S}\}}_{n+1}, \mathbf{y}^{\{\mathrm{F}\}}_{n+1}), \tag{4.10a}$$

$$\mathbf{y}^{\{\mathrm{F}\}}_{n+\ell/\mathtt{m}} = \mathbf{y}^{\{\mathrm{F}\}}_{n+(\ell-1)/\mathtt{m}} + h\,\mathbf{f}^{\{\mathrm{F}\}}(\mathrm{y}^{\{\mathrm{S}\}}_{n+\ell/\mathtt{m}}, \mathbf{y}^{\{\mathrm{F}\}}_{n+\ell+/\mathtt{m}}), \quad \ell = 1, \dots, \mathtt{m}; \tag{4.10b}$$

$$\mathrm{y}^{\{\mathrm{S}\}}_{n+\ell/\mathtt{m}} = \begin{cases} \mathbf{y}^{\{\mathrm{S}\}}_{n} \text{ or } \mathbf{y}^{\{\mathrm{S}\}}_{n+1}, & \text{constant}, \\ \frac{\mathtt{m}-\ell}{\mathtt{m}}\mathbf{y}^{\{\mathrm{S}\}}_{n} + \frac{\ell}{\mathtt{m}}\mathbf{y}^{\{\mathrm{S}\}}_{n+1}, & \text{linear}, \end{cases} \quad \ell = 1, \dots, \mathtt{m}. \tag{4.10c}$$

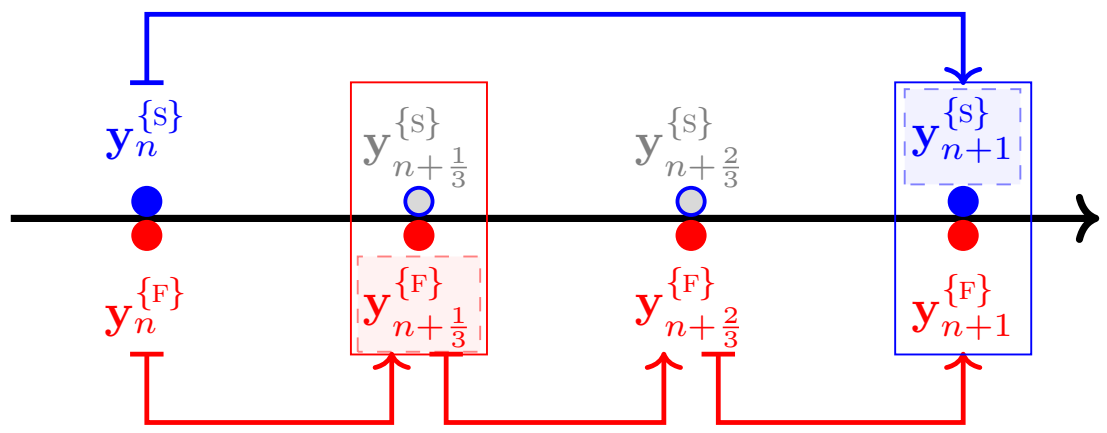

Fig. 6: Cartoon of MRBE method for $\mathtt{m} = 3$. The values computed by backward Euler steps are represented by filled discs. Empty circles represent the slow values not computed by Euler steps, and that need to be approximated.

Note that the value of the argument $\mathbf{y}_{n+1}^{\{F\}}$ is the solution after $\mathtt{m}$ small steps. The argument $\mathbf{y}_{n+(\ell+1)/\mathtt{m}}^{\{S\}}$ obtained by constant interpolation is sufficient for first order convergence.

The fully coupled approach computational process couples all the slow and the fast substeps, leading to a very large system of $d^{\{S\}} + \mathtt{m} \times d^{\{F\}}$ nonlinear equations for the slow solution $\mathbf{y}_{n+1}^{\{S\}}$ and for all intermediate fast solutions $\mathbf{y}_{n+\ell/\mathtt{m}}^{\{F\}}$, $\ell = 1, \dots, \mathtt{m}$. The total cost

$$\text{cost of fully coupled MRBE} \sim \mathcal{O}((d^{\{S\}} + \mathtt{m}\, d^{\{F\}})^3) \tag{4.11}$$

exceeds that of taking small steps with the classical backward Euler method (4.9). *Therefore the fully coupled backward multirate approach is impractical due to its very large computational cost, and will not be considered further.*

*The decoupled slowest-first approach.* In order to reduce computational costs a simple idea is to apply the backward Euler method to solve the slow and fast variables in a decoupled way. The decoupling is realized by using in the slow solution only known past values of the fast variable, and vice-versa.

We use a large step size $H$ to solve the slow components, which leads to a system of nonlinear equations in $\mathbf{y}_{n+1}^{\{S\}}$:

$$\mathbf{y}_{n+1}^{\{S\}} \;=\; \mathbf{y}_n^{\{S\}} + H\, \mathbf{f}^{\{S\}}(\mathbf{y}_{n+1}^{\{S\}}, \mathbf{y}_n^{\{F\}}). \tag{4.12a}$$

This formula is implicit in the slow variables and explicit in the fast variables. The fast solution is obtained by applying the backward Euler method with a small step size $h = H/\mathtt{m}$:

$$\mathbf{y}_{n+(\ell+1)/\mathtt{m}}^{\{F\}} = \mathbf{y}_{n+\ell/\mathtt{m}}^{\{F\}} + h\, \mathbf{f}^{\{F\}}(\mathrm{y}_{n+(\ell+1)/\mathtt{m}}^{\{S\}}, \mathbf{y}_{n+(\ell+1)/\mathtt{m}}^{\{F\}}), \quad \ell = 0, \dots, \mathtt{m}-1. \tag{4.12b}$$

The intermediate slow values $\mathrm{y}_{n+(\ell+1)/\mathtt{m}}^{\{S\}}$ are obtained from the known $\mathbf{y}_n^{\{S\}}$ and $\mathbf{y}_{n+1}^{\{S\}}$ by applying either constant or linear interpolation. Since the fast components are treated explicitly in the slow formula, the decoupled slowest-first approach raises stability concerns. The total cost of the slowest-first approach is that of one backward slow step plus $\mathtt{m}$ backward fast steps:

$$\text{cost of decoupled slowest-first MRBE} \sim \mathcal{O}(\mathtt{m}\, d^{\{F\}3} + d^{\{S\}3}). \tag{4.13}$$

*The decoupled fastest-first approach.* This approach proceeds with solving the fast variable as in the decoupled slowest-first approach, since (4.12b) depends only on the past value of the slow variable, i.e., we set $\mathbf{y}^{\{\mathrm{S}\}}_{n+(\ell+1)/\mathtt{m}} := \mathbf{y}^{\{\mathrm{S}\}}_{n}$:

$$\mathbf{y}^{\{\mathrm{F}\}}_{n+(\ell+1)/\mathtt{m}} = \mathbf{y}^{\{\mathrm{F}\}}_{n+\ell/\mathtt{m}} + h\,\mathbf{f}^{\{\mathrm{F}\}}(\mathbf{y}^{\{\mathrm{S}\}}_{n}, \mathbf{y}^{\{\mathrm{F}\}}_{n+(\ell+1)/\mathtt{m}}), \quad \ell = 0, \ldots, \mathtt{m}-1. \tag{4.14a}$$

The intermediate values of the slow variable are computed by constant interpolation (4.10c) with the value at $t_n$. The slow variable is treated explicitly, while the fast variable is treated implicitly in (4.14a).

The slow variable is then computed using:

$$\mathbf{y}^{\{\mathrm{S}\}}_{n+1} = \mathbf{y}^{\{\mathrm{S}\}}_{n} + H\,\mathbf{f}^{\{\mathrm{S}\}}(\mathbf{y}^{\{\mathrm{S}\}}_{n+1}, \mathbf{y}^{\{\mathrm{F}\}}_{n+1}). \tag{4.14b}$$

The fast variable value is the one obtained from the integration (4.14a). The total cost of the fastest-first approach is that of $\mathtt{m}$ backward fast steps (4.14a) plus one backward slow step (4.14b):

$$\text{cost of decoupled fastest-first MRBE} \sim \mathcal{O}(\mathtt{m}\, d^{\{\mathrm{F}\}3} + d^{\{\mathrm{S}\}3}). \tag{4.15}$$

*The coupled slowest-first approach.* In the coupled slowest–first approach both components are solved together:

$$\begin{aligned} \mathbf{y}^{\{\mathrm{S}\}}_{n+1} &= \mathbf{y}^{\{\mathrm{S}\}}_{n} + H\,\mathbf{f}^{\{\mathrm{S}\}}(\mathbf{y}^{\{\mathrm{S}\}}_{n+1}, \hat{\mathbf{y}}^{\{\mathrm{F}\}}_{n+1}), \\ \hat{\mathbf{y}}^{\{\mathrm{F}\}}_{n+1} &= \mathbf{y}^{\{\mathrm{F}\}}_{n} + H\,\mathbf{f}^{\{\mathrm{F}\}}(\mathbf{y}^{\{\mathrm{S}\}}_{n+1}, \hat{\mathbf{y}}^{\{\mathrm{F}\}}_{n+1}). \end{aligned} \tag{4.16a}$$

The predicted fast component $\hat{\mathbf{y}}^{\{\mathrm{F}\}}_{n+1}$ is inaccurate for large $H$ and is discarded. A corrected fast solution is obtained by applying the backward Euler method with a small step size $h = H/\mathtt{m}$

$$\mathbf{y}^{\{\mathrm{F}\}}_{n+(\ell+1)/\mathtt{m}} = \mathbf{y}^{\{\mathrm{F}\}}_{n+\ell/\mathtt{m}} + h\,\mathbf{f}^{\{\mathrm{F}\}}(\mathbf{y}^{\{\mathrm{S}\}}_{n+(\ell+1)/\mathtt{m}}, \mathbf{y}^{\{\mathrm{F}\}}_{n+(\ell+1)/\mathtt{m}}), \quad \ell = 0, \ldots, \mathtt{m}-1. \tag{4.16b}$$

The intermediate slow variables $\mathbf{y}^{\{\mathrm{S}\}}_{n+(\ell+1)/\mathtt{m}}$ for $\ell = 1, \ldots, \mathtt{m}-1$ can be approximated by any interpolation formula in (4.10c).

The cost is given by one large Euler step with the full system (4.16a) plus $\mathtt{m}$ small steps with the fast subsystem (4.16b):

$$\text{cost of coupled slowest-first MRBE} \sim \mathcal{O}((d^{\{\mathrm{F}\}} + d^{\{\mathrm{S}\}})^3 + \mathtt{m}\, d^{\{\mathrm{F}\}3}). \tag{4.17}$$

This idea to recompute components, which are not accurate enough when approximated with a large step size, by using smaller step sizes while accepting those components which are accurate enough has been introduced in [58] for stiff systems within a Rosenbrock-Wanner aproach.

*The coupled-first-step approach.* In order to avoid computing and discarding a fast solution in (4.16a) we can couple the slow backward Euler step with the first fast backward Euler step, and use zeroth order interpolation in both formulas:

$$\begin{aligned} \mathbf{y}^{\{\mathrm{S}\}}_{n+1} &= \mathbf{y}^{\{\mathrm{S}\}}_{n} + H\,\mathbf{f}^{\{\mathrm{S}\}}(\mathbf{y}^{\{\mathrm{S}\}}_{n+1}, \mathbf{y}^{\{\mathrm{F}\}}_{n+1/\mathtt{m}}), \\ \mathbf{y}^{\{\mathrm{F}\}}_{n+1/\mathtt{m}} &= \mathbf{y}^{\{\mathrm{F}\}}_{n} + h\,\mathbf{f}^{\{\mathrm{F}\}}(\mathbf{y}^{\{\mathrm{S}\}}_{n+1}, \mathbf{y}^{\{\mathrm{F}\}}_{n+1/\mathtt{m}}). \end{aligned} \tag{4.18a}$$

The remaining fast steps are then carried out in a decoupled manner:

$$\mathbf{y}^{\{\mathrm{F}\}}_{n+(\ell+1)/\mathtt{m}} = \mathbf{y}^{\{\mathrm{F}\}}_{n+\ell/\mathtt{m}} + h\, \mathbf{f}^{\{\mathrm{F}\}}(\mathbf{y}^{\{\mathrm{S}\}}_{n+1}, \mathbf{y}^{\{\mathrm{F}\}}_{n+(\ell+1)/\mathtt{m}}), \quad \ell = 1, \ldots, \mathtt{m}-1. \tag{4.18b}$$

The cost is given by:

$$\text{cost of coupled-first-step MRBE} \sim \mathcal{O}((d^{\{\mathrm{F}\}} + d^{\{\mathrm{S}\}})^3 + (\mathtt{m}-1)\, d^{\{\mathrm{F}\}3}). \tag{4.19}$$

The coupled-first step approach has been introduced in [37] for stiff systems within a Runge-Kutta approach.

As for the forward case, all instances of multirate backward Euler schemes have convergence order one provided that at least constant extrapolation is used (see Section 3).

**4.2.2. Linear stability analysis of multirate backward Euler methods.** As linear interpolation will not increase the convergence order of an backward multirate Euler scheme, we will only consider the case of constant interpolation in the following analysis.

The linear stability analysis follows the one developed in Section 4.1.2. To this end we apply the implicit schemes to solve the linear test problem (4.3). This gives an iteration of the form:

$$\begin{bmatrix} \mathbf{y}^{\{\mathrm{S}\}}_{n+1} \\ \mathbf{y}^{\{\mathrm{F}\}}_{n+1} \end{bmatrix} = \mathrm{R}_{\mathrm{MRBE}} \cdot \begin{bmatrix} \mathbf{y}^{\{\mathrm{S}\}}_{n} \\ \mathbf{y}^{\{\mathrm{F}\}}_{n} \end{bmatrix}, \quad \mathrm{R}_{\mathrm{MRBE}} \in \mathbb{R}^{2\times 2} \tag{4.20}$$

with $\mathrm{R}_{\mathrm{MRBE}}$ depending on the respective approach used. The multirate backward Euler method is linearly stable if both eigenvalues of the matrix $\mathrm{R}^{\mathrm{DFFC}}_{\mathrm{MRBE}}$ have absolute values smaller than or equal to one.

Definition 4.4 (Unconditional stability). *A multirate method is unconditionally stable if it is stable for any step sizes $H > 0$ and $h > 0$.*

Remark 9 (Linear stability for a decoupled system). *When the backward Euler method is applied to a decoupled test problem (4.3) where $z^{\{\mathrm{F},\mathrm{S}\}} = z^{\{\mathrm{S},\mathrm{F}\}} = 0$ one takes one step with the slow system and $\mathtt{m}$ steps with the fast system, and (4.6) becomes:*

$$\mathbf{y}^{\{\mathrm{S}\}}_{n+1} = \mathrm{R}^{\mathrm{S}}_{\mathrm{BE}} \cdot \mathbf{y}^{\{\mathrm{S}\}}_{n}, \quad \mathbf{y}^{\{\mathrm{F}\}}_{n+1} = \mathrm{R}^{\mathrm{F}}_{\mathrm{BE}} \cdot \mathbf{y}^{\{\mathrm{F}\}}_{n},$$

*where the slow and fast stability functions of the backward Euler method over a macro-step are defined as:*

$$\mathrm{R}^{\mathrm{S}}_{\mathrm{BE}} := (1 - z^{\{\mathrm{S},\mathrm{S}\}})^{-1} \in (0,1] \quad \text{and} \quad \mathrm{R}^{\mathrm{F}}_{\mathrm{BE}} := \left(1 - \frac{z^{\{\mathrm{F},\mathrm{F}\}}}{\mathtt{m}}\right)^{-\mathtt{m}} \in (0,1], \tag{4.21}$$

*respectively. The decoupled schemes are stable for any $z^{\{\mathrm{F},\mathrm{F}\}}, z^{\{\mathrm{S},\mathrm{S}\}} < 0$.*

*The decoupled slowest-first and fastest-first approach.* One can show that both the decoupled slowest-first and fastest-first approach with constant interpolation (4.10c) is unconditionally stable if the system is weakly coupled, $-1 \le \mathtt{k} < 1$; it becomes unstable for $\mathtt{k} \to -\infty$.

For unconditional stability we need some coupling of active and latent parts, i.e., the slow part has to be computed together with some part of the fast one. We will show this for

*The coupled slowest-first approach.* The macro-step (4.16) reads

$$\begin{bmatrix} 1 - z^{\{\text{S},\text{S}\}} & -z^{\{\text{F},\text{S}\}} \\ -z^{\{\text{S},\text{F}\}} & 1 - z^{\{\text{F},\text{F}\}} \end{bmatrix} \cdot \begin{bmatrix} \mathbf{y}_{n+1}^{\{\text{S}\}} \\ \hat{\mathbf{y}}_{n+1}^{\{\text{F}\}} \end{bmatrix} = \begin{bmatrix} \mathbf{y}_{n}^{\{\text{S}\}} \\ \mathbf{y}_{n}^{\{\text{F}\}} \end{bmatrix}, \tag{4.22}$$

which yields

$$\mathbf{y}_{n+1}^{\{\text{S}\}} = \frac{(1 - z^{\{\text{F},\text{F}\}})\mathbf{y}_{n}^{\{\text{S}\}} + z^{\{\text{F},\text{S}\}}\mathbf{y}_{n}^{\{\text{F}\}}}{(1 - z^{\{\text{S},\text{S}\}})(1 - z^{\{\text{F},\text{F}\}}) - z^{\{\text{S},\text{F}\}} z^{\{\text{F},\text{S}\}}}.$$

The micro-step solution reads:

$$\mathbf{y}_{n+1}^{\{\text{F}\}} = (\mathrm{R}_{\text{BE}}^{\text{F}} - 1)\frac{z^{\{\text{S},\text{F}\}}}{z^{\{\text{F},\text{F}\}}}\mathbf{y}_{n+1}^{\{\text{S}\}} + \mathrm{R}_{\text{BE}}^{\text{F}}\mathbf{y}_{n}^{\{\text{F}\}}$$

The overall solution can be written as:

$$\begin{bmatrix} \mathbf{y}_{n+1}^{\{\text{S}\}} \\ \mathbf{y}_{n+1}^{\{\text{F}\}} \end{bmatrix} = \mathrm{R}_{\text{MRBE}}^{\text{CSFC}} \cdot \begin{bmatrix} \mathbf{y}_{n}^{\{\text{S}\}} \\ \mathbf{y}_{n}^{\{\text{F}\}} \end{bmatrix},$$

with

$$\det(\mathrm{R}_{\text{MRBE}}^{\text{CSFC}}) = \mathrm{R}_{\text{BE}}^{\text{F}} \frac{1 - z^{\{\text{F},\text{F}\}}}{1 - z^{\{\text{S},\text{S}\}} - z^{\{\text{F},\text{F}\}} + (1 - \mathtt{k}) z^{\{\text{S},\text{S}\}} z^{\{\text{F},\text{F}\}}}, \tag{4.23a}$$

$$\mathrm{tr}(\mathrm{R}_{\text{MRBE}}^{\text{CSFC}}) = \mathrm{R}_{\text{BE}}^{\text{F}} + \frac{1 - z^{\{\text{F},\text{F}\}} + (\mathrm{R}_{\text{BE}}^{\text{F}} - 1) z^{\{\text{F},\text{S}\}} z^{\{\text{S},\text{F}\}} / z^{\{\text{F},\text{F}\}}}{1 - z^{\{\text{S},\text{S}\}} - z^{\{\text{F},\text{F}\}} + (1 - \mathtt{k}) z^{\{\text{S},\text{S}\}} z^{\{\text{F},\text{F}\}}}. \tag{4.23b}$$

THEOREM 4.5 (Stability of coupled slowest-first approach). *The coupled slowest–first multirate backward Euler method* (4.16) *with constant interpolation of the slow variable at $t_{n+1}$, i.e., $\mathbf{y}_{n+(\ell+1)/\mathtt{m}}^{\{\text{S}\}} = \mathbf{y}_{n+1}^{\{\text{S}\}}$ in* (4.16b)*, is unconditionally stable.*

*Proof.* One verifies the three Routh-Hurwitz criteria (4.8) using the method's stability matrix (4.23).

REMARK 10. *Following the lines of arguments above, one can show that the coupled-first step approach* (4.18) *with constant interpolation of the slow variable at $t_{n+1}$ is also unconditionally stable.*

**5. Multirate Runge–Kutta methods.** Classical Runge-Kutta methods for solving initial value problems (1.1) have a long history [7,35,46] and are now pervasive in computations. They achieve high order by computing, within each step, several internal stage approximations of the solution [8,30,31].

**5.1. Generalized-structure additive Runge-Kutta schemes.** Generalized-structure additive Runge-Kutta (GARK) methods were introduced in [52] to solve IVPs for *additively* partitioned systems of ODEs (2.2). For systems (2.1) split into slow (S) and fast (F) components, these schemes advance the slow/fast numerical solutions as follows:

$$\mathbf{Y}_i^{\{\sigma,\mu\}} = \mathbf{y}_n^{\{\sigma\}} + H \sum_{j=1}^{s^{\{\mu\}}} a_{i,j}^{\{\sigma,\mu\}} \mathbf{f}^{\{\mu\}}\left(\mathbf{Y}_j^{\{\mu,\text{S}\}}, \mathbf{Y}_j^{\{\mu,\text{F}\}}\right), \quad i = 1, \ldots, s^{\{\sigma\}}, \quad \forall \sigma, \mu \in \{\text{S}, \text{F}\}, \tag{5.1a}$$

$$\text{(5.1b)} \qquad \mathbf{y}_{n+1}^{\{\sigma\}} = \mathbf{y}_n^{\{\sigma\}} + H \sum_{i=1}^{s^{\{\mu\}}} b_i^{\{\sigma\}} \, \mathbf{f}^{\{\sigma\}} \left( \mathbf{Y}_i^{\{\sigma,\text{S}\}}, \mathbf{Y}_i^{\{\sigma,\text{F}\}} \right), \quad \forall \sigma \in \{\text{S},\text{F}\},$$

and are represented compactly by the extended Butcher tableau:

$$\text{(5.2)} \qquad \begin{array}{c} \mathbf{A}_{\text{GARK}} \\ \hline \mathbf{b}_{\text{GARK}}^{T} \end{array} = \begin{array}{c|c} \mathbf{A}^{\{\text{F},\text{F}\}} & \mathbf{A}^{\{\text{F},\text{S}\}} \\ \hline \mathbf{A}^{\{\text{S},\text{F}\}} & \mathbf{A}^{\{\text{S},\text{S}\}} \\ \hline \mathbf{b}^{\{\text{F}\}T} & \mathbf{b}^{\{\text{S}\}T} \end{array}.$$

The internal stages (5.1a) approximate solution components at internal time points, $\mathbf{Y}_i^{\{\sigma,\mu\}} = \mathbf{y}^{\{\sigma\}}(t_n + c_i^{\{\sigma,\mu\}} H) + \mathcal{O}\big(H^2\big)$ with $c_i^{\{\sigma,\mu\}} = \sum_j a_{i,j}^{\{\sigma,\mu\}}$. A GARK scheme (5.1) is called internally consistent if both the slow and fast arguments of each function call $\mathbf{f}^{\{\sigma\}}(\mathbf{Y}_j^{\{\sigma,\text{S}\}}, \mathbf{Y}_j^{\{\sigma,\text{F}\}})$ approximate the corresponding solution components at the same time (i.e., $c_i^{\{\sigma,\text{S}\}} = c_i^{\{\sigma,\text{F}\}} = c_i^{\{\sigma\}}$). This requirement is equivalent to the following algebraic condition:

$$\text{(5.3)} \qquad \mathbf{A}^{\{\sigma,\mu\}} \, \mathbf{1}^{\{\mu\}} =: \mathbf{c}^{\{\sigma,\mu\}} = \mathbf{c}^{\{\sigma\}}, \quad \forall \sigma, \mu \in \{\text{S},\text{F}\}.$$

**5.2. Multirate GARK schemes.** We construct multirate Runge-Kutta schemes using the GARK formalism [27,43,53]. The multirate GARK (MR-GARK) family includes many well-known multirate schemes as special cases. The order conditions theory follows directly from the GARK accuracy theory.

The MR-GARK method can be written as a GARK scheme (5.1) over the macro-step $H$ with the fast stage vectors $\mathbf{Y}^{\{\text{F},\text{F}\}} := [\mathbf{Y}^{\{\text{F},\text{F},1\}T}, \dots, \mathbf{Y}^{\{\text{F},\text{F},\mathtt{m}\}T}]^T$ and $\mathbf{Y}^{\{\text{F},\text{S}\}} := [\mathbf{Y}^{\{\text{F},\text{S},1\}T}, \dots, \mathbf{Y}^{\{\text{F},\text{S},\mathtt{m}\}T}]^T$ applied to the two-way component partitioned systems of ordinary differential equations (2.1). The corresponding Butcher tableau (5.2) reads:

$$\text{(5.4)} \qquad \begin{array}{c|c} \mathbf{A}^{\{\text{F},\text{F}\}} & \mathbf{A}^{\{\text{F},\text{S}\}} \\ \hline \mathbf{A}^{\{\text{S},\text{F}\}} & \mathbf{A}^{\{\text{S},\text{S}\}} \\ \hline \mathbf{b}^{\{\text{F}\}T} & \mathbf{b}^{\{\text{S}\}T} \end{array} := \begin{array}{ccc|c} \mathtt{m}^{-1} A^{\{\text{F},\text{F}\}} & \cdots & 0 & A^{\{\text{F},\text{S},1\}} \\ \vdots & \ddots & & \vdots \\ \mathtt{m}^{-1} \mathbf{1}^{\{\text{F}\}} \, b^{\{\text{F}\}T} & \cdots & \mathtt{m}^{-1} A^{\{\text{F},\text{F}\}} & A^{\{\text{F},\text{S},\mathtt{m}\}} \\ \hline \mathtt{m}^{-1} A^{\{\text{S},\text{F},1\}} & \cdots & \mathtt{m}^{-1} A^{\{\text{S},\text{F},\mathtt{m}\}} & A^{\{\text{S},\text{S}\}} \\ \hline \mathtt{m}^{-1} b^{\{\text{F}\}T} & \cdots & \mathtt{m}^{-1} b^{\{\text{F}\}T} & b^{\{\text{S}\}T} \end{array}$$

This method discretizes the slow component with a "base slow Runge-Kutta method" $\text{RK}^{\{\text{S}\}} \equiv (A^{\{\text{S},\text{S}\}}, b^{\{\text{S}\}})$. The slow variable is advanced from $t_n$ to $t_{n+1} = t_n + H$ using a single macro-step of size $H$:

$$
\begin{aligned}
&\underbrace{\mathbf{Y}_i^{\{\mathrm{S},\mathrm{S}\}} = \mathbf{y}_n^{\{\mathrm{S}\}} + H \sum_{j=1}^{s^{\{\mathrm{S}\}}} a_{i,j}^{\{\mathrm{S},\mathrm{S}\}}\, \mathbf{f}^{\{\mathrm{S}\}}\big(\mathbf{Y}_j^{\{\mathrm{S},\mathrm{S}\}}, \mathbf{Y}_j^{\{\mathrm{S},\mathrm{F}\}}\big)}_{\text{stage } i \text{ of RK}^{\{\mathrm{S}\}} \text{ solving slow component over } [t_n, t_n+H]}, \quad i = 1,\dots,s^{\{\mathrm{S}\}},\\
&\underbrace{\mathbf{Y}_i^{\{\mathrm{S},\mathrm{F}\}} = \mathbf{y}_n^{\{\mathrm{F}\}} + h \sum_{\ell=1}^{\mathtt{m}} \sum_{j=1}^{s^{\{\mathrm{F}\}}} a_{i,j}^{\{\mathrm{S},\mathrm{F},\ell\}}\, \mathbf{f}^{\{\mathrm{F}\}}\big(\mathbf{Y}_j^{\{\mathrm{F},\mathrm{S},\ell\}}, \mathbf{Y}_j^{\{\mathrm{F},\mathrm{F},\ell\}}\big)}_{\text{fast component argument of slow function value at RK}^{\{\mathrm{S}\}} \text{ stage } i}, \quad i = 1,\dots,s^{\{\mathrm{S}\}},\\
&\underbrace{\mathbf{y}_{n+1}^{\{\mathrm{S}\}} = \mathbf{y}_n^{\{\mathrm{S}\}} + H \sum_{i=1}^{s^{\{\mathrm{S}\}}} b_i^{\{\mathrm{S}\}}\, \mathbf{f}^{\{\mathrm{S}\}}\big(\mathbf{Y}_i^{\{\mathrm{S},\mathrm{S}\}}, \mathbf{Y}_i^{\{\mathrm{S},\mathrm{F}\}}\big)}_{\text{slow solution of RK}^{\{\mathrm{S}\}} \text{ over } [t_n, t_n+H]}.
\end{aligned}
\tag{5.5a}
$$

Remark 11 (Slow-fast coupling). *The slow component computation is illustrated by the blue arrows in Figure 7. At each stage RK*$^{\{\mathrm{S}\}}$ *evaluates the slow function* $\mathbf{f}^{\{\mathrm{S}\}}(\mathbf{Y}_i^{\{\mathrm{S},\mathrm{S}\}}, \mathbf{Y}_i^{\{\mathrm{S},\mathrm{F}\}})$*, where the slow component arguments* $\mathbf{Y}_i^{\{\mathrm{S},\mathrm{S}\}}$ *are the stage values computed by the base slow method. The fast component arguments* $\mathbf{Y}_i^{\{\mathrm{S},\mathrm{F}\}}$ *are obtained in a Runge-Kutta fashion by linear combinations of fast function values evaluated at all fast stages* $j$ *of each fast micro-step* $\ell$*. The GARK coefficients* $a_{i,j}^{\{\mathrm{S},\mathrm{F},\ell\}}$ *perform the coupling of fast dynamics in micro-step* $\ell$ *into the slow component solution.*

The MR-GARK scheme discretizes the fast component of (2.1) with a "base fast Runge-Kutta method" $\mathrm{RK}^{\{\mathrm{F}\}} \equiv \big(A^{\{\mathrm{F},\mathrm{F}\}}, b^{\{\mathrm{F}\}}\big)$ . The fast component is advanced from $t_n$ to $t_{n+1} = t_n + H$ using $\mathtt{m}$ fast micro-steps of size $h = H/\mathtt{m}$. The micro-steps are initialized with $\mathbf{y}_n^{\{\mathrm{F}\}}$ at $t_n$ and proceed as follows:

$$
\begin{aligned}
&\mathbf{y}_{n+0}^{\{\mathrm{F}\}} = \mathbf{y}_n^{\{\mathrm{F}\}}; \quad /\!/ \text{ micro-step initialization}\\
&\text{For} \quad \ell = 0,\dots,\mathtt{m}-1:\\
&\left\{
\begin{aligned}
&\text{For} \quad i = 1,\dots,s^{\{\mathrm{F}\}}:\\
&\underbrace{\mathbf{Y}_i^{\{\mathrm{F},\mathrm{S},\ell+1\}} = \mathbf{y}_n^{\{\mathrm{S}\}} + H \sum_{j=1}^{s^{\{\mathrm{S}\}}} a_{i,j}^{\{\mathrm{F},\mathrm{S},\ell+1\}}\, \mathbf{f}^{\{\mathrm{S}\}}\big(\mathbf{Y}_j^{\{\mathrm{S},\mathrm{S}\}}, \mathbf{Y}_j^{\{\mathrm{S},\mathrm{F}\}}\big),}_{\substack{\text{slow component argument of fast function value}\\ \text{at RK}^{\{\mathrm{F}\}} \text{ stage } i \text{ over } [t_n + \ell/\mathtt{m}H,\, t_n + (\ell+1)/\mathtt{m}H]}}\\
&\underbrace{\mathbf{Y}_i^{\{\mathrm{F},\mathrm{F},\ell+1\}} = \mathbf{y}_{n+\ell/\mathtt{m}}^{\{\mathrm{F}\}} + h \sum_{j=1}^{s^{\{\mathrm{F}\}}} a_{i,j}^{\{\mathrm{F},\mathrm{F}\}}\, \mathbf{f}^{\{\mathrm{F}\}}\big(\mathbf{Y}_j^{\{\mathrm{F},\mathrm{S},\ell+1\}}, \mathbf{Y}_j^{\{\mathrm{F},\mathrm{F},\ell+1\}}\big),}_{\text{stage } i \text{ of RK}^{\{\mathrm{F}\}} \text{ solving fast component over } [t_n+\ell/\mathtt{m}H,\, t_n+(\ell+1)/\mathtt{m}H]}\\
&\underbrace{\mathbf{y}_{n+(\ell+1)/\mathtt{m}}^{\{\mathrm{F}\}} = \mathbf{y}_{n+\ell/\mathtt{m}}^{\{\mathrm{F}\}} + h \sum_{i=1}^{s^{\{\mathrm{F}\}}} b_j^{\{\mathrm{F}\}}\, \mathbf{f}^{\{\mathrm{F}\}}\big(\mathbf{Y}_j^{\{\mathrm{F},\mathrm{S},\ell+1\}}, \mathbf{Y}_j^{\{\mathrm{F},\mathrm{F},\ell+1\}}\big);}_{\text{solution of RK}^{\{\mathrm{F}\}} \text{ over } [t_n+\ell/\mathtt{m}H,\, t_n+(\ell+1)/\mathtt{m}H]}
\end{aligned}
\right.\\
&\mathbf{y}_{n+1}^{\{\mathrm{F}\}} = \mathbf{y}_{n+\mathtt{m}/\mathtt{m}}^{\{\mathrm{F}\}}. \quad /\!/ \text{ fast component solution at the end of micro-steps}
\end{aligned}
\tag{5.5b}
$$

Remark 12 (Fast-slow coupling). *The fast component computation is illustrated by the red arrows in Figure 7. At stage $i$ in micro-step $\ell$ RK$^{\{\mathrm{F}\}}$ computes the fast stage solution value $\mathbf{Y}_i^{\{\mathrm{F},\mathrm{F},\ell\}}$. To evaluate the fast function $\mathbf{f}^{\{\mathrm{F}\}}(\mathbf{Y}_i^{\{\mathrm{F},\mathrm{S},\ell\}},\mathbf{Y}_i^{\{\mathrm{F},\mathrm{F},\ell\}})$, the corresponding slow component arguments $\mathbf{Y}_i^{\{\mathrm{F},\mathrm{S},\ell\}}$ are also needed. They are obtained in a Runge-Kutta fashion by linear combinations of slow function values evaluated at all slow stages $j$. The GARK coefficients $a_{i,j}^{\{\mathrm{F},\mathrm{S},\ell\}}$ perform the coupling of slow dynamics into the fast component solution.*

Remark 13 (Number of function evaluations). *One step of the scheme (5.5) requires $s^{\{\mathrm{S}\}}$ slow function evaluations plus $\mathtt{m}s^{\{\mathrm{F}\}}$ fast function evaluations. We note that the coupling stages $\mathbf{Y}_i^{\{\mathrm{F},\mathrm{S},\ell+1\}}$ in (5.5b) use the same $s^{\{\mathrm{S}\}}$ slow function values $\mathbf{f}^{\{\mathrm{S}\}}(\mathbf{Y}_j^{\{\mathrm{S},\mathrm{S}\}},\mathbf{Y}_j^{\{\mathrm{S},\mathrm{F}\}})$ as the slow stages $\mathbf{Y}_i^{\{\mathrm{S},\mathrm{S}\}}$ in (5.5a). Similarly, the coupling stages $\mathbf{Y}_i^{\{\mathrm{S},\mathrm{F}\}}$ in (5.5a) use the same $\mathtt{m}s^{\{\mathrm{F}\}}$ fast function values $\mathbf{f}^{\{\mathrm{F}\}}(\mathbf{Y}_j^{\{\mathrm{F},\mathrm{S},\ell\}},\mathbf{Y}_j^{\{\mathrm{F},\mathrm{F},\ell\}})$ as the fast stages $\mathbf{Y}_i^{\{\mathrm{F},\mathrm{F},\ell\}}$ in (5.5b).*

Remark 14 (Variable micro-steps). *One can integrate the fast component from $t_n$ to $t_{n+1}$ using $\mathtt{m}$ different time steps $h_\ell = \sigma_\ell H$, with $\sum_{\ell=1}^{\mathtt{m}} \sigma_\ell = 1$, and a different Runge-Kutta method $(A^{\{\mathrm{F},\mathrm{F},\ell\}}, b^{\{\mathrm{F},\ell\}})$ on each fast sub-interval $\ell$. The entire* MR-GARK *discussion can be immediately extended to this general case.*

Definition 5.1 (Telescopic MR-GARK). *Of particular interest are methods (5.5) which use the same base scheme for both the slow and the fast components,*

$$A^{\{\mathrm{F},\mathrm{F}\}} = A^{\{\mathrm{S},\mathrm{S}\}} = A, \quad b^{\{\mathrm{F}\}} = b^{\{\mathrm{S}\}} = b. \tag{5.6}$$

*Such methods can be easily extended to systems with more than two scales by applying them in a telescopic fashion.*

Remark 15 (Internal consistency). *Let $\mathbf{1}^{\{\mathrm{F}\}} \in \mathbb{R}^{s^{\{\mathrm{F}\}}}$ and $\mathbf{1}^{\{\mathrm{S}\}} \in \mathbb{R}^{s^{\{\mathrm{S}\}}}$ be vectors of ones. From (5.3), a* MR-GARK *scheme is internally consistent iff:*

$$\mathtt{m}\, c^{\{\mathrm{F},\mathrm{S},\ell\}} = c^{\{\mathrm{F},\mathrm{F}\}} + (\ell-1)\mathbf{1}^{\{\mathrm{F}\}}, \quad \ell = 1,\dots,\mathtt{m}, \quad \text{where} \quad c^{\{\mathrm{F},\mathrm{S},\ell\}} := A^{\{\mathrm{F},\mathrm{S},\ell\}}\,\mathbf{1}^{\{\mathrm{S}\}}; \tag{5.7a}$$

$$\sum_{\ell=1}^{\mathtt{m}} c^{\{\mathrm{S},\mathrm{F},\ell\}} = \mathtt{m}\, c^{\{\mathrm{S},\mathrm{S}\}}, \quad \text{where} \quad c^{\{\mathrm{S},\mathrm{F},\ell\}} := A^{\{\mathrm{S},\mathrm{F},\ell\}}\,\mathbf{1}^{\{\mathrm{F}\}}. \tag{5.7b}$$

### 5.3. Coupled and decoupled MR-GARK methods.

The general formulation of the method (5.5) leads to a coupled system of nonlinear equations where all slow and fast stages in a macro-step are solved simultaneously. Consequently, the resulting computational effort is larger, not smaller, than solving the coupled system with $\mathtt{m}$ small steps. For constructing practical MR-GARK method one seeks to build a computational process where slow stages and fast stages are evaluated separately from each other. This leads to the following definition.

Definition 5.2 (Decoupled MR-GARK methods). *An* MR-GARK *method is decoupled if the computation of its stages proceeds in sequence, such that each slow stage uses only information from other slow stages and the already computed fast stages, and vice-versa. There is no coupling that requires fast and slow stages to be solved together. Any form of implicitness is entirely within the fast or within the slow system.*

Remark 16 (Graph theory interpretation). *Consider a computation dependency graph where vertices are stages, and directed edges point to a specific stage from each stage required during its computation. If one thinks of $\mathbf{A}_{\mathrm{GARK}}$ as the adjacency matrix of the computation dependency graph [42], then decoupled means that any strongly*

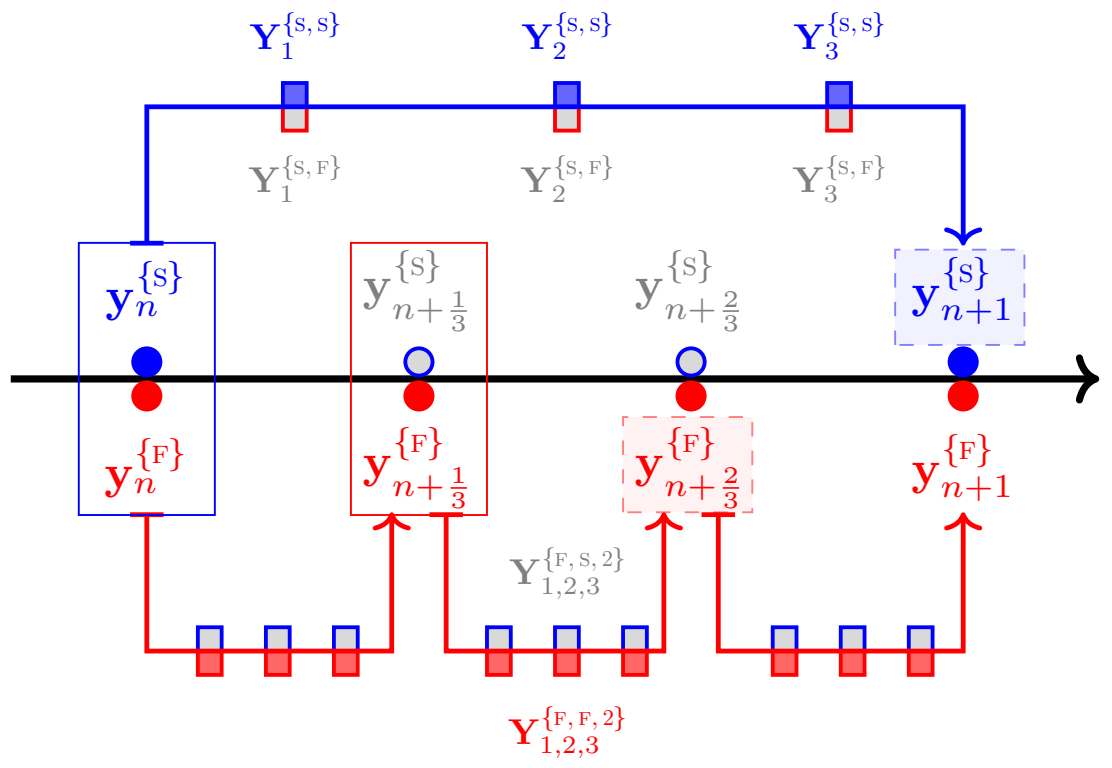

Fig. 7: Cartoon of explicit MR-GARK for $s = 3$ stages and $\mathtt{m} = 3$. The values computed by full steps are represented by filled discs. Empty circles represent the slow values not computed by full steps, and that need to be approximated. Stage values computed by the method are represented by filled squares, and stage values not computed and in need of approximation are represented by empty squares. GARK methods compute along the way both the slow and fast stage values used as arguments of the slow function, $\mathbf{f}^{\{S\}}(\mathbf{Y}_i^{\{S,S\}}, \mathbf{Y}_i^{\{S,F\}})$, $i = 1, \dots, s$. Similarly, GARK methods compute along the way both the slow and fast stage values used as arguments of the fast functions at each micro-step $\ell$, $\mathbf{f}^{\{F\}}(\mathbf{Y}_i^{\{F,F,\ell\}}, \mathbf{Y}_i^{\{F,S,\ell\}})$, $i = 1, \dots, s$.

*connected component has all vertices that are fast or all that are slow. A strongly connected component containing both types of vertices indicates a coupled multirate method.*

The need to solve for fast and slow stages together, or the possibility to avoid coupled fast-slow computations, are determined by the structure of the coupling matrices $\mathbf{A}^{\{S,F\}}$ and $\mathbf{A}^{\{F,S\}}$.

- For *decoupled* MR-GARK *methods* the sparsity structures of the coupling matrices $\mathbf{A}^{\{F,S\}}$ and $\mathbf{A}^{\{S,F\}}$ are complementary, in the sense that one must have zero entries where the other matrix can have nonzero elements:

$$\mathbf{A}^{\{S,F\}} \times \mathbf{A}^{\{F,S\}T} = \mathbf{0}_{s^{\{S\}} \times \mathtt{m}s^{\{F\}}}. \tag{5.8}$$

  This is schematically illustrated in Figure 8a. In the graph theory interpretation the directed, bipartite graph $\{\{0, \mathbf{A}^{\{F,S\}}\}, \{\mathbf{A}^{\{S,F\}}, 0\}\}$ must be acyclic [42].
- For *coupled* MR-GARK *methods* the sparsity structures can overlap, in the sense that they both can have non-zeros entries in the same location:

$$\mathbf{A}^{\{S,F\}} \times \mathbf{A}^{\{F,S\}T} \neq \mathbf{0}_{s^{\{S\}} \times \mathtt{m}s^{\{F\}}}. \tag{5.9}$$

  This is illustrated in Figure 8c. The overlapping non-zero coupling coefficients, indicated by dashed boxes in Figure 8c, imply that the corresponding fast and slow stages need to be computed together, in a step that involves the full non-partitioned system.

In the standard form of the GARK Butcher tableau (5.4) the first $\mathtt{m}s^{\{F\}}$ rows/columns correspond to the fast stages (5.5b), and the last $s^{\{S\}}$ rows/columns correspond

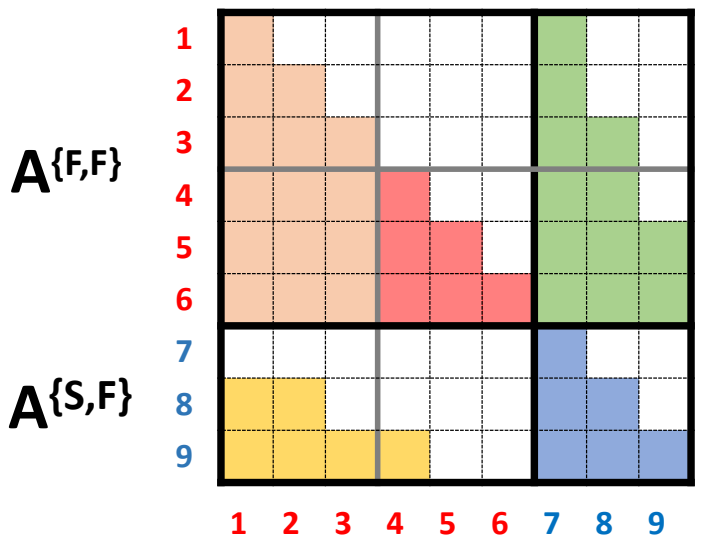

(a) Butcher tableau of a decoupled MR-GARK. The coupling matrices have complementary sparsity patterns.

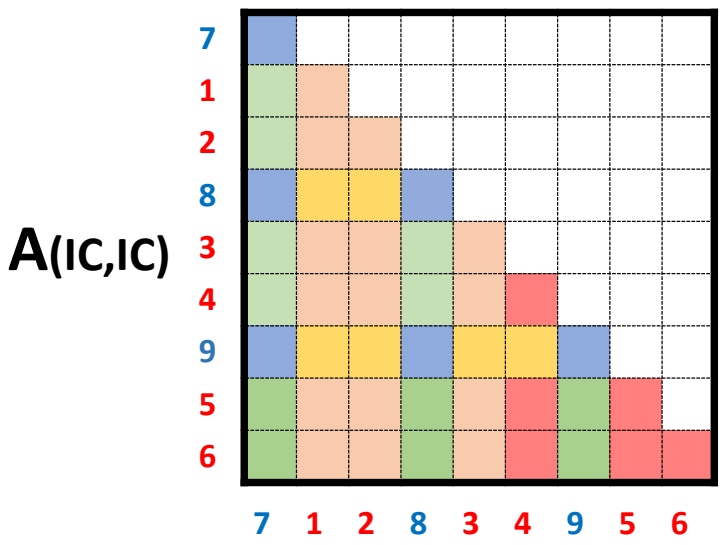

(b) Permuted Butcher tableau of a decoupled MR-GARK. The slow and fast stages can be computed in succession.

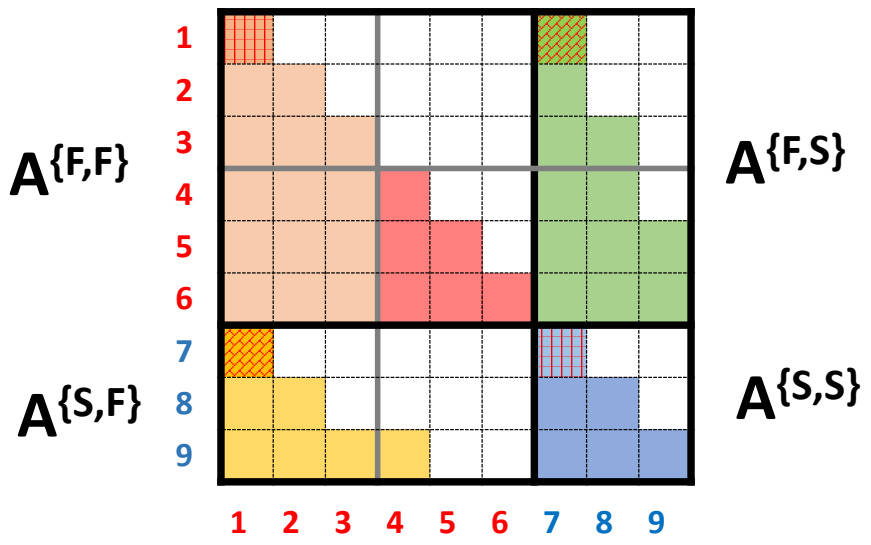

(c) Butcher tableau of a coupled MR-GARK. The brick-dashed entries of the coupling matrices violate the sparsity complementarity.

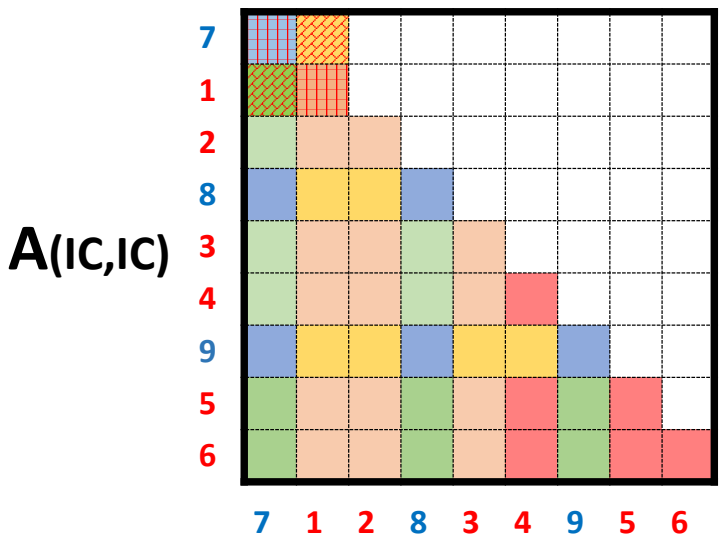

(d) Permuted Butcher tableau of a coupled MR-GARK. The first slow stage and the first fast stage of the first micro-step, dashed, need to be computed together in a coupled manner.

Fig. 8: Example of decoupled and coupled MR-GARK with $\mathtt{m} = 2$ and $s^{\{\mathrm{F}\}} = s^{\{\mathrm{S}\}} = 3$. Blue is the slow method, pink the first fast step, dark pink the second fast step, and green and yellow are the couplings. Rows 1:3 in the Butcher tableaus (a) and (c) correspond to first micro-step fast stages $\mathbf{Y}_{1:3}^{\{\mathrm{F},\cdot,1\}}$,rows 4:6 to second micro-step fast stages $\mathbf{Y}_{1:3}^{\{\mathrm{F},\cdot,2\}}$,and rows 7:9 to slow stages $\mathbf{Y}_{1:3}^{\{\mathrm{S},\cdot\}}$. The permuted versions of the tableaus reflect the sequential order of stage computations. Note the entry above the diagonal for the coupled, permuted tableau (d) due to the non-complementary coupling sparsity structure in (c).

to the slow stages (5.5a). Let $\textsc{ic} \in \mathbb{N}^{\mathtt{m}s^{\{\mathrm{F}\}}+s^{\{\mathrm{S}\}}}$ be the vector of GARK tableau row indices sorted in the order in which the corresponding stages are actually computed; this is a topological ordering of the condensation (computation dependency) graph. A renumbering of stages leads to a row and column permutation of the Butcher tableau (5.4). The reordered Butcher matrix is $\mathbf{A}_{\mathrm{GARK}}(\textsc{ic}, \textsc{ic})$. The reordering is illustrated in Figure 8. We distinguish the following cases:

- If the reordered Butcher tableau is strictly lower triangular then the MR-GARK method is explicit, and each stage uses only previously computed information.
- If the reordered Butcher tableau is lower triangular, with some non-zero diagonal entries, then the MR-GARK method is diagonally implicit and decoupled. The non-zero diagonal entries correspond to implicit fast or slow stages. This case is illustrated in Figure 8b.

- In general a permuted Butcher tableau of a MR-GARK scheme has a block lower triangular structure. Similar to the case of a fully implicit Runge-Kutta scheme, diagonal blocks indicate which groups of stages need to be solved together. If the blocks involve both fast and slow stages then the MR-GARK method is coupled. This case is illustrated in Figure 8d.

**5.4.** MR-GARK **order conditions.** Order conditions of MR-GARK methods (5.5) follow directly from the general GARK order conditions theory, see [27]. The order conditions are formulated in terms of the coefficients of the base methods and of the coupling coefficients using the block structure of (5.4).

There are no coupling conditions for order one. Moreover, if the internal consistency conditions (5.7) hold, and each individual method $RK^{\{\mathrm{S}\}}$ and $RK^{\{\mathrm{F}\}}$ has at least order two, then the MR-GARK scheme (5.5) has order two [27].

Next, assume that the internal consistency conditions (5.7) hold, and that each individual method $RK^{\{\mathrm{S}\}}$ and $RK^{\{\mathrm{F}\}}$ has at least order three. Then the MR-GARK scheme (5.5) has order three if and only if the following coupling conditions hold:

$$b^{\{\mathrm{F}\}T}\left(\sum_{\ell=1}^{\mathtt{m}} A^{\{\mathrm{F},\mathrm{S},\ell\}}\right)c^{\{\mathrm{S}\}} = \frac{\mathtt{m}}{6}, \tag{5.10a}$$

$$b^{\{\mathrm{S}\}T}\left(\sum_{\ell=1}^{\mathtt{m}}(\ell-1)\,c^{\{\mathrm{S},\mathrm{F},\ell\}}\right) + b^{\{\mathrm{S}\}T}\left(\sum_{\ell=1}^{\mathtt{m}} A^{\{\mathrm{S},\mathrm{F},\ell\}}\right)c^{\{\mathrm{F}\}} = \frac{\mathtt{m}^2}{6}. \tag{5.10b}$$

For higher order coupling conditions see [27].

**5.5. Matrix stability analysis.** For component partitioned systems (2.1) we consider the following model problem [36,52]:

$$\begin{bmatrix}\dot{\mathbf{y}}^{\{\mathrm{F}\}}\\ \dot{\mathbf{y}}^{\{\mathrm{S}\}}\end{bmatrix} = \underbrace{\begin{bmatrix}\lambda^{\{\mathrm{F},\mathrm{F}\}} & \frac{1-\xi}{\alpha}\left(\lambda^{\{\mathrm{F},\mathrm{F}\}}-\lambda^{\{\mathrm{S},\mathrm{S}\}}\right)\\ -\alpha\,\xi\,\left(\lambda^{\{\mathrm{F},\mathrm{F}\}}-\lambda^{\{\mathrm{S},\mathrm{S}\}}\right) & \lambda^{\{\mathrm{S},\mathrm{S}\}}\end{bmatrix}}_{\boldsymbol{\Omega}}\begin{bmatrix}\mathbf{y}^{\{\mathrm{F}\}}\\ \mathbf{y}^{\{\mathrm{S}\}}\end{bmatrix}. \tag{5.11}$$

The eigenvalue-eigenvector pairs of $\boldsymbol{\Omega}$ are linear combinations of the slow and fast diagonal terms:

$$\left\{\xi\,\lambda^{\{\mathrm{F},\mathrm{F}\}}+(1-\xi)\,\lambda^{\{\mathrm{S},\mathrm{S}\}},\ \begin{bmatrix}-\frac{1}{\alpha}\\ 1\end{bmatrix}\right\},\quad \left\{(1-\xi)\,\lambda^{\{\mathrm{F},\mathrm{F}\}}+\xi\,\lambda^{\{\mathrm{S},\mathrm{S}\}},\ \begin{bmatrix}-\frac{1-\xi}{\alpha}\\ \xi\end{bmatrix}\right\}. \tag{5.12}$$

Several simplifications are proposed in [52] for analysis.

- We restrict analysis to the case where the coupling variable $\xi\in[0,1]$, such that the eigenvalues (5.12) are convex combinations of $\lambda^{\{\mathrm{F},\mathrm{F}\}}$ and $\lambda^{\{\mathrm{S},\mathrm{S}\}}$, and the coupled system (5.11) is stable for any $\lambda^{\{\mathrm{F},\mathrm{F}\}},\lambda^{\{\mathrm{S},\mathrm{S}\}}\in\mathbb{C}^-$. For $|\xi|\ll 1$ the fast sub-system has a weak influence on the slow one; the first eigenvalue (5.12) is slow and the second one is fast. For $|1-\xi|\ll 1$ the slow sub-system has a weak influence on the fast one, and the first eigenvalue (5.12) is fast.
- The scale analysis in [52] shows that, in order to have the slow and fast contributions to $\dot{\mathbf{y}}^{\{\mathrm{S}\}}$ of similar magnitude, and the contributions to $\dot{\mathbf{y}}^{\{\mathrm{F}\}}$ of similar magnitude as well, one needs a coupling coefficient inversely proportional to the scale ratio of the system, $\xi\sim(|\lambda^{\{\mathrm{F},\mathrm{F}\}}/\lambda^{\{\mathrm{S},\mathrm{S}\}}|+1)^{-1}$.
- The $\alpha$ factor in (5.11) represents a scaling of the fast variable, as can be seen from rewriting the system (5.11) in terms of the variables $[\alpha^{-1}\,\mathbf{y}^{\{\mathrm{F}\}},\mathbf{y}^{\{\mathrm{S}\}}]$. For simplicity we take $\alpha=1$, i.e., the multirate method is applied to the system (5.11) after rescaling the fast variables. In general, however, the stability of multirate schemes depends on the scaling.

Let $z^{\{\text{F},\text{F}\}} = H\lambda^{\{\text{F},\text{F}\}}$, $z^{\{\text{S},\text{S}\}} = H\lambda^{\{\text{S},\text{S}\}}$, $z^{\{\text{S},\text{F}\}} = (1-\xi)/\alpha \cdot \left(z^{\{\text{F},\text{F}\}} - z^{\{\text{S},\text{S}\}}\right)$, and $z^{\{\text{F},\text{S}\}} = -\alpha\,\xi \cdot \left(z^{\{\text{F},\text{F}\}} - z^{\{\text{S},\text{S}\}}\right)$. Application of MR-GARK method (5.4) to test problem (5.11) advances over one step $H$ via the recurrence:

$$\begin{bmatrix} \mathbf{y}_{n+1}^{\{\text{F}\}} \\ \mathbf{y}_{n+1}^{\{\text{S}\}} \end{bmatrix} = \mathbf{M}_{\alpha,\xi}(z^{\{\text{F},\text{F}\}}, z^{\{\text{S},\text{S}\}}) \begin{bmatrix} \mathbf{y}_{n}^{\{\text{F}\}} \\ \mathbf{y}_{n}^{\{\text{S}\}} \end{bmatrix},$$

with the stability matrix:

$$\begin{aligned} \mathbf{M}_{\alpha,\xi}(z^{\{\text{F},\text{F}\}}, z^{\{\text{S},\text{S}\}}) &= \mathbf{I}_{2\times 2} + \begin{bmatrix} \mathbf{b}^{\{\text{F}\}\text{T}} & \mathbf{0}_{1\times s^{\{\text{S}\}}} \\ \mathbf{0}_{1\times \mathtt{m}s^{\{\text{F}\}}} & \mathbf{b}^{\{\text{S}\}\text{T}} \end{bmatrix} \cdot \\ &\cdot \begin{bmatrix} \mathbf{I}_{s^{\{\text{F}\}}\times s^{\{\text{F}\}}} - z^{\{\text{F},\text{F}\}}\,\mathbf{A}^{\{\text{F},\text{F}\}} & -z^{\{\text{S},\text{F}\}}\,\mathbf{A}^{\{\text{F},\text{S}\}} \\ -z^{\{\text{F},\text{S}\}}\,\mathbf{A}^{\{\text{S},\text{F}\}} & \mathbf{I}_{s^{\{\text{S}\}}\times s^{\{\text{S}\}}} - z^{\{\text{S},\text{S}\}}\,\mathbf{A}^{\{\text{S},\text{S}\}} \end{bmatrix}^{-1} \cdot \\ &\cdot \begin{bmatrix} z^{\{\text{F},\text{F}\}}\,\mathbf{1}_{\mathtt{m}s^{\{\text{F}\}}\times 1} & z^{\{\text{S},\text{F}\}}\,\mathbf{1}_{\mathtt{m}s^{\{\text{F}\}}\times 1} \\ z^{\{\text{F},\text{S}\}}\,\mathbf{1}_{s^{\{\text{S}\}}\times 1} & z^{\{\text{S},\text{S}\}}\,\mathbf{1}_{s^{\{\text{S}\}}\times 1} \end{bmatrix}. \end{aligned} \tag{5.13}$$

Definition 5.3 (Matrix stability). *The matrix slow stability region is:*

$$\mathcal{S}^{2\text{D}}_{\rho,\theta} := \left\{ z^{\{\text{S},\text{S}\}} \in \mathbb{C} \;:\; \begin{array}{r} \sup_k \|\mathbf{M}_{\alpha,\xi}(z^{\{\text{F},\text{F}\}}, z^{\{\text{S},\text{S}\}})^k\| \le 1,\ \forall\, z^{\{\text{F},\text{F}\}} \in \mathbb{C} \text{ with} \\ |z^{\{\text{F},\text{F}\}}| \le \rho \text{ and } |\arg(z^{\{\text{F},\text{F}\}}) - \pi| \le \theta \end{array} \right\}. \tag{5.14}$$

*The error propagation matrix $\mathbf{M}_{\alpha,\xi}(z^{\{\text{F},\text{F}\}}, z^{\{\text{S},\text{S}\}})$ is stable for any slow component $z^{\{\text{S},\text{S}\}} \in \mathcal{S}^{2\text{D}}_{\rho,\theta}$ paired with any fast component $z^{\{\text{F},\text{F}\}}$ in a $\theta$-wedge in the negative complex half plane. The matrix fast stability region $\mathcal{F}^{2\text{D}}_{\rho,\theta}$ is defined in a similar manner. The stability regions depend on both the coupling variable $\xi$ and on the scaling variable $\alpha$.*

For practical purposes we consider the following parameters:

$$\alpha = 1, \quad z^{\{\text{F},\text{F}\}} = \mathtt{m}\,\exp^{i\theta^{\{\text{F}\}}}, \quad z^{\{\text{S},\text{S}\}} = \exp^{i\theta^{\{\text{S}\}}}, \quad \theta^{\{\text{F}\}}, \theta^{\{\text{S}\}} \in [\pi/2, 3\pi/2].$$

Both slow and fast eigenvalues have negative real parts, and the ratio of their magnitudes coincides with the step size ration. This choice allows to visualize the stability region for each coupling coefficient $\xi$ in the $(\theta^{\{\text{F}\}}, \theta^{\{\text{S}\}})$ plane.

Example 3. *Consider the second order implicit midpoint as the slow and the fast base method:*

$$\begin{array}{c|c} c^{\{\text{F},\text{F}\}} & A^{\{\text{F},\text{F}\}} \\ \hline & b^{\{\text{F}\}\text{T}} \end{array} = \begin{array}{c|c} c^{\{\text{S},\text{S}\}} & A^{\{\text{S},\text{S}\}} \\ \hline & b^{\{\text{S}\}\text{T}} \end{array} = \begin{array}{c|c} \frac{1}{2} & \frac{1}{2} \\ \hline & 1 \end{array}.$$

*For odd* $\mathtt{m}$ *the* MR-GARK *coupling coefficients*

$$A^{\{\text{F},\text{S},\ell\}} = \begin{cases} 0, & \ell < \lfloor \mathtt{m}/2 \rfloor \\ \frac{1}{2}, & \ell = \lfloor \mathtt{m}/2 \rfloor \\ 1, & \ell > \lfloor \mathtt{m}/2 \rfloor \end{cases}, \qquad A^{\{\text{S},\text{F},\ell\}} = \begin{cases} 1, & \ell < \lfloor \mathtt{m}/2 \rfloor \\ \frac{1}{2}, & \ell = \lfloor \mathtt{m}/2 \rfloor \\ 0, & \ell > \lfloor \mathtt{m}/2 \rfloor \end{cases},$$

*lead to a coupled multirate midpoint method. For even* $\mathtt{m}$ *the coupling coefficients*

$$A^{\{\text{F},\text{S},\ell\}} = \begin{cases} 0 & \ell \le \mathtt{m}/2 \\ 1 & \textit{otherwise} \end{cases}, \qquad A^{\{\text{S},\text{F},\ell\}} = \begin{cases} 1 & \ell \le \mathtt{m}/2 \\ 0 & \textit{otherwise} \end{cases}. \tag{5.15}$$

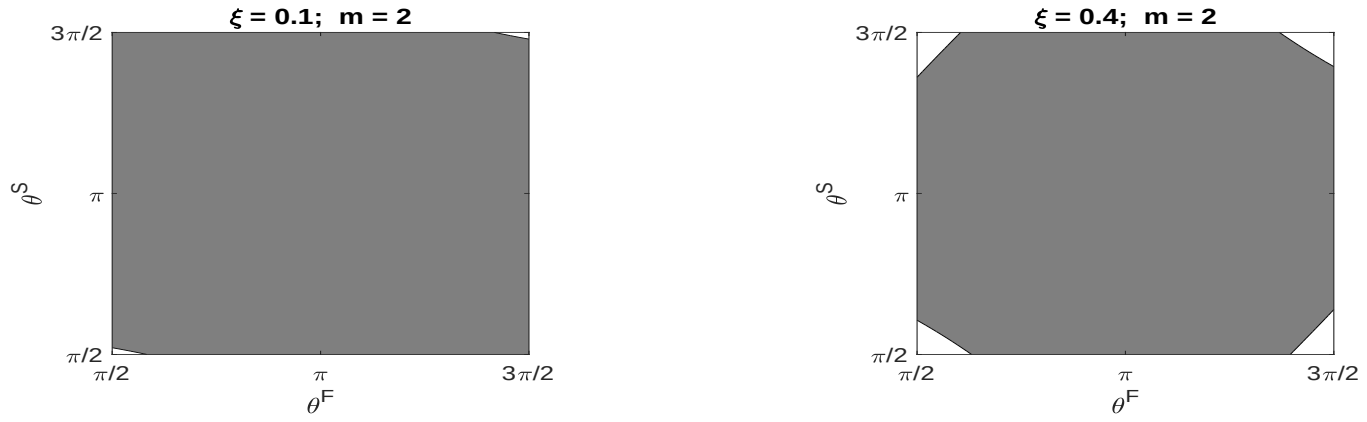

(a) $\mathtt{m} = 2$

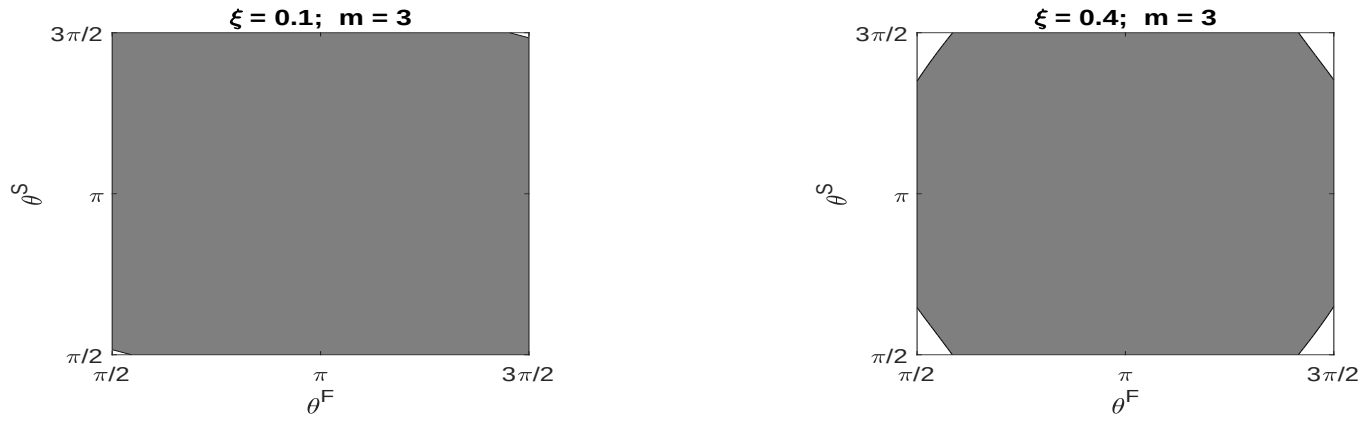

(b) $\mathtt{m} = 3$

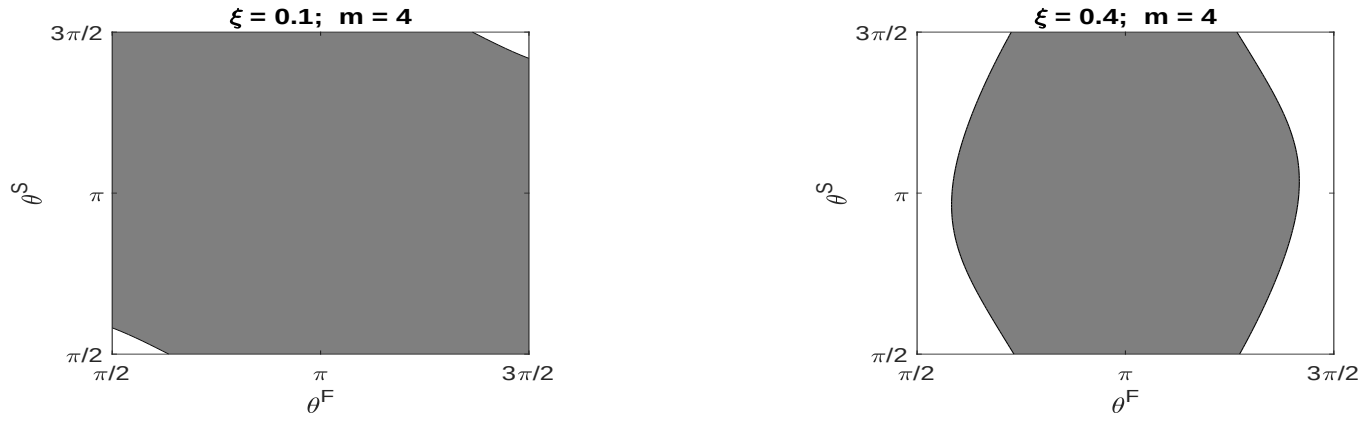

(c) $\mathtt{m} = 4$

Fig. 9: Stability plots for the multirate implicit midpoint method.

*lead to a decoupled multirate midpoint method. The stability plots are shown in Figure 9.*

**5.6. Examples of traditional multirate Runge-Kutta methods.** The traditional multirate Runge-Kutta methods proposed in the literature can be represented and analyzed in the MR-GARK framework. For a detailed discussion see [27].

**5.6.1. Predictor-corrector approach.** The predictor-corrector approach is based on the early work of Rice [41] and the later developments in [58, 69]. The predictor step solves the full system (2.1) with a standard Runge-Kutta method $(A, b, c)$ with the macro-step $H$, to obtain predicted stages $(\mathbf{Y}_j^{\{\mathrm{S}\}}, \mathbf{Y}_j^{\{\mathrm{F}\}*})$ and solution $(\mathbf{y}_{n+1}^{\{\mathrm{S}\}}, \mathbf{y}_{n+1}^{\{\mathrm{F}\}*})$. As the integration of the fast subsystem with large time step $H$ is inaccurate, the results of the fast component $\mathbf{y}_{n+1}^{\{\mathrm{F}\}*}$ are discarded. The corrector step repeats the integration of the fast partition using small timesteps $h$. Following Savcenco, Hundsdorfer, and co-workers [54–56, 58], one can employ the dense output of the scheme $(A, b, c)$ to interpolate the slow components, and fill in the slow values needed during the fast micro-step integration. A continuous approximation of the slow function values at intermediate points is

$$\mathbf{f}^{\{\mathrm{S}\}}(\mathbf{y}(t_n + \theta h)) = \sum_{j=1}^{s} \gamma_j(\theta)\, \mathbf{f}^{\{\mathrm{S}\}}(\mathbf{Y}_j^{\{\mathrm{S}\}}, \mathbf{Y}_j^{\{\mathrm{F}\}*}), \quad 0 \le \theta \le 1. \tag{5.16}$$

Approximations (5.16) are used to compute the slow function values at the micro-steps (5.5b); this leads to a multirate GARK approach with coupling coefficients

$$a_{i,j}^{\{\mathrm{F},\mathrm{S},\ell\}} = \gamma_j \left( \frac{\ell - 1 + c_i^{\{\mathrm{F},\mathrm{F}\}}}{\mathtt{m}} \right).$$

A third order predictor-corrector MR-GARK method, referred to as MR-GARK SDIRK3 [43], is built directly from the following base method of Alexander [2]:

$$\begin{array}{c|ccc}
\gamma & \gamma & 0 & 0 \\
\frac{\gamma}{2} + \frac{1}{2} & \frac{1}{2} - \frac{\gamma}{2} & \gamma & 0 \\
1 & -\frac{3\gamma^2}{2} + 4\gamma - \frac{1}{4} & \frac{3\gamma^2}{2} - 5\gamma + \frac{5}{4} & \gamma \\
\hline
 & -\frac{3\gamma^2}{2} + 4\gamma - \frac{1}{4} & \frac{3\gamma^2}{2} - 5\gamma + \frac{5}{4} & \gamma \\
\hdashline
 & -\frac{3\gamma^2}{2} + 3\gamma - \frac{1}{4} & \frac{3\gamma^2}{2} - 3\gamma + \frac{5}{4} & 0
\end{array}, \qquad \gamma \approx 0.43586\ldots. \tag{5.17}$$

The nonzero coupling coefficients are:

$$a_{1,1}^{\{\mathrm{F},\mathrm{S},\ell\}} = \frac{\left(6\gamma^2 - 24\gamma + 5\right)\left(2\gamma^2 + 2\gamma(\ell-1) + (\ell-1)^2\right) - 8\gamma^3\mathtt{m}^2}{(\gamma-1)\left(6\gamma^2 - 20\gamma + 5\right)\mathtt{m}^2} - \frac{\left(6\gamma^3 - 30\gamma^2 - 15\gamma + 5\right)\mathtt{m}(\gamma + \ell - 1)}{(\gamma-1)\left(6\gamma^2 - 20\gamma + 5\right)\mathtt{m}^2},$$

$$a_{1,2}^{\{\mathrm{F},\mathrm{S},\ell\}} = \frac{-5(\ell-1)^2 + 12\gamma^4(\mathtt{m}-1) + 4\gamma^3(\mathtt{m}-1)(3\ell + 4\mathtt{m} - 17)}{(\gamma-1)\left(6\gamma^2 - 20\gamma + 5\right)\mathtt{m}^2} + \frac{\gamma^2\left(-6\ell^2 + 68\ell + (82 - 72\ell)\mathtt{m} - 72\right) + 2\gamma(\ell-1)(14\ell + 5\mathtt{m} - 19)}{(\gamma-1)\left(6\gamma^2 - 20\gamma + 5\right)\mathtt{m}^2},$$

$$a_{1,3}^{\{\mathrm{F},\mathrm{S},\ell\}} = -\frac{4\gamma\left((\ell-1)^2 + 2\gamma^2(\mathtt{m}-1)^2 - 2\gamma(\ell-1)(2\mathtt{m}-1)\right)}{(\gamma-1)\left(6\gamma^2 - 20\gamma + 5\right)\mathtt{m}^2},$$

$$a_{2,1}^{\{\mathrm{F},\mathrm{S},\ell\}} = -\frac{-2\left(6\gamma^2 - 24\gamma + 5\right)\left(\gamma(\ell+1) + (\ell-1)\ell\right) + 16\gamma^3\mathtt{m}^2}{2(\gamma-1)\left(6\gamma^2 - 20\gamma + 5\right)\mathtt{m}^2} - \frac{\left(6\gamma^3 - 30\gamma^2 - 15\gamma + 5\right)\mathtt{m}(\gamma + 2\ell - 1)}{2(\gamma-1)\left(6\gamma^2 - 20\gamma + 5\right)\mathtt{m}^2},$$

$$a_{2,2}^{\{\mathrm{F},\mathrm{S},\ell\}} = \frac{\gamma\left(28\ell^2 - 33\ell + 5(2\ell-1)\mathtt{m} - 5\right) + 2\gamma^3\left(8\mathtt{m}^2 - 3(\ell+1) + 3(2\ell-7)\mathtt{m}\right)}{(\gamma-1)\left(6\gamma^2 - 20\gamma + 5\right)\mathtt{m}^2} + \frac{6\gamma^4\mathtt{m} - 5(\ell-1)\ell + \gamma^2\left(-6\ell^2 + 34\ell + (41 - 72\ell)\mathtt{m} + 28\right)}{(\gamma-1)\left(6\gamma^2 - 20\gamma + 5\right)\mathtt{m}^2},$$

$$a_{2,3}^{\{\mathrm{F},\mathrm{S},\ell\}} = -\frac{4\gamma\left((\ell-1)\ell + 2\gamma^2(\mathtt{m}-1)\mathtt{m} + \gamma(\ell + (2 - 4\ell)\mathtt{m} + 1)\right)}{(\gamma-1)\left(6\gamma^2 - 20\gamma + 5\right)\mathtt{m}^2},$$

$$\begin{aligned}
a_{3,1}^{\{\mathrm{F},\mathrm{S},\ell\}} &= \frac{-36\gamma^5 + 252\gamma^4 + 20\ell^2 - 4\gamma^3\left(8\mathtt{m}^2 + 6\ell\mathtt{m} + 129\right)}{4(\gamma-1)\left(6\gamma^2 - 20\gamma + 5\right)\mathtt{m}^2} \\
&\quad + \frac{24\gamma^2\left(\ell^2 + 5\ell\mathtt{m} + 13\right) + \gamma\left(-96\ell^2 + 60\ell\mathtt{m} - 69\right) - 20\ell\mathtt{m} + 5}{4(\gamma-1)\left(6\gamma^2 - 20\gamma + 5\right)\mathtt{m}^2},
\end{aligned}$$

$$\begin{aligned}
a_{3,2}^{\{\mathrm{F},\mathrm{S},\ell\}} &= \frac{36\gamma^5 - 276\gamma^4 - 5\left(4\ell^2 + 1\right) + \gamma^3\left(64\mathtt{m}^2 + 48\ell\mathtt{m} + 588\right)}{4(\gamma-1)\left(6\gamma^2 - 20\gamma + 5\right)\mathtt{m}^2} \\
&\quad + \frac{-12\gamma^2\left(2\ell^2 + 24\ell\mathtt{m} + 29\right) + \gamma\left(112\ell^2 + 40\ell\mathtt{m} + 73\right)}{4(\gamma-1)\left(6\gamma^2 - 20\gamma + 5\right)\mathtt{m}^2},
\end{aligned}$$

$$a_{3,3}^{\{\mathrm{F},\mathrm{S},\ell\}} = \frac{\gamma\left(6\gamma^3 - 4\ell^2 - 2\gamma^2\left(4\mathtt{m}^2 + 9\right) + \gamma(16\ell\mathtt{m} + 9) - 1\right)}{(\gamma-1)\left(6\gamma^2 - 20\gamma + 5\right)\mathtt{m}^2}.$$

**5.6.2. Multirate additive Runge-Kutta schemes.** Another special case of multirate GARK schemes (5.5) are the multirate additive Runge-Kutta (ARK) schemes, as discussed in [10]. For a second order multirate method one starts with a second order base explicit scheme $(A, b, c)$ that is strong stability preserving. The multirate MPRK-2 scheme proposed in [10] reads:

$$\begin{array}{c|c} \mathbf{c}^{\{\mathrm{F}\}} & \mathbf{A}^{\{\mathrm{F,F}\}} \\ \hline & \mathbf{b}^{\{\mathrm{F}\}} \end{array} = \begin{array}{c|cccc} \frac{1}{\mathtt{m}} c & \frac{1}{\mathtt{m}} A & & & \\ \frac{1}{\mathtt{m}}\mathbf{1} + \frac{1}{\mathtt{m}} c & \frac{1}{\mathtt{m}}\mathbf{1}\, b^T & \frac{1}{\mathtt{m}} A & & \\ \vdots & \vdots & & \ddots & \\ \frac{\mathtt{m}-1}{\mathtt{m}}\mathbf{1} + \frac{1}{\mathtt{m}} c & \frac{1}{\mathtt{m}}\mathbf{1}\, b^T & \dots & \frac{1}{\mathtt{m}}\mathbf{1}\, b^T & \frac{1}{\mathtt{m}} A \\ \hline & \frac{1}{\mathtt{m}} b^T & \frac{1}{\mathtt{m}} b^T & \dots & \frac{1}{\mathtt{m}} b^T \end{array}; \quad \begin{array}{c|c} \mathbf{c}^{\{\mathrm{S}\}} & \mathbf{A}^{\{\mathrm{S,S}\}} \\ \hline & \mathbf{b}^{\{\mathrm{S}\}} \end{array} = \begin{array}{c|cccc} c & A & & & \\ c & & A & & \\ \vdots & & & \ddots & \\ c & & & & A \\ \hline & \frac{1}{\mathtt{m}} b^T & \frac{1}{\mathtt{m}} b^T & \dots & \frac{1}{\mathtt{m}} b^T \end{array}.$$

The method is written here in ARK form with the fast and slow tableaus of coefficients having the same number of stages. For a GARK formulation one takes $\mathbf{A}^{\{\mathrm{F,S}\}} = \mathbf{A}^{\{\mathrm{S,S}\}}$ and $\mathbf{A}^{\{\mathrm{S,F}\}} = \mathbf{A}^{\{\mathrm{F,F}\}}$. Clearly the scheme is not internally consistent. It was shown in [10] that, under regular assumptions, the MPRK-2 scheme is second order accurate, preserves linear invariants, positivity, and a maximum principle, and is total variation bounded. Therefore MPRK-2 is suitable for multirate integration of hyperbolic conservation laws on non-uniform spatial grids, where shorter steps are used for smaller grid elements.

**5.6.3. Kvaerno and Rentrop methods.** The mRK class of multirate Runge-Kutta methods proposed by Kvaerno and Rentrop [37] can be formulated in the GARK framework. One chooses fast and slow base schemes $(A^{\{\mathrm{F,F}\}}, b^{\{\mathrm{F}\}})$ and $(A^{\{\mathrm{S,S}\}}, b^{\{\mathrm{S}\}})$, fixed coupling matrices $\bar{A}^{\{\mathrm{F,S}\}}$, $\bar{A}^{\{\mathrm{S,F}\}}$, and micro-step dependent coupling function $\eta(\ell)$. The coupling has the following structure.

a) The slow dynamics is impacted only by the first fast microstep:

$$A^{\{\mathrm{S,F},1\}} = \mathtt{m}\, \bar{A}^{\{\mathrm{S,F}\}}; \quad A^{\{\mathrm{S,F},\ell\}} = 0, \quad \ell = 2, \dots, \mathtt{m}.$$

b) The fast dynamics is influenced by the slow one via the coupling:

$$A^{\{\mathrm{F,S},1\}} = \frac{1}{\mathtt{m}}\, \bar{A}^{\{\mathrm{F,S}\}}; \quad A^{\{\mathrm{F,S},\ell\}} = \frac{1}{\mathtt{m}} \left( \bar{A}^{\{\mathrm{F,S}\}} + \mathbf{1}^{\{\mathrm{F}\}}\, \eta(\ell-1)^T \right);$$
$$\eta : \mathbb{N} \to \mathbb{R}^{s^{\{\mathrm{S}\}}}, \quad \eta(0) = \mathbf{0}^{\{\mathrm{S}\}}, \quad \eta^T(\ell)\, \mathbf{1}^{\{\mathrm{S}\}} = \ell.$$

This choice aids the fulfillment of the internal consistency conditions (5.7) by choosing the coupling matrices $\bar{A}^{\{\mathrm{F,S}\}}$ and $\bar{A}^{\{\mathrm{S,F}\}}$ such that:

$$\bar{A}^{\{\mathrm{S,F}\}}\, \mathbf{1}^{\{\mathrm{F}\}} = A^{\{\mathrm{S,S}\}} \mathbf{1}^{\{\mathrm{S}\}} = c^{\{\mathrm{S}\}}, \quad \bar{A}^{\{\mathrm{F,S}\}}\, \mathbf{1}^{\{\mathrm{S}\}} = A^{\{\mathrm{F,F}\}} \mathbf{1}^{\{\mathrm{F}\}} = c^{\{\mathrm{F}\}}. \tag{5.18}$$

Choose two order three schemes $(A^{\{\mathrm{F,F}\}}, b^{\{\mathrm{F}\}})$ and $(A^{\{\mathrm{S,S}\}}, b^{\{\mathrm{S}\}})$. To obtain a multirate method of order three choose the free parameters $\bar{A}^{\{\mathrm{F,S}\}}$, $\bar{A}^{\{\mathrm{S,F}\}}$, and $\eta(\ell)$ fulfill the simplifying conditions (5.18) together with the order three coupling conditions (5.10):

$$\frac{\mathtt{m}}{6} = b^{\{\mathrm{S}\}T}\, \bar{A}^{\{\mathrm{S,F}\}}\, c^{\{\mathrm{F}\}}, \tag{5.19a}$$

$$\frac{\mathtt{m}}{6} = b^{\{\mathrm{F}\}T}\, \bar{A}^{\{\mathrm{F,S}\}}\, c^{\{\mathrm{S}\}} + \frac{1}{\mathtt{m}} \sum_{\ell=1}^{\mathtt{m}} \eta^T(\ell-1)\, c^{\{\mathrm{S}\}}. \tag{5.19b}$$

Note that the additional condition imposed by Kvaerno and Rentrop [37]

$$\eta^T(\ell-1)\, c^{\{\mathrm{S}\}} = \frac{\ell(\ell+1)}{2\mathtt{m}}$$

ensures that the fast solution component ($\mathbf{y}^{\{\mathrm{F}\}}_{n+\ell/\mathtt{m}}$ in (5.5b)) has order three *at all micro-steps* $\ell = 1, \dots, \mathtt{m}$.

EXAMPLE 4 (Multirate Radau-IIA schemes of order 3). *Consider the case where both the fast and the slow schemes are the RADAU-IIA scheme with* $p = 3, s = 2$, $(A^{\{\mathrm{F},\mathrm{F}\}}, b^{\{\mathrm{F}\}}) = (A^{\{\mathrm{S},\mathrm{S}\}}, b^{\{\mathrm{S}\}}) = (A, b)$. *A Kvaerno-Rentrop third order mRK scheme is obtained by selecting* $\bar{A}^{\{\mathrm{F},\mathrm{S}\}} = \bar{A}^{\{\mathrm{S},\mathrm{F}\}} = \bar{A}$ *and:*

$$\begin{array}{c|c} c & A \\ \hline & b^T \end{array} = \begin{array}{c|cc} 0 & \frac{1}{4} & -\frac{1}{4} \\ \frac{2}{3} & \frac{1}{4} & \frac{5}{12} \\ \hline & \frac{1}{4} & \frac{3}{4} \end{array}, \quad \eta(\ell) = \begin{bmatrix} \frac{3}{2}\ell \\ -\frac{1}{2}\ell \end{bmatrix}, \quad \bar{A} = \begin{bmatrix} \frac{1}{3} & 0 \\ -\mathtt{m} & \mathtt{m}-1 \end{bmatrix}.$$

**5.7. Examples of decoupled** MR-GARK **methods.** This section presents the methods developed in [53] and their properties. We use the usual naming convention where $p$ is the method order, $\hat{p}$ is the embedded order, and $s$ is the number of stages of the base method.

**5.7.1.** MR-GARK **fast explicit** ($s = 2$)**, slow explicit** ($s = 2$)**,** $p = 2$**,** $\hat{p} = 1$**.** This explicit method uses a base method from [40] and is telescopic.

$$A^{\{\mathrm{F},\mathrm{F}\}} = A^{\{\mathrm{S},\mathrm{S}\}} = \begin{bmatrix} 0 & 0 \\ \frac{2}{3} & 0 \end{bmatrix}, \qquad A^{\{\mathrm{F},\mathrm{S},1\}} = \begin{bmatrix} 0 & 0 \\ \frac{2}{3\mathtt{m}} & 0 \end{bmatrix},$$

$$A^{\{\mathrm{F},\mathrm{S},\ell\}} = \begin{bmatrix} \frac{3\mathtt{m}^3-11\mathtt{m}^2+20\ell\mathtt{m}-20\mathtt{m}-20\ell+20}{20(\mathtt{m}-1)\mathtt{m}} & -\frac{\mathtt{m}(3\mathtt{m}-11)}{20(\mathtt{m}-1)} \\ \frac{-3\mathtt{m}^3-9\mathtt{m}^2+60\ell\mathtt{m}-20\mathtt{m}-60\ell+20}{60(\mathtt{m}-1)\mathtt{m}} & \frac{\mathtt{m}(\mathtt{m}+3)}{20(\mathtt{m}-1)} \end{bmatrix}, \quad \ell = 2, \dots, \mathtt{m},$$

$$A^{\{\mathrm{S},\mathrm{F},1\}} = \begin{bmatrix} 0 & 0 \\ -\frac{1}{3}(\mathtt{m}-2)\mathtt{m} & \frac{\mathtt{m}^2}{3} \end{bmatrix}, \qquad A^{\{\mathrm{S},\mathrm{F},\ell\}} = \begin{bmatrix} 0 & 0 \\ 0 & 0 \end{bmatrix}, \qquad \ell = 2, \dots, \mathtt{m},$$

$$b^{\{\mathrm{F}\}} = b^{\{\mathrm{S}\}} = \begin{bmatrix} \frac{1}{4} & \frac{3}{4} \end{bmatrix}^T, \qquad \hat{b}^{\{\mathrm{F}\}} = \hat{b}^{\{\mathrm{S}\}} = \begin{bmatrix} 1 & 0 \end{bmatrix}^T.$$

**5.7.2.** MR-GARK **fast explicit** ($s = 2$)**, slow implicit** ($s = 2$)**,** $p = 2$**,** $\hat{p} = 1$**.** This explicit-implicit method uses the fast method from [40] and slow method from [2]. The multirate scheme is stiffly accurate in the slow partition and the coupling error is independent of $\mathtt{m}$.

$$A^{\{\mathrm{F},\mathrm{F}\}} = \begin{bmatrix} 0 & 0 \\ \frac{2}{3} & 0 \end{bmatrix}, \qquad A^{\{\mathrm{S},\mathrm{S}\}} = \begin{bmatrix} 1-\frac{1}{\sqrt{2}} & 0 \\ \frac{1}{\sqrt{2}} & 1-\frac{1}{\sqrt{2}} \end{bmatrix}, \quad A^{\{\mathrm{S},\mathrm{F},1\}} = \begin{bmatrix} \mathtt{m}-\frac{\mathtt{m}}{\sqrt{2}} & 0 \\ \frac{1}{4} & \frac{3}{4} \end{bmatrix},$$

$$A^{\{\mathrm{F},\mathrm{S},\ell\}} = \begin{bmatrix} \frac{\ell-1}{\mathtt{m}} & 0 \\ \frac{3\ell-1}{3\mathtt{m}} & 0 \end{bmatrix}, \quad \ell = 1, \dots, \mathtt{m}, \qquad A^{\{\mathrm{S},\mathrm{F},\ell\}} = \begin{bmatrix} 0 & 0 \\ \frac{1}{4} & \frac{3}{4} \end{bmatrix}, \qquad \ell = 2, \dots, \mathtt{m},$$

$$b^{\{\mathrm{F}\}} = \begin{bmatrix} \frac{1}{4} & \frac{3}{4} \end{bmatrix}^T, \qquad b^{\{\mathrm{S}\}} = \begin{bmatrix} \frac{1}{\sqrt{2}} & 1-\frac{1}{\sqrt{2}} \end{bmatrix}^T, \qquad \hat{b}^{\{\mathrm{F}\}} = \begin{bmatrix} 1 & 0 \end{bmatrix}^T, \qquad \hat{b}^{\{\mathrm{S}\}} = \begin{bmatrix} \frac{3}{5} & \frac{2}{5} \end{bmatrix}^T.$$

**5.7.3.** MR-GARK **fast implicit** ($s = 2$)**, slow explicit** ($s = 2$)**,** $p = 2$**,** $\hat{p} = 1$**.** This implicit-explicit method uses the fast method from [2] and the slow method from [40]. The multirate scheme is stiffly accurate in the fast partition.

$$A^{\{\mathrm{F,F}\}} = \begin{bmatrix} 1-\frac{1}{\sqrt{2}} & 0 \\ \frac{1}{\sqrt{2}} & 1-\frac{1}{\sqrt{2}} \end{bmatrix}, \quad A^{\{\mathrm{S,S}\}} = \begin{bmatrix} 0 & 0 \\ \frac{2}{3} & 0 \end{bmatrix}, \quad A^{\{\mathrm{F,S},\mathtt{m}\}} = \begin{bmatrix} \frac{2\mathtt{m}-\sqrt{2}}{2\mathtt{m}} & 0 \\ \frac{1}{4} & \frac{3}{4} \end{bmatrix},$$

$$A^{\{\mathrm{F,S},\ell\}} = \begin{bmatrix} \frac{2\ell-\sqrt{2}}{2\mathtt{m}} & 0 \\ \frac{\ell}{\mathtt{m}} & 0 \end{bmatrix}, \qquad \ell = 1,\dots,\mathtt{m}-1, \qquad A^{\{\mathrm{S,F},\ell\}} = \begin{bmatrix} 0 & 0 \\ \frac{2}{3} & 0 \end{bmatrix}, \qquad \ell = 1,\dots,\mathtt{m},$$

$$b^{\{\mathrm{F}\}} = \begin{bmatrix} \frac{1}{\sqrt{2}} & 1-\frac{1}{\sqrt{2}} \end{bmatrix}^{\mathrm{T}}, \qquad b^{\{\mathrm{S}\}} = \begin{bmatrix} \frac{1}{4} & \frac{3}{4} \end{bmatrix}^{\mathrm{T}}, \qquad \hat{b}^{\{\mathrm{F}\}} = \begin{bmatrix} \frac{3}{5} & \frac{2}{5} \end{bmatrix}^{\mathrm{T}}, \qquad \hat{b}^{\{\mathrm{S}\}} = \begin{bmatrix} 1 & 0 \end{bmatrix}^{\mathrm{T}}.$$

**5.7.4.** MR-GARK **fast explicit ($s = 3$), slow explicit ($s = 3$), $p = 3$, $\hat{p} = 2$.** This explicit method uses a three stage base method and has a telescopic property. The free parameters are $c_2$, $L_2$, and $\hat{b}_2$. A practical choice is $L_2 = \text{floor}(c_2\,\mathtt{m})$.

$$A^{\{\mathrm{F,F}\}} = A^{\{\mathrm{S,S}\}} = \begin{bmatrix} 0 & 0 & 0 \\ c_2 & 0 & 0 \\ \frac{(3c_2^2-3c_2+1)}{c_2(3c_2-2)} & \frac{(c_2-1)}{c_2(3c_2-2)} & 0 \end{bmatrix},$$

$$A^{\{\mathrm{F,S},\ell\}} = \begin{bmatrix} \frac{\ell-1}{\mathtt{m}} & 0 & 0 \\ \frac{c_2+\ell-1}{\mathtt{m}} & 0 & 0 \\ \frac{\ell}{\mathtt{m}} & 0 & 0 \end{bmatrix}, \qquad \ell = 1,\dots,L_2,$$

$$A^{\{\mathrm{F,S},\ell\}} = \begin{bmatrix} \frac{2\ell-1}{12c_2(L_2-\mathtt{m})} + \frac{1-2\ell}{12c_2(L_2+\mathtt{m})} + \frac{2\ell-1}{2\mathtt{m}} & \frac{2\ell-1}{12c_2(L_2+\mathtt{m})} + \frac{1-2\ell}{12c_2(L_2-\mathtt{m})} - \frac{1}{2\mathtt{m}} & 0 \\ \frac{2\ell-1}{12c_2(L_2-\mathtt{m})} + \frac{1-2\ell}{12c_2(L_2+\mathtt{m})} + \frac{2\ell-1}{2\mathtt{m}} & \frac{2\ell-1}{12c_2(L_2+\mathtt{m})} + \frac{1-2\ell}{12c_2(L_2-\mathtt{m})} + \frac{2c_2-1}{2\mathtt{m}} & 0 \\ \frac{2\ell-1}{12c_2(L_2-\mathtt{m})} + \frac{1-2\ell}{12c_2(L_2+\mathtt{m})} + \frac{2\ell-1}{2\mathtt{m}} & \frac{2\ell-1}{12c_2(L_2+\mathtt{m})} + \frac{1-2\ell}{12c_2(L_2-\mathtt{m})} + \frac{1}{2\mathtt{m}} & 0 \end{bmatrix}, \qquad \ell = L_2+1,\dots,\mathtt{m},$$

$$A^{\{\mathrm{S,F},\ell\}} = \begin{bmatrix} 0 & 0 & 0 \\ c_2\left(\frac{2\ell(c_2(4L_2-3)-3L_2+3)}{(c_2-1)(3c_2^2+4c_2+1)(L_2+1)} + \frac{\mathtt{m}}{L_2}\right) & -\frac{c_2\ell(c_2(4L_2-3)-3L_2+3)}{(c_2-1)(3c_2^2+4c_2+1)(L_2+1)} & -\frac{c_2\ell(c_2(4L_2-3)-3L_2+3)}{(c_2-1)(3c_2^2+4c_2+1)(L_2+1)} \\ \frac{2c_2\ell(c_2(4L_2-3)-3L_2+3)}{(c_2-1)(3c_2^2+4c_2+1)(L_2+1)} & -\frac{c_2\ell(c_2(4L_2-3)-3L_2+3)}{(c_2-1)(3c_2^2+4c_2+1)(L_2+1)} & -\frac{c_2\ell(c_2(4L_2-3)-3L_2+3)}{(c_2-1)(3c_2^2+4c_2+1)(L_2+1)} \end{bmatrix}, \qquad \ell = 1,\dots,L_2,$$

$$A^{\{\mathrm{S,F},\ell\}} = \begin{bmatrix} 0 & 0 & 0 \\ 0 & 0 & 0 \\ \frac{\ell}{3c_2-2} + \frac{c_2(3L_2-4)-3L_2+3}{6c_2-4} + \frac{\mathtt{m}}{\mathtt{m}-L_2} & \frac{\ell}{2-3c_2} & \frac{c_2(4-3L_2)+3(L_2-1)}{6c_2-4} \end{bmatrix}, \qquad \ell = L_2+1,\dots,\mathtt{m},$$

$$b^{\{\mathrm{F}\}} = b^{\{\mathrm{S}\}} = \begin{bmatrix} \frac{3c_2-1}{6c_2} & -\frac{1}{6(c_2-1)c_2} & \frac{3c_2-2}{6(c_2-1)} \end{bmatrix}^{\mathrm{T}}, \quad \hat{b}^{\{\mathrm{F}\}} = \hat{b}^{\{\mathrm{S}\}} = \begin{bmatrix} \hat{b}_2\,(c_2-1)+\frac{1}{2} & \hat{b}_2 & \frac{1-2\hat{b}_2c_2}{2} \end{bmatrix}^{\mathrm{T}}.$$

**5.8. Examples of coupled** MR-GARK **methods.** In this subsection, we present implicit MR-GARK methods of orders one to four. All methods are telescopic and based on singly diagonally implicit Runge–Kutta (SDIRK) methods. At high order, coupling coefficients can become complicated rational functions of $\ell$ and $\mathtt{m}$.

Consider the L-stable, order two SDIRK base method from [2]

$$\begin{array}{c|c} c & A \\ \hline & b^{\mathrm{T}} \\ & \hat{b}^{\mathrm{T}} \end{array} = \begin{array}{c|cc} \gamma & \gamma & 0 \\ 1 & 1-\gamma & \gamma \\ \hline & 1-\gamma & \gamma \\ & \frac{3}{5} & \frac{2}{5} \end{array}, \qquad \gamma = 1 - \frac{\sqrt{2}}{2}. \tag{5.20}$$

For this base method, an internally consistent standard MR-GARK method must have at least one coupled stage. Enforcing stiff accuracy for both partitions uniquely determines a lightly coupled method:

$$
(5.21)\qquad
\begin{aligned}
A^{\{\text{F},\text{S},\ell\}} &= \begin{bmatrix} \frac{\ell-1+\gamma}{\mathtt{m}} & 0 \\ \begin{cases} 1-\gamma, & \ell=\mathtt{m} \\ \frac{\ell}{\mathtt{m}}, & \text{otherwise} \end{cases} & \begin{cases} \gamma, & \ell=\mathtt{m} \\ 0, & \text{otherwise} \end{cases} \end{bmatrix}, \\
A^{\{\text{S},\text{F},\ell\}} &= \begin{bmatrix} \begin{cases} \mathtt{m}\gamma, & \ell=1 \\ 0, & \text{otherwise} \end{cases} & 0 \\ 1-\gamma & \gamma \end{bmatrix}.
\end{aligned}
$$

For this method, the first slow and fast stages are coupled, an d so are the last slow and fast stages. All other stages are decoupled.

**5.9. Design of** MR-GARK **methods.** In this section, we discuss several desirable properties of practical MR-GARK methods and the design methodologies to incorporate them in the construction of high-order schemes.

The following set of design principles incorporates properties that are important for ensuring the practicality of MR-GARK schemes:

1. Methods where the fast and slow base schemes are both either explicit or implicit should be *telescopic* (5.6) in order to be easily applicable to multi-partitioned systems with multiple time scales.
2. In order to maintain computational efficiency, a complex coupling between multiple fast and slow stages needs to be avoided. *Decoupled* multirate methods satisfying equation (5.8) are less computationally demanding than coupled methods. However, the overall stability of the scheme may be affected by decoupling.
3. Practical use of MR-GARK methods demands an error control mechanism that is capable of adapting both the step size $H$ and the multirate ratio $\mathtt{m}$.
4. It is desirable to derive multirate methods with reduced coupling errors, e.g., by requiring that both the main and the embedded methods satisfy coupling conditions to a higher order than fast and slow local truncation errors. This strategy isolates the dominant local truncation errors to the solution of the slow and fast components and greatly simplifies the task of adaptively choosing $H$ and $\mathtt{m}$, as discussed in subsection 5.9.1.

The order conditions of subsection 5.4 and the additional constraints associated with the above design principles lead to large, nonlinear systems of equations that need to be solved for the method coefficients. The resulting coupling coefficients $\mathbf{A}^{\{\text{S},\text{F}\}}$ and $\mathbf{A}^{\{\text{F},\text{S}\}}$ are *functions* of the micro-step number $\ell$ and the multi-rate step size ratio $\mathtt{m}$. The proposed approach for designing MR-GARK methods of order $p$ is a three-step process, as follows.

1. The first step is to construct optimized base slow $\left(A^{\{\text{S},\text{S}\}}, b^{\{\text{S}\}}\right)$ and fast $\left(A^{\{\text{F},\text{F}\}}, b^{\{\text{F}\}}\right)$ schemes of order $p$.
2. The second step is to define the sparsity patterns of the coupling matrices $\mathbf{A}^{\{\text{S},\text{F}\}}$ and $\mathbf{A}^{\{\text{F},\text{S}\}}$. These patterns determine the computational flow of the method and its implementation complexity, and influence the overall stability and accuracy properties.
3. The third step computes the coupling coefficients $\mathbf{A}^{\{\text{S},\text{F}\}}$ and $\mathbf{A}^{\{\text{F},\text{S}\}}$ such as to satisfy the coupling conditions of subsection 5.4 up to order $p$. Any free parameters are used to minimize the Euclidean norm of the residuals of the order $p+1$ coupling conditions.

An alternative to this derivation strategy is a monolithic constrained optimization procedure that minimizes the residuals of the order $p+1$ coupling conditions, subject to solving the conditions up to order $p$, and subject to structural constraints such as decoupling (5.8) and stiff accuracy. The proposed three-step procedure is preferable for two reasons. First, we expect MR-GARK methods to be applied to problems with a rather weak coupling between partitions, where the primary sources of error are the base methods and not the coupling. This approach gives precedence to the base errors first. Second, the procedure is practical as it reduces the number of nonlinear equations that need to be solved together during the design.

**5.9.1. Error estimation and adaptive** MR-GARK **methods.** Adaptivity of traditional (single rate) methods adjusts the step size such as to ensure the desired accuracy of the solution at a minimal computational effort. In the context of Runge–Kutta schemes a second "embedded" method is used to provide an aposteriori estimate of the local truncation error [30, Chapter II.4]. The step size is adjusted to bring the local error estimate to the user prescribed level using the asymptotic relation that this error is proportional to $H^{p+1}$. Adaptivity of multirate methods is more complex, as there are two independent parameters that control the solution accuracy and efficiency: the macro-step size $H$ and the multirate ratio $\mathtt{m}$ (or, equivalently, the macro-step $H$ and the micro-step $h$.)

**5.9.2. Structure of the** MR-GARK **local truncation error.** The local truncation error of the MR-GARK method (5.5) has three components associated with the slow integration, the fast integration, and with the coupling [53]:

$$\text{LTE} = \text{LTE}^{\{\text{S}\}} + \text{LTE}^{\{\text{F}\}} + \text{LTE}^{\{\text{C}\}}. \tag{5.22}$$

For a method of order $p$ the slow local truncation error is that of applying one step with the base slow Runge–Kutta method with a step size $H$:

$$\|\text{LTE}^{\{\text{S}\}}\| = C^{\{\text{S}\}} \cdot H^{p+1} + \mathcal{O}(H^{p+2}). \tag{5.23a}$$

The fast truncation error is that of applying $\mathtt{m}$ consecutive steps with the base fast Runge–Kutta method with a step size $h = H/\mathtt{m}$:

$$\|\text{LTE}^{\{\text{F}\}}\| = C^{\{\text{F}\}} \, \frac{H^{p+1}}{\mathtt{m}^p} + \mathcal{O}(H^{p+2}) = C^{\{\text{F}\}} \, h^p \, H + \mathcal{O}(H^{p+2}). \tag{5.23b}$$

The principal terms of the coupling errors correspond to derivatives of order $p+1$ involving both $\mathbf{f}^{\{\text{S}\}}$ and $\mathbf{f}^{\{\text{F}\}}$. Each slow term scales with $H$, and each fast term with $H/\mathtt{m}$, leading to the following structure of the local coupling error:

$$\|\text{LTE}^{\{\text{C}\}}\| = \left(\textstyle\sum_{i=-p}^{p} C_i^{\{\text{C}\}} \, \mathtt{m}^{-i}\right) \cdot H^{p+1} + \mathcal{O}(H^{p+2}). \tag{5.23c}$$

While the slow and fast errors are simple functions of $H$ and $\mathtt{m}$, the dependency of the coupling error on $\mathtt{m}$ is more difficult to describe. For example, decreasing $h$ can increase the coupling error terms with negative $\mathtt{m}$ powers. Therefore, in an adaptive method, the effect of changing $\mathtt{m}$ on the coupling error is more difficult to quantify.

**5.9.3. Use of multiple embedded methods to estimate the local truncation error.** Following the traditional Runge–Kutta strategy, we look to use embedded methods to obtain estimates of different components of the local truncation error. Specifically, consider an MR-GARK scheme that produces a main solution $\mathbf{y}_{n+1}$ of order $p$, i.e., cancels all residuals for up to order $p$. A generic embedded solution $\hat{\mathbf{y}}_{n+1}$ cancels all residuals for up to order $p-1$. We seek to build an additional embedded solution $\hat{\mathbf{y}}^{\{\mathrm{S}\}}_{n+1}$ that cancels all slow residuals up to order $p-1$, as well as all fast and all coupling residuals of order $p$. Similarly, we seek to build another embedded solution $\hat{\mathbf{y}}^{\{\mathrm{F}\}}_{n+1}$ that cancels all fast residuals up to order $p-1$, as well as all slow and all coupling residuals of order $p$. The three components of the local truncation error (5.22) are then estimated as follows:

$$
\begin{aligned}
&\mathrm{LTE}_{n+1} = \mathbf{y}_{n+1} - \hat{\mathbf{y}}_{n+1}, \quad \mathrm{LTE}^{\{\mathrm{S}\}}_{n+1} = \mathbf{y}_{n+1} - \hat{\mathbf{y}}^{\{\mathrm{S}\}}_{n+1}, \quad \mathrm{LTE}^{\{\mathrm{F}\}}_{n+1} = \mathbf{y}_{n+1} - \hat{\mathbf{y}}^{\{\mathrm{F}\}}_{n+1}, \\
&\mathrm{LTE}^{\{\mathrm{C}\}}_{n+1} = \mathrm{LTE}_{n+1} - \mathrm{LTE}^{\{\mathrm{S}\}}_{n+1} - \mathrm{LTE}^{\{\mathrm{F}\}}_{n+1} = \hat{\mathbf{y}}^{\{\mathrm{S}\}}_{n+1} + \hat{\mathbf{y}}^{\{\mathrm{F}\}}_{n+1} - \hat{\mathbf{y}}_{n+1} - \mathbf{y}_{n+1}.
\end{aligned}
\tag{5.24}
$$

The construction of three different embedded methods can be very difficult to achieve for high-order schemes. A simplified error estimation strategy utilizes only two embedded methods. Consider an MR-GARK scheme with primary weights $(b^{\{\mathrm{F}\}}, b^{\{\mathrm{S}\}})$ that produces a main solution $\mathbf{y}_{n+1}$ of order $p$. The residuals corresponding to two-trees with up to $p$ nodes are zero. We construct three different embedded solutions, as follows: a generic embedded solution $\hat{\mathbf{y}}_{n+1}$ of order $p-1$, obtained with weights $(\hat{b}^{\{\mathrm{F}\}}, \hat{b}^{\{\mathrm{S}\}})$; an embedded solution $\hat{\mathbf{y}}^{\{\mathrm{S}\}}_{n+1}$, generated with the weights $(b^{\{\mathrm{F}\}}, \hat{b}^{\{\mathrm{S}\}})$; and an embedded solution $\hat{\mathbf{y}}^{\{\mathrm{F}\}}_{n+1}$, generated with $(\hat{b}^{\{\mathrm{F}\}}, b^{\{\mathrm{S}\}})$. In general, these mixed embeddings do not exactly isolate the slow, fast, and coupling error, but they can serve as approximations in (5.24).

**5.9.4. Controlling errors by adapting both the macro-step and the micro-step.** Based on our understanding of the structure of errors in the MR-GARK framework, we now look into practical approaches to adaptivity of multirate methods. The following error estimates are available via the set of embedded methods:

$$
\hat{\varepsilon}_{n+1} := \|\mathbf{y}_{n+1} - \hat{\mathbf{y}}_{n+1}\|_{\mathrm{err}}, \quad \hat{\varepsilon}^{\{\mathrm{S}\}}_{n+1} := \|\mathbf{y}_{n+1} - \hat{\mathbf{y}}^{\{\mathrm{S}\}}_{n+1}\|_{\mathrm{err}}, \quad \hat{\varepsilon}^{\{\mathrm{F}\}}_{n+1} := \|\mathbf{y}_{n+1} - \hat{\mathbf{y}}^{\{\mathrm{F}\}}_{n+1}\|_{\mathrm{err}},
\tag{5.25}
$$

where the error is measured by the traditional specific relative error norm [30, Chapter II.4]. Several strategies that use these error estimates to adapt both $H$ and $\mathtt{m}$ are discussed below.

- *Balancing error strategy.* In this approach, we first use the estimated total truncation error $\hat{\varepsilon}_{n+1}$ to control the macro-step size using the traditional mechanism based on the asymptotic error behavior [30, Chapter II.4]:

$$
\hat{\varepsilon}_{n+2} \le 1 \qquad \Rightarrow \qquad H_{\mathrm{new}} = \mathrm{fac} \cdot H \cdot (\hat{\varepsilon}_{n+1})^{-\frac{1}{p}},
$$

  where $\mathrm{fac} < 1$ is a safety factor. Next, the multirate ratio $\mathtt{m}$ is adjusted such that the estimated slow and fast error components over the next step are equal to each other:

$$
\mathtt{m}_{\mathrm{new}} = \mathtt{m} \cdot \left(\hat{\varepsilon}^{\{\mathrm{F}\}}_{n+1} / \hat{\varepsilon}^{\{\mathrm{S}\}}_{n+1}\right)^{\frac{1}{p}}.
\tag{5.26}
$$

- *Optimizing efficiency strategy.* This approach focuses on reducing the overall cost of multirate integration for decoupled MR-GARK schemes. Let $t^{\{\mathrm{S}\}}$ and $t^{\{\mathrm{F}\}}$ represent the computational costs of a slow step and a fast step, respectively; these can be evaluated online by timing the slow macro-steps and the fast micro-steps during their execution. We define the computational efficiency of a step as the progress made during the step ($H$) divided by the total cost of executing step ($t^{\{S\}} + \mathtt{m}_{\rm new} t^{\{F\}}$). The new values of $H$ and $\mathtt{m}$ are selected such as to achieve the desired accuracy while maximizing the computational efficiency [53]:

$$\text{(5.27a)} \quad \mathtt{m}_{\rm new} = \arg\min_{\mathtt{m}_*} \frac{t^{\{S\}} + \mathtt{m}_* \, t^{\{F\}}}{H} \left( \hat{\varepsilon}^{\{S\}}_{n+1} + \hat{\varepsilon}^{\{F\}}_{n+1} \cdot \mathtt{m}^q / \mathtt{m}_*^q \right)^{\frac{1}{q+1}},$$

$$\text{(5.27b)} \quad H_{\rm new} = H \cdot \left( \hat{\varepsilon}^{\{S\}}_{n+1} + \hat{\varepsilon}^{\{F\}}_{n+1} \cdot \mathtt{m}^q / \mathtt{m}_{\rm new}^q \right)^{-\frac{1}{q+1}}.$$

**6. Linear Multistep Methods.** In the following we discuss the general framework for constructing multirate methods using linear multistep methods (LMMs) as base schemes. Small step sizes $h$ are used for the fast components, and large step sizes $H$ for the slow components, with an integer multirate factor $\mathtt{m} = H/h$.

**6.1. Construction of multirate linear multistep schemes.** For this, we need the following ingredients:

- The base scheme to solve (1.1) is a $k$-step linear multistep scheme, defined by its coefficients $\alpha_0 = 1, \alpha_1, \dots, \alpha_r$, and $\beta_0, \dots, \beta_r$:

$$\sum_{r=0}^{k} \alpha_r \, \mathbf{y}_{n+1-r} = H \sum_{r=0}^{k} \beta_r \, \mathbf{f}_{n+1-r}, \quad \mathbf{f}_i := \mathbf{f}(\mathbf{y}_i);$$

- Slow variables $\mathbf{y}^{\{S\}}$ are discretized on the coarse grid $t_{n-k+1}, t_{n-k+2}, \dots, t_{n+1}$ with $t_n := t_0 + nH$. Fast variables $\mathbf{y}^{\{F\}}$ are discretized on the fine grid $t_n, t_{n+1/\mathtt{m}}, t_{n+2/\mathtt{m}}, \dots, t_{n+1}$ with $t_{n+\ell/\mathtt{m}} := t_n + \ell h = t_0 + (n + \ell/\mathtt{m})\, H$.
- Assume we have computed the solutions up to time $t_n$, and have the history of the slow variables on the coarse grid, and the history of the fast variables on the fine grid. A compound multirate step to advance the solution to $t_{n+1}$ is defined as:

$$\text{(6.1a)} \quad \sum_{r=0}^{k} \alpha_r \, \mathbf{y}^{\{S\}}_{n+1-r} = H \sum_{r=1}^{k} \beta_r \, \mathbf{f}^{\{S\}}(\mathbf{y}^{\{S\}}_{n+1-r}, \mathbf{y}^{\{F\}}_{n+1-r}) + H \, \beta_0 \, \mathbf{f}^{\{S\}}(\mathbf{y}^{\{S\}}_{n+1}, \mathrm{y}^{\{F\}}_{n+1}),$$

$$\text{(6.1b)} \quad \sum_{r=0}^{k} \alpha_r \, \mathbf{y}^{\{F\}}_{n+(\ell-r)/\mathtt{m}} = h \sum_{r=0}^{k} \beta_r \, \mathbf{f}^{\{F\}}(\mathrm{y}^{\{S\}}_{n+(\ell-r)/\mathtt{m}}, \mathbf{y}^{\{F\}}_{n+(\ell-r)/\mathtt{m}}), \quad \text{for } \ell = 0, \dots, \mathtt{m}-1.$$

  Note that, in the most general formulation, different base schemes can be used for the slow discretization (6.1a) and for the $\mathtt{m}$ fast discretizations (6.1b).
- The definition of the extrapolated/interpolated values $\mathrm{y}^{\{S\}}_{n+(\ell-r)/\mathtt{m}}$ and $\mathrm{y}^{\{F\}}_{n+1}$ depends on the sequence of computations, i.e., first fast, then slow (*fastest-first strategy*) or first slow, then fast (*slowest-first strategy*).

The analysis of extrapolation/interpolation based multirate schemes in Section 3 has shown that the multirate linear multistep method (6.1) inherits the convergence order $p$ of the underlying linear multistep scheme if a) the interpolation and extrapolation schemes are at least of order $p-1$ and b) the overall multirate scheme is stable. To check 0-stability of the numerical integration scheme (6.1) we apply it to the test

problem

$$\dot{\mathbf{y}}^{\{\mathrm{S}\}} = 0, \qquad \dot{\mathbf{y}}^{\{\mathrm{F}\}} = 0.$$

For this problem the fast integration (6.1a) and the slow integration (6.1a) proceed independently. Consequently, zero stability of the general multirate scheme (6.1) is equivalent to zero-stability of the fast and the slow base multistep schemes.

Summing up, we have established the following result.

THEOREM 6.1 (Convergence of multirate linear multistep schemes). *The multirate linear multistep scheme* (6.1) *is convergent with convergence order $p$, if the base linear multistep scheme is convergent with $p$ for each individual component system, and the interpolation/extrapolation schemes used for evaluating the coupling variables $\mathrm{y}^{\{\mathrm{S}\}}_{n+(\ell-r)/\mathtt{m}}$ and $\mathrm{y}^{\{\mathrm{F}\}}_{n+1}$ are of order at least $p-1$.*

**6.2. Explicit schemes.** As for explicit multirate Euler schemes discussed in Section 4, we distinguish slowest-first and fastest-first implementations of explicit mulitrate linear multistep methods. The computation of the slow components in (6.1a) requires only past fast approximations $\mathbf{y}^{\{\mathrm{F}\}}_{n+1-\ell}$, $\ell = 1, \dots, k$, as $\beta_0 = 0$ holds for an explicit scheme. Hence $\mathbf{y}^{\{\mathrm{S}\}}_{n+1}$ can be computed first in both the slowest-first and fastest-first approach. These approaches only differ in the interpolation or extrapolation procedures for the slow components in (6.1b).

*Slowest-first approach.* One first computes $\mathbf{y}^{\{\mathrm{S}\}}_{n+1}$. The fast components are computed using (6.1b); the slow components $\mathrm{y}^{\{\mathrm{S}\}}_{n+(\ell-r)/\mathtt{m}}$ on the fine grid are obtained by interpolation from $\mathbf{y}^{\{\mathrm{S}\}}_{n+1-\ell}$, $\ell = 0, 1, \dots, k-1$.

*Fastest-first approach.* Here one assumes that $\mathbf{y}^{\{\mathrm{S}\}}_{n+1}$ is not available yet, and thus $\mathrm{y}^{\{\mathrm{S}\}}_{n+(\ell-r)/\mathtt{m}}$ can only be obtained by extrapolation based on $\mathbf{y}^{\{\mathrm{S}\}}_{n+1-\ell}$, $\ell = 1, \dots, k$.

In both approaches the analysis of Section 3 guarantees convergence order $p$ for explicit multirate linear multistep methods if the base method has order $p$ and the interpolation/extrapolation has order at last $p-1$.

EXAMPLE 5. *Consider the case with $k = 2$ and $\mathtt{m} = 3$. The compound step*

$$\begin{pmatrix} \mathbf{y}^{\{\mathrm{S}\}}_{n-1} \\ \mathbf{y}^{\{\mathrm{S}\}}_{n} \\ \mathbf{y}^{\{\mathrm{F}\}}_{n-1} \\ \mathbf{y}^{\{\mathrm{F}\}}_{n-1+1/3} \\ \mathbf{y}^{\{\mathrm{F}\}}_{n-1+2/3} \\ \mathbf{y}^{\{\mathrm{F}\}}_{n} \end{pmatrix} \mapsto \begin{pmatrix} \mathbf{y}^{\{\mathrm{S}\}}_{n} \\ \mathbf{y}^{\{\mathrm{S}\}}_{n+1} \\ \mathbf{y}^{\{\mathrm{F}\}}_{n} \\ \mathbf{y}^{\{\mathrm{F}\}}_{n+1/3} \\ \mathbf{y}^{\{\mathrm{F}\}}_{n+2/3} \\ \mathbf{y}^{\{\mathrm{F}\}}_{n+1} \end{pmatrix}$$

*is defined by*

$$\begin{aligned}
\mathbf{y}^{\{\mathrm{S}\}}_{n+1} + \alpha_1\, \mathbf{y}^{\{\mathrm{S}\}}_{n} + \alpha_2\, \mathbf{y}^{\{\mathrm{S}\}}_{n-1} &= H\left(\beta_1\, \mathbf{f}^{\{\mathrm{S}\}}(\mathbf{y}^{\{\mathrm{S}\}}_{n}, \mathbf{y}^{\{\mathrm{F}\}}_{n}) + \beta_2\, \mathbf{f}^{\{\mathrm{S}\}}(\mathbf{y}^{\{\mathrm{S}\}}_{n-1}, \mathbf{y}^{\{\mathrm{F}\}}_{n-1})\right), \\
\mathbf{y}^{\{\mathrm{F}\}}_{n+1/3} + \alpha_1\, \mathbf{y}^{\{\mathrm{F}\}}_{n} + \alpha_2\, \mathbf{y}^{\{\mathrm{F}\}}_{n-1/3} &= h\left(\beta_1\, \mathbf{f}^{\{\mathrm{F}\}}(\mathbf{y}^{\{\mathrm{S}\}}_{n}, \mathbf{y}^{\{\mathrm{F}\}}_{n}) + \beta_2\, \mathbf{f}^{\{\mathrm{F}\}}(\mathrm{y}^{\{\mathrm{S}\}}_{n-1/3}, \mathbf{y}^{\{\mathrm{F}\}}_{n-1/3})\right), \\
\mathbf{y}^{\{\mathrm{F}\}}_{n+2/3} + \alpha_1\, \mathbf{y}^{\{\mathrm{F}\}}_{n+1/3} + \alpha_2\, \mathbf{y}^{\{\mathrm{F}\}}_{n} &= h\left(\beta_1\, \mathbf{f}^{\{\mathrm{F}\}}(\mathrm{y}^{\{\mathrm{S}\}}_{n+1/3}, \mathbf{y}^{\{\mathrm{F}\}}_{n+1/3}) + \beta_2\, \mathbf{f}^{\{\mathrm{F}\}}(\mathbf{y}^{\{\mathrm{S}\}}_{n}, \mathbf{y}^{\{\mathrm{F}\}}_{n})\right), \\
\mathbf{y}^{\{\mathrm{F}\}}_{n+1} + \alpha_1\, \mathbf{y}^{\{\mathrm{F}\}}_{n+2/3} + \alpha_2\, \mathbf{y}^{\{\mathrm{F}\}}_{n+1/3} &= h\left(\beta_1\, \mathbf{f}^{\{\mathrm{F}\}}(\mathrm{y}^{\{\mathrm{S}\}}_{n+2/3}, \mathbf{y}^{\{\mathrm{F}\}}_{n+2/3}) + \beta_2\, \mathbf{f}^{\{\mathrm{F}\}}(\mathrm{y}^{\{\mathrm{S}\}}_{n+1/3}, \mathbf{y}^{\{\mathrm{F}\}}_{n+1/3})\right),
\end{aligned}$$

*in the case of an explicit base linear multistep method with $\alpha_0 = 1$ and $\beta_0 = 0$. This case is illustrated in Fig. 10.*

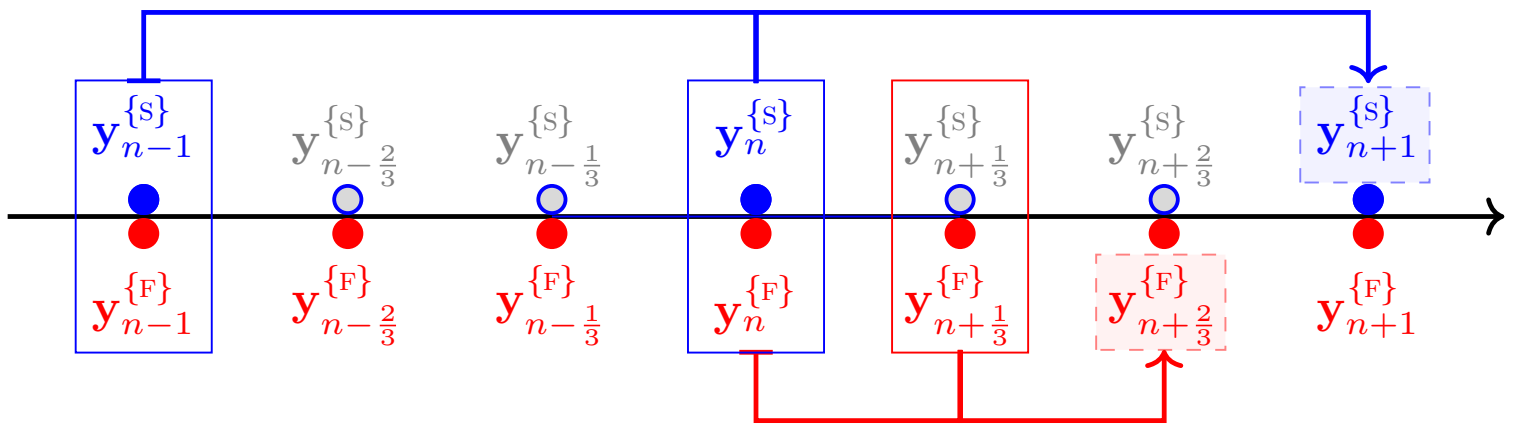

Fig. 10: Cartoon of explicit multirate linear multistep method (MR-LMM) for $k = 2$ steps and $\mathtt{m} = 3$. The values computed by MR-LMM steps are represented by filled discs. Empty circles represent the slow values not computed by MR-LMM steps, and that need to be approximated. The computation of $\mathbf{y}^{\{\mathrm{S}\}}_{n+1}$ uses the known values $(\mathbf{y}^{\{\mathrm{S}\}}_{n}, \mathbf{y}^{\{\mathrm{F}\}}_{n})$ and $(\mathbf{y}^{\{\mathrm{S}\}}_{n-1}, \mathbf{y}^{\{\mathrm{F}\}}_{n-1})$. The computation of $\mathbf{y}^{\{\mathrm{F}\}}_{n+\frac{2}{3}}$ uses the known values $(\mathbf{y}^{\{\mathrm{S}\}}_{n}, \mathbf{y}^{\{\mathrm{F}\}}_{n})$, as well as $(\mathbf{y}^{\{\mathrm{S}\}}_{n+\frac{1}{3}}, \mathbf{y}^{\{\mathrm{F}\}}_{n+\frac{1}{3}})$. While $\mathbf{y}^{\{\mathrm{F}\}}_{n+\frac{1}{3}}$ has been computed by the previous small fast step, $\mathbf{y}^{\{\mathrm{S}\}}_{n+\frac{1}{3}}$ is not computed by the slow method and needs to be approximated.

For implicit schemes the slow-fast coupling is more complex. As for implicit multirate Euler methods discussed in Section 4, we consider semi-implicit formulations, implicit formulations, and compound step formulations implicit MR-LMMs. For each case, we distinguish between the fastest-first and slowest-first first strategies.

### 6.3. Decoupled implicit formulations.

**Decoupled fastest-first strategy.** We first solve the fast part (6.1b) with $\mathtt{m}$ steps of step size $h$ to obtain $\mathbf{y}^{\{\mathrm{F}\}}_{n+1/\mathtt{m}}, \dots, \mathbf{y}^{\{\mathrm{F}\}}_{n+1}$:

$$\text{(6.2a)} \qquad \sum_{r=0}^{k} \alpha_r\, \mathbf{y}^{\{\mathrm{F}\}}_{n+(\ell-r)/\mathtt{m}} = h \sum_{r=0}^{k} \beta_r\, \mathbf{f}^{\{\mathrm{F}\}}\big(\mathbf{P}^{\{\mathrm{S}\}}_{\mathrm{ext}}(t_{n+(\ell-r)/\mathtt{m}}), \mathbf{y}^{\{\mathrm{F}\}}_{n+(\ell-r)/\mathtt{m}}\big), \quad \ell = 1, \dots, \mathtt{m}.$$

The slow components on the fine grid are obtained by evaluating the extrapolation polynomial $\mathbf{P}^{\{\mathrm{S}\}}_{\mathrm{ext}}$ at time point $t_{n+(\ell-r)/\mathtt{m}}$. Using Lagrangian bases $\ell_q(t)$ defined by $\ell_q(t_{n+1-r}) = \delta_{q,r}, r = 1, \dots, k$, the extrapolation polynomial reads [64]:

$$\mathbf{P}^{\{\mathrm{S}\}}_{\mathrm{ext}}(t) = \sum_{q=1}^{k} \ell_q(t) \cdot \mathbf{y}^{\{\mathrm{S}\}}_{n+1-q},$$

and provides approximations of order $k-1$.

Next, we compute the slow part $\mathbf{y}^{\{\mathrm{S}\}}_{n+1}$ by integrating (6.1a) from $t_n$ to $t_{n+1}$ with one large step $H = \mathtt{m}h$:

$$\text{(6.2b)} \qquad \sum_{r=0}^{k} \alpha_r\, \mathbf{y}^{\{\mathrm{S}\}}_{n+1-r} = H \sum_{r=0}^{k} \beta_r\, \mathbf{f}^{\{\mathrm{S}\}}\big(\mathbf{y}^{\{\mathrm{S}\}}_{n+1-r}, \mathbf{y}^{\{\mathrm{F}\}}_{n+1-r}\big).$$

The fast components on the coarse grid $\mathbf{y}^{\{F\}}_{n+1-k}, \mathbf{y}^{\{F\}}_{n+2-k}, \dots, \mathbf{y}^{\{F\}}_{n+1}$ have already been computed and do not require additional interpolations. If the basis $k$-step scheme has at least order $k$, the multirate version is convergent of order $k$.

**Decoupled slowest-first strategy.** One first solves the slow part (6.1a) with one step of size $H = \mathtt{m}h$ to obtain $\mathbf{y}^{\{S\}}_{n+1}$:

$$
\begin{aligned}
\sum_{r=0}^{k} \alpha_r \, \mathbf{y}^{\{S\}}_{n+1-r} \;&=\; H \sum_{r=1}^{k} \beta_r \, \mathbf{f}^{\{S\}}(\mathbf{y}^{\{S\}}_{n+1-r}, \mathbf{y}^{\{F\}}_{n+1-r}) \\
&\quad + H\, \beta_0 \, \mathbf{f}^{\{S\}}(\mathbf{y}^{\{S\}}_{n+1}, \mathbf{P}^{\{F\}}_{\mathbf{ext}}(t_{n+1})).
\end{aligned} \tag{6.3a}
$$

The fast components on the past coarse grid points $\mathbf{y}^{\{F\}}_{n-k+1}, \dots, \mathbf{y}^{\{F\}}_{n}$ have already been computed. The fast component on the next coarse point $\mathbf{y}^{\{F\}}_{n+1}$ is approximated by extrapolating the interpolation polynomial $\mathbf{P}^{\{F\}}_{\mathbf{ext}}(t)$ through the points $(t_{n-(k-1)}, \mathbf{y}^{\{F\}}_{n-(k-1)}), \dots, (t_n, \mathbf{y}^{\{F\}}_{n})$ at time point $t_{n+1}$:

$$
\mathbf{P}^{\{F\}}_{\mathbf{ext}}(t_{n+1}) = \sum_{r=1}^{k} \ell_r(t_{n+1}) \, \mathbf{y}^{\{F\}}_{n-r+1},
$$

which can be written with some $h$-dependent coefficients $\hat{\alpha}_r = \ell_r(t_{n+1})$ as a linear combination of the interpolation points [64]. The approximation is again of order $k-1$.

Next, we compute the fast part (6.1b) with $\mathtt{m}$ steps of step size $h$ to obtain $\mathbf{y}^{\{F\}}_{n+1/\mathtt{m}}, \dots, \mathbf{y}^{\{F\}}_{n+1}$.

$$
\begin{aligned}
\sum_{r=0}^{k} \alpha_r \, \mathbf{y}^{\{F\}}_{n+(\ell-r)/\mathtt{m}} \;&=\; h \sum_{r=0}^{k} \beta_r \, \mathbf{f}^{\{F\}}(\mathbf{P}^{\{S\}}_{\mathbf{int}}(t_{n+(\ell-r)/\mathtt{m}}), \mathbf{y}^{\{F\}}_{n+(\ell-r)/\mathtt{m}}), \\
&\quad \ell = 0, \dots, \mathtt{m}.
\end{aligned} \tag{6.3b}
$$

For this we need, as in the fastest-first strategy, the slow solutions on the fine grid $\mathbf{y}^{\{S\}}_{n+\ell/\mathtt{m}}$, $\ell = -k, -k+1, \dots, \mathtt{m}$, which are obtained by evaluating the interpolation polynomial $\mathbf{P}^{\{S\}}_{\mathbf{int}}(t)$ through the points $(t_{n-k+2}, \mathbf{y}^{\{S\}}_{n-k+2}), \dots, (t_{n+1}, \mathbf{y}^{\{S\}}_{n+1})$ at each time point $t_{n-(\ell-r)/\mathtt{m}}$:

$$
\mathbf{y}^{\{S\}}_{n-(\ell-r)/\mathtt{m}} := \mathbf{P}^{\{S\}}_{\mathbf{int}}(t_{n-(\ell-r)/\mathtt{m}}) = \sum_{q=0}^{k-1} \ell_q(t_{n-(\ell-r)/\mathtt{m}}) \mathbf{y}^{\{S\}}_{n+1-q},
$$

where we have used again the representation of the extrapolated values by the Lagrangian basis polynomials $\ell_q$ defined by $\ell_q(t_{n+1-m}) = \delta_{q,\ell}$, $m = 0, \dots, k-1$. Consequently, $\mathbf{y}^{\{S\}}_{n-(k-r)/\mathtt{m}}$ can be written with some $h$-dependent coefficients as a linear combination of the interpolation points [64]. The approximation is of order $k-1$. Again, if the basis $k$-step scheme has at least order $k$, the multirate version will be convergent of order $k$.

The slowest-first strategy was introduced by Gear [21]. In contrast to the fastest-first strategy (6.2), the slowest-first strategy (6.3) allows for a second update of $\mathbf{y}^{\{S\}}_{n+1}$ using (6.3a) with $\mathbf{y}^{\{F\}}_{n+1}$ replacing $\mathbf{P}^{\{F\}}_{\mathbf{ext}}(t_{n+1})$ [64].

**6.4. Coupled implicit formulations.** Both fastest and slowest-first strategies discussed above use a semi-implicit formulation of the joint step, i.e., the steps computing $\mathbf{y}^{\{S\}}_{n+1}$ and $\mathbf{y}^{\{F\}}_{n+1}$ are done sequentially: either first $\mathbf{y}^{\{F\}}_{n+1}$, and then $\mathbf{y}^{\{S\}}_{n+1}$ (fastest-first strategy), or vice-versa (slowest-first strategy). One may also compute both approximations together in a coupled manner, defining a large nonlinear system in both $\mathbf{y}^{\{S\}}_{n+1}$ and $\mathbf{y}^{\{F\}}_{n+1}$.

**Coupled fastest-first strategy.** The sequential computations (6.2) are replaced by

$$
\begin{aligned}
&\text{(6.4a)} \quad \sum_{r=0}^{k} \alpha_r\, \mathbf{y}^{\{F\}}_{n+(\ell-r)/\mathtt{m}} = h \sum_{r=0}^{k} \beta_r\, \mathbf{f}^{\{F\}}\big(\mathbf{P}^{\{S\}}_{\mathrm{ext}}(t_{n+(\ell-r)/\mathtt{m}}), \mathbf{y}^{\{F\}}_{n+(\ell-r)/\mathtt{m}}\big), \\
&\qquad\qquad \ell = 1, \dots, \mathtt{m}-1, \\
&\sum_{r=0}^{k} \alpha_r\, \mathbf{y}^{\{F\}}_{n+1-r/\mathtt{m}} = h \sum_{r=1}^{k} \beta_r\, \mathbf{f}^{\{F\}}\big(\mathbf{P}^{\{S\}}_{\mathrm{ext}}(t_{n+1-r/\mathtt{m}}), \mathbf{y}^{\{F\}}_{n+1-r/\mathtt{m}}\big) \\
&\text{(6.4b)} \qquad + h\, \beta_0\, \mathbf{f}^{\{F\}}\big(\mathbf{y}^{\{S\}}_{n+1}, \mathbf{y}^{\{F\}}_{n+1}\big), \\
&\sum_{r=0}^{k} \alpha_r\, \mathbf{y}^{\{S\}}_{n+1-r} = H \sum_{r=1}^{k} \beta_r\, \mathbf{f}^{\{S\}}\big(\mathbf{y}^{\{S\}}_{n+1-r}, \mathbf{y}^{\{F\}}_{n+1-r}\big) \\
&\text{(6.4c)} \qquad + H\, \beta_0\, \mathbf{f}^{\{S\}}\big(\mathbf{y}^{\{S\}}_{n+1}, \mathbf{y}^{\{F\}}_{n+1}\big),
\end{aligned}
$$

where (6.4b)–(6.4c) couple the last fast step with the slow step and define $\mathbf{y}^{\{S\}}_{n+1}$ and $\mathbf{y}^{\{F\}}_{n+1}$ simultaneously.

**Coupled slowest-first strategy.** The sequential computations (6.3) are replaced by

$$
\begin{aligned}
&\text{(6.5a)} \sum_{r=0}^{k} \alpha_r\, \mathbf{y}^{\{S\}}_{n+1-r} = H \sum_{r=1}^{k} \beta_r\, \mathbf{f}^{\{S\}}\big(\mathbf{y}^{\{S\}}_{n+1-r}, \mathbf{y}^{\{F\}}_{n+1-r}\big) + H\, \beta_0\, \mathbf{f}^{\{S\}}\big(\mathbf{y}^{\{S\}}_{n+1}, \mathbf{y}^{\{F\}}_{n+1}\big), \\
&\sum_{r=0}^{k} \alpha_r\, \mathbf{y}^{\{F\}}_{n+(\ell-r)/\mathtt{m}} = h \sum_{r=0}^{k} \beta_r\, \mathbf{f}^{\{F\}}\big(\mathbf{P}^{\{S\}}_{\mathrm{int}}(t_{n+(\ell-r)/\mathtt{m}}), \mathbf{y}^{\{F\}}_{n+(\ell-r)/\mathtt{m}}\big), \\
&\text{(6.5b)} \qquad \ell = 0, \dots, \mathtt{m}-1, \\
&\sum_{r=0}^{k} \alpha_r\, \mathbf{y}^{\{F\}}_{n+1-r/\mathtt{m}} = h \sum_{r=1}^{k} \beta_r\, \mathbf{f}^{\{F\}}\big(\mathbf{P}^{\{S\}}_{\mathrm{int}}(t_{n+1-r/\mathtt{m}}), \mathbf{y}^{\{F\}}_{n+1-r/\mathtt{m}}\big) + \\
&\text{(6.5c)} \qquad + h\, \beta_0\, \mathbf{f}^{\{F\}}\big(\mathbf{y}^{\{S\}}_{n+1}, \mathbf{y}^{\{F\}}_{n+1}\big),
\end{aligned}
$$

where (6.5a) and (6.5c) couple the last fast step with the slow step, defining $\mathbf{y}^{\{S\}}_{n+1}$ and $\mathbf{y}^{\{F\}}_{n+1}$ simultaneously. Note that the interpolated values $\mathbf{P}^{\{S\}}_{\mathrm{int}}$ depend on $\mathbf{y}^{\{S\}}_{n+1}$. Such a strategy was first introduced by Skelboe [64], and is the basis for dynamic iteration schemes [47]: one can replace the interpolation with an extrapolation $\mathbf{P}^{\{S\}}_{\mathrm{ext}}$ to obtain a first approximation of $\mathbf{y}^{\{S\}}_{n+1}$, then (6.5) with $\mathbf{P}^{\{S\}}_{\mathrm{int}}$ are iterated to obtain increasingly accurate solutions $\mathbf{y}^{\{S\}}_{n+1}$ and $\mathbf{y}^{\{F\}}_{n+1}$.

**6.5. Compound-step approach.** All approaches discussed so far – slowest-first and fastest-first, decoupled and coupled formulations – use extrapolation, which may cause some instability. One idea to overcome this problem and to avoid extrapolation was first introduced by Verhoeven [25]: first, in a compound predictor step, one computes at the same time $\mathbf{y}_{n+1}^{\{\mathrm{S}\}}$ and $\mathbf{y}_{n+1}^{\{\mathrm{F}\}}$ with step size $H$:

$$\text{(6.6)} \quad \sum_{r=0}^{k} \alpha_r \, \mathbf{y}_{n+1-r}^{\{\mathrm{S}\}} = H \sum_{r=1}^{k} \beta_r \, \mathbf{f}^{\{\mathrm{S}\}}(\mathbf{y}_{n+1-r}^{\{\mathrm{S}\}}, \mathbf{y}_{n+1-r}^{\{\mathrm{F}\}}) + H\, \beta_0 \, \mathbf{f}^{\{\mathrm{S}\}}(\mathbf{y}_{n+1}^{\{\mathrm{S}\}}, \mathbf{y}_{n+1}^{\{\mathrm{F}\}}),$$

$$\text{(6.7)} \quad \sum_{r=0}^{k} \alpha_r \, \mathbf{y}_{n+1-r}^{\{\mathrm{F}\}} = H \sum_{r=1}^{k} \beta_r \, \mathbf{f}^{\{\mathrm{F}\}}(\mathbf{y}_{n+1-r}^{\{\mathrm{S}\}}, \mathbf{y}_{n+1-r}^{\{\mathrm{F}\}}) + H\, \beta_0 \, \mathbf{f}^{\{\mathrm{F}\}}(\mathbf{y}_{n+1}^{\{\mathrm{S}\}}, \mathbf{y}_{n+1}^{\{\mathrm{F}\}}).$$

Improved solutions are obtained by a corrector approach, that applies the coupled or decoupled implicit strategy with modified polynomials $\mathbf{P}_{\mathbf{ext}}^{\{\mathrm{S}\}}$ and $\mathbf{P}_{\mathbf{ext}}^{\{\mathrm{F}\}}$:

- *Fastest-first approach (semi-implict and implicit)*: The slow extrapolation polynomial $\mathbf{P}_{\mathbf{ext}}^{\{\mathrm{S}\}}$ in (6.2a), (6.4a), (6.4b) is replaced by an interpolation polynomial that is constructed using the predicted slow values.
- *Slowest-first approach (semi-implicit and implicit)*: In (6.3a) the evaluation of the fast interpolation polynomial $\mathbf{P}_{\mathbf{ext}}^{\{\mathrm{F}\}}(t_{n+1})$ is replaced by the predicted value $\mathbf{y}_{n+1}^{\{\mathrm{F}\}}$. In (6.5a) the coupled $\mathbf{y}_{n+1}^{\{\mathrm{F}\}}$ is replaced by the predicted fast value.

One notes that $\mathbf{y}_{n+1}^{\{\mathrm{S}\}}$ and $\mathbf{y}_{n+1}^{\{\mathrm{F}\}}$ have to be computed twice: the first computation of $\mathbf{y}_{n+1}^{\{\mathrm{S}\}}$ and $\mathbf{y}_{n+1}^{\{\mathrm{F}\}}$ is used to replace extrapolation by interpolation. In the slowest-first approach, however, one may not change the predicted $\mathbf{y}_{n+1}^{\{\mathrm{S}\}}$, and compute the slow part only once.

One criticism of this compound step is to use the macro step size $H$ also for the fast component $\mathbf{y}^{\{\mathrm{F}\}}$. One may overcome this problem by combining one slow step of step size $H$ with one fast step of step size $h$ in the following generalized compound step, computing $\mathbf{y}_{n+1}^{\{\mathrm{S}\}}$ and $\mathbf{y}_{n+1/\mathtt{m}}^{\{\mathrm{F}\}}$ together:

$$\sum_{r=0}^{k} \alpha_r \, \mathbf{y}_{n+1-r}^{\{\mathrm{S}\}} = H \sum_{r=1}^{k} \beta_r \, \mathbf{f}^{\{\mathrm{S}\}}(\mathbf{y}_{n+1-r}^{\{\mathrm{S}\}}, \mathbf{y}_{n+1-r}^{\{\mathrm{F}\}}) + H\, \beta_0 \, \mathbf{f}^{\{\mathrm{S}\}}(\mathbf{y}_{n+1}^{\{\mathrm{S}\}}, \mathbf{P}_{\mathbf{ext}}^{\{\mathrm{F}\}}(t_{n+1})),$$

$$\sum_{r=0}^{k} \alpha_r \, \mathbf{y}_{n+(1-r)/\mathtt{m}}^{\{\mathrm{F}\}} = h \sum_{r=0}^{k} \beta_r \, \mathbf{f}^{\{\mathrm{F}\}}(\mathbf{P}_{\mathbf{ext}}^{\{\mathrm{S}\}}(t_{n+(1-r)/\mathtt{m}}), \mathbf{y}_{n+(1-r)/\mathtt{m}}^{\{\mathrm{F}\}}).$$

The polynomial $\mathbf{P}_{\mathbf{int}}^{\{\mathrm{S}\}}$ is based on the points $(t_{n-k+1}, \mathbf{y}_{n-k+1}^{\{\mathrm{S}\}}), \ldots, (t_n, \mathbf{y}_n^{\{\mathrm{S}\}})$, and thus $\mathbf{P}_{\mathbf{ext}}^{\{\mathrm{S}\}}(\mathbf{t}_{\mathbf{n}+(\mathbf{1}-\mathbf{r})/\mathtt{m}})$ is interpolating instead of extrapolating for $r > 1$. For the slow combined step, however, $\mathbf{P}_{\mathbf{ext}}^{\{\mathrm{F}\}}(t_{n+1})$ based on the points $(t_{n-k+1}, \mathbf{y}_{n-k+1}^{\{\mathrm{F}\}}), \ldots, (t_n, \mathbf{y}_n^{\{\mathrm{F}\}})$ is extrapolating instead of interpolating.

**6.6. Comparison of different MR-LMM strategies.** At a first glance, the fastest-first strategy seems to be advantageous. The extrapolation error introduced in (6.2a) necessary to compute the fast part is acceptable, as one is extrapolating over an interval smaller than the step size $H$ tailored to the activity level of the slow part. In the slowest-first strategy, however, we are extrapolating the fast variables over many step sizes in (6.3a) to compute the slow variables, which may only be tolerable if the coupling of the fast into the slow part, measured by $\|\partial \mathbf{f}^{\{\mathrm{S}\}} / \partial \mathbf{y}^{\{\mathrm{F}\}}\|$, is small. As a rule of thumb, this quantity will be usually small, as otherwise a high level of activity would be transferred from the fast into the slow part, and the slow part would not be slow anymore. This is, however, only a rule of thumb, and there exist academic counterexamples [22].

However, the slowest-first strategy becomes advantageous, if we leave the case of constant step sizes $h$ and adapt the step sizes according to the step size predictions gained from an error control. If the last integration step (6.3b) from $t$ to $t+H$ fails in the fastest-first strategy, the macro step size $H$ has to be lowered and must be repeated with a smaller step size. However, the values of $\mathbf{y}^{\{F\}}$ before $t+h$ are no longer available if back-up costs should be kept small, and the value of $\mathbf{y}^{\{F\}}$ at $t+H$ is not reliable, as it is based on predicted values of $\mathbf{y}^{\{S\}}$ which turned out not to be sufficiently accurate. In contrast, the slowest-first strategy does not introduce any problems: if an integration of the fast part fails, one only has to repeat it with a smaller step size. The necessary information to interpolate the slow part does not change. See [21, 22] for a detailed discussion on the pros and cons of both strategies.

As the implicit strategy, the compound step has to solve a nonlinear system of equations of dimension $2n$, but avoids extrapolation. However, one criticism of this compound step is to use the macro step size $H$ also for the fast component $\mathbf{y}^{\{F\}}$.

The implicit strategy is the workhorse for dynamic iteration schemes. It defines a first approximation of $\mathbf{y}^{\{S\}}, \mathbf{y}^{\{F\}}$ in $[t_n, t_{n+1}]$, which is iteratively improved by recomputing $\mathbf{y}^{\{S\}}, \mathbf{y}^{\{F\}}$ in $[t_n, t_{n+1}]$ until convergence is reached [64].

**6.7. A first discussion of A-stability.** Following the analysis done in the Euler chapter, we investigate stability by applying the methods to the linear test problem (4.3).

Before considering the general case, we start with an example.

Example 6 (BDF-2: Implicit formulation. slowest-first approach, linear interpolation, $\mathtt{m}=2$). *We have*

$$\begin{aligned}
\mathbf{y}^{\{S\}}_{n+1} &= \frac{4}{3}\mathbf{y}^{\{S\}}_n - \frac{1}{3}\mathbf{y}^{\{S\}}_{n-1} + \frac{2}{3}H\mathbf{f}^{\{S\}}(\mathbf{y}^{\{S\}}_{n+1}, \mathbf{y}^{\{F\}}_{n+1}),\\
\mathbf{y}^{\{F\}}_{n+1} &= \frac{4}{3}\mathbf{y}^{\{F\}}_{n+1/2} - \frac{1}{3}\mathbf{y}^{\{F\}}_n + \frac{2}{3}h\mathbf{f}^{\{F\}}(\mathbf{y}^{\{S\}}_{n+1}, \mathbf{y}^{\{F\}}_{n+1}),\\
\mathbf{y}^{\{F\}}_{n+1/2} &= \frac{4}{3}\mathbf{y}^{\{F\}}_n - \frac{1}{3}\mathbf{y}^{\{F\}}_{n-1/2} + \frac{2}{3}h\mathbf{f}^{\{F\}}(\frac{1}{2}\mathbf{y}^{\{S\}}_{n+1} + \frac{1}{2}\mathbf{y}^{\{S\}}_n, \mathbf{y}^{\{F\}}_{n+1/2}),
\end{aligned}$$

*which yields applied to the test equation* (4.3)

$$\underbrace{\begin{bmatrix} a & 0 & b & 0\\ 0 & 1 & 0 & 0\\ -2c & 0 & d & e\\ -c & 0 & 0 & d \end{bmatrix}}_{A:=} \begin{bmatrix} \mathbf{y}^{\{S\}}_{n+1}\\ \mathbf{y}^{\{S\}}_n\\ \mathbf{y}^{\{F\}}_{n+1}\\ \mathbf{y}^{\{F\}}_{n+1/2} \end{bmatrix} = \underbrace{\begin{bmatrix} \frac{4}{3} & -\frac{1}{3} & 0 & 0\\ 1 & 0 & 0 & 0\\ 0 & 0 & -\frac{1}{3} & 0\\ c & 0 & \frac{4}{3} & -\frac{1}{3} \end{bmatrix}}_{B:=} \begin{bmatrix} \mathbf{y}^{\{S\}}_n\\ \mathbf{y}^{\{S\}}_{n-1}\\ \mathbf{y}^{\{F\}}_n\\ \mathbf{y}^{\{F\}}_{n-1/2} \end{bmatrix}.$$

*using the abbreviations* $a = 1 - \frac{2}{3}z^{\{S,S\}}, b = -\frac{2}{3}z^{\{F,S\}}, c = \frac{1}{3}\frac{z^{\{S,F\}}}{2}, d = 1 - \frac{2}{3}\frac{z^{\{F,F\}}}{\mathtt{m}}$ *and* $e = -\frac{4}{3}$. *The method is A-stable, if the spectral radius of* $A^{-1}B < 1$ *holds. With* $f = 1/(d(d + \frac{2cb}{a}) - \frac{bce}{a})$ *we get*

$$A^{-1} = f\begin{bmatrix} \frac{1}{af} - \frac{b}{a}\left(\frac{2cd-ce}{a}\right) & 0 & -\frac{bd}{a} & \frac{be}{a}\\ 0 & \frac{1}{f} & 0 & 0\\ \frac{2cd-ce}{a} & 0 & d & -e\\ \frac{cd}{a} & 0 & -\frac{cb}{a} & d + \frac{2bc}{a} \end{bmatrix}$$

*and $A^{-1}B$ is given by*

$$f\begin{bmatrix} \frac{4}{3}\left(\frac{1}{af}-\frac{b}{a}\left(\frac{2cd-ce}{a}\right)\right)+\frac{bce}{a} & -\frac{1}{3}\left(\frac{1}{af}-\frac{b}{a}\left(\frac{2cd-ce}{a}\right)\right) & \frac{1}{3}\frac{bd}{a}+\frac{4}{3}\frac{be}{a} & -\frac{1}{3}\frac{be}{a} \\ \frac{1}{f} & 0 & 0 & 0 \\ \frac{4}{3}\frac{2cd-ce}{a}-ce & -\frac{1}{3}\frac{2cd-ce}{a} & -\frac{1}{3}d-\frac{4}{3}e & \frac{1}{3}e \\ \frac{4}{3}\frac{cd}{a}+c(d+\frac{2bc}{a}) & -\frac{1}{3}\frac{cd}{a} & \frac{1}{3}\frac{cb}{a}+\frac{4}{3}(d+\frac{2bc}{a}) & -\frac{1}{3}(d+\frac{2bc}{a}) \end{bmatrix}$$

*For $|z^{\{\mathrm{S,F}\}}|, |z^{\{\mathrm{F,S}\}}|$ sufficiently small, the method is unconditionally stable, as the spectral radius of*

$$\begin{bmatrix} \frac{4}{3a} & -\frac{1}{3a} & 0 & 0 \\ 1 & 0 & 0 & 0 \\ 0 & 0 & -\frac{\frac{1}{3}d+\frac{4}{3}e}{d^2} & \frac{1}{3}\frac{e}{d^2} \\ 0 & 0 & \frac{4}{3}\frac{1}{d} & -\frac{1}{3}\frac{1}{d} \end{bmatrix}$$

*is bounded by one:*

$$\max\left\{\frac{z^{\{\mathrm{F,F}\}}+5\pm 4\sqrt{z^{\{\mathrm{F,F}\}}+1}}{{z^{\{\mathrm{F,F}\}}}^2-6z^{\{\mathrm{F,F}\}}+9}, \frac{-2\pm\sqrt{z^{\{\mathrm{S,S}\}}+1}}{-3+2z^{\{\mathrm{S,S}\}}}\right\}\leq 1.$$

*For larger values of $|z^{\{\mathrm{S,F}\}}|, |z^{\{\mathrm{F,S}\}}|$ we see that the scheme is only stable in the wedges given by $|\arg(z^{\{\mathrm{F,F}\}})-\pi|\leq\theta$ and $|\arg(z^{\{\mathrm{S,S}\}})-\pi|\leq\theta$ with $\theta\approx 87.6°$. This is depicted in Fig. 11. Here we have set $|z^{\{\mathrm{S,S}\}}/z^{\{\mathrm{F,F}\}}| = \mathtt{m} = 2$ according to the multirate factor. The couplings $z^{\{\mathrm{S,F}\}}$ and $z^{\{\mathrm{F,S}\}}$ are chosen such that the eigenvalues of the test equation (4.3) are a convex combination of $z^{\{\mathrm{S,S}\}}$ and $z^{\{\mathrm{F,F}\}}$, and thus (4.3) is stable. For $\xi = 0$ or $\xi = 1$ we get a one-sided coupling and thus unconditional stability, as shown above. See 5.5 for more details.*

The generic multirate schemes (6.1) can be reformulated as

$$\begin{aligned} \mathbf{y}^{\{\mathrm{S}\}}_{n+1} &= \sum_{r=1}^{k} -\alpha_r \mathbf{y}^{\{\mathrm{S}\}}_{n+1-r} + H\sum_{r=1}^{k}\beta_r \underbrace{\mathbf{f}^{\{\mathrm{S}\}}(\mathbf{y}^{\{\mathrm{S}\}}_{n+1-r}, \mathbf{y}^{\{\mathrm{F}\}}_{n+1-r})}_{\mathbf{f}^{\{\mathrm{S}\}}_{n+1-r}} + \\ &\quad + H\,\beta_0 \underbrace{\mathbf{f}^{\{\mathrm{S}\}}(\mathbf{y}^{\{\mathrm{S}\}}_{n+1}, \tilde{\mathbf{y}}^{\{\mathrm{F}\}}_{n+1})}_{=:\tilde{\mathbf{f}}^{\{\mathrm{S}\}}_{n+1}}, \\ \mathbf{y}^{\{\mathrm{F}\}}_{n+\ell/\mathtt{m}} &= \sum_{r=1}^{k}\alpha_r \mathbf{y}^{\{\mathrm{F}\}}_{n+(\ell-r)/\mathtt{m}} + h\sum_{r=0}^{k}\beta_r \underbrace{\mathbf{f}^{\{\mathrm{F}\}}(\tilde{\mathbf{y}}^{\{\mathrm{S}\}}_{n+(\ell-r)/\mathtt{m}}, \mathbf{y}^{\{\mathrm{F}\}}_{n+(\ell-r)/\mathtt{m}})}_{\tilde{\mathbf{f}}^{\{\mathrm{F}\}}_{n+(\ell-r)/\mathtt{m}}}, \\ &\quad \ell = 1,\dots,\mathtt{m}. \end{aligned}$$

To derive a compact notation of this scheme, we use the notation introduced by [19]:

$$\begin{aligned} \mathbf{Y}_{\mathrm{S},n+1} &:= (\mathbf{y}^{\{\mathrm{S}\}}_{n+1}, H\,\mathbf{f}^{\{\mathrm{S}\}}_{n+1},\dots,\mathbf{y}^{\{\mathrm{S}\}}_{n-k+2}, H\,\mathbf{f}^{\{\mathrm{S}\}}_{n-k+2})^\top, \\ \mathbf{Y}_{\mathrm{F},\ell+1} &:= (\mathbf{y}^{\{\mathrm{F}\}}_{\ell+1}, h\tilde{\mathbf{f}}^{\{\mathrm{F}\}}_{\ell+1},\dots,\mathbf{y}^{\{\mathrm{F}\}}_{\ell+1-(k-1)/\mathtt{m}}, h\tilde{\mathbf{f}}^{\{\mathrm{F}\}}_{\ell+1-(k-1)/\mathtt{m}})^\top \\ \mathbf{a}^\top &:= (\alpha_1,\beta_1,\dots,\alpha_k,\beta_k), \\ \mathbf{\Phi} &:= (\mathbf{a}^\top, \mathbf{0}, \mathbf{e}_1^\top,\dots,\mathbf{e}_{2k-2}^\top)^\top, \\ \mathbf{\Gamma} &:= \beta_0\mathbf{e}_1 + \mathbf{e}_2, \end{aligned}$$

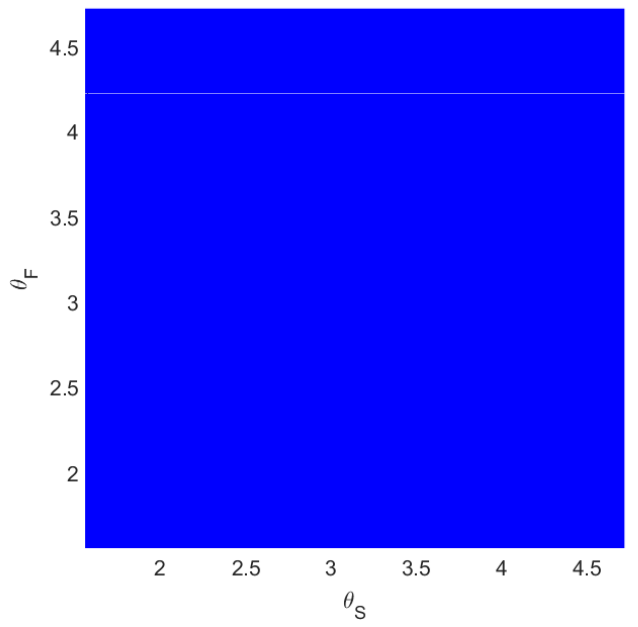

(a) Stability region for $\xi = 0, \alpha = 1$.

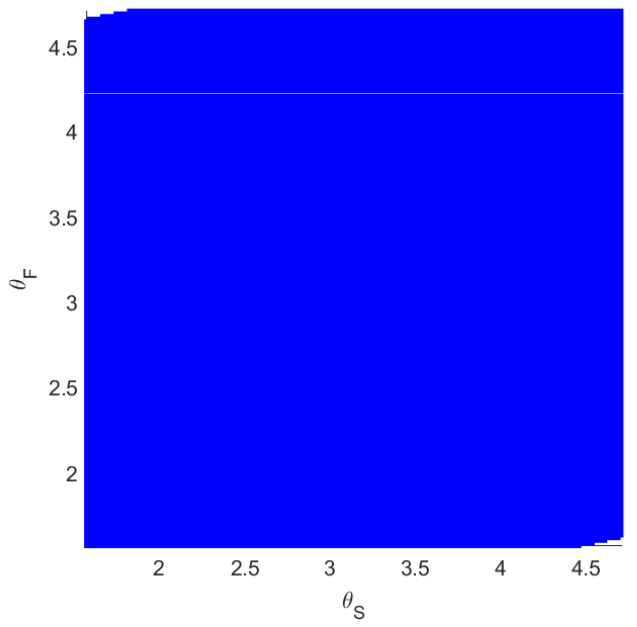

(b) Stability region for $\xi = \frac{1}{3}, \alpha = 1$.

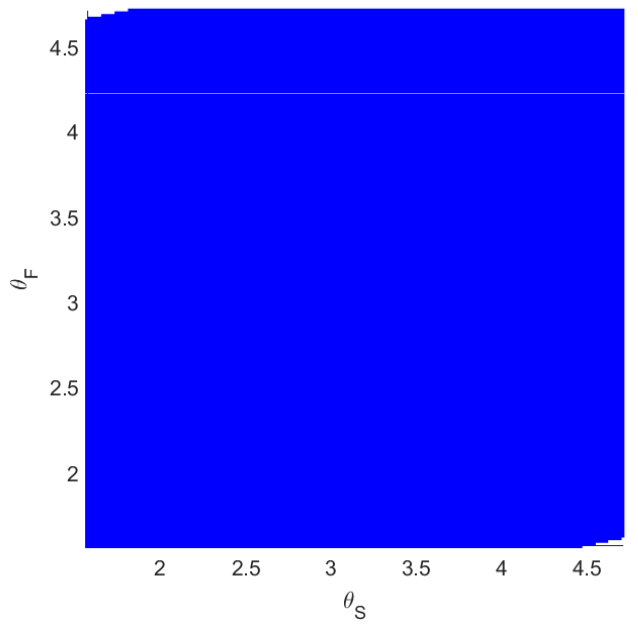

(c) Stability region for $\xi = \frac{2}{3}, \alpha = 1$.

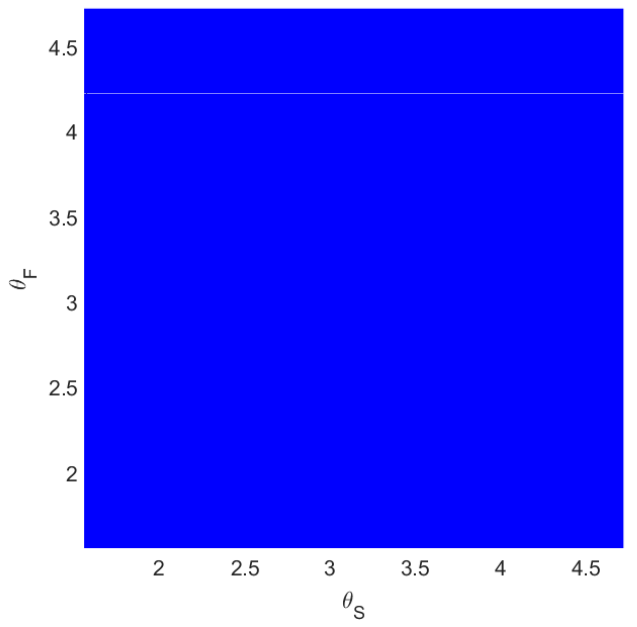

(d) Stability region for $\xi = 1, \alpha = 1$.

Fig. 11: Stability region for $\frac{1}{2}\pi \leq \theta_S, \theta_F \leq \frac{3}{2}\pi$ setting $z^{\{\mathrm{S,S}\}} = \exp(i\theta_S), z^{\{\mathrm{F,F}\}} = \mathtt{m}\exp(i\theta_F), z^{\{\mathrm{S,F}\}} = -\alpha\xi(z^{\{\mathrm{S,S}\}} - z^{\{\mathrm{F,F}\}})$ and $z^{\{\mathrm{F,S}\}} = \frac{1-\xi}{\alpha}(z^{\{\mathrm{S,S}\}} - z^{\{\mathrm{F,F}\}})$.

and get

$$
\begin{aligned}
&(6.8a) \qquad \mathbf{Y}_{\mathrm{S},n+1} = \boldsymbol{\Phi}\mathbf{Y}_{\mathrm{S},n} + H\boldsymbol{\Gamma}\tilde{\mathbf{f}}^{\{\mathrm{S}\}}_{n+1} + \mathbf{e}_2 H(\mathbf{f}^{\{\mathrm{S}\}}_{n+1} - \tilde{\mathbf{f}}^{\{\mathrm{S}\}}_{n+1}),\\
&(6.8b) \qquad \mathbf{Y}_{\mathrm{F},\ell+1} = \boldsymbol{\Phi}\mathbf{Y}_{\mathrm{F},\ell+1-1/\mathtt{m}} + h\boldsymbol{\Gamma}\tilde{\mathbf{f}}^{\{\mathrm{F}\}}_{\ell+1}, \qquad \ell = n-1+\frac{1}{\mathtt{m}} : \frac{1}{\mathtt{m}} : n.
\end{aligned}
$$

Applied to (4.3) this reads

$$
\begin{aligned}
(6.9a) \qquad \mathbf{Y}_{\mathrm{S},n+1} &= \boldsymbol{\Phi}\mathbf{Y}_{\mathrm{S},n} + \boldsymbol{\Gamma}(z^{\{\mathrm{S,S}\}}\mathbf{e}_1^\top\mathbf{Y}_{\mathrm{S},n+1} + z^{\{\mathrm{F,S}\}}\tilde{\mathbf{y}}^{\{\mathrm{F}\}}_{n+1})\\
&\quad + \mathbf{e}_2 z^{\{\mathrm{F,S}\}}(\mathbf{e}_1^\top\mathbf{Y}_{\mathrm{F},n+1} - \tilde{\mathbf{y}}^{\{\mathrm{F}\}}_{n+1}),\\
(6.9b) \qquad \mathbf{Y}_{\mathrm{F},\ell+1} &= \boldsymbol{\Phi}\mathbf{Y}_{\mathrm{F},\ell+1-1/\mathtt{m}} + \boldsymbol{\Gamma}\left(\frac{z^{\{\mathrm{S,F}\}}}{\mathtt{m}}\tilde{\mathbf{y}}^{\{\mathrm{S}\}}_{\ell+1} + \frac{z^{\{\mathrm{F,F}\}}}{\mathtt{m}}\mathbf{e}_1^\top\mathbf{Y}_{\mathrm{F},\ell+1}\right),\\
&\quad \ell = n-1+\frac{1}{\mathtt{m}} : \frac{1}{\mathtt{m}} : n.
\end{aligned}
$$

Applying the last equation $\mathtt{m}-1$ times recursively, one gets

$$
\mathbf{Y}_{\mathrm{F},n+1} = \mathbf{M}^{\mathtt{m}}(z^{\{\mathrm{F,F}\}}/\mathtt{m})\mathbf{Y}_{\mathrm{F},n} + \sum_{i=0}^{\mathtt{m}-1}\mathbf{M}^i(z^{\{\mathrm{F,F}\}}/\mathtt{m})\tilde{\mathbf{M}}(z/\mathtt{m})\boldsymbol{\Gamma}\frac{z^{\{\mathrm{S,F}\}}}{\mathtt{m}}\tilde{\mathbf{y}}^{\{\mathrm{S}\}}_{n+1-i}
$$

with

$$\mathbf{M}(z) := \underbrace{\left(\mathbf{I} + (z/(1-z\beta_0))\mathbf{\Gamma}\mathbf{e}_1^\top\right)^{-1}}_{=:\,\tilde{\mathbf{M}}(z)} \mathbf{\Phi}$$

being the companion matrix of single rate linear multistep schemes. When the multirate scheme is applied to a system without coupling ($z^{\{\text{S,F}\}} = z^{\{\text{F,S}\}} = 0$), we get

$$\mathbf{Y}_{\text{S},n+1} = \mathbf{M}(z^{\{\text{S,S}\}})\mathbf{Y}_{\text{S},n}, \tag{6.10a}$$

$$\mathbf{Y}_{\text{F},n+1} = \mathbf{M}^{\mathtt{m}}(z^{\{\text{F,F}\}}/\mathtt{m})\mathbf{Y}_{\text{F},n}. \tag{6.10b}$$

Thus the stability of a decoupled system is given by the stability of the basic multirate schemes with macro and micro step size $H$ and $h = H/\mathtt{m}$, resp.

In the case of non-vanishing coupling parameters $\eta_{\text{F}}$ and $\eta_{\text{S}}$, the coupling variables $\tilde{\mathbf{y}}_{n+1}^{\{\text{F}\}}$ and $\tilde{\mathbf{y}}_{m+1}^{\{\text{S}\}}$ have to be taken into consideration:

- for the semi-implicit and implicit formulations these are based on extrapolation or interpolation, i.e., they are given as a linear combination of $\mathbf{y}_n^{\{\text{F}\}}$ or $\mathbf{y}_{n+1}^{\{\text{F}\}}$, and $\mathbf{y}_n^{\{\text{S}\}}$ or $\mathbf{y}_{n+1}^{\{\text{S}\}}$, resp.
- in the compound-step approach, $\tilde{\mathbf{y}}_{n+1}^{\{\text{F}\}}$ depends on both $\mathbf{y}_n^{\{\text{S}\}}$ and $\mathbf{y}_n^{\{\text{F}\}}$, and $\tilde{\mathbf{y}}_{\ell+1}^{\{\text{S}\}}$ can be obtained via interpolation.
- in the generalized compound-step approach, $\tilde{\mathbf{y}}_{n+1}^{\{\text{F}\}}$ depends on both $\mathbf{y}_n^{\{\text{F}\}}$ and $\mathbf{y}_{n+1}^{\{\text{F}\}}$, and $\tilde{\mathbf{y}}_{\ell+1}^{\{\text{S}\}}$ can be obtained via interpolation.

Inserting these interpolated/extrapolated values into (6.9) leads to

$$\begin{bmatrix} I & \mathcal{O}(z^{\{\text{F,S}\}}) \\ \mathcal{O}(z^{\{\text{S,F}\}}) & I \end{bmatrix} \begin{bmatrix} \mathbf{Y}_{\text{S},n+1} \\ \mathbf{Y}_{\text{F},n+1} \end{bmatrix} = \begin{bmatrix} M(z_s) + \mathcal{O}(z^{\{\text{F,S}\}}) & \mathcal{O}(z^{\{\text{F,S}\}}) \\ \mathcal{O}(z^{\{\text{S,F}\}}) & \left(M(\frac{z_f}{\mathtt{m}})\right)^{\mathtt{m}} \end{bmatrix} \begin{bmatrix} \mathbf{Y}_{\text{S},n} \\ \mathbf{Y}_{\text{F},n} \end{bmatrix} \quad \Rightarrow$$

$$\begin{bmatrix} \mathbf{Y}_{\text{S},n+1} \\ \mathbf{Y}_{\text{F},n+1} \end{bmatrix} = \mathbf{M}_{\text{MRLMM}} \begin{bmatrix} \mathbf{Y}_{\text{S},n+1} \\ \mathbf{Y}_{\text{F},n+1} \end{bmatrix} \quad \text{with}$$

$$\mathbf{M}_{\text{MRLMM}} = \begin{bmatrix} M(z_s) & 0 \\ 0 & M(\frac{z_f}{\mathtt{m}})^{\mathtt{m}} \end{bmatrix} + \begin{bmatrix} \mathcal{O}(w_f) & \mathcal{O}(w_f) \\ \mathcal{O}(w_s) & \mathcal{O}(w_s w_f) \end{bmatrix}.$$

Consequently, as the spectral radius of $\mathbf{M}_{\text{MRLMM}}$ is a continuous function in both $w_s$ and $w_f$, the A-stability of the multirate linear multistep methods follows from the A-stability of the underlying multistep schemes provided that $z^{\{\text{S,F}\}}$ *or* $z^{\{\text{F,S}\}}$ is small enough. In addition, the same results holds for one-way couplings, i.e., either $z^{\{\text{S,F}\}} = 0$ or $z^{\{\text{F,S}\}}$ =0, with $z^{\{\text{F,S}\}}$ or $z^{\{\text{S,F}\}}$ arbitrarily.

Remark 17. *As A-stability is not possible for linear multistep schemes with order larger than two, we can only hope for a generalized $A(\alpha)$-stability in the general case. For example, for fixed $z^{\{\text{S,F}\}}$ and $z^{\{\text{F,S}\}}$, we can consider the slow stability region, which is defined as the set of all $z^{\{\text{S,S}\}} \in \mathbb{C}$, for which the spectral radius of $M$ is bounded by one for all $z^{\{\text{F,F}\}}$ in a section of the left half plane defined by an angle $\alpha \le \pi/2$ (see Sect. 5.5 for more details).*

**6.8. Remarks.** Based on preliminary work by Wells [17], Gear [21] and Orailoglu [39], the first comprehensive study and still the fundamental work on multirate linear multistep methods was published in 1984 by Gear and Wells [22]. It ranges from efficiency considerations, convergence and error analysis to absolute stability and numerical test results for an experimental code called MRATE, which is based on the Adams-type algorithm DIFSUB. For the first time, the concept of fastest-first and slowest-first methods was introduced, provided with a discussion on the advantages and disadvantages of both approaches.

So far, multirate linear multistep schemes have been generalized in two directions: a parallel version, based on the fastest-first approach, has been suggested by Chang [9]; for applications in circuit simulation, Verhoeven et al. have introduced the concept of (general) compound-step schemes into the class of multirate linear multistep schemes: in a first step for the implicit Euler scheme [25], and then for BDF schemes [70]. These schemes generalize the idea of compound-step schemes to schemes which jointly compute in the beginning the slow part at the end of the macro step $t_n$ and the fast part arbitrarily in the interval $[t_{n-1}, t_n]$. In a series of papers, Verhoeven et al. have considered different aspects of general compound-step multirate BDF schemes: error analysis and convergence [68], step size control [72] and automatic partitioning [71], finally summarized in a survey paper [72].

Though there has obviously been only few activity with respect to method derivation of multirate linear multistep schemes, stability considerations attained much more interest. In general, linear stability is based on $2 \times 2$ linear block test system of differential equations. In [22] it was shown that multirate linear multistep methods preserve the A-stability properties of the underlying linear multistep schemes, if the linear system of ordinary differential equations to be solved is nearly block triangular, with other words, if the the coupling of fast to slow components is only weak. A more general discussion on A-stability for linear multistep methods, for both semi-implicit and implicit, slowest-first and fastest-first algorithm, can be found in [65] for a scalar test problem described by a $2 \times 2$-matrix $A$. Absolute stability is not inherited from A-stable basic schemes for unbounded coefficients in $A$. However, Skelboe and Andersen conjecture that the implicit approach is A-stable for linear one- and two-step methods if $A$ is real with eigenvalues in the left-hand half-plane and, in addition, $a_{11} \leq 0$. The second condition prevents a match of A-stable multirate schemes and stable linear differential equations. For the semi-implicit approach with linear and zero-order interpolation, Gomez et al. [19] were the first to derive necessary and sufficient stability conditions for the $2 \times 2$-scalar test equation by a general companion matrix approach. First stability conditions for the general compound step approach are presented in [69] for the multirate backward Euler scheme, and later generalized to BDF slowest-first multirate schemes. A comprehensive stability analysis of the latter class of schemes can be found in the survey paper [70].

Nonlinear stability analysis results of multirate linear multistep methods are reported by Sand and Skelboe [47]: for monotonically max-norm stable systems of ordinary differential equations, semi-implicit multirate backward Euler schemes yield stable numerical approximations. This result can be transferred to implicit schemes, which are iteratively solved by waveform relaxation schemes.

## 7. Multirate linearly-implicit GARK schemes.

Rosenbrock-Wanner (ROS) and W-methods can be easily generalized into a linearly-implicit generalization of GARK schemes, the GARK-ROS/ROW methods:

Definition 7.1 (GARK-ROS/ROW method). *One macro-step of a GARK-ROS/ROW method applied to* (2.1) *with stepsize* $H$ *reads:*

$$\begin{aligned}
\mathbf{k}^{\{F\}} &= H\,\mathbf{f}^{\{F\}}\Big(\mathbf{1}_{s^{\{F\}}} \otimes \mathbf{y}_n^{\{F\}} + \boldsymbol{\alpha}^{\{F,F\}} \circledast_{F} \mathbf{k}^{\{F\}},\, \mathbf{1}_{s^{\{F\}}} \otimes \mathbf{y}_n^{\{S\}} + \boldsymbol{\alpha}^{\{F,S\}} \circledast_{S} \mathbf{k}^{\{S\}}\Big) \\
&\quad + \Big(\mathbf{I}_{s^{\{F\}} \times s^{\{F\}}} \otimes h\mathbf{L}^{\{F,F\}}\Big)\Big(\boldsymbol{\gamma}^{\{F,F\}} \circledast_{F} \mathbf{k}^{\{F\}}\Big) \qquad (7.1a)
\end{aligned}$$

$$+\left(\mathbf{I}_{s^{\{F\}}\times s^{\{F\}}} \otimes h\mathbf{L}^{\{F,S\}}\right)\left(\boldsymbol{\gamma}^{\{F,S\}} \circledast \mathbf{k}^{\{S\}}\right);$$

$$\mathbf{k}^{\{S\}} = H\,\mathbf{f}^{\{S\}}\left(\mathbf{1}_{s^{\{S\}}} \otimes \mathbf{y}_n^{\{F\}} + \boldsymbol{\alpha}^{\{S,F\}} \circledast \mathbf{k}^{\{F\}}, \mathbf{1}_{s^{\{S\}}} \otimes \mathbf{y}_n^{\{S\}} + \boldsymbol{\alpha}^{\{S,S\}} \circledast \mathbf{k}^{\{S\}}\right)$$

$$+(\mathbf{I}_{s^{\{S\}}\times s^{\{S\}}} \otimes h\mathbf{L}^{\{S,F\}})\left(\boldsymbol{\gamma}^{\{S,F\}} \circledast \mathbf{k}^{\{F,\ell\}}\right) \tag{7.1b}$$

$$+(\mathbf{I}_{s^{\{S\}}\times s^{\{S\}}} \otimes h\mathbf{L}^{\{S,S\}})\left(\boldsymbol{\gamma}^{\{S,S\}} \circledast \mathbf{k}^{\{S\}}\right);$$

$$\mathbf{y}_{n+1}^{\{F\}} = \mathbf{y}_n^{\{F\}} + \mathsf{b}^{\{F\}T} \circledast \mathbf{k}^{\{F\}}; \tag{7.1c}$$

$$\mathbf{y}_{n+1}^{\{S\}} = \mathbf{y}_n^{\{S\}} + \mathsf{b}^{\{S\}T} \circledast \mathbf{k}^{\{S\}}. \tag{7.1d}$$

Here $\mathbf{1}_s \in \mathbb{R}^s$ is a vector of ones, $\mathbf{I}_{s\times s} \in \mathbb{R}^{s\times s}$ is the identity matrix, $\otimes$ is the Kronecker product, and the following matrix notation is used for $m, q \in \{F, S\}$:

$$\mathbf{k}^{\{q\}} := \begin{bmatrix} \mathbf{k}_1^{\{q\}} \\ \vdots \\ \mathbf{k}_{s^{\{q\}}}^{\{q\}} \end{bmatrix} \in \mathbb{R}^{d^{\{q\}}s^{\{q\}}}, \qquad \begin{aligned} \boldsymbol{\alpha}^{\{m,q\}} \circledast \mathbf{k}^{\{q\}} &:= \left(\boldsymbol{\alpha}^{\{m,q\}} \otimes \mathbf{I}_{d\times d}\right)\mathbf{k}^{\{q\}}, \\ \boldsymbol{\gamma}^{\{m,q\}} \circledast \mathbf{k}^{\{q\}} &:= \left(\boldsymbol{\gamma}^{\{m,q\}} \otimes \mathbf{I}_{d\times d}\right)\mathbf{k}^{\{q\}}, \end{aligned} \tag{7.2a}$$

with $\boldsymbol{\alpha}^{\{m,q\}}, \boldsymbol{\gamma}^{\{m,q\}} \in \mathbb{R}^{s^{\{m\}}\times s^{\{q\}}}$ and $\mathbf{f}^{\{q\}}\left(\mathbf{1}_s \otimes \mathbf{y}_n + \boldsymbol{\alpha}^{\{q,F\}} \circledast \mathbf{k}^{\{F\}} + \boldsymbol{\alpha}^{\{q,S\}} \circledast \mathbf{k}^{\{S\}}\right) :=$

$$\begin{bmatrix} \mathbf{f}^{\{q\}}\left(\mathbf{y}_n + \sum_{j=1}^{s^{\{F\}}} \alpha_{1,j}^{\{q,F\}}\mathbf{k}_j^{\{F\}} + \sum_{j=1}^{s^{\{S\}}} \alpha_{1,j}^{\{q,S\}}\mathbf{k}_j^{\{S\}}\right) \\ \mathbf{f}^{\{q\}}\left(\mathbf{y}_n + \sum_{j=1}^{s^{\{F\}}} \alpha_{2,j}^{\{q,F\}}\mathbf{k}_j^{\{F\}} + \sum_{j=1}^{s^{\{S\}}} \alpha_{2,j}^{\{q,S\}}\mathbf{k}_j^{\{S\}}\right) \\ \vdots \\ \mathbf{f}^{\{q\}}\left(\mathbf{y}_n + \sum_{j=1}^{s^{\{F\}}} \alpha_{s^{\{q\}},j}^{\{q,F\}}\mathbf{k}_j^{\{F\}} + \sum_{j=1}^{s^{\{S\}}} \alpha_{s^{\{q\}},j}^{\{q,S\}}\mathbf{k}_j^{\{S\}}\right) \end{bmatrix} \tag{7.2b}$$

The matrices $\boldsymbol{\alpha}^{\{q,m\}}$ are strictly lower triangular and $\boldsymbol{\gamma}^{\{q,m\}}$ are lower triangular. Depending on the choice of matrices $\mathbf{L}^{\{m,q\}} \in \mathbb{R}^{d^{\{m\}}\times d^{\{q\}}}$ one distinguishes several types of methods, as follows:

- GARK-ROS schemes use the exact Jacobian information, i.e., $\mathbf{L}^{\{m,q\}} := \partial\mathbf{f}^{\{m\}}/\partial\mathbf{y}^{\{q\}}(\mathbf{y}_n^{\{F\}}, \mathbf{y}_n^{\{S\}})$ are the Jacobians of the component functions evaluated at the current solution;
- GARK-ROW schemes allow any approximation of the Jacobian, i.e., $\mathbf{L}^{\{m,q\}}$ may be arbitrary;
- In the case of GARK-ROS schemes with time-lagged Jacobians one has $\mathbf{L}^{\{m,q\}} := \partial\mathbf{f}^{\{m\}}/\partial\mathbf{y}^{\{q\}}(\mathbf{y}_n^{\{F\}}, \mathbf{y}_n^{\{S\}}) + \mathcal{O}(h)$.

The scheme (7.1)is characterized by the extended Butcher tableau (with $\boldsymbol{\beta}^{\{q,m\}} := \boldsymbol{\alpha}^{\{q,m\}} + \boldsymbol{\gamma}^{\{q,m\}}$)

$$\begin{array}{c|c} \mathbf{A} & \mathbf{G} \\ \hline \mathbf{b}^\top & \end{array} = \begin{array}{cc|cc} \boldsymbol{\alpha}^{\{F,F\}} & \boldsymbol{\alpha}^{\{F,S\}} & \boldsymbol{\gamma}^{\{F,F\}} & \boldsymbol{\gamma}^{\{F,S\}} \\ \boldsymbol{\alpha}^{\{S,F\}} & \boldsymbol{\alpha}^{\{S,S\}} & \boldsymbol{\gamma}^{\{S,F\}} & \boldsymbol{\gamma}^{\{S,S\}} \\ \hline \mathbf{b}^{\{F\}T} & \mathbf{b}^{\{S\}T} & & \end{array}. \tag{7.3}$$

REMARK 18 (GARK-ROS and GARK-ROW scheme structure).

- *Similar to GARK, the scheme* (7.1) *uses only one function* $\mathbf{f}^{\{q\}}$ *evaluation for the increment* $\mathbf{k}^{\{q\}}$*, and increments* $\mathbf{k}^{\{m\}}$ *both as function arguments and as additive terms.*
- *For all* $\boldsymbol{\gamma}^{\{q,m\}} = 0$ (7.1) *is an explicit GARK scheme.*

As in the GARK case, the split system (2.1) directly allows to generalize a multirate version of GARK-ROS/ROW schemes by solving the slow expensive component with a large macro-step $H$, and the fast inexpensive one with $\mathtt{m}$ micro-steps $h = H/\mathtt{m}$:

Definition 7.2 (Multirate GARK-ROS/ROW method). *One step of a multirate GARK-ROS/ROW method (MR-GARK-ROS/ROW for short) applied to* (2.1) *computes the solution as follows. The slow component is discretized with a ROS/ROW method* $(\mathtt{b}^{\{\mathrm{S}\}}, \boldsymbol{\alpha}^{\{\mathrm{S},\mathrm{S}\}}, \boldsymbol{\gamma}^{\{\mathrm{S},\mathrm{S}\}})$ *and macro-step* $H$*. The fast component uses* $\mathtt{m}$ *micro-steps* $h = H/\mathtt{m}$*, and at each micro-step* $\ell$ *a (possibly different) ROS/ROW method* $(\mathtt{b}^{\{\mathrm{F},\ell\}}, \boldsymbol{\alpha}^{\{\mathrm{F},\mathrm{F},\ell\}}, \boldsymbol{\gamma}^{\{\mathrm{F},\mathrm{F},\ell\}})$ *is applied. The computational process reads:*

$$
\begin{aligned}
\mathbf{k}^{\{\mathrm{F},\ell\}} &= h\,\mathbf{f}^{\{\mathrm{F}\}}\Big(\mathbf{1}_{s^{\{\mathrm{F}\}}} \otimes \mathbf{y}_n^{\{\mathrm{F}\}} + \sum_{\ell=1}^{\ell-1} \mathtt{b}^{\{\mathrm{F},\ell\}T} \otimes \mathbf{k}^{\{\mathrm{F},\ell\}} + \boldsymbol{\alpha}^{\{\mathrm{F},\mathrm{F},\ell\}} \otimes \mathbf{k}^{\{\mathrm{F},\ell\}}, \\
&\qquad \mathbf{1}_{s^{\{\mathrm{F}\}}} \otimes \mathbf{y}_n^{\{\mathrm{S}\}} + \boldsymbol{\alpha}^{\{\mathrm{F},\mathrm{S},\ell\}} \otimes \mathbf{k}^{\{\mathrm{S}\}}\Big) \qquad (7.4a)\\
&\quad + \Big(\mathbf{I}_{s^{\{\mathrm{F}\}}\times s^{\{\mathrm{F}\}}} \otimes h\mathbf{L}^{\{\mathrm{F},\mathrm{F}\}}\Big)\Big(\boldsymbol{\gamma}^{\{\mathrm{F},\mathrm{F},\ell\}} \otimes \mathbf{k}^{\{\mathrm{F},\ell\}}\Big) \\
&\quad + \Big(\mathbf{I}_{s^{\{\mathrm{F}\}}\times s^{\{\mathrm{F}\}}} \otimes h\mathbf{L}^{\{\mathrm{F},\mathrm{S}\}}\Big)\Big(\boldsymbol{\gamma}^{\{\mathrm{F},\mathrm{S},\ell\}} \otimes \mathbf{k}^{\{\mathrm{S}\}}\Big) \\
&\quad \text{for} \quad \ell = 1, \dots, \mathtt{m}; \\
\mathbf{k}^{\{\mathrm{S}\}} &= H\,\mathbf{f}^{\{\mathrm{S}\}}\Big(\mathbf{1}_{s^{\{\mathrm{S}\}}} \otimes \mathbf{y}_n^{\{\mathrm{F}\}} + \sum_{\ell=1}^{\mathtt{m}} \boldsymbol{\alpha}^{\{\mathrm{S},\mathrm{F},\ell\}} \otimes \mathbf{k}^{\{\mathrm{F},\ell\}}, \\
&\qquad \mathbf{1}_{s^{\{\mathrm{S}\}}} \otimes \mathbf{y}_n^{\{\mathrm{S}\}} + \boldsymbol{\alpha}^{\{\mathrm{S},\mathrm{S}\}} \otimes \mathbf{k}^{\{\mathrm{S}\}}\Big) \qquad (7.4b)\\
&\quad + (\mathbf{I}_{s^{\{\mathrm{S}\}}\times s^{\{\mathrm{S}\}}} \otimes H\mathbf{L}^{\{\mathrm{S},\mathrm{F}\}})\Big(\sum_{\ell=1}^{\mathtt{m}} \boldsymbol{\gamma}^{\{\mathrm{S},\mathrm{F},\ell\}} \otimes \mathbf{k}^{\{\mathrm{F},\ell\}}\Big) \\
&\quad + (\mathbf{I}_{s^{\{\mathrm{S}\}}\times s^{\{\mathrm{S}\}}} \otimes H\mathbf{L}^{\{\mathrm{S},\mathrm{S}\}})\Big(\boldsymbol{\gamma}^{\{\mathrm{S},\mathrm{S}\}} \otimes \mathbf{k}^{\{\mathrm{S}\}}\Big);
\end{aligned}
$$

$$
\mathbf{y}_{n+1}^{\{\mathrm{F}\}} = \mathbf{y}_n^{\{\mathrm{F}\}} + \sum_{\ell=1}^{\mathtt{m}} \mathtt{b}^{\{\mathrm{F},\ell\}T} \otimes \mathbf{k}^{\{\mathrm{F},\ell\}}; \qquad (7.4c)
$$

$$
\mathbf{y}_{n+1}^{\{\mathrm{S}\}} = \mathbf{y}_n^{\{\mathrm{S}\}} + \mathtt{b}^{\{\mathrm{S}\}T} \otimes \mathbf{k}^{\{\mathrm{S}\}}. \qquad (7.4d)
$$

The matrices $\mathbf{L}^{\{m,q\}}$, $m, q \in \{\mathrm{S}, \mathrm{F}\}$, can be chosen as exact Jacobians (ROS schemes) or as arbitrary approximations to the Jacobians (ROW schemes).[1] Method (7.4) is a GARK ROS/ROW scheme with step $H$ defined by the Butcher

[1]Note that in the GARK-ROS case, all Jacobians are evaluated at the current macro grid point $t_n$ for two reasons: first of all, one saves computational and overhead costs. Secondly, we only aim at achieving the order of the method in the macro grid points, and not in addition in the micro grid points, as we regard our scheme as a partitioned scheme on the macro step size. The same point of view was taken in [4].

tableau (7.3):

(7.5a)
$$
\begin{array}{c|c}
\boldsymbol{\alpha}^{\{\mathrm{F},\mathrm{F}\}} & \boldsymbol{\alpha}^{\{\mathrm{F},\mathrm{S}\}} \\ \hline
\boldsymbol{\alpha}^{\{\mathrm{S},\mathrm{F}\}} & \boldsymbol{\alpha}^{\{\mathrm{S},\mathrm{S}\}} \\ \hline\hline
\mathbf{b}^{\{\mathrm{F}\}\top} & \mathbf{b}^{\{\mathrm{S}\}\top}
\end{array}
:=
\begin{array}{ccc|c}
\frac{1}{\mathtt{m}}\boldsymbol{\alpha}^{\{\mathrm{F},\mathrm{F},1\}} & \cdots & 0 & \boldsymbol{\alpha}^{\{\mathrm{F},\mathrm{S},1\}} \\
\vdots & \ddots & \vdots & \vdots \\
\frac{1}{\mathtt{m}}\mathbf{1}^{\{\mathrm{F}\}}\mathtt{b}^{\{\mathrm{F},1\}T} & \cdots & \frac{1}{\mathtt{m}}\boldsymbol{\alpha}^{\{\mathrm{F},\mathrm{F},\mathtt{m}\}} & \boldsymbol{\alpha}^{\{\mathrm{F},\mathrm{S},\mathtt{m}\}} \\ \hline
\frac{1}{\mathtt{m}}\boldsymbol{\alpha}^{\{\mathrm{S},\mathrm{F},1\}} & \cdots & \frac{1}{\mathtt{m}}\boldsymbol{\alpha}^{\{\mathrm{S},\mathrm{F},\mathtt{m}\}} & \boldsymbol{\alpha}^{\{\mathrm{S},\mathrm{S}\}} \\ \hline\hline
\frac{1}{\mathtt{m}}\mathtt{b}^{\{\mathrm{F},1\}T} & \cdots & \frac{1}{\mathtt{m}}\mathtt{b}^{\{\mathrm{F},\mathtt{m}\}T} & \mathtt{b}^{\{\mathrm{S}\}T}
\end{array},
$$

(7.5b)
$$
\begin{array}{c|c}
\boldsymbol{\gamma}^{\{\mathrm{F},\mathrm{F}\}} & \boldsymbol{\gamma}^{\{\mathrm{F},\mathrm{S}\}} \\ \hline
\boldsymbol{\gamma}^{\{\mathrm{S},\mathrm{F}\}} & \boldsymbol{\gamma}^{\{\mathrm{S},\mathrm{S}\}}
\end{array}
:=
\begin{array}{ccc|c}
\frac{1}{\mathtt{m}}\boldsymbol{\gamma}^{\{\mathrm{F},\mathrm{F},1\}} & \cdots & 0 & \boldsymbol{\gamma}^{\{\mathrm{F},\mathrm{S},1\}} \\
\vdots & \ddots & \vdots & \vdots \\
0 & \cdots & \frac{1}{\mathtt{m}}\boldsymbol{\gamma}^{\{\mathrm{F},\mathrm{F},\mathtt{m}\}} & \boldsymbol{\gamma}^{\{\mathrm{F},\mathrm{S},\mathtt{m}\}} \\ \hline
\frac{1}{\mathtt{m}}\boldsymbol{\gamma}^{\{\mathrm{S},\mathrm{F},1\}} & \cdots & \frac{1}{\mathtt{m}}\boldsymbol{\gamma}^{\{\mathrm{S},\mathrm{F},\mathtt{m}\}} & \boldsymbol{\gamma}^{\{\mathrm{S},\mathrm{S}\}}
\end{array},
$$

and $\boldsymbol{\beta}^{\{m,q\}} := \boldsymbol{\alpha}^{\{m,q\}} + \boldsymbol{\gamma}^{\{m,q\}}$ for $m, q \in \{\mathrm{S}, \mathrm{F}\}$.

REMARK 19 (Non-uniform micro-steps). *Definition 7.2 can be immediately extended to accommodate nonuniform micro-steps $h_\ell = \sigma_e llH$ with $\sum_{\ell=1}^{\mathtt{m}} \sigma_\ell = 1$ by using $h_\ell$ in each fast step (7.4a), and scaling columns $\ell$ of the Butcher tableaus (7.5a) and (7.5b) by $\sigma_\ell$ instead of $1/\mathtt{m}$.*

REMARK 20 (Pure multirate approach). *In this paper we focus on "pure" multirate methods where all micro-steps use the same fast base method, i.e.,*

(7.6)
$$
\left(\mathtt{b}^{\{\mathrm{F},\ell\}}, \boldsymbol{\alpha}^{\{\mathrm{F},\mathrm{F},\ell\}}, \boldsymbol{\gamma}^{\{\mathrm{F},\mathrm{F},\ell\}}\right) \equiv \left(\mathtt{b}^{\{\mathrm{F}\}}, \boldsymbol{\alpha}^{\{\mathrm{F}\}}, \boldsymbol{\gamma}^{\{\mathrm{F}\}}\right), \quad \ell = 1, \dots, \mathtt{m}.
$$

*Nevertheless, the general structure of Definition 7.2 remains of interest as it allows to build special types of fast-slow couplings.*

REMARK 21 (Order conditions). *As for GARK schemes, the order conditions of Multirate GARK-ROS/ROW schemes can be easily derived from the order conditions of GARK-ROS/ROW schemes. See for more details*

REMARK 22 (Linear stability).

*As in Section 5.5, we consider the model problem (5.11) with $z^{\{\mathrm{F},\mathrm{F}\}} = H\lambda^{\{\mathrm{F},\mathrm{F}\}}$, $z^{\{\mathrm{S},\mathrm{S}\}} = H\lambda^{\{\mathrm{S},\mathrm{S}\}}$, $z^{\{\mathrm{S},\mathrm{F}\}} = (1-\xi)/\alpha \cdot \left(z^{\{\mathrm{F},\mathrm{F}\}} - z^{\{\mathrm{S},\mathrm{S}\}}\right)$, and $z^{\{\mathrm{F},\mathrm{S}\}} = -\alpha\,\xi \cdot \left(z^{\{\mathrm{F},\mathrm{F}\}} - z^{\{\mathrm{S},\mathrm{S}\}}\right)$, see [36]. Application of the MGARK-ROS method (7.4), regarded as a partitioned GARK-ROS scheme according to the Butcher tableau (7.5) over one step $H$, advances now via the recurrence:*

$$
\begin{bmatrix} \mathbf{y}_{n+1}^{\{\mathrm{F}\}} \\ \mathbf{y}_{n+1}^{\{\mathrm{S}\}} \end{bmatrix} = \mathbf{M}(z^{\{\mathrm{F},\mathrm{F}\}}, z^{\{\mathrm{S},\mathrm{S}\}}, z^{\{\mathrm{S},\mathrm{F}\}}, z^{\{\mathrm{F},\mathrm{S}\}}) \begin{bmatrix} \mathbf{y}_{n}^{\{\mathrm{F}\}} \\ \mathbf{y}_{n}^{\{\mathrm{S}\}} \end{bmatrix},
$$

*with the stability matrix:*

$$
\mathbf{M}(z^{\{\mathrm{F,F}\}}, z^{\{\mathrm{S,S}\}}, z^{\{\mathrm{S,F}\}}, z^{\{\mathrm{F,S}\}}) = \mathbf{I}_{2\times 2} + \begin{bmatrix} \mathbf{b}^{\{\mathrm{F}\}} & \mathbf{0}_{\mathtt{m}s^{\{\mathrm{F}\}}} \\ \mathbf{0}_{s^{\{\mathrm{S}\}}} & \mathbf{b}^{\{\mathrm{S}\}} \end{bmatrix}^{\top} \cdot
$$

$$
\begin{bmatrix} \mathbf{I}_{Ms^{\{\mathrm{F}\}}\times Ms^{\{\mathrm{F}\}}} - z^{\{\mathrm{F,F}\}}\,\boldsymbol{\beta}^{\{\mathrm{F,F}\}} & -z^{\{\mathrm{S,F}\}}\,\boldsymbol{\beta}^{\{\mathrm{F,S}\}} \\ -z^{\{\mathrm{F,S}\}}\,\boldsymbol{\beta}^{\{\mathrm{S,F}\}} & \mathbf{I}_{s^{\{\mathrm{S}\}}\times s^{\{\mathrm{S}\}}} - z^{\{\mathrm{S,S}\}}\,\boldsymbol{\beta}^{\{\mathrm{S,S}\}} \end{bmatrix}^{-1} \cdot \begin{bmatrix} z^{\{\mathrm{F,F}\}}\,\mathbf{1}_{\mathtt{m}s^{\{\mathrm{F}\}}} & z^{\{\mathrm{S,F}\}}\,\mathbf{1}_{\mathtt{m}s^{\{\mathrm{F}\}}} \\ z^{\{\mathrm{F,S}\}}\,\mathbf{1}_{s^{\{\mathrm{S}\}}} & z^{\{\mathrm{S,S}\}}\,\mathbf{1}_{s^{\{\mathrm{S}\}}} \end{bmatrix}.
$$

*One immediately sees that for one-sided coupled problems with* $z^{\{\mathrm{S,F}\}} = 0$ *or* $z^{\{\mathrm{F,S}\}} = 0$ *the stability of the base schemes guarantees the stability of the multirate schemes.*

**7.1. Coupling the fast and slow systems in a computationally-efficient manner.** In traditional Rosenbrock methods the coefficient matrix $\boldsymbol{\alpha}$ is strictly lower triangular, and the matrix $\boldsymbol{\gamma}$ lower triangular with equal diagonal entries. Due to this structure the stages $\mathbf{k}_i$ are evaluated sequentially in a decoupled manner, each stage computation is only implicit in the current stage.

Multirate GARK ROS/ROW schemes compute both slow and fast stage vectors. We call a stage computation "decoupled" if it is implicit in only the current stage $\mathbf{k}_i^{\{\mathrm{S}\}}$ or $\mathbf{k}_i^{\{\mathrm{F},\ell\}}$. We call computations "coupled" if one (or more) slow stages, and one (or more) fast stages are computed together by solving a single large system of linear equations. Computational efficiency of multirate methods relies on evaluating less frequently the expensive slow part. Consequently, an efficient multirate method keeps the coupling at a minimum.

In the following we consider multirate GARK ROS/ROW schemes where the base methods are Rosenbrock(-W) schemes with matrices $\boldsymbol{\alpha}^{\{\mathrm{F,F},\ell\}}$, $\boldsymbol{\alpha}^{\{\mathrm{S,S}\}}$ strictly lower triangular, and matrices $\boldsymbol{\gamma}^{\{\mathrm{F,F},\ell\}}$ and $\boldsymbol{\gamma}^{\{\mathrm{S,S}\}}$ lower triangular. From the Butcher tableau (7.5) compute the coupling structure matrix

$$
\mathbf{S} := \Big(|\boldsymbol{\alpha}^{\{\mathrm{S,F}\}}| + |\boldsymbol{\gamma}^{\{\mathrm{S,F}\}}|\Big)^{\top} \times \Big(|\boldsymbol{\alpha}^{\{\mathrm{F,S}\}}| + |\boldsymbol{\gamma}^{\{\mathrm{F,S}\}}|\Big) \in \mathbb{R}^{\mathtt{m}s^{\{\mathrm{F}\}}\times s^{\{\mathrm{S}\}}}, \tag{7.7}
$$

where $|\cdots|$ takes element-wise absolute values, and $\times$ is the element-wise product. We make the following observations:

- In order to compute stages sequentially, in a completely decoupled manner, it must hold that $\mathbf{S} = \mathbf{0}$.
- The non-zero entries in this matrix $\mathbf{S}$ correspond to slow and fast stages that are computed together, in a coupled manner. Specifically, if element in row $(\ell, i)$ and column $j$ is non-zero then stages $\mathbf{k}_i^{\{\mathrm{F},\ell\}}$ and $\mathbf{k}_j^{\{\mathrm{S}\}}$ are computed by solving a joint linear system.
- Defining $\mathbf{1}^{\{q\}} \in \mathbb{R}^{s^{\{q\}}}$ as vectors of ones and $\mathbf{g}^{\{m,q\}} := \boldsymbol{\gamma}^{\{m,q\}}\,\mathbf{1}^{\{q\}}$ for $m, q \in \{\mathrm{F,S}\}$, the internal consistency condition $\mathbf{g}^{\{m,q\}} = \mathbf{g}^{\{m\}}$ for $\mathbf{g}$ requires that at least one slow and one fast stage are computed together in a coupled manner.

Example 7 (Second order, two-rate method). *Consider the following example using* $\mathtt{m} = 2$ *and two-stage base methods:*

$$
\begin{array}{cc:cc|cc}
0 & 0 & 0 & 0 & 0 & 0 \\
\frac{1}{2}\,\alpha_{2,1}^{\{\mathrm{F,F},1\}} & 0 & 0 & 0 & \alpha_{2,1}^{\{\mathrm{F,S},1\}} & \alpha_{2,2}^{\{\mathrm{F,S},1\}} \\ \hdashline
\frac{1}{2}\,b_1^{\{\mathrm{F},1\}} & \frac{1}{2}\,b_2^{\{\mathrm{F},1\}} & 0 & 0 & \alpha_{1,1}^{\{\mathrm{F,S},2\}} & \alpha_{1,2}^{\{\mathrm{F,S},2\}} \\
\frac{1}{2}\,b_1^{\{\mathrm{F},1\}} & \frac{1}{2}\,b_2^{\{\mathrm{F},1\}} & \frac{1}{2}\,\alpha_{2,1}^{\{\mathrm{F,F},2\}} & 0 & \alpha_{2,1}^{\{\mathrm{F,S},2\}} & \alpha_{2,2}^{\{\mathrm{F,S},2\}} \\ \hline
0 & 0 & 0 & 0 & 0 & 0 \\
\frac{1}{2}\,\alpha_{2,1}^{\{\mathrm{S,F},1\}} & \frac{1}{2}\,\alpha_{2,2}^{\{\mathrm{S,F},1\}} & \frac{1}{2}\,\alpha_{2,1}^{\{\mathrm{S,F},2\}} & \frac{1}{2}\,\alpha_{2,2}^{\{\mathrm{S,F},2\}} & \alpha_{2,1}^{\{\mathrm{S,S}\}} & 0 \\ \hline
\frac{1}{2}\,b_1^{\{\mathrm{F},1\}} & \frac{1}{2}\,b_2^{\{\mathrm{F},1\}} & \frac{1}{2}\,b_1^{\{\mathrm{F},2\}} & \frac{1}{2}\,b_2^{\{\mathrm{F},2\}} & b_1^{\{\mathrm{S}\}} & b_2^{\{\mathrm{S}\}}
\end{array}\,,
\qquad
\begin{array}{cc:cc|cc}
\frac{1}{2}\,\gamma^{\{\mathrm{F,F},1\}} & 0 & 0 & 0 & \gamma_{1,1}^{\{\mathrm{F,S},1\}} & \gamma_{1,2}^{\{\mathrm{F,S},1\}} \\
\frac{1}{2}\,\gamma_{2,1}^{\{\mathrm{F,F},1\}} & \frac{1}{2}\,\gamma^{\{\mathrm{F,F},1\}} & 0 & 0 & \gamma_{2,1}^{\{\mathrm{F,S},1\}} & \gamma_{2,2}^{\{\mathrm{F,S},1\}} \\ \hdashline
0 & 0 & \frac{1}{2}\,\gamma^{\{\mathrm{F,F},2\}} & 0 & \gamma_{1,1}^{\{\mathrm{F,S},2\}} & \gamma_{1,2}^{\{\mathrm{F,S},2\}} \\
0 & 0 & \frac{1}{2}\,\gamma_{2,1}^{\{\mathrm{F,F},2\}} & \frac{1}{2}\,\gamma^{\{\mathrm{F,F},2\}} & \gamma_{2,1}^{\{\mathrm{F,S},2\}} & \gamma_{2,2}^{\{\mathrm{F,S},2\}} \\ \hline
\frac{1}{2}\,\gamma_{1,1}^{\{\mathrm{S,F},1\}} & \frac{1}{2}\,\gamma_{1,2}^{\{\mathrm{S,F},1\}} & \frac{1}{2}\,\gamma_{1,1}^{\{\mathrm{S,F},2\}} & \frac{1}{2}\,\gamma_{1,2}^{\{\mathrm{S,F},2\}} & \gamma^{\{\mathrm{S,S}\}} & 0 \\
\frac{1}{2}\,\gamma_{2,1}^{\{\mathrm{S,F},1\}} & \frac{1}{2}\,\gamma_{2,2}^{\{\mathrm{S,F},1\}} & \frac{1}{2}\,\gamma_{2,1}^{\{\mathrm{S,F},2\}} & \frac{1}{2}\,\gamma_{2,2}^{\{\mathrm{S,F},2\}} & \gamma_{2,1}^{\{\mathrm{S,S}\}} & \gamma^{\{\mathrm{S,S}\}}
\end{array}\,.
$$

*We conveniently choose $\alpha^{\{\mathrm{S,F},\ell\}}_{1,j} = \alpha^{\{\mathrm{F,S},\ell\}}_{1,j} = 0$ for all $j$ and $\ell$ such as to satisfy the first internal consistency conditions. The coupling structure matrix* (7.7) *is:*

$$
\mathbf{S} = \tfrac{1}{2}\begin{bmatrix}
|\gamma^{\{\mathrm{S,F},1\}}_{1,1}|\cdot|\gamma^{\{\mathrm{F,S},1\}}_{1,1}| & (|\alpha^{\{\mathrm{S,F},1\}}_{2,1}|+|\gamma^{\{\mathrm{S,F},1\}}_{2,1}|)\cdot|\gamma^{\{\mathrm{F,S},1\}}_{1,2}| \\
|\gamma^{\{\mathrm{S,F},1\}}_{1,2}|\cdot(|\alpha^{\{\mathrm{F,S},1\}}_{2,1}|+|\gamma^{\{\mathrm{F,S},1\}}_{2,1}|) & (|\alpha^{\{\mathrm{S,F},1\}}_{2,2}|+|\gamma^{\{\mathrm{S,F},1\}}_{2,2}|)\cdot(|\alpha^{\{\mathrm{F,S},1\}}_{2,2}|+|\gamma^{\{\mathrm{F,S},1\}}_{2,2}|) \\
|\gamma^{\{\mathrm{S,F},2\}}_{1,1}|\cdot(|\alpha^{\{\mathrm{F,S},2\}}_{1,1}|+|\gamma^{\{\mathrm{F,S},2\}}_{1,1}|) & (|\gamma^{\{\mathrm{S,F},2\}}_{2,1}|+|\alpha^{\{\mathrm{S,F},2\}}_{2,1}|)\cdot(|\alpha^{\{\mathrm{F,S},2\}}_{1,2}|+|\gamma^{\{\mathrm{F,S},2\}}_{1,2}|) \\
|\gamma^{\{\mathrm{S,F},2\}}_{1,2}|\cdot(|\alpha^{\{\mathrm{F,S},2\}}_{2,1}|+|\gamma^{\{\mathrm{F,S},2\}}_{2,1}|) & (|\gamma^{\{\mathrm{S,F},2\}}_{2,2}|+|\alpha^{\{\mathrm{S,F},2\}}_{2,2}|)\cdot(|\alpha^{\{\mathrm{F,S},2\}}_{2,2}|+|\gamma^{\{\mathrm{F,S},2\}}_{2,2}|)
\end{bmatrix}. \tag{7.8}
$$

In the following we briefly discuss some coupling approaches for an efficient implementation of GARK-ROS/ROW schemes. For more details we refer to [28].

**7.2. IMEX approach.** If one chooses $\boldsymbol{\gamma}^{\{\mathrm{F,F},\ell\}} = \mathbf{0}$ and $\boldsymbol{\gamma}^{\{\mathrm{F,S},\ell\}} = \mathbf{0}$ then the fast component is integrated with a Runge-Kutta method; this method is explicit if $\boldsymbol{\alpha}^{\{\mathrm{F,F},\ell\}}$ are strictly lower triangular matrices. For a decoupled computation one needs to select the coupling coefficients $\boldsymbol{\alpha}^{\{\mathrm{F,S},\ell\}}$, $\boldsymbol{\alpha}^{\{\mathrm{S,F},\ell\}}$, and $\boldsymbol{\gamma}^{\{\mathrm{S,F},\ell\}}$ such that the matrix $\mathbf{S} = \mathbf{0}$. The fast stages are computed as:

$$
\begin{aligned}
\mathbf{y}\widetilde{\mathbf{y}}^{\{\mathrm{F}\}}_{n+(\ell-1)/\mathtt{m}} &= \mathbf{y}_{n-1} + \sum_{\ell=1}^{\ell-1}\sum_{j=1}^{s^{\{\mathrm{F}\}}} b^{\{\mathrm{F},\ell\}}_{i,j}\,\mathbf{k}^{\{\mathrm{F},\ell\}}_j, \\
\mathbf{k}^{\{\mathrm{F},\ell\}}_i &= h\,\mathbf{f}^{\{\mathrm{F}\}}\Big(\mathbf{y}\widetilde{\mathbf{y}}^{\{\mathrm{F}\}}_{n-1+(\ell-1)/\mathtt{m}} + \sum_{j=1}^{s^{\{\mathrm{F}\}}} \alpha^{\{\mathrm{F,F},\ell\}}_{i,j}\,\mathbf{k}^{\{\mathrm{F},\ell\}}_j, \mathbf{y}^{\{\mathrm{S}\}}_n + \sum_{j=1}^{s^{\{\mathrm{S}\}}} \alpha^{\{\mathrm{F,S},\ell\}}_{i,j}\,\mathbf{k}^{\{\mathrm{S}\}}_j\Big).
\end{aligned}
$$

The corresponding Butcher tableau (7.5) is:

$$
\begin{array}{cccc|c}
\frac{1}{\mathtt{m}}\boldsymbol{\alpha}^{\{\mathrm{F,F},1\}} & 0 & \cdots & 0 & \boldsymbol{\alpha}^{\{\mathrm{F,S},1\}} \\
\frac{1}{\mathtt{m}}\mathbf{1b}^{\{\mathrm{F},1\}T} & \frac{1}{\mathtt{m}}\boldsymbol{\alpha}^{\{\mathrm{F,F},2\}} & \cdots & 0 & \boldsymbol{\alpha}^{\{\mathrm{F,S},2\}} \\
\vdots & & \ddots & & \vdots \\
\frac{1}{\mathtt{m}}\mathbf{1b}^{\{\mathrm{F},1\}T} & \frac{1}{\mathtt{m}}\mathbf{1b}^{\{\mathrm{F},2\}T} & \cdots & \frac{1}{\mathtt{m}}\boldsymbol{\alpha}^{\{\mathrm{F,F},\mathtt{m}\}} & \boldsymbol{\alpha}^{\{\mathrm{F,S},\mathtt{m}\}} \\
\hline
\frac{1}{\mathtt{m}}\boldsymbol{\alpha}^{\{\mathrm{S,F},1\}} & \frac{1}{\mathtt{m}}\boldsymbol{\alpha}^{\{\mathrm{S,F},2\}} & \cdots & \frac{1}{\mathtt{m}}\boldsymbol{\alpha}^{\{\mathrm{S,F},\mathtt{m}\}} & \boldsymbol{\alpha}^{\{\mathrm{S,S}\}}
\end{array}, \quad
\begin{array}{cccc|c}
0 & 0 & \cdots & 0 & 0 \\
0 & 0 & \cdots & 0 & 0 \\
\vdots & \vdots & \ddots & \vdots & \vdots \\
0 & 0 & \cdots & 0 & 0 \\
\hline
\frac{1}{\mathtt{m}}\boldsymbol{\gamma}^{\{\mathrm{S,F},1\}} & \frac{1}{\mathtt{m}}\boldsymbol{\gamma}^{\{\mathrm{S,F},2\}} & \cdots & \frac{1}{\mathtt{m}}\boldsymbol{\gamma}^{\{\mathrm{S,F},\mathtt{m}\}} & \boldsymbol{\gamma}^{\{\mathrm{S,S}\}}
\end{array}. \tag{7.9}
$$

**7.3. Compound-first-step approach: coupling the macro-step with the first micro-step.** We consider base methods with the same number of stages $s^{\{\mathrm{F}\}} = s^{\{\mathrm{S}\}} = s$. Moreover, we set the coupling coefficients $\boldsymbol{\alpha}^{\{\mathrm{S,F},\ell\}} = \boldsymbol{\gamma}^{\{\mathrm{S,F},\ell\}} = 0$ for $\ell = 2, \dots, \mathtt{m}$. The resulting Butcher tableau (7.5) is:

$$
\begin{array}{cccc|c}
\frac{1}{\mathtt{m}}\boldsymbol{\alpha}^{\{\mathrm{F,F},1\}} & 0 & \cdots & 0 & \boldsymbol{\alpha}^{\{\mathrm{F,S},1\}} \\
\frac{1}{\mathtt{m}}\mathbf{1b}^{\{\mathrm{F},1\}T} & \frac{1}{\mathtt{m}}\boldsymbol{\alpha}^{\{\mathrm{F,F},2\}} & \cdots & 0 & \boldsymbol{\alpha}^{\{\mathrm{F,S},2\}} \\
\vdots & \vdots & \ddots & \vdots & \vdots \\
\frac{1}{\mathtt{m}}\mathbf{1b}^{\{\mathrm{F},1\}T} & \frac{1}{\mathtt{m}}\mathbf{1b}^{\{\mathrm{F},2\}T} & \cdots & \frac{1}{\mathtt{m}}\boldsymbol{\alpha}^{\{\mathrm{F,F},\mathtt{m}\}} & \boldsymbol{\alpha}^{\{\mathrm{F,S},\mathtt{m}\}} \\
\hline
\frac{1}{\mathtt{m}}\boldsymbol{\alpha}^{\{\mathrm{S,F},1\}} & 0 & \cdots & 0 & \boldsymbol{\alpha}^{\{\mathrm{S,S}\}}
\end{array}, \quad
\begin{array}{cccc|c}
\frac{1}{\mathtt{m}}\boldsymbol{\gamma}^{\{\mathrm{F,F},1\}} & 0 & \cdots & 0 & \boldsymbol{\gamma}^{\{\mathrm{F,S},1\}} \\
0 & \frac{1}{\mathtt{m}}\boldsymbol{\gamma}^{\{\mathrm{F,F},2\}} & \cdots & 0 & \boldsymbol{\gamma}^{\{\mathrm{F,S},2\}} \\
\vdots & \vdots & \ddots & \vdots & \vdots \\
0 & 0 & \cdots & \frac{1}{\mathtt{m}}\boldsymbol{\gamma}^{\{\mathrm{F,F},\mathtt{m}\}} & \boldsymbol{\gamma}^{\{\mathrm{F,S},\mathtt{m}\}} \\
\hline
\frac{1}{\mathtt{m}}\boldsymbol{\gamma}^{\{\mathrm{S,F},1\}} & 0 & \cdots & 0 & \boldsymbol{\gamma}^{\{\mathrm{S,S}\}}
\end{array}. \tag{7.10}
$$

The coefficient matrices $\alpha^{\{\mathrm{F,F},1\}}$, $\alpha^{\{\mathrm{F,S},1\}}$, $\alpha^{\{\mathrm{S,F},1\}}$, $\alpha^{\{\mathrm{S,S}\}}$ are chosen strictly lower triangular. The coefficient matrices $\gamma^{\{\mathrm{F,F},1\}}$, $\gamma^{\{\mathrm{F,S},1\}}$, $\gamma^{\{\mathrm{S,F},1\}}$, and $\gamma^{\{\mathrm{S,S}\}}$ are chosen lower triangular, with equal diagonal entries: $\gamma^{\{\mathrm{S,S}\}}_{i,i} := \gamma^{\{\mathrm{S,S}\}}, \gamma^{\{\mathrm{F,F},1\}}_{i,i} := \gamma^{\{\mathrm{F,F},1\}}, \gamma^{\{\mathrm{S,F},1\}}_{i,i} := \gamma^{\{\mathrm{S,F},1\}}, \gamma^{\{\mathrm{F,S},1\}}_{i,i} := \gamma^{\{\mathrm{F,S},1\}}$ for $i = 1, \dots, s$.

In the structure matrix (7.7) the entries $\mathbf{S}_{i,i} \neq 0$ for $i = 1, \dots, s^{\{\mathrm{S}\}}$. This means that each slow stage $\mathbf{k}_i^{\{\mathrm{S}\}}$ and the corresponding fast stage of the first micro-step $\mathbf{k}_i^{\{\mathrm{F},1\}}$ are computed in a coupled manner, by solving a full coupled system of linear equations. Since all diagonal entries are equal to each other, only one LU decomposition of the compound matrix is necessary for computing all stage vectors $\mathbf{k}_i^{\{\mathrm{S}\}}$ and $\mathbf{k}_i^{\{\mathrm{F},1\}}$ for $i = 1, \dots, s$.

Since all slow stages are known after the first micro-step, $\boldsymbol{\alpha}^{\{\mathrm{F},\mathrm{S},\ell\}}$ and $\boldsymbol{\gamma}^{\{\mathrm{F},\mathrm{S},\ell\}}$ can be full matrices for $\ell = 2, \dots, \mathtt{m}$. For all remaining micro steps a single additional LU decomposition is necessary if $\gamma_{i,i}^{\{\mathrm{F},\mathrm{F},\ell\}} := \gamma^{\{\mathrm{F},\mathrm{F}\}}$ is constant for all $i = 1, \dots, s^{\{\mathrm{F}\}}$ and $\ell = 2, \dots, \mathtt{m}$.

A simple choice of coefficients for compound-first-step coupling is $\alpha^{\{\mathrm{S},\mathrm{F},1\}} = \mathtt{m}\,\alpha^{\{\mathrm{S},\mathrm{S}\}}$, $\gamma^{\{\mathrm{S},\mathrm{F},1\}} = \mathtt{m}\,\gamma^{\{\mathrm{S},\mathrm{S}\}}$, $\alpha^{\{\mathrm{F},\mathrm{S},1\}} = (1/\mathtt{m})\,\alpha^{\{\mathrm{F},\mathrm{F},1\}}$, $\gamma^{\{\mathrm{F},\mathrm{S},1\}} = (1/\mathtt{m})\,\gamma^{\{\mathrm{F},\mathrm{F},1\}}$.

Example 8. *Consider the scheme of Example 7 with the following coefficients:*

$$
\begin{array}{cc:cc|cc}
0 & 0 & 0 & 0 & 0 & 0 \\
\frac{1}{2}\,\alpha_{2,1}^{\{\mathrm{F},\mathrm{F},1\}} & 0 & 0 & 0 & \alpha_{2,1}^{\{\mathrm{F},\mathrm{S},1\}} & 0 \\
\hdashline
\frac{1}{2}\,b_1^{\{\mathrm{F},1\}} & \frac{1}{2}\,b_2^{\{\mathrm{F},1\}} & 0 & 0 & \alpha_{1,1}^{\{\mathrm{F},\mathrm{S},2\}} & \alpha_{1,2}^{\{\mathrm{F},\mathrm{S},2\}} \\
\frac{1}{2}\,b_1^{\{\mathrm{F},1\}} & \frac{1}{2}\,b_2^{\{\mathrm{F},1\}} & \frac{1}{2}\,\alpha_{2,1}^{\{\mathrm{F},\mathrm{F},2\}} & 0 & \alpha_{2,1}^{\{\mathrm{F},\mathrm{S},2\}} & \alpha_{2,2}^{\{\mathrm{F},\mathrm{S},2\}} \\
\hline
0 & 0 & 0 & 0 & 0 & 0 \\
\frac{1}{2}\,\alpha_{2,1}^{\{\mathrm{S},\mathrm{F},1\}} & 0 & 0 & 0 & \alpha_{2,1}^{\{\mathrm{S},\mathrm{S}\}} & 0
\end{array}\,,
\qquad
\begin{array}{cc:cc|cc}
\frac{1}{2}\,\gamma^{\{\mathrm{F},\mathrm{F},1\}} & 0 & 0 & 0 & \gamma^{\{\mathrm{F},\mathrm{S},1\}} & 0 \\
\frac{1}{2}\,\gamma_{2,1}^{\{\mathrm{F},\mathrm{F},1\}} & \frac{1}{2}\,\gamma^{\{\mathrm{F},\mathrm{F},1\}} & 0 & 0 & \gamma_{2,1}^{\{\mathrm{F},\mathrm{S},1\}} & \gamma^{\{\mathrm{F},\mathrm{S},1\}} \\
\hdashline
0 & 0 & \frac{1}{2}\,\gamma^{\{\mathrm{F},\mathrm{F},2\}} & 0 & \gamma_{1,1}^{\{\mathrm{F},\mathrm{S},2\}} & \gamma_{1,2}^{\{\mathrm{F},\mathrm{S},2\}} \\
0 & 0 & \frac{1}{2}\,\gamma_{2,1}^{\{\mathrm{F},\mathrm{F},2\}} & \frac{1}{2}\,\gamma^{\{\mathrm{F},\mathrm{F},2\}} & \gamma_{2,1}^{\{\mathrm{F},\mathrm{S},2\}} & \gamma_{2,2}^{\{\mathrm{F},\mathrm{S},2\}} \\
\hline
\frac{1}{2}\,\gamma^{\{\mathrm{S},\mathrm{F},1\}} & 0 & 0 & 0 & \gamma^{\{\mathrm{S},\mathrm{S}\}} & 0 \\
\frac{1}{2}\,\gamma_{2,1}^{\{\mathrm{S},\mathrm{F},1\}} & \frac{1}{2}\,\gamma^{\{\mathrm{S},\mathrm{F},1\}} & 0 & 0 & \gamma_{2,1}^{\{\mathrm{S},\mathrm{S}\}} & \gamma^{\{\mathrm{S},\mathrm{S}\}}
\end{array}\,.
$$

*The coupling matrix (7.8) reads:*

$$
\mathbf{S} = \frac{1}{2}\begin{bmatrix}
|\gamma^{\{\mathrm{S},\mathrm{F},1\}}|\cdot|\gamma^{\{\mathrm{F},\mathrm{S},1\}}| & 0 \\
0 & |\gamma^{\{\mathrm{S},\mathrm{F},1\}}|\cdot|\gamma^{\{\mathrm{F},\mathrm{S},1\}}| \\
0 & 0 \\
0 & 0
\end{bmatrix},
$$

*which indicates that $\mathbf{k}_1^{\{\mathrm{F},1\}}$ and $\mathbf{k}_1^{\{\mathrm{S}\}}$ are computed together, and so are $\mathbf{k}_2^{\{\mathrm{F},1\}}$ and $\mathbf{k}_2^{\{\mathrm{S}\}}$.*

Remark 23. *The multirate ROW schemes introduced by Bartel and Günther [4] fall into the class of multirate GARK-ROW schemes. They consider the case of time-lagged Jacobians (which differ by a term of magnitude $\mathcal{O}(H)$ from the exact Jacobian). In addition, the same order $p$ within the micro steps is demanded.*

**7.4. Coupling only the first fast and the first slow stage computations.** The smallest amount of coupling that allows the construction of internally consistent implicit schemes is a lighter version of the strategy discussed in section 7.3, where only the first fast stage $\mathbf{k}_1^{\{\mathrm{F},1\}}$ and the first slow stage $\mathbf{k}_1^{\{\mathrm{S}\}}$ are computed together. It is possible to select coefficients such that all subsequent stage computations are implicit in either fast or slow stages, and are computed in a decoupled manner.

Example 9. *Consider the scheme from Example 7 with the following coefficients:*

$$
\begin{array}{cc:cc|cc}
0 & 0 & 0 & 0 & 0 & 0 \\
\frac{1}{2}\,\alpha_{2,1}^{\{\mathrm{F},\mathrm{F},1\}} & 0 & 0 & 0 & \alpha_{2,1}^{\{\mathrm{F},\mathrm{S},1\}} & 0 \\
\hdashline
\frac{1}{2}\,b_1^{\{\mathrm{F},1\}} & \frac{1}{2}\,b_2^{\{\mathrm{F},1\}} & 0 & 0 & \alpha_{1,1}^{\{\mathrm{F},\mathrm{S},2\}} & 0 \\
\frac{1}{2}\,b_1^{\{\mathrm{F},1\}} & \frac{1}{2}\,b_2^{\{\mathrm{F},1\}} & \frac{1}{2}\,\alpha_{2,1}^{\{\mathrm{F},\mathrm{F},2\}} & 0 & \alpha_{2,1}^{\{\mathrm{F},\mathrm{S},2\}} & \alpha_{2,2}^{\{\mathrm{F},\mathrm{S},2\}} \\
\hline
0 & 0 & 0 & 0 & 0 & 0 \\
\frac{1}{2}\,\alpha_{2,1}^{\{\mathrm{S},\mathrm{F},1\}} & \frac{1}{2}\,\alpha_{2,2}^{\{\mathrm{S},\mathrm{F},1\}} & \frac{1}{2}\,\alpha_{2,1}^{\{\mathrm{S},\mathrm{F},2\}} & 0 & \alpha_{2,1}^{\{\mathrm{S},\mathrm{S}\}} & 0
\end{array}\,,
\qquad
\begin{array}{cc:cc|cc}
\frac{1}{2}\,\gamma^{\{\mathrm{F},\mathrm{F},1\}} & 0 & 0 & 0 & \gamma_{1,1}^{\{\mathrm{F},\mathrm{S},1\}} & 0 \\
\frac{1}{2}\,\gamma_{2,1}^{\{\mathrm{F},\mathrm{F},1\}} & \frac{1}{2}\,\gamma^{\{\mathrm{F},\mathrm{F},1\}} & 0 & 0 & \gamma_{2,1}^{\{\mathrm{F},\mathrm{S},1\}} & 0 \\
\hdashline
0 & 0 & \frac{1}{2}\,\gamma^{\{\mathrm{F},\mathrm{F},2\}} & 0 & \gamma_{1,1}^{\{\mathrm{F},\mathrm{S},2\}} & 0 \\
0 & 0 & \frac{1}{2}\,\gamma_{2,1}^{\{\mathrm{F},\mathrm{F},2\}} & \frac{1}{2}\,\gamma^{\{\mathrm{F},\mathrm{F},2\}} & \gamma_{2,1}^{\{\mathrm{F},\mathrm{S},2\}} & \gamma_{2,2}^{\{\mathrm{F},\mathrm{S},2\}} \\
\hline
\frac{1}{2}\,\gamma_{1,1}^{\{\mathrm{S},\mathrm{F},1\}} & 0 & 0 & 0 & \gamma^{\{\mathrm{S},\mathrm{S}\}} & 0 \\
\frac{1}{2}\,\gamma_{2,1}^{\{\mathrm{S},\mathrm{F},1\}} & \frac{1}{2}\,\gamma_{2,2}^{\{\mathrm{S},\mathrm{F},1\}} & \frac{1}{2}\,\gamma_{2,1}^{\{\mathrm{S},\mathrm{F},2\}} & 0 & \gamma_{2,1}^{\{\mathrm{S},\mathrm{S}\}} & \gamma^{\{\mathrm{S},\mathrm{S}\}}
\end{array}\,.
$$

*The coupling matrix (7.8) has single non-zero element, $\frac{1}{2}\,|\gamma_{1,1}^{\{\mathrm{S},\mathrm{F},1\}}|\,|\gamma_{1,1}^{\{\mathrm{F},\mathrm{S},1\}}|$, corresponding to first computing stages $\mathbf{k}_1^{\{\mathrm{F},1\}}$ and $\mathbf{k}_1^{\{\mathrm{S}\}}$ in a coupled manner. Next, $\mathbf{k}_2^{\{\mathrm{F},1\}}$*

*and* $\mathbf{k}_1^{\{\text{F},2\}}$ *are computed in a decoupled manner, since they only depend on the known slow stage* $\mathbf{k}_1^{\{\text{S}\}}$. *After this,* $\mathbf{k}_2^{\{\text{S}\}}$ *is computed in a decoupled manner as it does not depend on the (yet unknown) last fast stage* $\mathbf{k}_2^{\{\text{F},2\}}$. *Finally,* $\mathbf{k}_2^{\{\text{F},2\}}$ *is evaluated using both slow stages.*

REMARK 24. *Note that the coupling approaches in Section 7.2 and 7.3 consider only the influence of the fast scale within the first micro step onto the slow variables. If the evolution of the fast part changes drastically after the first micro-step, this would only affect the slow part in the next macro step. However, the order conditions show that the overall scheme has order p at the macro step grid points. In addition, as a rule of thumb, the fast part has typically only a weak influence on the slow part, as otherwise the slow part might incorporate some fast dynamics, and require to be treated as fast, too. The compound-step approach idea traces back to Kvaerno and Rentrop [37] in 1999, and to our knowledge, no time-lag effect has been observed since then in any multirate scheme based on the compound-step approach.*

**7.5. Fully decoupled approach.** In the completely decoupled approach each stage follows a regular Rosenbrock computation, implicit in either the fast or the slow stages, but *not* in both at the same time. In this case the second internal consistency conditions $\gamma^{\{\text{F},\text{F},1\}} = \sum_j \gamma_{1,j}^{\{\text{F},\text{S},1\}}$ do not hold unless this first stage is explicit. Consequently, the coupling order conditions for the entire method become more complex, but such methods are possible to construct.

EXAMPLE 10. *Consider the scheme from Example 7 with the following coefficients:*

$$
\begin{array}{cc:cc|cc}
0 & 0 & 0 & 0 & 0 & 0 \\
\frac{1}{2}\,\alpha_{2,1}^{\{\text{F},\text{F},1\}} & 0 & 0 & 0 & \alpha_{2,1}^{\{\text{F},\text{S},1\}} & 0 \\
\hdashline
\frac{1}{2}\,b_1^{\{\text{F},1\}} & \frac{1}{2}\,b_2^{\{\text{F},1\}} & 0 & 0 & \alpha_{1,1}^{\{\text{F},\text{S},2\}} & 0 \\
\frac{1}{2}\,b_1^{\{\text{F},1\}} & \frac{1}{2}\,b_2^{\{\text{F},1\}} & \frac{1}{2}\,\alpha_{2,1}^{\{\text{F},\text{F},2\}} & 0 & \alpha_{2,1}^{\{\text{F},\text{S},2\}} & \alpha_{2,2}^{\{\text{F},\text{S},2\}} \\
\hline
\frac{1}{2}\,\alpha_{1,1}^{\{\text{S},\text{F},1\}} & 0 & 0 & 0 & 0 & 0 \\
\frac{1}{2}\,\alpha_{2,1}^{\{\text{S},\text{F},1\}} & \frac{1}{2}\,\alpha_{2,2}^{\{\text{S},\text{F},1\}} & \frac{1}{2}\,\alpha_{2,1}^{\{\text{S},\text{F},2\}} & 0 & \alpha_{2,1}^{\{\text{S},\text{S}\}} & 0
\end{array},
\quad
\begin{array}{cc:cc|cc}
\frac{1}{2}\,\gamma^{\{\text{F},\text{F},1\}} & 0 & 0 & 0 & 0 & 0 \\
\frac{1}{2}\,\gamma_{2,1}^{\{\text{F},\text{F},1\}} & \frac{1}{2}\,\gamma^{\{\text{F},\text{F},1\}} & 0 & 0 & \gamma_{2,1}^{\{\text{F},\text{S},1\}} & 0 \\
\hdashline
0 & 0 & \frac{1}{2}\,\gamma^{\{\text{F},\text{F},2\}} & 0 & \gamma_{1,1}^{\{\text{F},\text{S},2\}} & 0 \\
0 & 0 & \frac{1}{2}\,\gamma_{2,1}^{\{\text{F},\text{F},2\}} & \frac{1}{2}\,\gamma^{\{\text{F},\text{F},2\}} & \gamma_{2,1}^{\{\text{F},\text{S},2\}} & \gamma_{2,2}^{\{\text{F},\text{S},2\}} \\
\hline
\frac{1}{2}\,\gamma_{1,1}^{\{\text{S},\text{F},1\}} & 0 & 0 & 0 & \gamma^{\{\text{S},\text{S}\}} & 0 \\
\frac{1}{2}\,\gamma_{2,1}^{\{\text{S},\text{F},1\}} & \frac{1}{2}\,\gamma_{2,2}^{\{\text{S},\text{F},1\}} & \frac{1}{2}\,\gamma_{2,1}^{\{\text{S},\text{F},2\}} & 0 & \gamma_{2,1}^{\{\text{S},\text{S}\}} & \gamma^{\{\text{S},\text{S}\}}
\end{array}.
$$

*Note the complementary sparsity structure of the off-diagonal coupling blocks. The first fast stage is that of a classical Rosenbrock method:*

$$
\left(\mathbf{I} - h\,\gamma^{\{\text{F},\text{F},1\}}\,\mathbf{L}^{\{\text{F},\text{F}\}}\right)\mathbf{k}_1^{\{\text{F},1\}} \quad = \quad h\,\mathbf{f}^{\{\text{F}\}}(\mathbf{y}_n^{\{\text{F}\}}, \mathbf{y}_n^{\{\text{S}\}}).
$$

*Similarly, the first slow stage is computed in a decoupled manner:*

$$
\left(\mathbf{I} - H\,\gamma^{\{\text{S},\text{S}\}}\,\mathbf{L}^{\{\text{S},\text{S}\}}\right)\mathbf{k}_1^{\{\text{S}\}} \quad = \quad H\,\mathbf{f}^{\{\text{S}\}}\left(\mathbf{y}_n^{\{\text{F}\}} + \alpha_{1,1}^{\{\text{S},\text{F},1\}}\,\mathbf{k}_1^{\{\text{F},1\}}, \mathbf{y}_n^{\{\text{S}\}}\right) + H\,\gamma_{1,1}^{\{\text{S},\text{F},1\}}\,\mathbf{L}^{\{\text{S},\text{F}\}}\,\mathbf{k}_1^{\{\text{F},1\}},
$$

*and the decoupled computations continue alternating fast and slow stages.*

**7.6. Step-predictor-corrector approach.** This approach starts with a "predictor" step where the slow Rosenbrock method is applied with step size $H$ to the entire system, in a classical fashion. The slow components are sufficiently accurate, but the fast components are not; for this reason we keep only the computed $\mathbf{k}^{\{\text{S}\}}$, but discard $\mathbf{k}^{\{\text{F}\}}$. The "corrector" re-computes $\mathbf{k}^{\{\text{F},\ell\}}$ for all sub-steps $\ell$, with the small steps sizes $h$, and uses these values to construct the final solution. The Butcher

tableaus (7.5) read:

$$
(7.11)\qquad
\begin{array}{cccc|c}
\boldsymbol{\alpha}^{\{\mathrm{S,S}\}} & 0 & \cdots & 0 & \boldsymbol{\alpha}^{\{\mathrm{S,S}\}} \\
0 & \frac{1}{\mathtt{m}}\boldsymbol{\alpha}^{\{\mathrm{F,F}\}} & \cdots & 0 & \boldsymbol{\alpha}^{\{\mathrm{F,S},1\}} \\
\vdots & \vdots & \ddots & \vdots & \vdots \\
0 & \frac{1}{\mathtt{m}}\mathbf{1}\mathbf{b}^{\{\mathrm{F}\}T} & \cdots & \frac{1}{\mathtt{m}}\boldsymbol{\alpha}^{\{\mathrm{F,F}\}} & \boldsymbol{\alpha}^{\{\mathrm{F,S},\mathtt{m}\}} \\
\hline
\boldsymbol{\alpha}^{\{\mathrm{S,S}\}} & 0 & \cdots & 0 & \boldsymbol{\alpha}^{\{\mathrm{S,S}\}} \\
\hline
0 & \frac{1}{\mathtt{m}}\mathbf{b}^{\{\mathrm{F}\}T} & \cdots & \frac{1}{\mathtt{m}}\mathbf{b}^{\{\mathrm{F}\}T} & \mathbf{b}^{\{\mathrm{S}\}T}
\end{array}\,,
\qquad
\begin{array}{cccc|c}
\boldsymbol{\gamma}^{\{\mathrm{S,S}\}} & 0 & \cdots & 0 & \boldsymbol{\gamma}^{\{\mathrm{S,S}\}} \\
0 & \frac{1}{\mathtt{m}}\boldsymbol{\gamma}^{\{\mathrm{F,F}\}} & \cdots & 0 & \boldsymbol{\gamma}^{\{\mathrm{F,S},1\}} \\
\vdots & \vdots & \ddots & \vdots & \vdots \\
0 & 0 & \cdots & \frac{1}{\mathtt{m}}\boldsymbol{\gamma}^{\{\mathrm{F,F}\}} & \boldsymbol{\gamma}^{\{\mathrm{F,S},\mathtt{m}\}} \\
\hline
\boldsymbol{\gamma}^{\{\mathrm{S,S}\}} & 0 & \cdots & 0 & \boldsymbol{\gamma}^{\{\mathrm{S,S}\}}
\end{array}\,.
$$

**8. Extrapolation Methods.** Here we are concerned with solving the partitioned system (2.1). We assume that the slow and fast partitions are relatively easy to obtain algorithmically.

We now construct multirate time stepping methods using extrapolation. The approach is to apply an Euler type discretization with a large step to $\mathbf{f}^{\{\mathrm{S}\}}$, an Euler discretization with a small to $\mathbf{f}^{\{\mathrm{F}\}}$, and achieve high orders of consistency by extrapolation [15, 16, 30, 31]. The discussion follows the work of Constantinescu and Sandu [11,12,12,13,13,49–51,51]. Extrapolation allows to easily construct multirate schemes of arbitrarily high orders, which are easy to implement and straightforward to parallelize.

**8.1. Global error expansions.** Consider a smooth system (1.1) and solve it numerically with a one-step "base" method, written in Henrici notation as follows:

$$\mathbf{y}_h(t+h) = \mathbf{y}_h(t) + h\,\Phi(\mathbf{y}_h(t), h). \tag{8.1}$$

Here $\mathbf{y}_h(t)$ is the numerical solution at time $t$ obtained by applying the method with a constant step size $h$ (where $t$ is an integer multiple of $h$). The method (8.1) is consistent with the differential equation (1.1) if the increment function $\Phi$ satisfies

$$\Phi(\mathbf{y}, 0) \equiv \mathbf{f}(\mathbf{y}) \quad \forall\, \mathbf{y}. \tag{8.2}$$

The following result from [23, 24, 29] establishes a fundamental property that allows extrapolation to work.

Theorem 8.1 (Asymptotic expansion of the global error [23,24,29]). *Consider a one-step method (8.1) with an increment function $\Phi$ that is smooth with respect to all of its arguments, and satisfies the consistency condition (8.2). The global error of the method (8.1) applied to a smooth system over the interval $[t_0, T]$ has an asymptotic expansion of the form [23, 24, 29, 30]*

$$\mathbf{y}(t) - \mathbf{y}_h(t) = e_p(t)\,h^p + \cdots + e_N(t)\,h^N + E_N(t,h)\,h^{N+1}, \tag{8.3}$$

*where $e_k(t)$ are smooth functions with $e_k(t_0) = 0$, and $E_N(t,h)$ is bounded for $t_0 \le t \le T$ and $0 \le h \le h_{\max}$.*

**8.2. Extrapolation procedure.** The extrapolation idea is to use multiple solutions computed with multiple time steps, and employ Richardson extrapolation to eliminate the dominant terms $e_k(t)h^k$ in the global error expansion (8.3). Specifically, consider a sequence of positive integers $\{n_j\}_{1\le j\le K}$, with $n_j < n_{j+1}$, and define a sequence of step sizes $\{h_j\}_{1\le j\le K}$ by $h_j = H/n_j$. Different numerical solutions $\{T_{j,1}\}_{1\le j\le K}$ at the final time $T = t_0 + H$ are obtained by applying $n_j$ steps of the base method with step size $h_j$:

$$T_{j,1} \coloneqq \mathbf{y}_{h_j}(t_0 + H), \quad 1 \le j \le K. \qquad \text{[Solutions of base method]} \tag{8.4}$$

Using the $K$ approximations (8.4) obtained with different $h_j$'s, one can eliminate the dominant $K$ error terms in the global error asymptotic expansion (8.3) via Richardson extrapolation (see [30, Ch. II.9]).

If the base method has order $p = 1$ the most economical approach to eliminate error terms is given by the Aitken-Neville formula [1, 20, 38]:

$$T_{j,k+1} = T_{j,k} + \frac{T_{j,k} - T_{j-1,k}}{(n_j/n_{j-k}) - 1}, \quad j \le K, \quad k < j. \tag{8.5}$$

The numerical scheme (8.4)–(8.5) is called the extrapolation method. Each application of (8.5) cancels the next term in the asymptotic error expansion (8.3); consequently, $T_{j,k}$ (8.5) are solutions of order $k$. It is customary to represent the solutions $T_{j,k}$ in a tableau. The workings of extrapolation are illustrated in Figure 12.

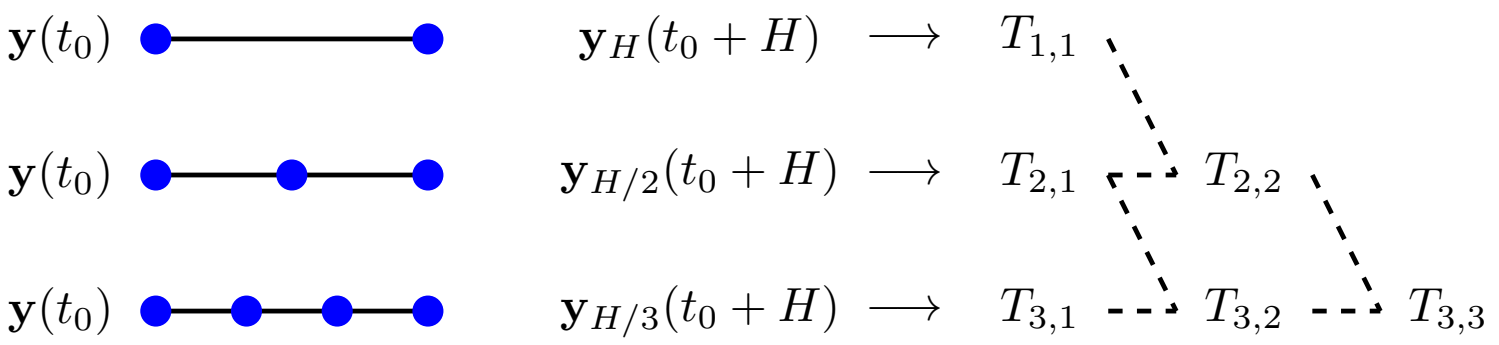

Fig. 12: Cartoon of the workings of extrapolation procedure for the sequence of steps $h_i = H/i$. Dashed lines represent the Aitken-Neville formula (8.5). We notice that, by design, extrapolation uses multiple time steps in the solution process.

Remark 25 (Extrapolation for symmetric base methods). *The asymptotic expansion of the global error* (8.3) *of a symmetric method [31] contains only even powers of the step size, and extrapolation yields numerical solutions* $T_{j,k}$ *of order* $2k$.

Remark 26 (Linear stability). *Let* $R_\Phi(z)$ *be the stability function/matrix of the base method* (8.1). *The stability functions of the extrapolation methods are calculated recursively from the extrapolation formula* (8.5) *[31, Chap. IV] where each* $T_{j,k}$ *is replaced with* $R_{j,k}(z)$.

Remark 27 (Extrapolation for stiff differential equations). *When the underlying differential equations are very stiff the global error does not enjoy an expansion with smooth coefficients* $e_p(t)$. *Rather,* (8.3) *is replaced by a perturbed asymptotic expansion [31, Section VI.5], and the order of the numerical solutions* $T_{j,k}$ *affected by these perturbations is smaller than* $k$.

**8.3. Multirate schemes by selective extrapolation: the MURX approach.** The first multirate extrapolation approach was proposed by Engstler and Lubich [18], and is based on adaptively selecting the level of extrapolation for different components. They call the resulting <u>mul</u>ti<u>r</u>ate e<u>x</u>trapolation implementation "MURX", and we will use this name for their approach.

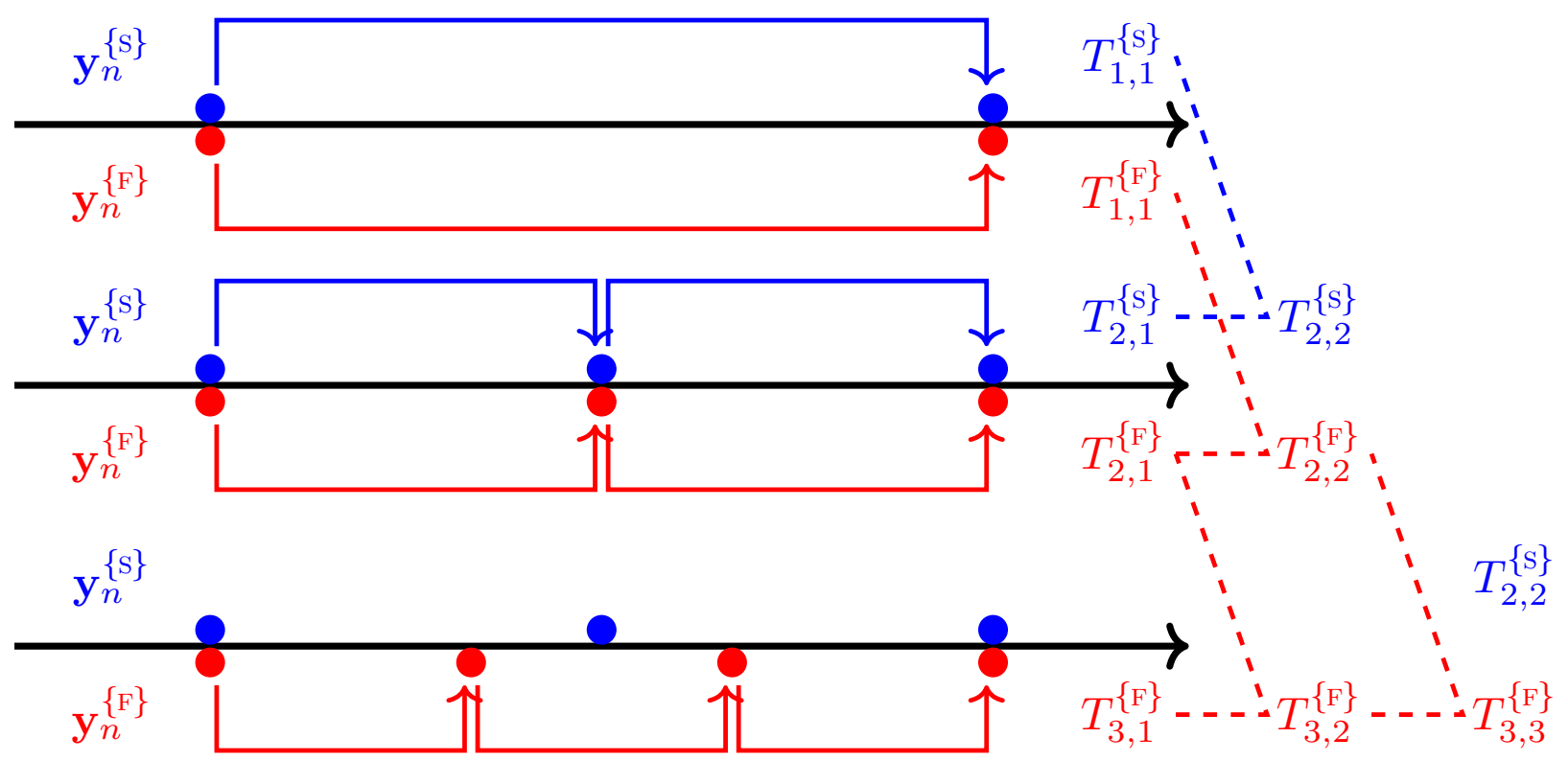

Fig. 13: Cartoon of MURX extrapolation method for $\mathtt{m} = 3$. Both the slow and the fast components are computed with step sizes $H$ and $H/2$, and are extrapolated to obtain the second column values $T_{2,2}^{\{\text{S}\}}$ and $T_{2,2}^{\{\text{F}\}}$, respectively. Dashed lines represent extrapolation operations. At this point the slow component $T_{2,2}^{\{\text{S}\}}$ is sufficiently accurate, and the slow components are not solved for any longer. The fast component is not sufficiently accurate and it is solved again with $H/3$ to obtain $T_{3,1}^{\{\text{F}\}}$; extrapolation is applied to fast component only to obtain $T_{3,2}^{\{\text{F}\}}$ and $T_{3,3}^{\{\text{F}\}}$. The final solution is $[\mathbf{y}_{n+1}^{\{\text{S}\}}, \mathbf{y}_{n+1}^{\{\text{F}\}}] = [T_{2,2}^{\{\text{S}\}}, T_{3,3}^{\{\text{F}\}}]$. The final slow component is computed using a smallest step size $H/2$, and with a second order discretization scheme. The final fast component is computed using a smallest step size $H/3$ (multirate), and with a third order discretization scheme (multiorder).

MURX works with the explicit Euler base method, and builds differently the extrapolation tableaus for fast and slow components. The workings of MURX are illustrated in Figure 13.

1. In the beginning all components are integrated simultaneously over a macro-step, with different step sizes, and the first rows of traditional extrapolation tableau for the entire system are constructed.
2. Once the slow components have been solved with sufficient accuracy no more integrations of the slow subsystem are carried out; the construction of the extrapolation tableau for the slow components is frozen. Specifically, after the construction of row $\ell$ check the error in the numerical solution of component $i$; if it is smaller than the desired tolerance, then component $i$ is slow, was solved sufficiently accurately, and it is inactivated:

$$\frac{|T_{\ell,\ell}^{\{\text{S}\}} - T_{\ell,\ell-1}^{\{\text{S}\}}|}{\text{Atol} + \text{Rtol}\,|T_{\ell,\ell}^{\{\text{S}\}}|} \le \text{tol} \quad \Rightarrow \quad \text{inactivate } \mathbf{y}_n^{\{\text{S}\}} \text{ in rows } \ge \ell + 1.$$

   The user does not need to apriori split the system into slow and fast components. This is done dynamically as the integration progresses.
3. The components that were not yet solved sufficiently accurately are fast, and they require smaller step sizes. The fast subsystem continues to be integrated with smaller step sizes.
   The first column $T_{\ell+1,1}^{\{\text{F}\}}$ for the next row of the extrapolation table is obtained by applying the forward Euler method with step $h_{\ell+1}$ to the fast subsystem:

$$\mathbf{y}^{\{F\}}_{\lambda/(\ell+1)} = \mathbf{y}^{\{F\}}_{(\lambda-1)/(\ell+1)} + h_{\ell+1}\,\mathbf{f}^{\{F\}}(\mathbf{y}^{\{S\}}_{(\lambda-1)/(\ell+1)}, \mathbf{y}^{\{F\}}_{(\lambda-1)/(\ell+1)}), \tag{8.6}$$

where $\mathbf{y}^{\{F\}}_{\lambda/(\ell+1)} \approx \mathbf{y}^{\{F\}}(t_0 + \lambda\, h_{\ell+1})$.

4. Euler steps do not advance the slow components, which are frozen beyond extrapolation row $\ell$. Therefore the slow components are not available, and during the fast integration (8.6) one needs to use instead the approximations $\mathbf{y}^{\{S\}}_{\lambda/(\ell+1)} \approx \mathbf{y}^{\{S\}}(t_0 + \lambda h_{\ell+1})$.
   In order to allow for a correct continuation of the extrapolation tableau, the slow approximate solutions $\overline{\mathbf{y}}^{\{S\}}_{(\lambda-1)/(\ell+1)}$ in (8.6) are carefully constructed such as to have the same asymptotic expansions as the Euler solution they replace [18]; this special construction is expensive.
5. The construction of the extrapolation tableau for the fast components continues, as illustrated in Figure 13, until all components are solved with sufficient accuracy.

Remark 28 (Properties of the MURX algorithm). *By freezing extrapolation tableau of the slow components after $\ell$ rows, the MURX local truncation is formally of order $\mathcal{O}(H^{\ell+1})$ for both the fast and the slow variables. The MURX approach creates implicitly an arbitrary number* N *of partitions.*

**8.4. Extrapolation using first-order multirate base methods.** We next discuss the approach of Sandu and Constantinescu, which starts with multirate base methods, and then performs extrapolation on the results [12, 12, 13, 13, 49, 51]. A somewhat similar philosophy is employed by Zhang et al. [74] for coupled simulations. This approach is illustrated in Figure 14.

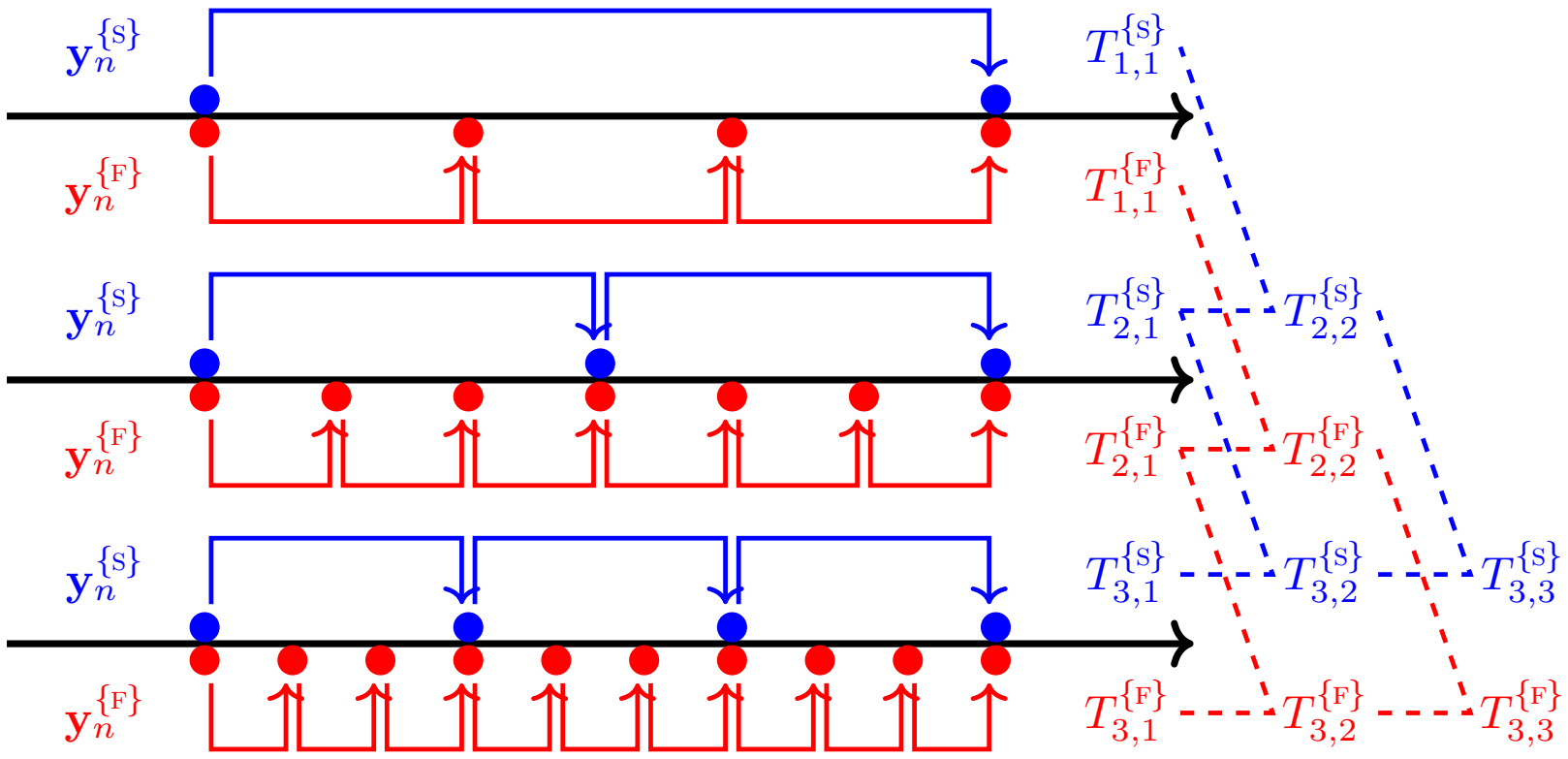

Fig. 14: Cartoon of Extrapolated MR-Euler method for $\mathtt{m} = 3$. One MR-Euler step with macro-step $H$ provides the solution $[T_{1,1}^{\{S\}}, T_{1,1}^{\{F\}}]$. Two successive MR-Euler steps with macro-steps $H/2$ yield the solution $[T_{2,1}^{\{S\}}, T_{2,1}^{\{F\}}]$, and three successive MR-Euler steps with macro-step $H/3$ yield the solution $[T_{3,1}^{\{S\}}, T_{3,1}^{\{F\}}]$. Dashed lines represent extrapolation operations. The final solution is $[\mathbf{y}_{n+1}^{\{S\}}, \mathbf{y}_{n+1}^{\{S\}}] = [T_{3,3}^{\{S\}}, T_{3,3}^{\{F\}}]$. In row $k$ of the extrapolation table, the slow component is solved with step size $H/k$, and the fast component is solved with a step size $h/k$ where $h = H/\mathtt{m}$.

The process is illustrated in Figure 14. While the discussion is based on the two-way partitioned system (2.1), the approach extends immediately to general N-way partitions.

1. One starts with a one-step, first order, base multirate method, and fixes the step size ratio $\mathtt{m}$. The method can be written in Henrici form (8.1) over a macro-step as follows:

$$\mathbf{y}_h(t+H) = \mathbf{y}_H(t) + h\,\Phi_{\mathtt{m}}(\mathbf{y}_H(t), H). \tag{8.7}$$

2. For a fixed step size ratio $\mathtt{m}$, the increment functions $\Phi_{\mathtt{m}}$ (8.7) are smooth, and satisfy the consistency condition (8.2). Application of Theorem 8.1 show that the multirate solutions have a global error expansion (8.3) in the macro-step $H$.
3. The construction of the extrapolation table now proceeds as in the traditional case. Solutions obtained by applying $n_j$ consecutive MR steps (with macro-step $H/n_j$ and microsteps $H/(\mathtt{m}n_j)$) build the first column entries $T_{j,1}$ (8.4). The other table columns are filled by applying the Aitken-Neville formula (8.5).

First order multirate schemes that form suitable base methods for MR extrapolation include:

1. Multirate explicit Euler method (4.2).
2. Multirate implicit Euler method (4.10).
3. Slowest-first multirate linearly-implicit Euler method, which starts with a regular linearly implicit Euler step applied to the entire system (2.1) with a large step size $H$:

(8.8a)

$$\begin{bmatrix} \mathbf{I} - H\,\mathbf{J}^{\{\mathrm{S},\mathrm{S}\}} & -H\,\mathbf{J}^{\{\mathrm{S},\mathrm{F}\}} \\ -H\,\mathbf{J}^{\{\mathrm{F},\mathrm{S}\}} & \mathbf{I} - H\,\mathbf{J}^{\{\mathrm{F},\mathrm{F}\}} \end{bmatrix} \cdot \begin{bmatrix} \mathbf{y}^{\{\mathrm{S}\}}_{n+1} - \mathbf{y}^{\{\mathrm{S}\}}_{n} \\ \bar{\mathbf{y}}^{\{\mathrm{F}\}}_{n+1} - \mathbf{y}^{\{\mathrm{F}\}}_{n} \end{bmatrix} = \begin{bmatrix} H\,\mathbf{f}^{\{\mathrm{S}\}}\left(\mathbf{y}^{\{\mathrm{S}\}}_{n}, \mathbf{y}^{\{\mathrm{F}\}}_{n}\right) \\ H\,\mathbf{f}^{\{\mathrm{F}\}}\left(\mathbf{y}^{\{\mathrm{S}\}}_{n}, \mathbf{y}^{\{\mathrm{F}\}}_{n}\right) \end{bmatrix},$$

where $\mathbf{J}^{\{\mathrm{S},\mathrm{S}\}} = \mathbf{f}^{\{\mathrm{S}\}}_{\mathbf{y}^{\{\mathrm{S}\}}|0}$, $\mathbf{J}^{\{\mathrm{S},\mathrm{F}\}} = \mathbf{f}^{\{\mathrm{S}\}}_{\mathbf{y}^{\{\mathrm{F}\}}|0}$, $\mathbf{J}^{\{\mathrm{F},\mathrm{F}\}} = \mathbf{f}^{\{\mathrm{F}\}}_{\mathbf{y}^{\{\mathrm{F}\}}|0}$, and $\mathbf{J}^{\{\mathrm{F},\mathrm{S}\}} = \mathbf{f}^{\{\mathrm{F}\}}_{\mathbf{y}^{\{\mathrm{S}\}}|0}$ denote the Jacobian matrices evaluated at $t_0$, i.e., at the beginning of the current extrapolation time step.

Next, the fast component $\bar{\mathbf{y}}^{\{\mathrm{F}\}}_{n+1}$ is discarded as it too inaccurate. To obtain an accurate fast solution the fast subsystem is solved by a sequence of small linearly implicit Euler steps:

$$\left(\mathbf{I} - \frac{H}{\mathtt{m}}\mathbf{J}^{\{\mathrm{F},\mathrm{F}\}}\right)\left(\mathbf{y}^{\{\mathrm{F}\}}_{n+\frac{i}{\mathtt{m}}} - \mathbf{y}^{\{\mathrm{F}\}}_{n+\frac{i-1}{\mathtt{m}}}\right) = \frac{H}{\mathtt{m}}\,\mathbf{f}^{\{\mathrm{F}\}}\left(\mathbf{y}^{\{\mathrm{S}\}}_{n+\frac{i-1}{\mathtt{m}}}, \mathbf{y}^{\{\mathrm{F}\}}_{n+\frac{i-1}{\mathtt{m}}}\right), \quad i = 1, \ldots, \mathtt{m}. \tag{8.8b}$$

Just as for the explicit case, $\mathbf{y}^{\{\mathrm{S}\}}_{n+\frac{i}{\mathtt{m}}}$ is approximated.

4. Compound multirate linearly-implicit Euler method. In order to avoid discarding the computed fast solution, one takes a linearly implicit compound step that computes the slow component with a step-size $H$, together with the first fast sub-step of size $H/\mathtt{m}$:

(8.9)

$$\begin{bmatrix} \mathbf{I} - H\,\mathbf{J}^{\{\mathrm{S,S}\}} & -H\,\mathbf{J}^{\{\mathrm{S,F}\}} \\ -\frac{H}{\mathtt{m}}\,\mathbf{J}^{\{\mathrm{F,S}\}} & \mathbf{I} - \frac{H}{\mathtt{m}}\,\mathbf{J}^{\{\mathrm{F,F}\}} \end{bmatrix} \cdot \begin{bmatrix} \mathbf{y}_{n+1}^{\{\mathrm{S}\}} - \mathbf{y}_n^{\{\mathrm{S}\}} \\ \mathbf{y}_{n+\frac{1}{\mathtt{m}}}^{\{\mathrm{F}\}} - \mathbf{y}_n^{\{\mathrm{F}\}} \end{bmatrix} = \begin{bmatrix} H\,\mathbf{f}^{\{\mathrm{S}\}}\left(\mathbf{y}_n^{\{\mathrm{S}\}}, \mathbf{y}_n^{\{\mathrm{F}\}}\right) \\ \frac{H}{\mathtt{m}}\,\mathbf{f}^{\{\mathrm{F}\}}\left(\mathbf{y}_n^{\{\mathrm{S}\}}, \mathbf{y}_n^{\{\mathrm{F}\}}\right) \end{bmatrix}.$$

The slow solution and the first fast substep solution are computed in a coupled manner. Next, the fast component solutions $\mathbf{y}_{n+2/\mathtt{m}}^{\{\mathrm{F}\}}, \ldots, \mathbf{y}_{n+1}^{\{\mathrm{F}\}}$, are obtained by advancing the fast system with small linearly implicit Euler steps.

REMARK 29. *The extrapolation of multirate base methods approach has (in principle) the full order associated with the last entry in the extrapolation Table. However, the user needs to provide a fast-slow partition, which remains fixed during an extrapolation step. The integration can be stopped adaptively, when the error estimate provided by the difference of two consecutive table entries is sufficiently small. Finally, the extrapolation of multirate base methods approach can be parallelized as easily as the standard extrapolation methods.*

**9. Multirate Infinitesimal Methods.** An interesting question is what happens in a multirate integration when the ratio of steps becomes very large, $\mathtt{m} \to \infty$? To study this question we revisit the MRFE scheme (4.2) with a slowest-first strategy and a constant interpolation of the slow variable. Using the notation $\tau := t_n + \frac{\ell}{\mathtt{m}} H$, $\mathbf{y}_\tau^{\{\mathrm{F}\}} := \mathbf{y}_{n+\ell/\mathtt{m}}^{\{\mathrm{F}\}}$, and $\mathbf{y}_{\tau+h}^{\{\mathrm{F}\}} := \mathbf{y}_{n+(\ell+1)/\mathtt{m}}^{\{\mathrm{F}\}}$ the MRFE method reads:

(9.1)
$$\begin{cases} \mathbf{y}_{n+1}^{\{\mathrm{S}\}} = \mathbf{y}_n^{\{\mathrm{S}\}} + H\,\mathbf{f}^{\{\mathrm{S}\}}(\mathbf{y}_n^{\{\mathrm{S}\}}, \mathbf{y}_n^{\{\mathrm{F}\}}), \\ \frac{\mathbf{y}_{\tau+h}^{\{\mathrm{F}\}} - \mathbf{y}_\tau^{\{\mathrm{F}\}}}{h} = \mathbf{f}^{\{\mathrm{F}\}}(\mathbf{y}_n^{\{\mathrm{S}\}}, \mathbf{y}_\tau^{\{\mathrm{F}\}}), \quad \tau = t_n + \frac{\ell}{\mathtt{m}} H, \quad 0 \le \ell \le \mathtt{m}. \end{cases}$$

We consider $H$ fixed; taking the limit $\mathtt{m} \to \infty$ is equivalent to taking the limit $h \to 0$ in (9.1), which leads to:

(9.2)
$$\begin{cases} \mathbf{y}_{n+1}^{\{\mathrm{S}\}} = \mathbf{y}_n^{\{\mathrm{S}\}} + H\,\mathbf{f}^{\{\mathrm{S}\}}(\mathbf{y}_n^{\{\mathrm{S}\}}, \mathbf{y}_n^{\{\mathrm{F}\}}), \\ \frac{d}{d\tau}\,\mathbf{y}_\tau^{\{\mathrm{F}\}} = \mathbf{f}^{\{\mathrm{F}\}}(\mathbf{y}_n^{\{\mathrm{S}\}}, \mathbf{y}_\tau^{\{\mathrm{F}\}}), \quad t_n \le \tau \le t_{n+1}. \end{cases}$$

The resulting scheme (9.2) is called the "multirate infinitesimal forward Euler" (MRI-FE) scheme. The "multirate infinitesimal" (MRI) approach:

1. solves the slow components discretely, and
2. advances the fast components continuously via a modified fast ODE (here the modification consists in keeping the argument $\mathbf{y}_n^{\{\mathrm{S}\}}$ constant throughout the integration).

The fast ODE is the limit of the fast solver when the substep $h \to 0$ becomes infinitesimally small, which justifies the name "multirate infinitesimal (fast step)" approach. The multirate infinitesimal step approach was first proposed in a seminal paper by Knoth and Wolke [34], and later extended in [14, 33, 44, 48, 59–61, 66, 73].

REMARK 30. *In practice, the modified fast ODE* (9.2) *is solved numerically; any discretization method and any appropriate sequence of small steps can be used as long as the fast ODE is solved sufficiently accurately. Consequently, the MRI philosophy provides extreme flexibility for treating the fast component and building multirate schemes. MRI provides a middle ground between splitting methods and traditional single rate methods.*

### 9.1. MRI Adams Methods.

**9.1.1. MRI Adams-Bashforth.** We start with Adams-Bashforth/Moulton schemes [14, 66] and make use of the interpolatory derivation of these methods.

Consider the degree $k-1$ Lagrange polynomial basis $\{\ell_j^{[\mathrm{AB}]}(\tau)\}_{j=0,\dots,k-1}$ for the scaled time grid $\{(t_n - t_{n-j})/H\}_{j=0,\dots,k-1}$ and define

$$\tilde{\ell}_{n-j}^{[\mathrm{AB}]}(\tau) := \int_0^\tau \ell_{n-j}^{[\mathrm{AB}]}(\tau')d\tau', \quad \beta_j^{[\mathrm{AB}]} := \tilde{\ell}_{n-j}^{[\mathrm{AB}]}(1). \tag{9.3}$$

Build the degree $k-1$ polynomial interpolant of the right hand side functions (2.1):

$$F^{\{\mathrm{S}\}}(t_n + \tau H) = \sum_{j=0}^{k-1} \mathbf{f}_{n-j}^{\{\mathrm{S}\}}\, \ell_{n-j}^{[\mathrm{AB}]}(\tau), \quad F^{\{\mathrm{F}\}}(t_n + \tau H) = \sum_{j=0}^{k-1} \mathbf{f}_{n-j}^{\{\mathrm{F}\}}\, \ell_{n-j}^{[\mathrm{AB}]}(\tau). \tag{9.4}$$

The standard Adams-Bashforth (AB) scheme for (2.1) is derived by replacing the both right hand side functions with their interpolants (9.4), and integrating the polynomials analytically. Using standard theory, the local truncation errors of the AB solution with step size $H$ are $\boldsymbol{\delta}_{n+1}^{\{\mathrm{S}\}} \propto \mathbf{y}^{\{\mathrm{S}\}(k+1)}\, H^{k+1}$ and $\boldsymbol{\delta}_{n+1}^{\{\mathrm{F}\}} \propto \mathbf{y}^{\{\mathrm{F}\}(k+1)}\, H^{k+1}$, for the slow and fast components, respectively. Here we denote $\mathbf{y}^{\{\mathrm{F}\}(\ell)} =: d^\ell \mathbf{y}^{\{\mathrm{F}\}}/dt^\ell$ and $\mathbf{y}^{\{\mathrm{S}\}(\ell)} =: d^\ell \mathbf{y}^{\{\mathrm{S}\}}/dt^\ell$. Since the time derivatives of the fast component can be quite large, $\|\mathbf{y}^{\{\mathrm{F}\}(k+1)}\| \gg \|\mathbf{y}^{\{\mathrm{S}\}(k+1)}\|$, the step size is restricted by the accuracy requirements for the fast component: $H \propto (\mathbf{y}^{\{\mathrm{F}\}(k+1)})^{-k-1}$.

To obtain a multirate infinitesimal version of the AB scheme, we modify this approach as follows. We replace the slow function by its interpolant, but use the exact fast function in the integration:

$$\begin{bmatrix} \mathbf{y}_{n+1}^{\{\mathrm{S}\}} \\ \mathbf{y}_{n+1}^{\{\mathrm{F}\}} \end{bmatrix} = \begin{bmatrix} \mathbf{y}_{n}^{\{\mathrm{S}\}} \\ \mathbf{y}_{n}^{\{\mathrm{F}\}} \end{bmatrix} + \int_0^H \begin{bmatrix} F^{\{\mathrm{S}\}}(t_n+\theta) \\ \mathbf{f}^{\{\mathrm{F}\}}(\mathbf{y}^{\{\mathrm{S}\}}(t_n+\theta), \mathbf{y}^{\{\mathrm{F}\}}(t_n+\theta)) \end{bmatrix} d\theta. \tag{9.5}$$

The numerical solutions (9.5) are obtained by solving the underlying ODEs. Working out the (exact) integration of the slow polynomial right hand side in (9.5), and using (9.3), we arrive at the multirate infinitesimal Adams-Bashforth (MRI-AB) scheme:

$$v^{\{\mathrm{F}\}}(0) = \mathbf{y}_n^{\{\mathrm{F}\}}; \tag{9.6a}$$

$$v^{\{\mathrm{S}\}}(\theta) = \mathbf{y}_n^{\{\mathrm{S}\}} + H \sum_{j=0}^{k-1} \tilde{\ell}_{n-j}^{[\mathrm{AB}]}(\theta/H)\, \mathbf{f}_{n-j}^{\{\mathrm{S}\}}, \tag{9.6b}$$

$$\dot{v}^{\{\mathrm{F}\}}(\theta) = \mathbf{f}^{\{\mathrm{F}\}}\left(v^{\{\mathrm{S}\}}(\theta), v^{\{\mathrm{F}\}}(\theta)\right), \quad \theta \in [0, H]; \tag{9.6c}$$

$$\mathbf{y}_{n+1}^{\{\mathrm{F}\}} = v^{\{\mathrm{F}\}}(H), \tag{9.6d}$$

$$\mathbf{y}_{n+1}^{\{\mathrm{S}\}} = \mathbf{y}_n^{\{\mathrm{S}\}} + H \sum_{j=0}^{k-1} \beta_j^{[\mathrm{AB}]}\, \mathbf{f}_{n-j}^{\{\mathrm{S}\}}. \tag{9.6e}$$

The MRI-AB method (9.6) has the following properties:

- The slow component is discretized with the standard Adams-Bashforth method with a macro-step size $H$ (9.6e). When no fast component is present, $\mathbf{f}^{\{\mathrm{F}\}} \equiv \mathbf{0}$, the method (9.6) reduces to the standard Adams-Bashforth scheme applied to the slow system.

- The fast component is advanced by integrating the modified fast ODE (9.6a)–(9.6d). This can be numerically solved with any scheme and any step size, as long as the solution is sufficiently accurate. When no slow component is present, $\mathbf{f}^{\{\mathrm{S}\}} \equiv \mathbf{0}$, and for $\mathbf{y}_n = \mathbf{y}(t_n)$, the method (9.6) performs an exact integration of the fast system.

We compute the local truncation error of (9.6) by carrying out the numerical solution starting from the exact past solution values. Direct calculations show that:

$$
\begin{aligned}
\boldsymbol{\delta}_{n+1}^{\{\mathrm{S}\}} &= C\, \frac{d^k \mathbf{f}^{\{\mathrm{S}\}}}{dt^k}(\mathbf{y}(t_n))\, \frac{H^{k+1}}{k!} + \mathcal{O}(H^{k+2}), \\
\boldsymbol{\delta}_{n+1}^{\{\mathrm{F}\}} &= C\, e^{H\, \mathbf{f}^{\{\mathrm{F}\}}_{\mathbf{y}^{\{\mathrm{F}\}}}(t_n)}\, \mathbf{f}^{\{\mathrm{F}\}}_{\mathbf{y}^{\{\mathrm{S}\}}}(t_n)\, \frac{d^k \mathbf{f}^{\{\mathrm{S}\}}}{dt^k}(\mathbf{y}(t_n))\, \frac{H^{k+2}}{k!} + \mathcal{O}(H^{k+3}).
\end{aligned}
\tag{9.7}
$$

We see that the dominant $\mathcal{O}(H^{k+1})$ term of the local truncation error only affects the slow component, and the error constant is proportional to the high derivative of slow function $\mathbf{y}^{\{\mathrm{S}\}(k+1)}(\mathbf{y}(t_n))$. The large term coming from the fast function derivative $\mathbf{y}^{\{\mathrm{F}\}(k+1)}$ is not present in the constant of the dominant term, thanks to the exact integration of the fast component (specifically, the slow error leads to an $\mathcal{O}(H^{k+1})$ perturbation of the right hand side of the modified ODE, and upon integration over a time interval of length $H$ this gives an $\mathcal{O}(H^{k+2})$ solution error).

The dominant $\mathcal{O}(H^{k+2})$ term of the fast component local truncation error has a constant that depends on the high derivative of slow function $\mathbf{y}^{\{\mathrm{S}\}(k+1)}$, but also on the Jacobians $\mathbf{f}^{\{\mathrm{F}\}}_{\mathbf{y}^{\{\mathrm{S}\}}}$ and $\mathbf{f}^{\{\mathrm{F}\}}_{\mathbf{y}^{\{\mathrm{F}\}}}$ of the fast function. These Jacobians represent the propagation of the slow interpolation error through the fast dynamics. Note that these Jacobians are associated with the first derivative of the fast function $\mathbf{f}^{\{\mathrm{F}\}\prime}$, but the term does not contain higher time derivatives.

*In summary, the MRI approach reduces the large constants (that depend on the fast dynamics) scaling the dominant error terms.*

Example 11. *Consider the decoupled system*

$$
\mathbf{y}^{\{\mathrm{S}\}\prime} = i\, \mathbf{y}^{\{\mathrm{S}\}}, \quad \mathbf{y}^{\{\mathrm{F}\}\prime} = i\, \mathbf{y}^{\{\mathrm{S}\}} + i\, \mathtt{m}\, \mathbf{y}^{\{\mathrm{F}\}}, \quad \mathbf{y}^{\{\mathrm{S}\}}(t_n) = \mathbf{y}^{\{\mathrm{F}\}}(t_n) = 1,
$$

*with* $\mathtt{m} > 1$, *where*

$$
\mathbf{y}^{\{\mathrm{S}\}(k+1)}(t_n) = i^k \quad \mathbf{y}^{\{\mathrm{F}\}(k+1)}(t_n) = i^k\, \frac{\mathtt{m}^k - 1}{\mathtt{m} - 1} + (i\, \mathtt{m})^k \;= i^k\, \frac{\mathtt{m}^{k+1} - 1}{\mathtt{m} - 1}.
$$

*For the single-rate AB scheme the local truncation errors are:*

$$
|\boldsymbol{\delta}^{\{\mathrm{S}\}}| \approx C\, H^{k+1}/k!\,, \quad |\boldsymbol{\delta}^{\{\mathrm{F}\}}| \approx C\, \mathtt{m}^k\, H^{k+1}/k!\,.
$$

*Keeping the slow local truncation error below a fixed threshold* $\mathtt{tol}$ *requires step sizes* $H^{\{\mathrm{S}\}} \propto \mathtt{tol}^{1/(k+1)}$, *and keeping the fast truncation error below threshold requires step sizes* $H^{\{\mathrm{F}\}} \propto (\mathtt{tol}/\mathtt{m}^k)^{1/(k+1)} = H^{\{\mathrm{S}\}}\, \mathtt{m}^{-k/(k+1)}$, *i.e., the fast dynamics causes step sizes to decrease (almost) linearly with increasing* $\mathtt{m}$.

*For the MRI-AB scheme the local truncation errors* (9.7) *are:*

$$
|\boldsymbol{\delta}^{\{\mathrm{S}\}}| \approx C\, H^{k+1}/k!\,, \quad |\boldsymbol{\delta}^{\{\mathrm{F}\}}| \approx C\, \mathtt{m}\, H^{k+2}/k!\,.
$$

*Keeping the fast truncation error below threshold requires step sizes* $H^{\{\mathrm{F}\}} \propto H^{\{\mathrm{S}\}}\, \mathtt{tol}^{(k+1)/(k+2)}\, \mathtt{m}^{-1/(k+2)}$, *which decrease very slowly with increasing* $\mathtt{m}$.

Remark 31 (MRI Adams-Bashforth with direct solution interpolation). *One can replace (9.6b) with the direct interpolant of the slow solution, and use it in the formulation of the modified fast ODE:*

$$v^{\{\mathrm{s}\}}(\theta) = \sum_{j=0}^{k-1} \mathbf{y}^{\{\mathrm{s}\}}_{n-j}\, \ell_{n-j}((\theta - t_n)/H). \tag{9.8}$$

*In this case, however,* $v^{\{\mathrm{s}\}}(H) \neq \mathbf{y}^{\{\mathrm{s}\}}_{n+1}$*, the solution computed via de AB scheme. One still needs to compute the slow solution* $\mathbf{y}^{\{\mathrm{s}\}}_{n+1}$ *using (9.6e)*

**9.1.2. MRI Adams-Moulton.** We start with an Adams-Bashforth predictor step with step size $H$ to obtain the (inaccurate) predicted solution $[\mathbf{y}^{\{\mathrm{s}*\}}_{n+1}, \mathbf{y}^{\{\mathrm{F}*\}}_{n+1}]$, and compute the corresponding function value $\mathbf{f}^{\{\mathrm{s}*\}}_{n+1}$. Choosing the degree $k$ Lagrange polynomial basis $\{\ell^{[\mathrm{AM}]}_j(t)\}_{j=-1,\dots,k-1}$ for the scaled time grid $\{(t_n - t_{n-j})/H\}_{j=-1,\dots,k-1}$, we build the following degree $k$ interpolant for the slow function values:

$$F^{\{\mathrm{s}\}}(t_n + \tau H) = \sum_{j=0}^{k-1} \mathbf{f}^{\{\mathrm{s}\}}_{n-j}\, \ell^{[\mathrm{AM}]}_{n-j}(\tau) + \mathbf{f}^{\{\mathrm{s}*\}}_{n+1}\, \ell^{[\mathrm{AM}]}_{n+1}(\tau). \tag{9.9}$$

The component partitioned MRI-AM scheme follows (9.6), except that we use the slow function interpolant (9.9) constructed using the predicted values in place of (9.6b). As in Remark 31, the slow component can also be obtained from a direct interpolation of the slow solution. We note that, if $\mathbf{f}^{\{\mathrm{F}\}} \equiv 0$, the slow component is solved using an AB prediction step followed by an AM corrector step.

**9.2. MRI GARK Methods.** We revisit the partitioned system (2.1). Consider the MR-GARK scheme (5.5) which solves the slow component with a Runge-Kutta method $(A^{\{\mathrm{s,s}\}}, b^{\{\mathrm{s}\}}, c^{\{\mathrm{s}\}})$ and step size $H$. The slow method has non-decreasing abscissae $\Delta c^{\{\mathrm{s}\}}_i = c^{\{\mathrm{s}\}}_{i+1} - c^{\{\mathrm{s}\}}_i \geq 0$ for $i = 1, \dots, s^{\{\mathrm{s}\}}$, where $c^{\{\mathrm{s}\}}_{s^{\{\mathrm{s}\}}+1} = 1$.

We consider the infinitesimal fast step limit where $H$ is fixed and $h \to 0$. The fast stages (5.5b) approach the solution of an ordinary differential equation, that advances the fast component in between consecutive slow abscissae. The fast numerical scheme $(A^{\{\mathrm{F,F}\}}, b^{\{\mathrm{F}\}}, c^{\{\mathrm{F}\}})$ is replaced by an exact integration of a modified fast ODE. The discrete fast-to-slow coupling coefficients $a^{\{\mathrm{F,s},\ell\}}_{i,j}$, that bring slow dynamic information into the fast integration at time $\tau = t_n + \ell\, h$ become, in the limit $h \to 0$, continuous functions of time $\gamma_{i,j}(\tau)$. This limit approach leads to multirate infinitesimal GARK (MRI-GARK) schemes [44, 48].

**9.3. Decoupled MRI-GARK.** Decoupled schemes compute the slow and fast components separately; even of each component is computed implicitly, there is no nonlinear system that requires solving for both fast and slow components in a coupled manner. A decoupled MRI-GARK scheme advances the solution as follows [44, 48]:

$$\mathbf{Y}^{\{\mathrm{s}\}}_1 = \mathbf{y}^{\{\mathrm{s}\}}_n, \quad \mathbf{Y}^{\{\mathrm{F}\}}_1 = \mathbf{y}^{\{\mathrm{F}\}}_n; \tag{9.10a}$$

(9.10b) For each stage $i = 1, \dots, s^{\{\mathrm{s}\}}$ :

$$\begin{cases} v_i(0) = \mathbf{Y}^{\{\mathrm{F}\}}_i, \\ v_i' = \Delta c^{\{\mathrm{s}\}}_i\, \mathbf{f}^{\{\mathrm{F}\}}\left( v_i, \mathbf{Y}_i + H \sum_{j=1}^{i+1} \tilde{\gamma}_{i,j}\left(\frac{\theta}{H}\right) \mathbf{f}^{\{\mathrm{s}\}}_j \right), \quad \theta \in [0, H]; \\ \mathbf{Y}^{\{\mathrm{F}\}}_{i+1} = v_i(H); \end{cases} \tag{9.10c}$$

$$
\text{(9.10d)} \qquad \begin{cases} \mathbf{Y}_{i+1}^{\{\mathrm{S}\}} = \mathbf{Y}_{i}^{\{\mathrm{S}\}} + H \displaystyle\sum_{j=1}^{i+1} a_{i+1,j}^{\{\mathrm{S}\}} \, \mathbf{f}_j^{\{\mathrm{S}\}}, \quad a_{i+1,j}^{\{\mathrm{S}\}} \equiv \widetilde{\gamma}_{i,j}\,(1)\,; \\ \mathbf{f}_{i+1}^{\{\mathrm{S}\}} = \mathbf{f}^{\{\mathrm{S}\}}\big(\,\mathbf{Y}_{i+1}^{\{\mathrm{F}\}},\, \mathbf{Y}_{i+1}^{\{\mathrm{S}\}}\,\big); \end{cases}
$$

$$
\text{(9.10e)} \qquad \mathbf{y}_{n+1}^{\{\mathrm{S}\}} = \mathbf{Y}_{s^{\{\mathrm{S}\}}+1}, \quad \mathbf{y}_{n+1}^{\{\mathrm{F}\}} = \mathbf{Y}^{\{\mathrm{F}\}}_{s^{\{\mathrm{S}\}}+1}.
$$

The workings of the MRI-GARK scheme (9.10) are summarized in Figure 15. The method proceeds as follows:

- The fast component is advanced at each stage by the exact integration of a modified fast ODE (9.10c). The slow components in the modified fast ODE are reconstructed from slow stage information. The slow-to-fast coupling coefficients $\widetilde{\gamma}_{i,j}\,(\tau)$, $\tau \in [0,1]$, are functions of the fast ODE integration time.
- The slow component is advanced at each stage by a Runge-Kutta stage computation (9.10d). The slow method coefficients are $a_{i+1,j}^{\{\mathrm{S}\}} \equiv \widetilde{\gamma}_{i,j}\,(1)$.
- The slow Runge-Kutta method is explicit if $a_{i+1,i+1}^{\{\mathrm{S}\}} \equiv \widetilde{\gamma}_{i,i+1}\,(1) = 0$. In this case the fast and slow stages $\mathbf{Y}_{i+1}^{\{\mathrm{F}\}}$, $\mathbf{Y}_{i+1}^{\{\mathrm{S}\}}$ are computed independently, in succession, in a decoupled fashion. This justifies the name "decoupled MRI-GARK" for scheme (9.10).
- The last stage values provide the next step solutions (9.10e).

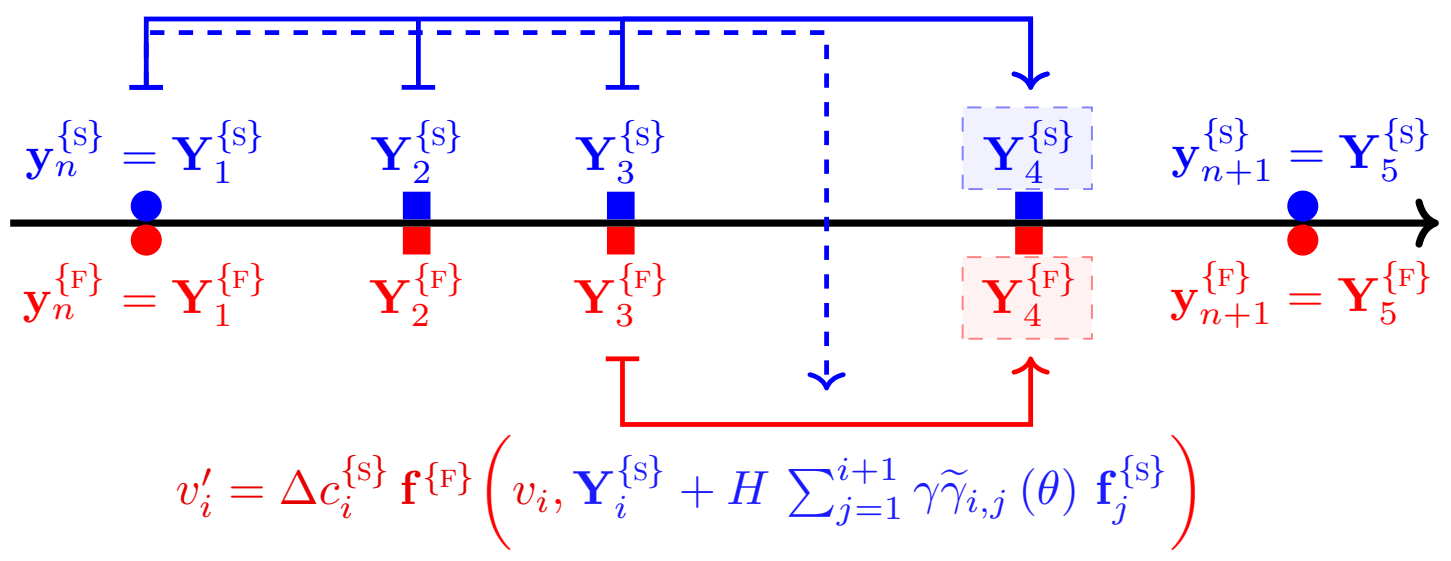

Fig. 15: Cartoon of a decoupled MRI-GARK scheme (9.10) with $s = 5$ for component partitioning. In this example, fast stage $\mathbf{Y}_4^{\{\mathrm{F}\}}$ is obtained by solving a modified fast equation over the interval $[c_3, c_4]$, started from $\mathbf{Y}_3^{\{\mathrm{F}\}}$. The fast equation is added slow forcing terms computed from $\mathbf{f}_{1,2,3}^{\{\mathrm{S}\}}$. Slow stage $\mathbf{Y}_4^{\{\mathrm{S}\}}$ is obtained through a regular Runge-Kutta step started with $\mathbf{Y}_3^{\{\mathrm{S}\}}$ plus a linear combination of $\mathbf{f}_{1,2,3}^{\{\mathrm{S}\}}$. Full step values are marked by circles, and stage values by squares.

Remark 32 (Decoupled computations for diagonally implicit slow method). *The slow Runge-Kutta method is diagonally implicit if* $a_{i+1,i+1}^{\{\mathrm{S}\}} \equiv \widetilde{\gamma}_{i,i+1}\,(1) > 0$. *In this case the slow stage* (9.10d) *requires solving a nonlinear system of equations:*

$$
\text{(9.11)} \qquad \mathbf{Y}_{i+1}^{\{\mathrm{S}\}} = \mathbf{Y}_{i}^{\{\mathrm{S}\}} + H \sum_{j=1}^{i} a_{i+1,j}^{\{\mathrm{S}\}} \, \mathbf{f}_j^{\{\mathrm{S}\}} + H\, a_{i+1,i+1}^{\{\mathrm{S}\}}\, \mathbf{f}^{\{\mathrm{S}\}}(\mathbf{Y}_{i+1}^{\{\mathrm{S}\}}, \mathbf{Y}_{i+1}^{\{\mathrm{F}\}}).
$$

*In order to keep computations decoupled, one chooses a slow method with abscissae* $c_{i+1}^{\{\mathrm{S}\}} = c_i^{\{\mathrm{S}\}}$ *whenever* $a_{i+1,i+1}^{\{\mathrm{S}\}} > 0$. *This makes* $\Delta c_i^{\{\mathrm{S}\}} = 0$, *and the resulting trivial fast ODE* (9.10c) *gives the fast stage solution* $\mathbf{Y}_{i+1}^{\{\mathrm{F}\}} = \mathbf{Y}_i^{\{\mathrm{F}\}}$. *The nonlinear system* (9.11) *is solved independently for the slow component $\mathbf{Y}^{\{\mathrm{s}\}}_{i+1}$. This makes the MRI-GARK method decoupled even for a diagonally implicit slow solution.*

**9.3.1. MRI-GARK order conditions.** We start with a slow Runge-Kutta method $(A^{\{\mathrm{s},\mathrm{s}\}}, b^{\{\mathrm{s}\}}, c^{\{\mathrm{s}\}})$ of order $p$. To obtain an order $p$ MRI-GARK method (9.10) the coupling coefficients $\tilde{\gamma}(\tau)$ need to satisfy certain order conditions. It is convenient to define the time-dependent combination coefficients as polynomials in time:

$$\tilde{\gamma}_{i,j}(\tau) := \sum_{k\geq 0} \gamma^k_{i,j}\, \frac{\tau^{k+1}}{k+1}\,, \quad a^{\{\mathrm{s}\}}_{i+1,k} = \tilde{\gamma}_{i,j}(1)\,, \tag{9.12}$$

and to gather the gamma coefficients in the matrices $\mathbf{\Gamma}^k := [\gamma^k_{i,j}]_{i,j} \in \mathbb{R}^{s^{\{\mathrm{s}\}} \times s^{\{\mathrm{s}\}}}$.

The MRI-GARK method is internally consistent iff the following conditions hold:

$$\mathbf{\Gamma}^0 \cdot \mathbf{1}^{\{\mathrm{s}\}} = \Delta c^{\{\mathrm{s}\}} \quad \text{and} \quad \mathbf{\Gamma}^k \cdot \mathbf{1}^{\{\mathrm{s}\}} = 0 \quad \forall k \geq 1. \tag{9.13}$$

Like with any GARK scheme, an internally consistent MRI-GARK scheme has order at least two. For higher order conditions, it is convenient to define:

$$\zeta_k := \frac{1}{(k+1)(k+2)}; \quad \omega_k := \frac{1}{(k+1)(k+3)}; \quad \xi_k := \frac{1}{(k+1)(k+2)(k+3)}, \tag{9.14}$$

and

$$\mathfrak{A}^{\{\mathrm{s}\}} := A^{\{\mathrm{s}\}} + \sum_{k\geq 0} \zeta_k\, \mathbf{\Gamma}^k. \tag{9.15}$$

Theorem 9.1 (Third order coupling condition [48]). *An internally consistent MRI-GARK method* (9.10) *has order three iff the slow base scheme has order at least three, and the following coupling condition holds:*

$$\Delta c^{\{\mathrm{s}\}\mathrm{T}}\, \mathfrak{A}^{\{\mathrm{s}\}}\, c^{\{\mathrm{s}\}} = \frac{1}{6}. \tag{9.16}$$

Theorem 9.2 (Fourth order coupling conditions [48]). *An internally consistent MRI-GARK method* (9.10) *has order four iff the slow base scheme has order at least four, and the following coupling conditions hold:*

$$\frac{1}{2}\,\zeta^{\{\mathrm{s}\}\mathrm{T}}\, A^{\{\mathrm{s}\}}\, c^{\{\mathrm{s}\}} + \sum_{k\geq 0} \Big(\Delta c^{\{\mathrm{s}\}} \times (\zeta_k\, c^{\{\mathrm{s}\}} + \omega_k\, \Delta c^{\{\mathrm{s}\}})\Big)^{\mathrm{T}} \mathbf{\Gamma}^k\, c^{\{\mathrm{s}\}} = \frac{1}{8}, \tag{9.17a}$$

$$\Delta c^{\{\mathrm{s}\}\mathrm{T}}\, \mathfrak{A}^{\{\mathrm{s}\}}\, c^{\{\mathrm{s}\}\times 2} = \frac{1}{12}, \tag{9.17b}$$

$$d^{\{\mathrm{s}\}\mathrm{T}}\, \mathfrak{A}^{\{\mathrm{s}\}}\, c^{\{\mathrm{s}\}} = \frac{1}{24}, \tag{9.17c}$$

$$(\Delta c^{\{\mathrm{s}\}\times 2})^{\mathrm{T}} \Big(\frac{1}{2} A^{\{\mathrm{s}\}} + \sum_{k\geq 0} \xi_k\, \mathbf{\Gamma}^k\Big)\, c^{\{\mathrm{s}\}} + t^{\{\mathrm{s}\}\mathrm{T}}\, \mathfrak{A}^{\{\mathrm{s}\}}\, c^{\{\mathrm{s}\}} = \frac{1}{24}, \tag{9.17d}$$

$$\Delta c^{\{\mathrm{s}\}\mathrm{T}}\, \mathfrak{A}^{\{\mathrm{s}\}}\, A^{\{\mathrm{s}\}}\, c^{\{\mathrm{s}\}} = \frac{1}{24}, \tag{9.17e}$$

*where $\times$ denotes componentwise operatons (multiplication, power), the coefficients $\zeta_k, \omega_k, \xi_k$ defined in* (9.14), *$\mathfrak{A}^{\{\mathrm{s}\}}$ defined in* (9.15), *and:*

$$\zeta^{\{\mathrm{s}\}} \quad := \quad \Big[c^{\{\mathrm{s}\}2}_{i+1} - c^{\{\mathrm{s}\}2}_{i}\Big]_{1\leq i \leq s^{\{\mathrm{s}\}}}, \tag{9.17f}$$

$$d^{\{S\}} := \left[ \Delta c_i^{\{S\}} \left( 1 - \sum_{\ell=1}^{i} b_\ell^{\{S\}} \right) \right]_{1 \le i \le s^{\{S\}}}, \tag{9.17g}$$

$$t^{\{S\}} := \left[ \sum_{j=i+1}^{s^{\{S\}}} (\Delta c_j^{\{S\}})^2 \right]_{1 \le i \le s^{\{S\}}}. \tag{9.17h}$$

### 9.3.2. Examples of explicit and implicit decoupled MRI-GARK methods.

*Explicit midpoint method.* Consider the Butcher tableau of the explicit midpoint method:

$$A_{\text{EMIDP}}^{\{S\}} = \begin{array}{c|cc} 0 & 0 & 0 \\ \frac{1}{2} & \frac{1}{2} & 0 \\ \hline 1 & 0 & 1 \end{array}, \qquad \Delta c_{\text{EMIDP}}^{\{S\}} = \begin{bmatrix} \frac{1}{2} \\ \frac{1}{2} \end{bmatrix}.$$

The MRI explicit midpoint method in component-wise splitting:

$$\begin{aligned}
\mathbf{k}_1^{\{S\}} &:= \mathbf{f}^{\{S\}}(\mathbf{y}_n^{\{F\}}, \mathbf{y}_n^{\{S\}}); \\
v_1(0) &= \mathbf{y}_n^{\{F\}}; \quad v_1' = \tfrac{1}{2} \mathbf{f}^{\{F\}} \left( v_1, \mathbf{y}_n^{\{S\}} + \tfrac{\theta}{2} \mathbf{k}_1^{\{S\}} \right), \quad \theta \in [0, H]; \\
\mathbf{Y}_2^{\{F\}} &= v_1(H); \\
\mathbf{Y}_2^{\{S\}} &= \mathbf{y}_n^{\{S\}} + \tfrac{H}{2} \mathbf{f}_n^{\{S\}}; \\
\mathbf{k}_2^{\{S\}} &:= \mathbf{f}^{\{S\}}(\mathbf{Y}_2^{\{F\}}, \mathbf{Y}_2^{\{S\}}); \\
v_2(0) &= \mathbf{Y}_2^{\{F\}}; \quad v_2' = \tfrac{1}{2} \mathbf{f}^{\{F\}} \left( v_2, \tfrac{H-\theta}{2} \mathbf{k}_1^{\{S\}} + \theta \, \mathbf{k}_2^{\{S\}} \right), \quad \theta \in [0, H]; \\
\mathbf{y}_{n+1}^{\{F\}} &= v_2(H); \\
\mathbf{y}_{n+1}^{\{S\}} &= \mathbf{y}_n^{\{S\}} + H \, \mathbf{k}_2^{\{S\}}.
\end{aligned}$$

*Implicit trapezoidal method.* Consider the Butcher tableau of the implicit trapezoidal Runge-Kutta method:

$$A_{\text{ITRAP}}^{\{S\}} = \begin{array}{c|cc} 0 & 0 & 0 \\ 1 & \frac{1}{2} & \frac{1}{2} \\ \hline 1 & \frac{1}{2} & \frac{1}{2} \end{array}, \qquad \Delta c_{\text{ITRAP}}^{\{S\}} = \begin{bmatrix} 1 \\ 0 \end{bmatrix}.$$

The MRI implicit trapezoidal method in component-wise splitting reads:

$$\begin{aligned}
\mathbf{k}_1^{\{S\}} &:= \mathbf{f}^{\{S\}}(\mathbf{y}_n^{\{F\}}, \mathbf{y}_n^{\{S\}}); \\
v^{\{F\}}(0) &= \mathbf{y}_n^{\{F\}}; \quad v^{\{F\}\prime} = \mathbf{f}^{\{F\}} \left( v^{\{F\}}, \mathbf{y}_n^{\{S\}} + \theta \, \mathbf{k}_1^{\{S\}} \right), \quad \theta \in [0, H]; \\
\mathbf{y}_{n+1}^{\{F\}} &= v^{\{F\}}(H); \\
\mathbf{y}_{n+1}^{\{S\}} &= \mathbf{y}_n^{\{S\}} + \tfrac{H}{2} \mathbf{k}_1^{\{S\}} + \tfrac{H}{2} \mathbf{f}^{\{S\}}(\mathbf{y}_{n+1}^{\{F\}}, \mathbf{y}_{n+1}^{\{S\}}).
\end{aligned}$$

*Explicit third order methods.* The MRI-GARK-ERK33 family of methods has the following base slow scheme and coupling coefficients:

$$A^{\{\mathrm{S}\}} = \left[\begin{array}{c|ccc} 0 & 0 & 0 & 0 \\ \frac{1}{3} & \frac{1}{3} & 0 & 0 \\ \frac{2}{3} & 0 & \frac{2}{3} & 0 \\ \hline 1 & \frac{1}{4} & 0 & \frac{3}{4} \\ \hline 1 & \frac{1}{12} & \frac{1}{3} & \frac{7}{12} \end{array}\right], \quad \mathbf{\Gamma}^{0} = \left[\begin{array}{ccc} \frac{1}{3} & 0 & 0 \\ -\frac{1}{3} & \frac{2}{3} & 0 \\ 0 & -\frac{2}{3} & 1 \\ \hline \frac{1}{12} & -\frac{1}{3} & \frac{7}{12} \end{array}\right], \quad \mathbf{\Gamma}^{1} = \left[\begin{array}{ccc} 0 & 0 & 0 \\ 0 & 0 & 0 \\ \frac{1}{2} & 0 & -\frac{1}{2} \\ \hline 0 & 0 & 0 \end{array}\right].$$

**9.4. Stability of MRI GARK.**

*Additively partitioned systems.* We consider the scalar model problem for additively partitioned systems:

$$y' = \lambda^{\{\mathrm{F,F}\}}\, y + \lambda^{\{\mathrm{S,S}\}}\, y, \qquad \lambda^{\{\mathrm{F,F}\}}, \lambda^{\{\mathrm{S,S}\}} \in \mathbb{C}^{-}. \tag{9.18}$$

Application of the MRI-GARK scheme (9.10) to (9.18) leads to stability function in the variables $z^{\{\mathrm{F,F}\}} = H\,\lambda^{\{\mathrm{F,F}\}}$ and $z^{\{\mathrm{S,S}\}} = H\,\lambda^{\{\mathrm{S,S}\}}$. The stability functions are rational in $z^{\{\mathrm{F,F}\}}$, a consequence of the discrete computation of slow variables. The fast variables, however, enter the stability as arguments of the special functions:

$$\varphi_0(z) \coloneqq e^{z}, \quad \varphi_k(z) \coloneqq \int_0^1 e^{z\,(1-t)}\, t^{k-1}\, dt,$$

which is a consequence of the exact integration of the fast components. For example, for the second order methods the stability functions read:

$$\mathrm{R}_{\mathrm{EMIDP}}(z^{\{\mathrm{F,F}\}}, z^{\{\mathrm{S,S}\}}) = \varphi_0(z^{\{\mathrm{F,F}\}}) + \left(\tfrac{3}{2}\,\varphi_0(\tfrac{1}{2} z^{\{\mathrm{F,F}\}}) - \tfrac{1}{2}\right) \varphi_1(\tfrac{1}{2} z^{\{\mathrm{F,F}\}})\, z^{\{\mathrm{S,S}\}} + \tfrac{1}{2}\,\varphi_1^2(\tfrac{1}{2} z^{\{\mathrm{F,F}\}})\, z^{\{\mathrm{S,S}\}\,2},$$

$$\mathrm{R}_{\mathrm{ITRAP}}(z^{\{\mathrm{F,F}\}}, z^{\{\mathrm{S,S}\}}) = \frac{\varphi_0(z^{\{\mathrm{F,F}\}}) + \left(\varphi_1(z^{\{\mathrm{F,F}\}}) - \tfrac{1}{2}\right)\, z^{\{\mathrm{S,S}\}}}{1 - \tfrac{1}{2} z^{\{\mathrm{S,S}\}}}.$$

We note that there is no numerical damping of the fast components, which are integrated exactly. The damping of the fast error components is done solely by the fast dynamics, and when this is not sufficiently dissipative the stability of the entire MRI scheme suffers. Consider for example the case with an imaginary fast eigenvalue. In the stiff limit of the slow variable, the MRI trapezoidal stability function

$$\lim_{z^{\{\mathrm{S,S}\}} \to -\infty} |\mathrm{R}_{\mathrm{ITRAP}}(z^{\{\mathrm{F,F}\}} = i\,\omega, z^{\{\mathrm{S,S}\}})| = 1 - 2\,\varphi_1(z^{\{\mathrm{F,F}\}} = i\,\omega)$$

can take values that smaller or larger than one, depending on $\omega$, as seen in Figure 16.

*Component partitioned systems.* For component partitioned systems we consider the matrix model problem (5.11). Application of the MRI-GARK scheme (9.10) to (5.11) leads to the stability matrix in the variables $z^{\{\mathrm{F,F}\}} = H\,\lambda^{\{\mathrm{F,F}\}}$ and $z^{\{\mathrm{S,S}\}} = H\,\lambda^{\{\mathrm{S,S}\}}$. For example, the stability matrix of the implicit trapezoidal method is

$$\mathbf{R}^{\mathrm{ITRAP}} = \begin{bmatrix} \varphi_0(z^{\{\mathrm{F,F}\}}) + w\,\varphi_2(z^{\{\mathrm{F,F}\}}) & -\frac{1-\xi}{\alpha}\left(z^{\{\mathrm{F,F}\}} - z^{\{\mathrm{S,S}\}}\right)\left(\varphi_1(z^{\{\mathrm{F,F}\}}) + z^{\{\mathrm{F,F}\}}\,\varphi_2(z^{\{\mathrm{F,F}\}})\right) \\ \frac{\alpha\,\xi\left(z^{\{\mathrm{F,F}\}} - z^{\{\mathrm{S,S}\}}\right)}{2 - z^{\{\mathrm{S,S}\}}} & \frac{2 + z^{\{\mathrm{S,S}\}}}{2 - z^{\{\mathrm{S,S}\}}} \end{bmatrix}.$$

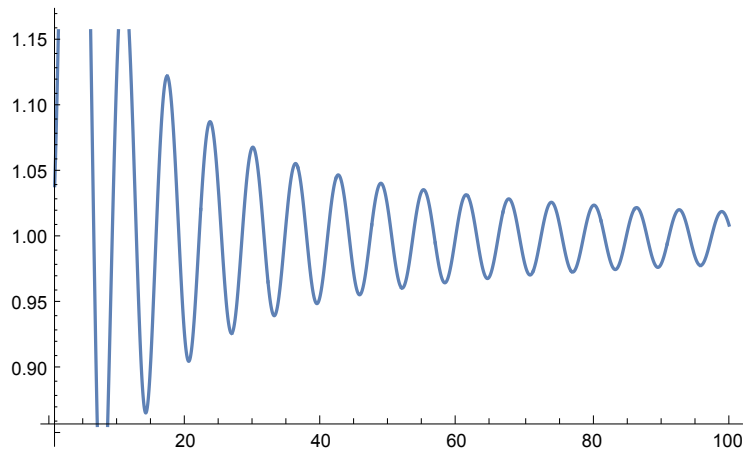

Fig. 16: $|\mathrm{R}_{\mathrm{ITRAP}}(z^{\{\mathrm{F,F}\}} = i\,\omega, z^{\{\mathrm{S,S}\}} \to -\infty)|$

The entries are rational in $z^{\{\mathrm{S,S}\}}$. The fast eigenvalues $z^{\{\mathrm{F,F}\}}$ appear both as arguments of the bounded $\varphi$ functions, and also as rational variables, due to the presence of $\lambda^{\{\mathrm{F,F}\}}$ in the coupling terms.

**9.5. Coupled implicit methods: Step predictor-corrector MRI GARK.** The discussion of coupled methods is based on [44].

We start with a "slow" Runge–Kutta base method $(A^{\{\mathrm{S}\}}, b^{\{\mathrm{S}\}}, c^{\{\mathrm{S}\}})$ with $s^{\{\mathrm{S}\}}$ stages. Unlike other multirate infinitesimal strategies, the base method is not restricted to be explicit or diagonally implicit.

One coupling strategy commonly used in discrete multirate methods is a predictor-corrector approach, where the predictor evolves the entire system, while the corrector is only applied to the fast partition whose solution was "predicted" inaccurately (see [41, 54, 57]). First, a combined Runge–Kutta macro-step is taken which serves as the predictor. The fast parts of the predicted stages are inaccurate and are refined by sub-stepping the fast component only. Approximations of the slow values needed during the micro-steps are obtained from interpolating the slow predicted values. In contrast, the decoupled schemes discussed in Subsection 9.3 advance the slow and fast components successively, and never compute a solution of the full system. The SPC-MRI-GARK methods, as depicted in Figure 17, can be viewed as an extreme case of this coupling strategy where the multirate ratio is infinite, i.e., the corrector takes infinitely many steps to refine the fast solution.

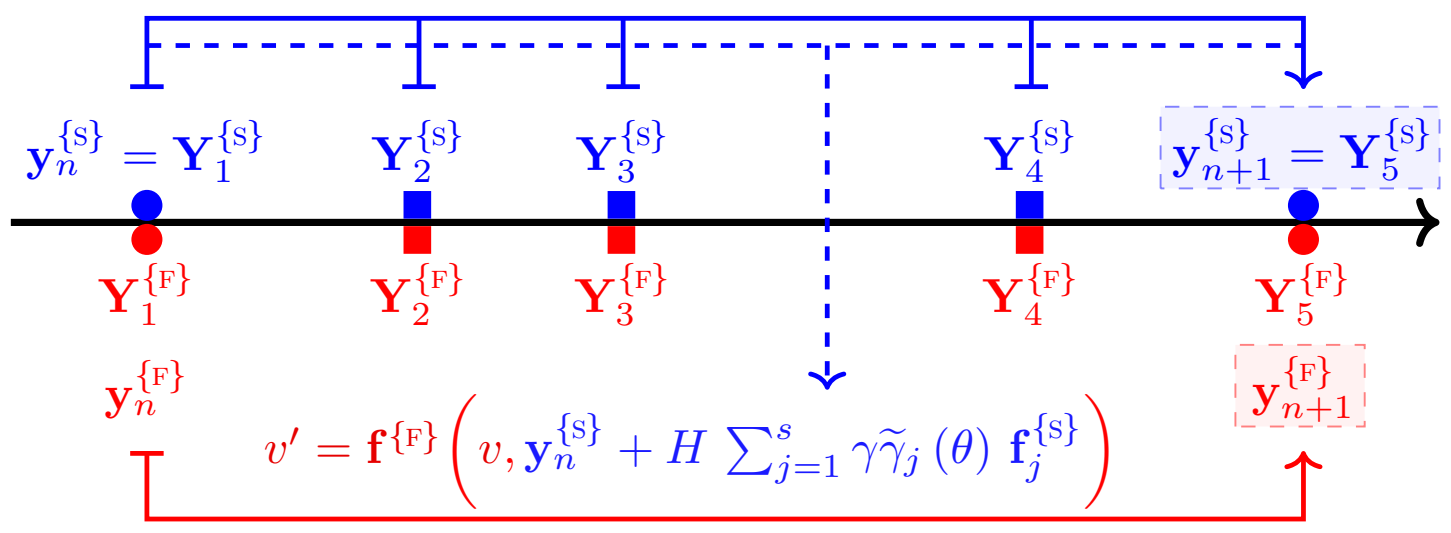

Fig. 17: Cartoon of the workings of a step predictor corrector MRI GARK method.

One step of an SPC-MRI-GARK method applied to (2.1) reads:

$$\begin{bmatrix} \mathbf{Y}_i^{\{F\}} \\ \mathbf{Y}_i^{\{S\}} \end{bmatrix} = \begin{bmatrix} \mathbf{y}_n^{\{F\}} + H \sum_{j=1}^{s^{\{S\}}} a_{i,j}^{\{S\}} \, \mathbf{f}_j^{\{F\}} \\ \mathbf{y}_n^{\{S\}} + H \sum_{j=1}^{s^{\{S\}}} a_{i,j}^{\{S\}} \, \mathbf{f}_j^{\{S\}} \end{bmatrix}, \tag{9.19a}$$

$$\mathbf{y}_{n+1}^{\{S\}} = \mathbf{y}_n^{\{S\}} + H \sum_{j=1}^{s^{\{S\}}} b_j^{\{S\}} \, \mathbf{f}_j^{\{S\}}, \tag{9.19b}$$

$$\begin{cases} v^{\{F\}}(0) = \mathbf{y}_n^{\{F\}}, \\ v^{\{F\}\prime} = \mathbf{f}^{\{F\}} \left( v^{\{F\}}, \mathbf{y}_n^{\{S\}} + H \displaystyle\sum_{j=1}^{s^{\{S\}}} \tilde{\gamma}_j \left( \frac{\theta}{H} \right) \mathbf{f}_j^{\{S\}} \right), \quad \theta \in [0, H], \\ \mathbf{y}_{n+1}^{\{F\}} = v^{\{F\}}(H), \end{cases} \tag{9.19c}$$

where $\mathbf{f}_j^{\{F\}} := \mathbf{f}^{\{F\}}(\mathbf{Y}_j^{\{F\}}, \mathbf{Y}_j^{\{S\}})$ and $\mathbf{f}_j^{\{S\}} := \mathbf{f}^{\{S\}}(\mathbf{Y}_j^{\{F\}}, \mathbf{Y}_j^{\{S\}})$. With (9.19) the the slow component is solved with the slow (base) Runge–Kutta method. The internal fast ODE (9.19c) integrates (and corrects) only the fast component.

Similar to (9.12), it is convenient to define the time-dependent slow-to-fast coupling coefficients $\tilde{\gamma}(\tau) \in \mathbb{R}^{s^{\{S\}}}$ as polynomials in time:

$$\tilde{\gamma}(\tau) := \sum_{k \geq 0} \gamma^k \, \frac{\tau^{k+1}}{k+1}, \quad \gamma^k \in \mathbb{R}^{s^{\{S\}}}. \tag{9.20}$$

Internal consistency conditions are:

$$(\gamma^0)^{\mathrm{T}} \, \mathbf{1}^{\{S\}} = 1 \quad \text{and} \quad (\gamma^k)^{\mathrm{T}} \, \mathbf{1}^{\{S\}} = 0 \quad \forall k \geq 1. \tag{9.21}$$

An internally consistent method with a second order base scheme is second order [44].

Moreover, it was shown in [44] that an internally consistent SPC-MRI-GARK method (9.19) satisfying (9.21) has order three/four iff the slow base scheme $(A^{\{S\}}, b^{\{S\}}, c^{\{S\}})$ has order at least three/four, and the following coupling conditions hold:

$$\frac{1}{6} = \Big( \sum_{k \geq 0} \zeta_k \, \gamma^k \Big)^{\mathrm{T}} c^{\{S\}}, \qquad \text{(order 3)} \tag{9.22a}$$

$$\frac{1}{8} = \Big( \sum_{k \geq 0} \omega_k \, \gamma^k \Big)^{\mathrm{T}} c^{\{S\}}, \qquad \text{(order 4)} \tag{9.22b}$$

$$\frac{1}{12} = \Big( \sum_{k \geq 0} \zeta_k \, \gamma^k \Big)^{\mathrm{T}} c^{\{S\} \times 2}, \qquad \text{(order 4)} \tag{9.22c}$$

$$\frac{1}{24} = \Big( \sum_{k \geq 0} \zeta_k \, \gamma^k \Big)^{\mathrm{T}} A^{\{S\}} \, c^{\{S\}}, \qquad \text{(order 4)} \tag{9.22d}$$

where $c^{\times 2}$ denotes the component-wise power operation.

Remark 33. *A related family of schemes, named "internal stage predictor-corrector" (IPC-MRI-GARK), is also given in [44]. SPC-MRI-GARK (9.19) takes a full step with the base method, and then corrects the fast component solution at $t_{n+1}$. IPC-MRI-GARK interlace predictor and corrector stages: each stage slow solution* $\mathbf{Y}_i^{\{\mathrm{s}\}}$ *is computed with the base method, then the fast component* $\mathbf{Y}_i^{\{\mathrm{s}\}}$ *is corrected by integrating over* $[c_{i-1}, c_i]$. *One IPC-MRI-GARK step reads:*

$$\mathbf{Y}_0^{\{\mathrm{F}\}} = \mathbf{y}_n^{\{\mathrm{F}\}}, \quad \mathbf{Y}_0^{\{\mathrm{s}\}} = \mathbf{y}_n^{\{\mathrm{s}\}}, \quad c_0^{\{\mathrm{s}\}} = 0, \tag{9.23a}$$

$$\begin{cases} \begin{bmatrix} \mathbf{Y}_i^{\{\mathrm{F}\}*} \\ \mathbf{Y}_i^{\{\mathrm{s}\}} \end{bmatrix} = \begin{bmatrix} \mathbf{y}_n^{\{\mathrm{F}\}} + H \sum_{j=1}^{i} a_{i,j}^{\{\mathrm{s}\}} \mathbf{f}_j^{\{\mathrm{F}\}} \\ \mathbf{y}_n^{\{\mathrm{s}\}} + H \sum_{j=1}^{i} a_{i,j}^{\{\mathrm{s}\}} \mathbf{f}_j^{\{\mathrm{s}\}} \end{bmatrix}, \\ v^{\{\mathrm{F}\}}(0) = \mathbf{Y}_{i-1}^{\{\mathrm{F}\}}, \qquad T_{i-1} = t_n + c_{i-1}^{\{\mathrm{s}\}} H, \\ v_i^{\{\mathrm{F}\}\prime} = \Delta c_i^{\{\mathrm{s}\}} \mathbf{f}^{\{\mathrm{F}\}} \left( v_i^{\{\mathrm{F}\}}, \mathbf{Y}_{i-1}^{\{\mathrm{s}\}} + H \sum_{j=1}^{i} \tilde{\delta}_{i,j}\left(\frac{\theta}{H}\right) \mathbf{f}_j^{\{\mathrm{s}\}} \right), \\ \qquad \textit{for } \theta \in [0, H], \\ \mathbf{Y}_i^{\{\mathrm{F}\}} = v_i^{\{\mathrm{F}\}}(H), \quad i = 1, \dots, s^{\{\mathrm{s}\}}, \end{cases} \tag{9.23b}$$

$$\begin{bmatrix} \mathbf{y}_{n+1}^{\{\mathrm{F}\}} \\ \mathbf{y}_{n+1}^{\{\mathrm{s}\}} \end{bmatrix} = \begin{bmatrix} \mathbf{Y}_{s^{\{\mathrm{s}\}}}^{\{\mathrm{F}\}} \\ \mathbf{Y}_{s^{\{\mathrm{s}\}}}^{\{\mathrm{s}\}} \end{bmatrix}, \tag{9.23c}$$

*where* $\tilde{\delta}_{i,j}\left(\frac{\theta}{H}\right)$ *are time-dependent coupling coefficients. The advantage of the IPC strategy is that all subsequent slow stages* $\mathbf{Y}_j^{\{\mathrm{s}\}}$, $j > i$, *benefit from the more accurate fast components. The IPC order conditions are more complex than the SPC ones.*

**9.6. MRI-GARK-ROS/ROW methods.** Multirate infinitesimal MRI-GARK-ROS/ROW methods were developed in [28]. Consider a base slow ROS/ROW scheme $(\mathsf{b}^{\{\mathrm{s}\}}, \boldsymbol{\alpha}^{\{\mathrm{s,s}\}}, \boldsymbol{\gamma}^{\{\mathrm{s,s}\}})$ with non-decreasing abscissae $\mathsf{c}_1^{\{\mathrm{s,s}\}} \le \mathsf{c}_2^{\{\mathrm{s,s}\}} \le \cdots \le \mathsf{c}_{s^{\{\mathrm{s}\}}}^{\{\mathrm{s,s}\}}$, and denote:

$$\mathsf{c}_0^{\{\mathrm{s,s}\}} = 0; \qquad \Delta \mathsf{c}_i^{\{\mathrm{s,s}\}} := \mathsf{c}_i^{\{\mathrm{s,s}\}} - \mathsf{c}_{i-1}^{\{\mathrm{s,s}\}}, \quad i = 1, \dots, s.$$

For simplicity of presentation we discuss this scheme for the additively partitioned system (2.2). A multirate infinitesimal step GARK ROS/ROW method advances the solution of (2.2) via the following computational process:

$$\tilde{\mathbf{y}}_{n-1} = \mathbf{y}_{n-1}; \tag{9.24a}$$

For $\lambda = 1, \dots, s^{\{\mathrm{s}\}}$ :

$$\begin{cases} \mathbf{v}_\lambda(0) = \tilde{\mathbf{y}}_{n-1+(\lambda-1)/s^{\{\mathrm{s}\}}}, \\ \mathbf{v}_\lambda'(\theta) = \Delta \mathsf{c}_\lambda^{\{\mathrm{s,s}\}} \mathbf{f}^{\{\mathrm{F}\}} \left( \mathbf{v}_\lambda(\theta) + \sum_{j=1}^{\lambda-1} \boxed{\mathbf{r}_{\lambda,j}\left(\tfrac{\theta}{H}\right)} \mathbf{k}_j^{\{\mathrm{s}\}} \right), \qquad \theta \in [0, H], \\ \tilde{\mathbf{y}}_{n-1+\lambda/s^{\{\mathrm{s}\}}} = \mathbf{v}_\lambda(H); \end{cases} \tag{9.24b}$$

$$\begin{aligned} \mathbf{k}_\lambda^{\{\mathrm{s}\}} &= H\, \mathbf{f}^{\{\mathrm{s}\}} \left( \tilde{\mathbf{y}}_{n-1+\lambda/s^{\{\mathrm{s}\}}} + \sum_{j=1}^{\lambda-1} \boldsymbol{\alpha}_{\lambda,j}^{\{\mathrm{s,s}\}} \mathbf{k}_j^{\{\mathrm{s}\}} \right) \\ &\quad + H\, \mathbf{L}^{\{\mathrm{s}\}} \left( \sum_{\ell=1}^{\lambda} \boxed{\mathbf{q}_{\lambda,\ell}} (\tilde{\mathbf{y}}_{n-1+\ell/s^{\{\mathrm{s}\}}} - \tilde{\mathbf{y}}_{n-1+(\ell-1)/s^{\{\mathrm{s}\}}}) + \sum_{j=1}^{\lambda} \boldsymbol{\gamma}_{\lambda,j}^{\{\mathrm{s,s}\}} \mathbf{k}_j^{\{\mathrm{s}\}} \right); \end{aligned} \tag{9.24c}$$

$$\mathbf{y}_n = \tilde{\mathbf{y}}_n + (\mathsf{b}^{\{\mathrm{s}\}T} \otimes \mathbf{I}_{d \times d})\, \mathbf{k}^{\{\mathrm{s}\}}, \tag{9.24d}$$

A modified fast ODE (9.24b) is integrated between consecutive stages of the base slow method. The slow components influence the fast dynamics via the time dependent coefficients $\mathbf{r}_\lambda(\cdot) \in \mathbb{R}^{s^{\{s\}}}$, $\lambda = 1, \dots, s^{\{s\}}$. The fast solutions impact the computation of the slow stages (9.24c) via the coupling coefficients $\mathbf{q}_\lambda \in \mathbb{R}^{s^{\{s\}}}$, $\lambda = 1, \dots, s^{\{s\}}$. The next step solution (9.24d) combines the fast solution $\tilde{\mathbf{y}}_n$ and the slow solution increment given by the stages $\mathbf{k}^{\{s\}}$.

The internal consistency conditions read:

$$\mathbf{r}_\lambda(\theta) = \textstyle\sum_{k\ge 0} \mathbf{r}_{\lambda,k}\, \theta^k \in \mathbb{R}^{s^{\{s\}}}, \quad \bar{\mathbf{r}}_\lambda := \textstyle\sum_{k\ge 0} \frac{1}{k+1}\, \mathbf{r}_{\lambda,k},$$
$$\mathbf{r}_{\lambda,0}{}^T\, \mathbf{1}^{\{s\}} = \mathsf{c}_{\lambda-1}^{\{s,s\}}, \quad \mathbf{r}_{\lambda,1}{}^T\, \mathbf{1}^{\{s\}} = \Delta \mathsf{c}_{\lambda}^{\{s,s\}}, \quad \mathbf{r}_{\lambda,k}{}^T\, \mathbf{1}^{\{s\}} = 0,\ k \ge 2.$$

The order three MRI-GARK-ROW coupling conditions (see Remark 21) are:

$$\mathsf{b}^{\{s\}T} \sum_{\lambda=1}^{s} \begin{bmatrix} \mathbf{0}_{(\lambda-1)\times 1} \\ q_\lambda \end{bmatrix} \left( \mathsf{c}_{\lambda}^{\{s,s\}2} - \mathsf{c}_{\lambda-1}^{\{s,s\}2} \right) = 0,$$
$$\sum_{\lambda=1}^{s} \Delta \mathsf{c}_{\lambda}^{\{s,s\}}\, \bar{\mathbf{r}}_\lambda{}^T\, \mathsf{c}^{\{s,s\}} = \frac{1}{6}, \qquad \sum_{\lambda=1}^{s} \Delta \mathsf{c}_{\lambda}^{\{s,s\}}\, \bar{\mathbf{r}}_\lambda{}^T\, \mathsf{g}^{\{s,s\}} = 0.$$

**9.7. Step predictor-corrector MRI-GARK-ROS/ROW.** Step predictor-corrector MRI-GARK-ROS/ROW methods were also developed in [28]. Consider a step-predictor-corrector method (7.11) with the slow base method $(\mathsf{b}^{\{s\}}, \boldsymbol{\alpha}^{\{s,s\}}, \boldsymbol{\gamma}^{\{s,s\}})$.

The computation starts with a step predictor by applying the ROW method to the entire system, as in Section 7.6. This builds the ROW slopes $\mathbf{k}^{\{s\}}$ and $\mathbf{k}^{\{F\}}$. The slow component is advanced using the discrete ROW formula:

$$\mathbf{y}_{n+1}^{\{s\}} = \mathbf{y}_{n}^{\{s\}} + (\mathsf{b}^{\{s\}T} \otimes \mathbf{I}_{d\times d})\mathbf{k}^{\{s\}}. \tag{9.25}$$

For the fast component, a step corrector is applied. An accurate fast component is obtained solving a modified fast ODE over the step:

$$\begin{aligned} \mathbf{v}(0) &= \mathbf{y}_n^{\{F\}}, \\ \mathbf{v}' &= \mathbf{f}^{\{F\}}\Big(\mathbf{v}, \mathbf{y}_n^{\{s\}} + (\boldsymbol{\mu}^T(\tfrac{\theta}{H}) \otimes \mathbf{I}_{d\times d})\mathbf{k}^{\{s\}}\Big), \quad 0 \le \theta \le H; \\ \mathbf{y}_{n+1}^{\{F\}} &= \mathbf{v}(H). \end{aligned} \tag{9.26}$$

Define the $\boldsymbol{\mu}(t) \in \mathbb{R}^{s^{\{s\}}\times 1}$ coupling coefficients as polynomials in time:

$$\boldsymbol{\mu}(t) := \sum_{k\ge 0} \boldsymbol{\mu}_k\, \frac{t^k}{k+1}, \qquad \overline{\boldsymbol{\mu}} := \boldsymbol{\mu}(1).$$

The internal consistency conditions are satisfied with:

$$\boldsymbol{\mu}_1{}^T\, \mathbf{1}^{\{s\}} = 1; \quad \boldsymbol{\mu}_k{}^T\, \mathbf{1}^{\{s\}} = 0,\ k \ne 1.$$

*SPC-MRI-GARK-ROS methods.* The order three ROS condition reads:

$$\overline{\boldsymbol{\mu}}^T\, \mathsf{e}^{\{s,s\}} = \frac{1}{6}.$$

The order four ROS conditions are:

$$\sum_{k\ge 0} \frac{\boldsymbol{\mu}_k{}^T\, \mathsf{e}^{\{s,s\}}}{k+2} = \frac{1}{8}, \qquad \overline{\boldsymbol{\mu}}^T\, \mathsf{c}^{\{s,s\}\times 2} = \frac{1}{12},$$
$$\sum_{k\ge 0} \frac{\boldsymbol{\mu}_k{}^T\, \mathsf{e}^{\{s,s\}}}{(k+1)(k+2)} = \frac{1}{24}, \qquad \overline{\boldsymbol{\mu}}^T\, \boldsymbol{\beta}^{\{s,s\}}\, \mathsf{e}^{\{s,s\}} = \frac{1}{24}.$$

For example, an MRIl version of the celebrated order 4 RODAS method computes the slow component using the standard formula (9.25). The the fast correction (9.26) is computed with with the coefficients:

$$
\begin{aligned}
&\boldsymbol{\mu}(t) \coloneqq \boldsymbol{\mu}_0 + \boldsymbol{\mu}_1\, t,\\
&\boldsymbol{\mu}_0 = \begin{bmatrix} \theta_1 \\ \theta_2 \\ -1.923968128204745\,\theta_1 + 2.446324727549974\text{e-}01\,\theta_2 + 4.509689603795104\text{e-}02 \\ 1.405229246707428\,\theta_1 - 2.181782847233643\,\theta_2 + 1.289372580090594\text{e-}01 \\ -4.812611185026834\text{e-}01\,\theta_1 + 9.371503744786454\text{e-}01\,\theta_2 - 1.740341540470105\text{e-}01 - \theta_3 \\ \theta_3 \end{bmatrix},\\
&\boldsymbol{\mu}_1 = \begin{bmatrix} -2\,\theta_1 + 4.061438468864431\text{e-}01 \\ -2\,\theta_2 + 5.932358823451654\text{e-}01 \\ 3.847936256409489\,\theta_1 - 4.892649455099948\text{e-}01\,\theta_2 - 3.657016798231872\text{e-}01 \\ -2.810458493414855\,\theta_1 + 4.363565694467286\,\theta_2 + 5.003688760525202\text{e-}02 \\ 9.625222370053667\text{e-}01\,\theta_1 - 1.874300748957291\,\theta_2 + 3.162850629863266\text{e-}01 - \theta_4 \\ \theta_4 \end{bmatrix}.
\end{aligned}
$$

**10. Conclusions.** This paper gives an overview on multirate schemes for initial-value problems of ordinary differential equations. Multirate schemes exploit different system dynamics given by an ODE split either component-wise or additively into a slow and a fast part. The aim is to construct such schemes which inherit the convergence order (and maybe the good stability properties) of the underlying scheme, but at the same time save computational costs by sampling the slow and fast part at a different pace.

A straightforward way to define multirate schemes is to solve each of the subsystems independently with an interpolation and extrapolation procedure, respectively, taking care of the coupling of the different parts. We have seen that one gets an order $p$ scheme if the basic scheme has order $p$ and the interpolation/extrapolation procedure at least order $p-1$. The most basic approach here defines multirate Euler schemes, which may vary by the order the subsystems are solved and the level of implicitness in the coupling. We have seen that unconditional linear stability holds, if the implicit Euler scheme is used as basic scheme and the overall scheme contains some implicitness in coupling the slow and fast part. Another example for interpolation/extrapolation-based multirate schemes are multirate linear multistep schemes. Besides the interpolation/extrapolation based treatment of coupling variables, the coupling can be done by defining internal stages for the coupling parts. Multirate Runge-Kutta and linear-implicit Runge-Kutta schemes fall into this class, as well their extensions in the GARK framework. We finish the overview on multirate schemes by discussing multirate extensions of extrapolation schemes and infinitesimal schemes.

**Acknowledgments.** The work of A. Sandu was supported in part by the the National Science Foundation through awards DMS–2436357, DMS–2411069, by the U.S. Department of Energy, and by the Computational Science laboratory at Virginia Tech. The work of M. Günther was supported in part by the German Research foundation (DFG) through the DFG-research unit FOR 5269.

REFERENCES

[1] A. Aitken, On interpolation by iteration of proportional parts without the use of differences, Proc. Edinburgh Math. Soc., 3 (1932), pp. 56–76.

[2] R. Alexander, Diagonally implicit Runge–Kutta methods for stiff ODE's, SIAM Journal on Numerical Analysis, 14 (1977), pp. 1006–1021.

[3] M. Arnold, Multi-rate time integration for large scale multibody system models, in IUTAM Symposium on Multiscale Problems in Multibody System Contacts: Proceedings of the IUTAM Symposium held in Stuttgart, Germany, February 20–23, 2006, P. Eberhard, ed., Springer Netherlands, Dordrecht, 2007, pp. 1–10.

[4] A. Bartel and M. Günther, A multirate W-method for electrical networks in state-space formulation, Journal of Computational Applied Mathematics, 147 (2002), pp. 411–425.

[5] A. Bartel and M. Günther, Inter/extrapolation-based multirate schemes: A dynamic-iteration perspective, in Progress in Differential-Algebraic Equations II, T. Reis, S. Grundel, and S. Schöps, eds., Cham, 2020, Springer International Publishing, pp. 73–90.

[6] ———, Multirate Schemes — An Answer of Numerical Analysis to a Demand from Applications, Springer International Publishing, Cham, 2022, pp. 5–27.

[7] J. Butcher, A history of Runge–Kutta methods, Applied Numerical Mathematics, 20 (1996).

[8] ———, Numerical Methods for Ordinary Differential Equations, Wiley, 3 ed., 2016.

[9] J. T. Chang, S. R. P. Ploen, P. Sohl, Garett. A, and B. J. Martin, Parallel multi-step /multi-rate integration of two-time scale dynamic systems, tech. rep., Pasadena, CA : Jet Propulsion Laboratory, National Aeronautics and Space Administration, 2004.

[10] E. Constantinescu and A. Sandu, Multirate timestepping methods for hyperbolic conservation laws, Journal on Scientific Computing, 33 (2007), pp. 239–278.

[11] ———, Explicit time stepping methods with high stage order and monotonicity properties, in International Conference on Computational Science (ICCS-2009), G. Allen, J. Nabrzyski, E. Seidel, G. van Albada, J. Dongarra, and P. Sloot, eds., vol. 5545, International Conference on Computational Science, 2009, pp. 293–301.

[12] ———, On extrapolated multirate methods, in Progress in Industrial Mathematics at ECMI 2008, H.-G. Bock, F. Hoog, A. Friedman, A. Gupta, H. Neunzert, W. R. Pulleyblank, T. Rusten, F. Santosa, A.-K. Tornberg, V. Capasso, R. Mattheij, H. Neunzert, O. Scherzer, A. D. Fitt, J. Norbury, H. Ockendon, and E. Wilson, eds., vol. 15 of Mathematics in Industry, Springer Berlin Heidelberg, 2010, pp. 341–347.

[13] ———, Extrapolated multirate methods for differential equations with multiple time scales, Journal of Scientific Computing, 56 (2013), pp. 28–44.

[14] A. Demirel, J. Niegemann, K. Busch, and M. Hochbruck, Efficient multiple time-stepping algorithms of higher order, Journal of Computational Physics, 285 (2015), pp. 133–148.

[15] P. Deuflhard, Recent progress in extrapolation methods for ordinary differential equations, SIAM Review, 27 (1985), pp. 505–535.

[16] P. Deuflhard, E. Hairer, and J. Zugck, One-step and extrapolation methods for differential-algebraic systems, Numer. Math., 51 (1987), pp. 501–516.

[17] D.R. Wells, Multirate linear multistep methods for the solution of systems of ordinary differential equations, PhD thesis, University of Illinois at Urbana-Champaign, Dept. of Computer Science, 1982.

[18] C. Engstler and C. Lubich, Multirate extrapolation methods for differential equations with different time scales, Computing, 58 (1997), pp. 173–185.

[19] G. Rodríguez-Gómez, P. González-Casanova, P. Martínez-Carballido, J., Computing general companion matrices and stability regions of multirate methods, Int. J.for Numerical Methods in Engineering, (2004), pp. 255–273.

[20] M. Gasca and T. Sauer, Polynomial interpolation in several variables, Advances in Computational Mathematics, 12 (2000), pp. 377–410.

[21] C. Gear, Automatic multirate methods for ordinary differential equations, in Information Processing, NTIS, Springfield, Virginia, 1980, p. 16.

[22] C. Gear and D. Wells, Multirate linear multistep methods, BIT, 24 (1984), pp. 484–502.

[23] W. Gragg, Repeated extrapolation to the limit in the numerical solution of ordinary differential equations, PhD thesis, University of California, Los Angeles - Mathematics, 1964.

[24] ———, On extrapolation algorithms for ordinary initial value problems, Journal of the Society for Industrial and Applied Mathematics, Series B: Numerical Analysis, 2 (1965), pp. 384–403.

[25] A. Guennouni, A. Verhoeven, E. Maten, and T. Beelen, Aspects of multirate time integration methods in circuit simulation problems, in Progress in Industrial Mathematics at ECMI 2004, H. B. et al., ed., vol. 8 of Mathematics in Industry, Springer Berlin Heidelberg, 2006, pp. 579–584.

[26] M. Günther and P. Rentrop, Multirate ROW-methods and latency of electric circuits, Applied Numerical Mathematics, 13 (1993), pp. 83–102.

[27] M. Günther and A. Sandu, Multirate generalized additive Runge-Kutta methods, Numerische Mathematik, 133 (2016), pp. 497–524.

[28] M. Günther and A. Sandu, Multirate linearly-implicit GARK schemes, BIT Numerical Mathematics, 62 (2022), pp. 869–901.
[29] E. Hairer and C. Lubich, Asymptotic expansions of the global error of fixed-stepsize methods, Numerische Mathematik, 45 (1984), pp. 345–360.
[30] E. Hairer, S. Norsett, and G. Wanner, Solving ordinary differential equations I: Nonstiff problems, no. 8 in Springer Series in Computational Mathematics, Springer-Verlag Berlin Heidelberg, 1993.
[31] E. Hairer and G. Wanner, Solving ordinary differential equations II: Stiff and differential-algebraic problems, no. 14 in Springer Series in Computational Mathematics, Springer-Verlag Berlin Heidelberg, 2 ed., 1996.
[32] A. Hurwitz, Ueber die bedingungen, unter welchen eine gleichung nur wurzeln mit negativen reellen theilen besitzt, Mathematische Annalen, 46 (1895), pp. 273–284.
[33] O. Knoth and J. Wensch, Generalized split-explicit Runge–Kutta methods for the compressible euler equations, Monthly Weather Review, 142 (2014).
[34] O. Knoth and R. Wolke, Implicit-explicit Runge-Kutta methods for computing atmospheric reactive flows, Applied Numerical Mathematics, 28 (1998).
[35] M. Kutta, Beitrag zur naherungsweisen integration totaler differentialgleichungen, Z. Math. Phys., 46 (1901).
[36] A. Kværnø, Stability of multirate Runge-Kutta schemes, International Journal of Differential Equations and Applications, 1 (2000), pp. 97–105.
[37] A. Kværnø and P. Rentrop, Low order multirate Runge-Kutta methods in electric circuit simulation, 1999.
[38] E. Neville, Iterative interpolation, J. Indian Math. Soc., 20 (1934), pp. 87–120.
[39] Orailoglu, A., A multirate ordinary differential equation integrator, PhD thesis, University of Illinois at Urbana-Champaign, Dept. of Computer Science, 1979.
[40] A. Ralston, Runge–Kutta methods with minimum error bounds, Mathematics of Computation, 16 (1962), pp. 431–437.
[41] J. Rice, Split Runge-Kutta methods for simultaneous equations, Journal of Research of the National Institute of Standards and Technology, 64 (1960).
[42] S. Roberts, Graph interpretation of decoupled multirate methods. Personal communication, 2004.
[43] S. Roberts, J. Loffeld, A. Sarshar, C. Woodward, and A. Sandu, Implicit multirate GARK methods, Journal of Scientific Computing, 87 (2021), p. 4.
[44] S. Roberts, A. Sarshar, and A. Sandu, Coupled multirate infinitesimal GARK methods for stiff differential equations with multiple time scales, SIAM Journal on Scientific Computing, 42 (2020), pp. A1609–A1638.
[45] E. Routh, A treatise on the stability of a given state of motion: Particularly steady motion. Macmillan, 1877.
[46] C. Runge, Uber die numerische auflösung yon differentialgleichungen, Math. Ann., 46 (1895).
[47] Sand, J. and Skelboe, S., Stability of backward Euler multirate methods and convergence of waveform relaxation, BIT, 32 (1992), pp. 350–366.
[48] A. Sandu, A class of multirate infinitesimal GARK methods, SIAM Journal on Numerical Analysis, 57 (2019), pp. 2300–2327.
[49] A. Sandu and E. Constantinescu, Multirate, multimethod time stepping techniques for PDEs, in Tenth Copper Mountain conference on iterative methods, April 2008, p. 70.
[50] ———, Multirate Adams methods for hyperbolic equations, Journal of Scientific Computing, 38 (2009), pp. 229–249.
[51] ———, Multirate time discretizations for hyperbolic partial differential equations, in International Conference of Numerical Analysis and Applied Mathematics (ICNAAM 2009), vol. 1168-1 of American Institute of Physics (AIP) Conference Proceedings, 2009, pp. 1411–1414.
[52] A. Sandu and M. Günther, A generalized-structure approach to additive Runge-Kutta methods, SIAM Journal on Numerical Analysis, 53 (2015), pp. 17–42.
[53] A. Sarshar, S. Roberts, and A. Sandu, Design of high-order decoupled multirate GARK schemes, SIAM Journal on Scientific Computing, 41 (2019), pp. A816–A847.
[54] V. Savcenco, Comparison of the asymptotic stability properties for two multirate strategies, Journal of Computational and Applied Mathematics, 220 (2008), pp. 508–524.
[55] ———, Construction of a multirate rodas method for stiff odes, Journal of Computational and Applied Mathematics, 225 (2009), pp. 323 – 337.
[56] V. Savcenco, W. Hundsdorfer, and J. Verwer, A multirate time stepping strategy for parabolic PDE, Tech. Rep. MAS-E0516, Centrum voor Wiskundeen Informatica, 2005.
[57] V. Savcenko, W. Hundsdorfer, and J. Verwer, Multirate time–stepping methods for

hyperbolic conservation laws, BIT Numer. Math., 47 (2007), pp. 137–155.
[58] ———, A multirate time stepping strategy for stiff odes, BIT Numerical Mathematics, 47 (2007), pp. 137–155.
[59] M. Schlegel, O. Knoth, M. Arnold, and R. Wolke, Multirate Runge-Kutta schemes for advection equations, Journal of Computational and Applied Mathematics, 226 (2009), pp. 345–357.
[60] ———, Multirate implicit-explicit time integration schemes in atmospheric modelling, in AIP Conference Proceedings, vol. 1281, International Conference of Numerical Analysis and Applied Mathematics, 2010, pp. 1831–1834.
[61] M. Schlegel, O. Knoth, M. Arnold, and R. Wolke, Implementation of splitting methods for air pollution modeling, Geoscientific Model Development Discussions, 4 (2011), pp. 2937–2972.
[62] K. Schäfers, A. Bartel, M. Günther, and C. Hachtel, Spline-oriented inter/extrapolation-based multirate schemes of higher order, Applied Mathematics Letters, 136 (2023), p. 108464.
[63] J. M. Sexton and D. R. Reynolds, Relaxed multirate infinitesimal step methods, arXiv preprint arXiv:1808.03718, (2018).
[64] S. Skelboe, Stability properties of backward differentiation multirate formulas, Applied Numerical Mathematics, 5 (1989), pp. 151–160.
[65] S. Skelboe and U. Andersen, Stability properties of backward euler multirate formulas, SIAM J. Sci. Stat. Comput., 10 (1989), pp. 1000–1009.
[66] W. Throwe and S. Teukolsky, A high-order, conservative integrator with local time-stepping, SIAM Journal on Scientific Computing, 42 (2020), pp. A3730–A3760.
[67] U.M. Ascher and S.J. Ruuth and B.T.R. Wetton, Implicit-explicit methods for time-dependent partial differential equations, SIAM Journal on Numerical Analysis, 32 (1995), pp. 797–823.
[68] A. Verhoeven, T. Beelen, A. E. Guennouni, E. ter Maten, R. Mattheij, and B. Tasic, Error analysis of BDF compound-fast multirate method for differential-algebraic equations, tech. rep., Centre for analysis, scientific computing and applications, Technical University Eindhoven, 2006.
[69] A. Verhoeven, A. El Guennouni, E. Ter Maten, and R. Mattheij, A general compound multirate method for circuit simulation problems, in Scientific Computing in Electrical Engineering, Springer, 2006, pp. 143–149.
[70] A. Verhoeven, E. J. W. T. Maten, R. M. M. Mattheij, and B. Tasic, Stability analysis of the bdf slowest-first multirate methods, International Journal of Computer Mathematics, 84 (2007), pp. 895–923.
[71] A. Verhoeven, B. Tasic, T. Beelen, E. Maten, and R. Mattheij, Automatic partitioning for multirate methods, in Scientific Computing in Electrical Engineering, H.-G. Bock, F. Hoog, A. Friedman, A. Gupta, H. Neunzert, W. R. Pulleyblank, T. Rusten, F. Santosa, A.-K. Tornberg, V. Capasso, R. Mattheij, H. Neunzert, O. Scherzer, G. Ciuprina, and D. Ioan, eds., vol. 11 of Mathematics in Industry, Springer Berlin Heidelberg, 2007, pp. 229–236.
[72] A. Verhoeven, B. Tasiic, T. Beelen, E. ter Maten, and R. Mattheij, Bdf compound-fast multirate transient analysis with adaptive stepsize control, Journal of Numerical Analysis, Industrial and Applied Mathematics, 3 (2008), pp. 275–297.
[73] J. Wensch, O. Knoth, and A. Galant, Multirate infinitesimal step methods for atmospheric flow simulation, BIT Numerical Mathematics, 49 (2009), pp. 449–473.
[74] H. Zhang, S. Liang, S. Song, and H. Wanf, Truncation error calculation based on Richardson extrapolation for variable-step collaborative simulation, Science China Information Sciences, 54 (2011), pp. 1238–1250.